\documentclass{amsart}
\usepackage{amsmath,amsfonts,amssymb,amsthm}
\newcommand{\torusemptyset}{0}
\renewcommand{\epsilon}{\varepsilon}
\newcommand{\cx}{{\mathbb C}}
\newcommand{\integers}{{\mathbb Z}}
\newcommand{\Z}{{\mathbb Z}}
\newcommand{\natls}{{\mathbb N}}
\newcommand{\reals}{{\mathbb R}}
\newcommand\R[1]{{\mathbb R}^{#1}}
\newcommand{\st}{\;\: : \;\:}         
\newcommand{\nsp}{\hspace*{-0.7pt}}

\renewcommand{\mod}{\operatorname{mod}}
\newcommand{\Vol}{\operatorname{Vol}}
\newcommand{\vol}{\operatorname{vol}}
\newcommand{\area}{\operatorname{area}}

\newtheorem{theorem}{Theorem}[section]
\newtheorem*{Th}{Theorem}
\newtheorem{proposition}[theorem]{Proposition}
\newtheorem*{StarProposition}{Proposition}
\newtheorem{lemma}[theorem]{Lemma}
\newtheorem{corollary}[theorem]{Corollary}
\newtheorem{formula}{Formula}[section]

\theoremstyle{definition}

\newtheorem{Convention}{Convention}

\theoremstyle{remark}
\newtheorem{remark}[theorem]{Remark}
\newtheorem{example}[theorem]{Example}

\newcommand{\cC}{{\mathcal C}}
\newcommand{\cF}{{\mathcal F}}
\newcommand{\cH}{{\mathcal H}}
\newcommand{\cU}{{\mathcal U}}

\renewcommand{\setminus}{-}
\newlength{\halfbls}

\begin{document}
\ifx\href\undefined\else\hypersetup{linktocpage=true}\fi 
\begin{picture}(0,0)(0,0)
\put(-20,60){Publications de l'IHES,}
\put(-20,50){Vol. {\bf 97} no.1 (2003), pp. 61--179.}
\end{picture}

\date{November 11, 2003}

\title
[MODULI SPACES OF ABELIAN DIFFERENTIALS: THE PRINCIPAL BOUNDARY]
{Moduli Spaces of Abelian Differentials: \\ The Principal Boundary,
Counting Problems, \\ and the Siegel--Veech Constants}

\author{Alex Eskin}
\thanks{Research of the first author is partially supported by
                 NSF grant DMS-9704845, the  Sloan
                Foundation and the Packard Foundation}

%
\author{Howard Masur}
\thanks{Research of the second author is partially supported by
NSF grant 9803497}

%
\author{Anton Zorich}

\begin{abstract}

A holomorphic 1-form on  a  compact Riemann surface $S$ naturally
defines a flat metric on $S$  with  cone-type  singularities.  We
present the  following  surprising  phenomenon:  having  found  a
geodesic  segment  (saddle connection) joining a pair of  conical
points one can find with a  nonzero  probability  another  saddle
connection on  $S$ having the  same direction and the same length
as  the  initial one.  A  similar  phenomenon  is  valid  for the
families of parallel closed geodesics.

We give a  complete description of all possible configurations of
parallel saddle connections  (and  of families of parallel closed
geodesics) which might be found on a generic flat surface $S$. We
count the number of saddle connections of length less than $L$ on
a  generic  flat  surface  $S$;  we  also  count  the  number  of
admissible  configurations  of  pairs  (triples,...)  of   saddle
connections; we count the analogous numbers  of configurations of
families  of   closed   geodesics.   By   the   previous   result
of~\cite{Eskin:Masur}  these  numbers have  quadratic asymptotics
$c\cdot (\pi L^2)$. Here  we  explicitly compute the constant $c$
for a configuration of every type. The constant $c$ is found from
a Siegel---Veech formula.

To perform this computation we elaborate the detailed description
of  the  principal part of the  boundary  of the moduli space  of
holomorphic  1-forms  and  we  find  the  numerical value of  the
normalized volume of the tubular neighborhood of the boundary. We
use this for evaluation of integrals over the moduli space.
\end{abstract}

\vspace*{-4pt}

\maketitle
\tableofcontents

\addcontentsline{toc}{part}{Introduction}
{\Large\bf Introduction}
%
\section{Flat Surfaces and Geodesics on Flat Surfaces}
\label{sec:flat:surfaces:and:geodesics}

In  this  paper  we  shall  consider  flat metrics with  isolated
conical singularities on a closed  orientable  surface  of  genus
$g$. An  important class of  such flat metrics corresponds to the
{\it translation surfaces}: those surfaces  for  which  the  flat
metric on the surface has trivial linear holonomy. These surfaces
have been  studied by several authors  in various guises;  in the
context of  Abelian  differentials  (and more generally quadratic
differentials) on  compact Riemann surface in ~\cite{Strebel}, in
the  context  of  rational  billiards  in~\cite{Katok:Zemlyakov},
\cite{Kerckhoff:Masur:Smillie},           \cite{Veech:billiards},
and~\cite{Gutkin:Judge}. Moduli spaces of  these  structures have
been  studied  in   ~\cite{Veech:Siegel},   \cite{Masur:Smillie},
\cite{Kontsevich},                      \cite{Kontsevich:Zorich},
and~\cite{Eskin:Masur}.       A       survey       is       given
in~\cite{Masur:Tabachnikov}.

Triviality  of  holonomy means, in particular, that the  parallel
transport of  a vector along a small loop  going around a conical
point brings a vector back to itself. This implies that  all cone
angles of such surface are integer multiples of $2\pi$.

Choose  a  direction in  the  tangent space  to  some base  point
$x_0\in S$ of the flat surface and transport this direction using
the parallel transport  to all nonsingular points of the surface.
Since the monodromy is  trivial  this parallel transport does not
depend on  the path. Thus, any  direction is globally  defined on
the flat surface punctured at the  singularities. Throughout this
paper by a  {\it  flat  surface} we mean a  flat  surface  having
trivial holonomy  representation  in  $SO(2,\reals)$. Following a
tradition we often call the conical points the {\it saddles}.

\begin{figure}[ht]
%
\includegraphics{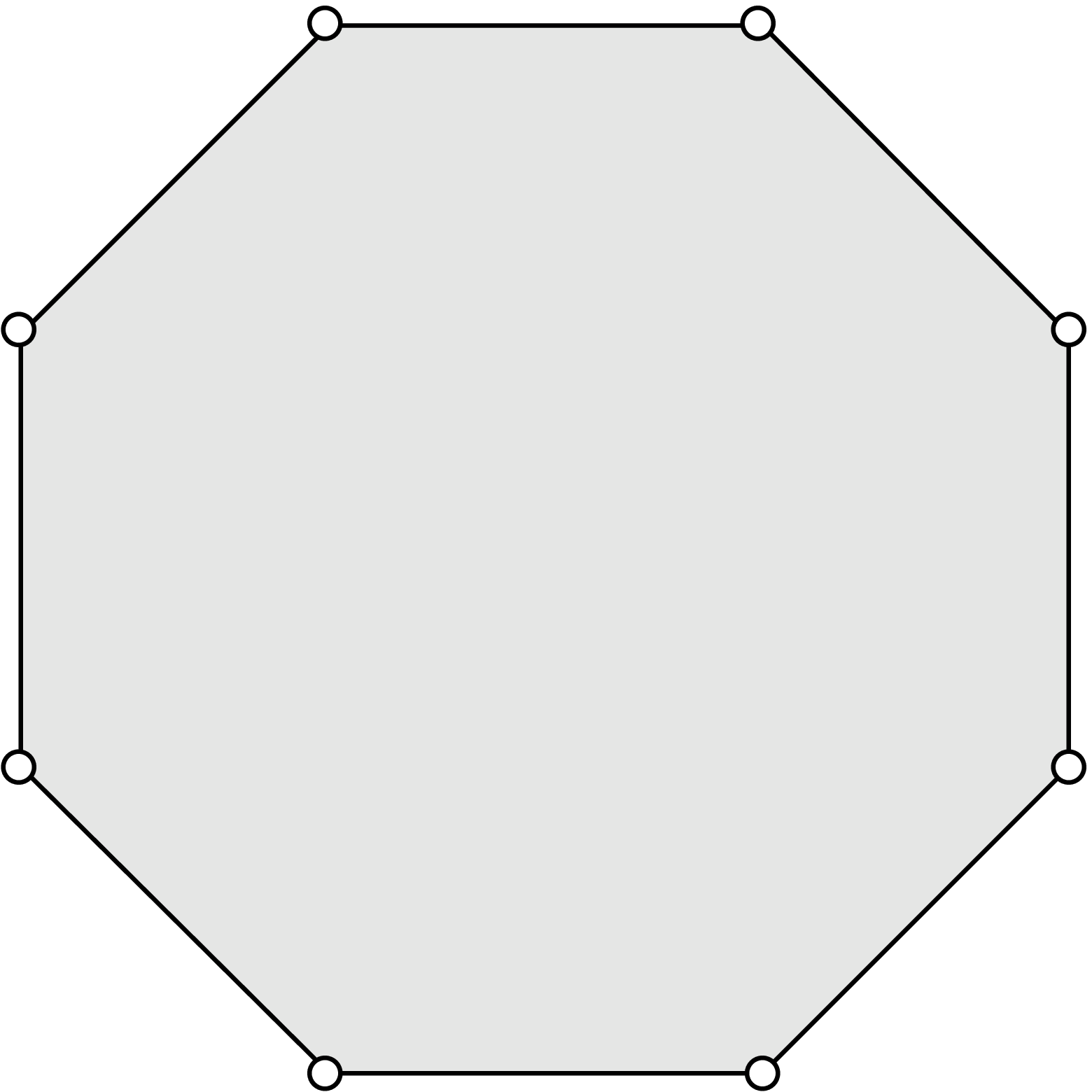}
\includegraphics{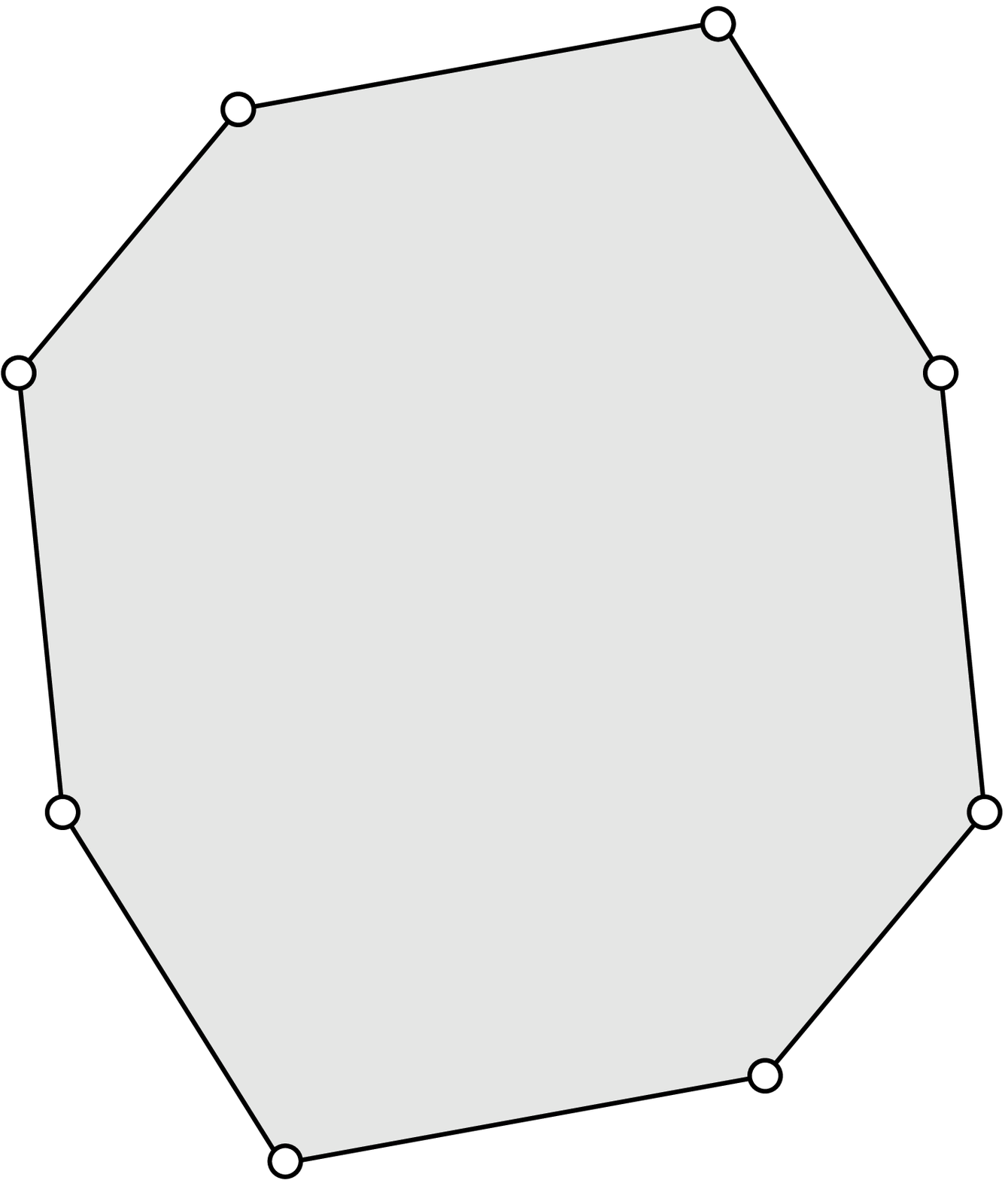}
\vspace{85bp} 
\caption{
\label{pic:octagon}
Identifying the  opposite  sides  of  these  octagons by parallel
translations we obtain a pair of flat surfaces of genus $g=2$.
}
\end{figure}

A flat torus  gives  an  example of a flat  surface  without  any
saddles.  One  can glue  a  flat torus  from  a parallelogram  by
identifying opposite  sides.  A  generalization  can  be found by
considering a  polygon in the  plane having the property that its
sides are  distributed into pairs  that are parallel and of equal
length. Identifying the pairs  by  translations one gets a closed
oriented flat  surface.  For  example,  identifying  the opposite
sides of a regular  octagon one gets a surface of genus  $2$, see
Figure~\ref{pic:octagon}.  All  vertices  of   the   octagon  are
identified  to a  single  saddle which has  an  angle of  $6\pi$.
Deforming the regular  octagon  in such  way  that the sides  are
organized into  pairs of parallel sides of equal  length we get a
{\it family} of flat surfaces of genus $g=2$ each having a single
conical singularity with the cone angle  $6\pi$.

We may assume that we have cartesian coordinates in the polygons.
In  these  coordinates away from the saddles,  geodesics  on  the
surface $S$ are straight lines. The behavior of geodesics on flat
surfaces is in many aspects similar to the  behavior of geodesics
on a flat  torus.  In particular,  a  geodesic cannot change  the
direction, and thus cannot have  self  intersections.  A  generic
geodesic is dense in the surface  analogously  to  an  irrational
winding line which is dense on the torus.

However,  a  flat  surface  different  from  the torus must  have
saddles and some geodesics  that  hit the saddles. Some geodesics
hit the saddles going both in forward and in backward directions.
Such a  geodesic segment joining a pair of  saddles and having no
saddles in its interior is called a {\it saddle connection}.

Note  that  the shortest representative in a  homotopy  class  of
curves joining two  saddles exists. Typically, even for the class
of a simple closed curve, one gets a broken line  containing many
geodesic  segments.  The segments  of  these  broken  lines  pass
through the saddles  changing directions there. We want to stress
that  throughout  this  paper  we  consider   only  those  saddle
connections which contain a  single  geodesic segment and thus no
saddles in its  interior. However, we consider the situation when
the  endpoints  of  a  saddle connection coincide. We  call  such
saddle connections the {\it closed saddle connections}.

A geodesic leaving  a  regular point $P$ may  return  back to $P$
without meeting any singularities. Since  it  cannot  change  its
direction, we necessarily  get a {\it closed regular geodesic} in
this case.  Note that  as in the case of  the torus, any geodesic
leaving  a nearby  point  in a parallel  direction  will be  also
closed and will  have  the same length as  the  initial one. Thus
closed regular geodesics appear in families of parallel geodesics
of the  same length. However, unlike the case  of the torus, such
parallel closed geodesics  do not typically fill all the surface,
but only a cylindrical subset. Each boundary component  of such a
cylinder  is   comprised  of  saddle  connections.  For  example,
choosing the vertical direction on the flat surface obtained from
a regular  octagon  (see  Figure~\ref{pic:octagon})  we  get  two
families  of  regular  closed  geodesics,  and  each  of  the two
families fills  a  cylinder. Generically, each boundary component
of a cylinder  filled with closed  geodesics is a  single  closed
saddle connection.

The converse, however,  is false. A closed saddle connection does
not necessarily bound  a cylinder of regular closed geodesics. In
fact, it bounds such a cylinder if and only if  the  angle at the
saddle between the outgoing and incoming segments is exactly $\pi
$. In this paper we  will  be interested in both counting  closed
saddle connections  and  counting  cylinders  filled with regular
closed geodesics.

Now fix a flat  surface $S$. Let $\hat S$ be the  universal cover
of  $S\setminus  \Sigma$,  where  $\Sigma$ is the  collection  of
saddles. There is  an  isometric map  of  $\hat S$ to  $\reals^2$
called the {\it developing map} which maps a lift $\hat\gamma$ of
an oriented curve $\gamma$ on  $S$  to a curve in $\reals^2$.  We
denote by $hol(\gamma)$ the difference  of  the  endpoints of the
image. The  vector  $hol(\gamma)\in\reals^2$  is  called the {\em
holonomy} vector of $\gamma$.

\begin{remark}
By definition a flat surface has  trivial holonomy representation
in the linear group
$$
\pi_1(S)\to H_1(S,\integers)\to SO(2,\reals).
$$
However, the holonomy representation in the {\it affine} group is
nontrivial: its image  belongs to the group of translations. From
now on we reserve the  notion  of {\it holonomy} for this  second
representation. It  matches  the notion of holonomy $hol(\gamma)$
defined by an  arbitrary  (not necessarily closed) oriented curve
$\gamma$ suggested in the paragraph above.
\end{remark}

Now fix a pair of saddles  on a flat surface $S$, and some length
$L$. Consider  all  those  saddle  connections  joining the fixed
saddles whose  length is shorter than  $L$. We are  interested in
the  asymptotics  of the number of saddle  connections  when  the
bound $L$  tends to  infinity. In the model case  of the torus of
unit area it  is  sufficient to  pass  to the universal  covering
plane to see that the  number  of geodesic segments of length  at
most $L$  joining a generic pair of distinct  points on the torus
grows quadratically as the number of lattice points in a  disc of
radius $L$, so we get asymptotics $\pi L^2$. Note that the number
of (homotopy  classes) of closed  geodesics of length at most $L$
has  different  asymptotics.  Since  we want to count  only  {\it
primitive} geodesics (those  which  do not repeat themselves) now
we have to count  only {\it coprime} lattice points in a  disc of
radius $L$,  considered up to a symmetry of  the torus (issues of
symmetry will be considered  in  details later). Therefore we get
the     asymptotics     $      \cfrac{1}{2\zeta(2)}\cdot      \pi
L^2=\cfrac{3}{\pi^2}\cdot \pi L^2$.

It is proved in~\cite{Eskin:Masur} that  the  growth  rate of the
number of saddle connections for a generic flat  surface also has
quadratic asymptotics  $c\cdot  (\pi L^2)$, and, moreover, almost
all flat surfaces in some natural families of flat surfaces share
the same  constant $c$ in the  asymptotics. One of  the principal
goals of this paper is to compute these constants for all natural
families  (which  are explicitly defined below). We also  compute
the  constants  in the quadratic asymptotics for  the  number  of
families of closed regular geodesics.

It is relatively easy to show that for a generic flat surface one
can  never find  a  pair of saddle  connections  having the  same
direction        but         different        lengths,        see
Proposition~\ref{pr:no:parallel:nonhomologous}  at   the  end  of
Section~\ref{s:volume:estimates}. However, somewhat surprisingly,
the configurations of pairs of parallel saddle connections of the
{\it same} length and same direction can be found on almost every
flat  surface,  and, what may seem even  more  surprising,  their
number also  has  quadratic  asymptotics  (with  another constant
$c$).

One can ask now the similar questions about the triples of saddle
connections, or about some more specific configurations of saddle
connections, or similar questions concerning closed geodesics. We
give a complete  description of those configurations which can be
found  on a generic  surface;  we  show  that each  of  them  has
quadratic  asymptotics  (this  was  essentially  already   proved
in~\cite{Eskin:Masur}),  and   in   every  case  we  compute  the
corresponding constant in the quadratic asymptotics.

\begin{remark}
A billiard in a polygon with rational angles gives rise to a flat
surface,       see       for      example~\cite{Katok:Zemlyakov},
\cite{Kerckhoff:Masur:Smillie},      or~\cite{Masur:Tabachnikov}.
However,  rational  billiards are  not  generic  in  the  natural
families of flat surfaces, which we consider. Thus, the constants
described in this  paper do not  apply to the  billiard  problem,
with an  exception  for  some  particular  classes off billiards,
see~\cite{Eskin:Zorich}.
\end{remark}

\subsection{Flat Surfaces and Abelian Differentials}

If we  cut a flat surface  $S$ successively along  an appropriate
collection  of  saddle  connections  we  can  decompose  it  into
polygons   contained   in   $\reals^2=\cx$,  see,  for   example,
Figure~\ref{pic:octagon}.  We  may then view $S$ as  a  union  of
polygons  glued  along  parallel  sides  by   means  of  parallel
translations.  Note  that we have endowed every  polygon  with  a
complex coordinate. Since our gluing rules are just translations,
the  transition  functions in these complex coordinates have  the
form
$$
z \to z + const.
$$
Thus any flat surface  with  the conical singularities removed is
endowed  with  a natural complex structure. Moreover, consider  a
holomorphic   1-form   $\omega=dz$  on   every   polygon.   Since
$dz=d(z+const)$ we  obtain  a globally defined holomorphic 1-form
on the surface with removed singularities. It is an easy exercise
to show that the complex structure and the holomorphic 1-form can
be extended to the singularities;  the  holomorphic  1-form  (or,
what is the same, the {\it Abelian differential})  $\omega$ has a
zero at every conical point.

Conversely,  given  a pair $(M,\omega)$ where $M$  is  a  Riemann
surface, $\omega$ is  a  holomorphic Abelian differential on $M$,
and a point $P \in M$ such that $\omega(P) \ne 0$, there exists a
local coordinate $z$  near  $P$ such that $\omega  =  dz$. Such a
local coordinate is  unique up to  change of coordinates  $z  \to
z+const$. Thus $|dz|^2$  is  a flat  metric  on $M$; this  metric
develops conical singularities at the  zeroes  of  $\omega$. At a
zero  of  order  $k$  the  total angle  is  $2  \pi  (k+1)$. Such
collection   of   coordinate  charts   determines   a   structure
$S=(M,\omega)$  of  a ``translation surface'' (which we agree  to
call  just  a ``flat  surface'')  on  $M$,  namely,  an  atlas of
coordinate  charts  which   cover   the  surface  away  from  the
singularities,   such   that   the   transition   functions   are
translations
$$
z \to z + const.
$$

When we consider  {\em moduli spaces}  of flat surfaces  we  will
wish  to  distinguish  between $\omega$ and  $e^{i\theta}\omega$.
This is equivalent to specifying a distinguished direction on the
flat surface.

\begin{Convention}
Depending on the context we shall use one of the synonyms
``saddle'', ``conical singularity'', ``zero of Abelian
differential'' or just ``zero''.
\end{Convention}

\section{Moduli Spaces of Abelian Differentials}
\label{ss:flat:sadd}
   %

\subsection{Stratification}

Let $\alpha$ be a  partition of $2 g - 2$ (i.e.  a representation
of  $2g  -  2$ as an  unordered  sum  of  positive integers). Let
$\cH(\alpha)$ denote the moduli space of pairs $(M,\omega)$ where
$M$ is a closed  Riemann surface of genus $g$, and $\omega$  is a
holomorphic Abelian differential on  $M$  such that the orders of
its  zeroes  is given by $\alpha$. Here  we  distinguish  between
$\omega$ and $e^{i\theta}\omega$. The zeroes on  each surface are
assumed to be {\em named}.  We  say that two such structures  are
{\em equivalent} if there is an isomorphism from one to the other
which takes  zeroes to zeroes,  preserving the naming. The set of
equivalence classes is denoted $\cH(\alpha)$ and is called a {\it
stratum}. This term is justified  by  the fact that the space  of
all Abelian differentials on Riemann surfaces  of  genus  $g$  is
stratified by the  spaces  $\cH(\alpha)$, as $\alpha$ varies over
the partitions of $2g-2$.  In the case of flat tori the  space is
$GL_+(2,\reals)/SL(2,\Z)$,  the  moduli space  of  lattices.  The
stratum corresponding  to  the  partition $(1,\dots,1)$ is called
the  {\it  principal  stratum}:  it  corresponds  to  holomorphic
differentials with simple zeroes. For  example,  the  surface  in
Figure~\ref{pic:m2simple}  belongs  to  $\cH(1,1)$.  The  stratum
might  be   nonconnected,   although  the  components  have  been
completely    classified    (see~\cite{Kontsevich:Zorich}).   The
classification  involves  the notion of parity of spin  structure
and  issues   of   hyperellipticity.  These  components  will  be
described later in Section~\ref{ss:con:comp}.

\subsection{Local Coordinates and Volume Element}

The stratum $\cH(\alpha)$  can be topologized and given a natural
``Lebesgue''   measure   as   follows.   For   a   flat   surface
$S_0\in\cH(\alpha)$ with conical points $P_1,\dots,P_k$, choose a
basis     of     cycles     for     the     relative     homology
$H_1(S_0,\{P_1,\dots,P_k\};\integers)$. This basis  may be chosen
in  such  a  way that each  element  is  represented  by a saddle
connection. Equivalently, the  saddle  connections cut $S$ into a
union of polygons. For any $S$ near $S_0$  the vectors associated
to these saddle  connections serve as local coordinates, see, for
example, Figure~\ref{pic:octagon}. By  construction these vectors
are  the  relative {\it  periods}  of  the  Abelian  differential
$\omega$ (i.e. the  integrals of $\omega$ along the paths joining
points  $P_i$ and  $P_j$,  where $i=j$ is  also  allowed). So  we
actually  use  a  domain  in  the  space  of  relative cohomology
$H^1(S_0,\{P_1,\dots,P_k\};\cx)$  as a  local  coordinate  chart.
Ignoring   the   complex   structure   we   get   a   domain   in
$\reals^n=H^1(S_0,\{P_1,\dots,P_k\};\cx)$,   where   $n=4g+2k-2$.
Note that we  have  a natural  cubic  lattice in our  coordinates
given by the ``integer'' cohomology
$$
H^1(S_0,\{P_1,\dots,P_k\};\integers\oplus i\integers)\subset
H^1(S_0,\{P_1,\dots,P_k\};\cx)=\reals^n
$$
We define a  measure or volume element $d\nu(S)$ on $\cH(\alpha)$
as Lebesgue measure  defined  by these coordinates, normalized so
that the volume of  a unit cube in $\reals^n$ is $1$.  One easily
checks that the volume element is well-defined: it is independent
of choice of basis.

Let  $\cH_1(\alpha)\subset\cH(\alpha)$  be  the  hypersurface  in
$\cH(\alpha)$ of unit area flat  surfaces.  Choose  a  symplectic
homology basis of closed curves $A_i,B_i$,  $i=1,\ldots,g$ on the
surface $S=(M,\omega)$. The  area of the flat surface $S$ defined
by  the  Abelian differential $\omega$ is given  by  the  Riemann
bilinear relation,
$$ 
\int_S |\omega|^2 dx dy=\frac{i}{2}\int_S
\omega\wedge\bar\omega=
\frac{i}{2}\sum_i\left(\int_{A_i}\omega\int_{B_i}\bar\omega-
\int_{A_i}\bar\omega\int_{B_i} \omega\right)
$$ 
Thus the unit  area surfaces represent a ``hyperboloid'' in terms
of the coordinates chosen above.

Throughout   this   paper   we   shall  mostly  work   with   the
``hyperboloid''  $\cH_1(\alpha)\subset\cH(\alpha)$.  The   volume
element in the embodying space $\cH(\alpha)$  induces naturally a
volume  element   on  the  hypersurface  $\cH_1(\alpha)$  in  the
following  way.   There   is   a   natural  $\cx^\ast$-action  on
$\cH(\alpha)$: having $\lambda\in\cx-\{0\}$ we  associate  to the
flat  surface  $S=(M,\omega)$  the  flat  surface   $\lambda\cdot
S=(M,\lambda\cdot\omega)$. In  particular,  we  can represent any
$S\in\cH(\alpha)$ as $S  = r S'$, where $r\in\reals_+$, and where
$S'$  belongs   to  the  ``hyperboloid'':   $S'\in\cH_1(\alpha)$.
Geometrically this means that the metric on $S$  is obtained from
the  metric  on  $S'$  by  rescaling  with  coefficient  $r$.  In
particular, vectors associated to saddle connections  on $S'$ are
multiplied by $r$ to give  vectors  associated  to  corresponding
saddle  connections  on  $S$.  It  means  also  that  $\area(S) =
r^2\cdot\area(S')=r^2$, since $\area(S') = 1$. We define the {\it
volume    element}    $d\vol(S')$    on    the    ``hyperboloid''
$\cH_1(\alpha)$ by disintegration of the volume element $d\nu(S)$
on $\cH(\alpha)$:
$$
d\nu(S) = r^{n-1} \, dr\, d\vol(S')\, ,
$$
where $n=\dim_\reals\cH(\alpha)$.  Using  this  volume element we
define the total {\it volume of the stratum} $\cH_1(\alpha)$:
\begin{equation}
\label{eq:int:cF:nu1}
\Vol(\cH_1(\alpha)):= \int_{\cH_1(\alpha)}d\vol(S')\,.
\end{equation}
Since the measure on $\cH_1(\alpha)$ is  induced  by  the  volume
element, the total volume of the stratum coincides with the total
measure.

For     a      subset     $E\subset\cH_1(\alpha)$     we      let
$C(E)\subset\cH(\alpha)$ denote the ``cone'' based on $E$:
\begin{equation}
\label{eq:cone}
 C(E):=\{S=rS'\,|\, S'\in E,\ 0<r\le 1\}\,.
\end{equation}
Our  definition  of the  volume  element  on  $\cH_1(\alpha)$  is
consistent with the following normalization:
\begin{equation}
\label{eq:normalization}
\Vol(\cH_1(\alpha)) = n\cdot\nu(C(\cH_1(\alpha))\,,
\end{equation}
where $n=\dim_\reals\cH(\alpha)$, and  $\nu(C(\cH_1(\alpha))$  is
the   total    volume    (total    measure)   of   the   ``cone''
$C(\cH_1(\alpha))\subset\cH(\alpha)$  measured by  means  of  the
volume element $d\nu(S)$ on $\cH(\alpha)$ defined above.

\begin{Convention}
\label{conv:vol:and:nu}
We need to use both volume elements. To distinguish them  we keep
the  notation  $d\nu(S)$  for  the  volume  element in the  whole
stratum $\cH(\alpha)$  and  notation  $d\vol(S)$  for  the volume
(hypersurface)    element    in    the    ``unit    hyperboloid''
$\cH_1(\alpha)$  which  is  a   hypersurface   in  $\cH(\alpha)$.
\end{Convention}

It  is   a   result  of  \cite{Masur},  \cite{Masur:Smillie}  and
\cite{Veech:moduli} that  these  volumes are finite. Their values
have     been     calculated    in~\cite{Eskin:Okounkov},     see
also~\cite{Zorich}  for  the  volumes  of  some   low-dimensional
strata.  They  will   be  needed  for  the  calculations  of  the
asymptotic constants. Tables of the values of the volumes for all
connected  components  of  the  strata  in  genera  $g\le  4$ are
presented in section~\ref{ss:examples:part1}.

\section{Principal  Boundary  of the Moduli Spaces and Method  of
Evaluation of the Siegel---Veech Constants}

\subsection{Holonomy and Homology}

Note that in terms of the local coordinates  $(x,y)$ defining the
flat structure  the  corresponding  Abelian  differential has the
form  $\omega  =  dz  =  dx  +  i\,dy$.  This  has  an  important
consequence:  for  any oriented curve $\beta$ on  $S$,  the  {\em
holonomy} of $\beta$ (as defined above in terms of the developing
map) coincides with the integral $\int_\beta  \omega$ of $\omega$
over $\beta$  (here  we  have  identified  $\reals^2$ and $\cx$).
Since the 1-form $\omega$ is closed,  it  means,  in  particular,
that if $\beta$  and $\gamma$ are homologous, then $hol(\gamma) =
hol(\beta)$. In particular, if we have  a  cylinder  filled  with
regular closed geodesics, then all these geodesics share the same
holonomy vector.

When a path $\beta$ joins a pair of distinct points we say that a
path $\gamma$ is homologous to $\beta$ if it joins the  same pair
of  points  and  if  the closed loop  $\beta\cdot\gamma^{-1}$  is
homologous  to  zero,   that  is  breaks  the  surface  into  two
components.  The  paths $\beta,\gamma$ which interest us in  this
paper are  saddle connections or  closed geodesics on $S$. We can
formalize the  above  observation  saying that saddle connections
$\beta,    \gamma$   representing    homologous    elements    in
$H_1(S,\{P_1,\dots,P_k\};\Z)$  (where $P_1,  \dots,P_k$  are  the
saddles) share the same holonomy vector $hol(\beta)=hol(\gamma)$.

\begin{figure}[ht]
%
\includegraphics{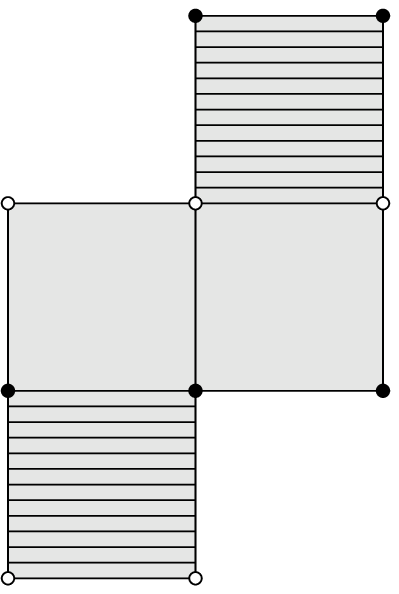}
\includegraphics{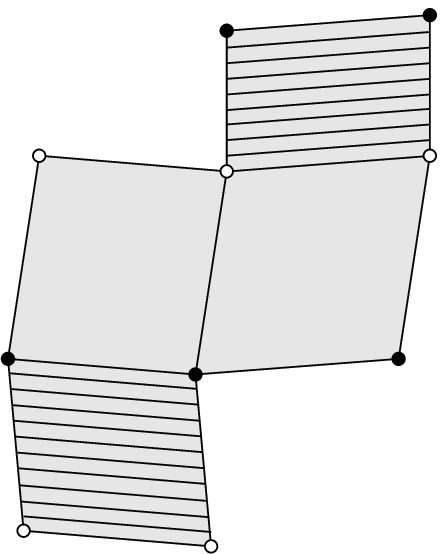}
\begin{picture}(0,0)(0,0)
\put(10,-10)
 {\begin{picture}(0,0)(0,0)
 \put(-73,-8){$\tau$}
 \put(-38,-40){$\eta$}
 \put(-105,-40){$\eta$}
 \put(-127,-60){$\rho$}
 \put(-160,-95){$\gamma$} 
 \put(-92,-95){$\gamma_1$}
 \put(-38,-95){$\gamma$} 
 \put(-73,-130){$\tau$}
 \put(-160,-150){$\beta$}
 \put(-92,-150){$\beta$}
 \put(-127,-185){$\rho$}
 \end{picture}}
\put(180,-10)
 {\begin{picture}(0,0)(0,0)
 \put(-68,-20){$\tau$} 
 \put(-30,-45){$\eta$} 
 \put(-100,-45){$\eta$}
 \put(-127,-60){$\rho$}   %
 \put(-160,-95){$\gamma$} %
 \put(-92,-95){$\gamma_1$}%
 \put(-35,-95){$\gamma$}  
 \put(-68,-133){$\tau$} 
 \put(-162,-150){$\beta$}
 \put(-95,-150){$\beta$}
 \put(-127,-185){$\rho$}
 \end{picture}}
\end{picture}
\vspace{200bp}
\caption{
\label{pic:m2simple}
Nonhomologous saddle  connections  which  have  the same holonomy
lose this property after a  generic  deformation  of the surface,
while homologous ones, $\gamma\sim\gamma_1$, keep the same holonomy.
}
\end{figure}

On the  other  hand  the  equality  $hol(\beta)=hol(\gamma)$  for
nonhomologous  $\beta,\gamma$   occurs   only   for  non  generic
surfaces. Consider the following example. Take  four unit squares
with sides glued as indicated  at  the  left-hand-side picture at
Figure~\ref{pic:m2simple}. The result is a closed  surface $S$ of
genus $2$. The five points indicated by the  filled circle symbol
are identified, as are the five points indicated  by the unfilled
circle symbol. Hence $S$  has  two conical singularities, each of
total angle $4\pi$. Each of the hatched squares  corresponds to a
cylinder filled  with  horizontal  regular  closed  geodesics  of
length $1$. Hence, for any closed  geodesics  $\sigma_1$  in  the
first cylinder and  any closed geodesics $\sigma_2$ in the second
cylinder we get $hol(\sigma_1) = hol(\sigma_2) = (1,0)$.

Note, however, that $\sigma_1$ and $\sigma_2$ are not homologous.
This  means that  the  relation above does  not  survive under  a
generic deformation of this flat surface, see the right-hand-side
picture at Figure~\ref{pic:m2simple}. On the other  hand the pair
of saddle connections $\gamma$ and $\gamma_1$  are homologous and
we   see   that    they    keep   sharing   the   same   holonomy
$hol(\gamma)=hol(\gamma_1)$ under  any  small  deformation of the
surface. Figure~\ref{pic:parallel:constr}  gives an example  when
such a relation between two distinct  families  of  {\em  closed}
geodesics  is  stable  under  any small deformation of  the  flat
surface.

This observation can be formalized as follows.

\begin{proposition}
\label{pr:nonhomologous:implies:nonholonomous}
Almost any flat  surface $S$ in  any connected component  of  any
stratum  does  not  have  a single pair of  nonhomologous  saddle
connections sharing the same holonomy vector.
\end{proposition}
\begin{proof}
Consider     a     basis     of     (relative)     cycles      in
$H_1(S,\{P_1,\dots,P_k\};\Z)$.  We know  that  the  corresponding
periods of $\omega$  (integrals over these basic cycles) serve as
local coordinates in the stratum  $\cH(\alpha)$.  Note  that  for
almost  all  points  of  the cohomology space these  periods  are
rationally independent:  no  linear combination of periods equals
to  zero  when  the  linear  combination  is taken with  rational
coefficients. Note  also, that if  the periods of a closed 1-form
$\omega$  are  rationally independent, then the same property  is
obviously   true   for  the   form   $\lambda\cdot\omega$   where
$\lambda\in\R{}$.  Hence,  this  is  true  for  almost  all  flat
surfaces in the ``unit hyperboloid'' $\cH_1(\alpha)$

A  pair  of saddle connections $\beta, \gamma$ represent  integer
cycles in  $H_1(S,\{P_1,\!\dots\!,P_k\};\Z)$.  Hence, if they are not
homologous, for almost all points $(S,\omega)$ of $\cH_1(\alpha)$
the       corresponding        periods       are       different:
$\int_{\beta}\omega\neq\int_{\gamma}\omega$.                Since
$hol(\beta)=\int_{\beta}\omega$                               and
$hol(\gamma)=\int_{\gamma}\omega$  the   Proposition  is  proved.
\end{proof}

\subsection{Configurations of Saddle Connections}
\label{sec:configurations}

Consider the set of all saddle connections on a flat  surface $S$
and consider its image $V_{sc}(S)\subset\R{2}$ under the holonomy
map.  (We  remind   the  reader  that  by  definition,  a  saddle
connection is represented  by a single geodesic segment joining a
pair  of  conical points, and not  by  a broken line of  geodesic
segments.)

Choose,  for  example,  a  flat torus $\R{2}/(\Z\oplus\Z)$  as  a
surface $S$ and mark a generic  pair of points $P_0, P_1$ on this
torus.      Let      $\vec{v}_0:=\vec{P_0 P_1}$.     The      set
$V_{sc}(S)\subset\R{2}$  is   the  square  lattice  of  the  form
$\{\vec{v}_0+\vec{u}\}_{\vec{u}\in\Z\oplus\Z}$.

Similarly to  the torus case, the  set $V_{sc}(S)$ is  a discrete
subset of $\reals^2$. However, when $S$  is  different  from  the
torus  the  map  from the set  of  saddle  connections  on $S$ to
$V_{sc}(S)$ is  not  injective:  different saddle connections may
have the  same holonomy. We define  the {\em multiplicity}  of an
element $\vec{v} \in V_{sc}(S)$ to be the number  $p$ of distinct
saddle  connections   $\gamma_1,   \dots,   \gamma_p$  such  that
$hol(\gamma_i)  =  \vec{v}$.  For  example, for the  surface  $S$
presented on  the  left  of  Figure~\ref{pic:m2simple} the vector
$(0,1)\in   V_{sc}(S)$   has  multiplicity   four:   the   saddle
connections $\beta,\gamma,\gamma_1,\eta$  chosen with appropriate
orientation have holonomy $(0,1)$.

In  this   example   the   saddle   connections   represented  by
$\beta,\gamma,\gamma_1,\eta$  are  not  pairwise  homologous.  As
indicated above, the fact that  their  holonomy  coincides is not
generic;  deforming  the   surface   slightly  we  see  that  the
corresponding saddle connections have different holonomy; see the
right-hand-side   picture   at   Figure~\ref{pic:m2simple}.   The
situation is different  with  the saddle connections $\gamma$ and
$\gamma_1$, for  they are homologous. The right-hand-side picture
at Figure~\ref{pic:m2simple}  confirms  that even after any small
deformation of  the  surface the corresponding saddle connections
share the same holonomy vector $hol(\gamma)=hol(\gamma_1)$.

From      the      measure-theoretical     point      of     view
Proposition~\ref{pr:nonhomologous:implies:nonholonomous}   allows
us to  assume from now on that if  two saddle connections $\beta,
\gamma$ have the same holonomy, then  $\gamma$  and  $\beta$  are
homologous. In  particular,  if  $\beta$  joins  distinct saddles
$P_1, P_2$,  then $\gamma$ joins  the same pair of saddles. Since
the   surfaces    as    in    the   left-hand-side   picture   at
Figure~\ref{pic:m2simple} form  a  set  of  measure  zero we will
ignore them in our considerations.

Somewhat surprisingly,  higher  multiplicity  is very common even
for  a  generic surface. For example,  we  will show that if  the
genus is  at least $3$, then  counting separately the  vectors of
multiplicity one in  $V_{sc}(S)$ of length  at most $L$  and  the
vectors of  multiplicity two, we  get quadratic growth in $L$ for
the number of elements of both multiplicities.

Multiplicity is not the only property  distinguishing elements in
$V_{sc}(S)$:  some   of   them   correspond   to   closed  saddle
connections,  other  elements correspond  to  saddle  connections
joining  distinct  zeroes.  We  may  also   consider  only  those
$\vec{v}\in V_{sc}(S)$,  which  correspond  to saddle connections
joining   some   particular  pair  of  saddles.  We  refine   our
consideration slightly more specifying the following data.

Suppose that we have precisely $p$  homologous saddle connections
$\gamma_1,\dots, \gamma_p$ joining a zero $z_1$ of order $m_1$ to
a zero  $z_2$  of order $m_2$, see Figure~\ref{pic:dist:sad:mult}
for a topological picture. All  the  $\gamma_i$, $1 \le i \le  p$
have the same  holonomy (where the orientation of each $\gamma_i$
is from $z_1$ to $z_2$).  By  convention the cyclic order of  the
$\gamma_i$ at $z_1$ is clockwise  in  the  orientation defined by
the  flat  structure.  Let  the  angle   between  $\gamma_i$  and
$\gamma_{i+1}$ at $z_1$ be $2\pi(a'_i+1)$; let  the angle between
$\gamma_i$ and  $\gamma_{i+1}$  at  $z_2$  be $2\pi(a''_i+1)$. We
call this data the {\em configuration}
$$
\cC= (m_1,m_2,a_i',a_i'')
$$
of the $p$  homologous  saddle connections.  If  there is just  a
single saddle connection ($p=1$) then $\cC=(m_1,m_2)$.

\begin{Convention}
We reserve the notion {\it ``configuration''} for geometric types
of possible {\it  collections} of saddle connections, and not for
the collections themselves.
\end{Convention}

Now  given  a  surface   $S$   and  a  configuration  $\cC$,  let
$V_{\cC}(S)\subset  V_{sc}(S)$  denote  the vectors $\vec{v}$  in
$\reals^2$ such that  there  are precisely $p$ saddle connections
$\gamma_1, \dots, \gamma_p$ forming the configuration of the type
$\cC$ and  having  holonomy  $hol(\gamma_i)=\vec{v}$.  We want to
compute the  number of collections  $\{\gamma_1,\dots,\gamma_p\}$
of the  type $\cC$ having holonomy vector of  length at most $L$.
In  other  words, we want to  compute  the asymptotics as $L  \to
\infty$ of the cardinality
$$
|V_{\cC}(S)\cap B(L)|
$$
of intersection of the discrete  set  $V_{\cC}(S)$  with the disc
$B(L)\subset\R{2}$ of radius $L$ centered at the origin.

Now consider the second problem of this paper: counting closed
saddle connections.

We have a surface $S$ and a saddle  connection $\gamma_1$ joining
a zero $z_1$  to itself. There  may be other  saddle  connections
$\gamma_2,\ldots,\gamma_m$  having  the  same holonomy vector  as
$\gamma_1$. By the same argument as above we  shall always assume
that all such $\gamma_i$ are homologous.

Some of the $\gamma_i$ may start  and end at the same zero $z_1$,
the others may start and end at the other zeroes. Each $\gamma_i$
returns  to  its zero at some  angle  $\theta_i$ which is an  odd
multiple of $\pi$. Some of the  $\gamma_i$  may  bound  cylinders
filled with regular closed geodesics. A {\em configuration} $\cC$
describes  all  these   geometric   data;  this  notion  will  be
formalized in Section~\ref{s:principal:boundary:II}.

We may  have  numerous  collections  of  precisely $p$ homologous
closed saddle connections corresponding to a configuration $\cC$,
i.e. defining  the  prescribed  number  of  cylinders, prescribed
angles $\theta_i$, etc. As in the first problem for a surface $S$
and configuration  $\cC$  let  $V_{\cC}(S)$  denote  the  vectors
$\vec{v}$ in $\reals^2$ such that the saddle connections in $\cC$
have  holonomy  $hol(\gamma_1)=\dots=hol(\gamma_p)=\vec{v}$ equal
to  that  vector.  Again  we  shall  compute  the  asymptotics as
$L\to\infty$ of
$$
|V_{\cC}(S)\cap B(L)|,
$$
where $B(L)$ is the disc of radius $L$.

\begin{Th}[A.~Eskin, H.~Masur; \cite{Eskin:Masur}]
For either of the problems, given a configuration $\cC$, there is
a constant $c=c(\alpha,\cC)$ such that for almost all  $S$ in any
connected component of any stratum $\cH_1(\alpha)$ one has
\begin{equation}
\label{eq:Siegel:Veech:constant}
\lim_{L\to\infty}\cfrac{|V_{\cC}(S)  \cap  B(L)   |}{\pi  L^2}  =
c(\alpha,\cC).
\end{equation}
The  constant  $c(\alpha,\cC)$  depends  only  on  the  connected
component of the stratum and on the configuration $\cC$.
\end{Th}

Versions of the above theorems where convergence a.e. is replaced
by convergence in $L^1$ are proved  in ~\cite{Veech:Siegel}. Thus
the  growth rate  is  quadratic for a  generic  surface. What  is
perhaps surprising is that this formula says that the growth rate
is quadratic for a  generic  surface even for multiple homologous
saddle connections. By  contrast, a generic surface does not have
{\em any} pairs of  saddle  connections that determine vectors in
the  same  direction  but  with  {\em   different}  lengths,  see
Proposition~\ref{pr:no:parallel:nonhomologous}  at   the  end  of
Section~\ref{s:volume:estimates}.

\begin{remark}
The issue of higher multiplicity and specifying the angles is not
explicitly  addressed  in~\cite{Eskin:Masur}.  However,  the  set
$V_{\cC}(S)$ defined above, satisfies  all  of the axioms of that
paper,  so  the  asymptotic  formula  from~\cite{Eskin:Masur}  is
applicable to the problems formulated above.
\end{remark}

One of  the main results of this  paper is  a description of  all
possible configurations $\cC$ for any connected  component of any
stratum,  and   an  evaluation  of  the  corresponding  constants
$c(\alpha,\cC)$.

We note that the lists  of  configurations and the values of  the
constants have  been  verified  {\em  numerically}. Beside direct
considerations  of  flat surfaces we used the following  implicit
experiments.  A  formula  based  on  the  formula  of  Kontsevich
(see~\cite{Kontsevich})  expresses  the  sum   of   the  Lyapunov
exponents  of  the Teichm\"uller geodesic flow on each  connected
component  of  every  stratum  $\cH(\alpha)$  in  terms  of  some
rational  function  of  $\alpha$  and  of  some  specific  linear
combination   of   the   constants   $c(\alpha,\cC)$   for    the
corresponding   component  of   $\cH(\alpha)$.   These   Lyapunov
exponents  have  been  computed  numerically  and  the  constants
achieved in this paper are  consistent  with  these exponents and
with Kontsevich's formula.

\subsection{Siegel---Veech Formula}
\label{ss:Siegel:Veech:Formula}
Fix a  topological type of a  configuration $\cC$. To  every flat
surface  $(S,\omega)$  in  $\cH_1(\alpha)$  we  can  associate  a
discrete  set  $V_{\cC}(S)\subset  \R{2}$.  This  enables  us  to
construct  an  operator  $f\mapsto  \hat  f$  from  the  space of
integrable functions with compact support in $\R{2}$ to functions
on $\cH_1(\alpha)$. We define $\hat f$ as
$$
\hat{f}(S) := \sum_{\vec{v} \in V_{\cC}(S)} f(\vec{v}).
$$
Averaging  $\hat  f$  over  the  entire  moduli space  we  get  a
functional  on  the space of functions with  compact  support  in
$\R{2}$.   In~\cite{Veech:Siegel}   it   was  proved  that   this
functional is $SL(2,\R{})$-invariant which implies the  following
Siegel-type  formula.  For each configuration $\cC$ and for  each
connected component  $\cH_1(\alpha)$  of each stratum the average
of $\hat f$  over  $\cH_1(\alpha)$ and  the  average of $f$  over
$\R{2}$ coincide up to a multiplicative factor $const(\cC)$ which
depends  only  on  the  configuration type and on  the  connected
component of the stratum. The formula
\begin{equation}
\label{eq:Siegel:Veech}
\frac{1}{\Vol(\cH_1(\alpha))}\int_{\cH_1 (\alpha)} \hat{f}(S) \,
d\vol(S) = const(\cC) \int_{\reals^2} f
\end{equation}
was   called  in~\cite{Eskin:Masur}   {\it   the   Siegel---Veech
formula}.  It   was   shown   in   ~\cite{Eskin:Masur}  that  the
Siegel-Veech constant $const(\cC)$  in  the formula above and the
$c(\cC)$  in  equation~\eqref{eq:Siegel:Veech:constant} coincide.
Thus, one of the major goals of this paper can  be  thought of as
finding      the      constant     in      the     Siegel---Veech
formula~\eqref{eq:Siegel:Veech}.

The  strategy  for  evaluating  $c(\cC)=const(\cC)$  consists  of
choosing a  convenient  function  $f$  for  which the complicated
integral on  the  left in~\eqref{eq:Siegel:Veech} can be computed
as explicitly as possible.  As  such a ``convenient function'' we
take  the   characteristic  function  $f_\epsilon$  of  the  disc
$B(\epsilon)$ centered at the origin of $\R{2}$, where $\epsilon$
is  the  radius  of  the  disc.  The  integral  on the  right  in
equation~\eqref{eq:Siegel:Veech} is equal to $\pi\epsilon^2$. The
integral on the left  gives  the average number of configurations
of the type $\cC$ having holonomy vector shorter than $\epsilon$.
We  cannot  compute   this   average  explicitly  for  any  given
$\epsilon$. However, we can compute the  leading  term  of  order
$\epsilon^2$  in   the   asymptotics  of  the  left-hand-side  of
(\ref{eq:Siegel:Veech})  as   $\epsilon  \to  0$.  This  will  be
sufficient for our purposes. To compute this leading  term we use
the  idea which  we  illustrate in the  simplest  case, when  the
configuration  $\cC$  consists  of  a  single  saddle  connection
joining two distinct saddles $z_1, z_2$.  The  general  case,  as
well  as  a  rigorous proof of  equation~\eqref{eq:thin:over:all}
below is treated in Section~\ref{s:volume:estimates}.

Denote     the     support    of     $\hat     f_\epsilon$     by
$\cH_1^\epsilon(\alpha,\cC)$. This is the subset of flat surfaces
$S$, which possess  at least one short saddle connection $\gamma$
of multiplicity one  joining the chosen zeroes $z_1, z_2$. Saying
``short''   we   mean   that   $|\vec{v}|\le   \epsilon$,   where
$\vec{v}=hol(\gamma)$. This space is not compact: for any compact
subset $K$ of $\cH_1(\alpha)$  the  length of the shortest saddle
connection on  any flat surface $S\in  K$ is bounded  from below.
The space  $\cH_1^\epsilon(\alpha,\cC)$  is  difficult to analyze
since a surface $S\in\cH_1^\epsilon(\alpha,\cC)$  may  have short
saddle connections  different  from  $\gamma$, possibly even many
saddle connections that intersect.

To simplify  our  task  we decompose $\cH_1^\epsilon(\alpha,\cC)$
into two complementary subsets: a ``thick'' and a ``thin'' part:
$
\cH_1^\epsilon(\alpha,\cC)=
\cH_1^{\epsilon, thick}(\alpha,\cC)\sqcup
\cH_1^{\epsilon,thin}(\alpha,\cC).
$
The  {\it   thick  part}  $\cH_1^{\epsilon,   thick}(\alpha,\cC)$
consists  of  surfaces  $S$  having  {\it   exactly  one}  saddle
connection $\gamma$ shorter than $\epsilon$; moreover, we require
that this  $\gamma$ joins the  two chosen zeroes $z_1$ and $z_2$.
The {\it thin part} $\cH_1^{\epsilon, thin}(\alpha,\cC)$ consists
of  surfaces  $S$  having  at least one short  saddle  connection
$\gamma$ as above and at least one other  short saddle connection
$\beta$ nonhomologous to  $\gamma$.  Here $\beta$ is an arbitrary
short saddle connection, not necessarily joining  $z_1$ to $z_2$,
possibly closed.

The  function  $\hat  f_\epsilon$  is  equal  to zero outside  of
$\cH_1^\epsilon(\alpha,\cC)$. By  definition  of  the thick part,
the value  of  $\hat  f_\epsilon$  on  any  $S\in\cH_1^{\epsilon,
thick}(\alpha,\cC)$    is    identically   $1$.    Finally,   for
$S\in\cH_1^{\epsilon,thin}(\alpha,\cC)$     we     have     $\hat
f_\epsilon(S)\ge 1$. Thus, we get the following representation of
the integral above:
\begin{equation*}
\int_{\cH_1 (\alpha)}
\hat{f_\epsilon}(S) \,d\vol(S) =
\int_{\cH^{\epsilon,thick}_1 (\alpha,\cC)}  \,d\vol(S) +
\int_{\cH^{\epsilon,thin}_1 (\alpha,\cC)}
\hat{f_\epsilon}(S) \,d\vol(S)
\end{equation*}
The first term in this sum  is just the volume of the thick part.
Even though $\hat f_\epsilon$ is unbounded on the  thin part, its
measure is so small  that the integral of $\hat f$ over  the thin
part is negligible: it is of the order  $o(\epsilon^2)$. We prove
this  statement   in  Section~\ref{s:volume:estimates}  using   a
nontrivial   upper    bound    for    $\hat   f_\epsilon$   found
in~\cite{Eskin:Masur}.  The  fact  that  the thin part  is  small
implies  also  that   $\Vol(\cH_1^{\epsilon}(\alpha,\cC))=   \Vol
(\cH_1^{\epsilon,thick}(\alpha,\cC))+o(\epsilon^2)$.  Hence,   we
get the following key statement:

\begin{proposition}
\label{prop:thick}
For any connected  component of any stratum $\cH(\alpha)$ and for
any configuration $\cC$  the following limit exists and is equals
to the corresponding Siegel---Veech constant:
\begin{equation}
\label{eq:thin:over:all}
c(\cC)=\lim_{\epsilon\to 0}\frac{1}{\pi\epsilon^2}\frac
{\Vol (\cH_1^{\epsilon}(\alpha,\cC))}
{\Vol(\cH_1(\alpha))}
\end{equation}
\end{proposition}

We  prove  this Proposition  in Section~\ref{s:volume:estimates}.
Formula~\eqref{eq:thin:over:all}   explains  why   one   of   the
principal  goals  of   this  paper  is  the  calculation  of  the
asymptotics of  the volume $\Vol  (\cH_1^{\epsilon}(\alpha,\cC))$
for  any  connected  component   of   any  stratum  and  for  any
configuration $\cC$ of homologous  saddle  connections admissible
for this component.

\subsection{Principal Boundary of the Moduli Spaces}
\label{ss:introduction:principal:boundary}
Now as  one lets $\epsilon\to 0$, the flat  surfaces in the thick
part degenerate to  simpler  surfaces. In  the  case of a  single
saddle connection joining a pair  of  distinct  zeroes the zeroes
coalesce to a higher order zero  giving a surface in a stratum in
the same genus; the complex  dimension  of  the resulting stratum
decreases by one  with respect to  the initial one.  In  general,
collapsing all saddle connections from a configuration $\cC$ to a
single point we obtain a  surface  from  more degenerate stratum;
the surface is disconnected if the configuration has multiplicity
two or  more. We shall  say that the resulting surface (connected
or not)  belongs to the  {\em principal boundary} of the original
stratum.  Each  configuration $\cC$ will determine a ``face''  of
this principal boundary. A major part of this paper then  will be
to describe this principal boundary,  an  object  of  independent
interest.

Conversely, given a  flat surface $S'$ in the principal boundary,
we will  describe a set  of surgeries  on $S'$ that  allow us  to
recover surfaces $S$ in the  thick  part. The first surgery is  a
splitting apart  of a higher  order zero into simpler zeroes that
is the opposite  of  coalescing.  The second one is  a  slit  and
gluing construction which  is the opposite of the degeneration of
multiple   saddle   connections.   To  construct  surfaces   with
homologous  closed  geodesics,  we  shall  also  need  a  pair of
constructions, called the figure eight and  creating  a  pair  of
holes construction.

The above constructions  lead to the following picture. Let $\cC$
be a configuration  of  homologous saddle connections joining two
distinct   points.   We   shall   prove  that  the   thick   part
$\cH^{\epsilon,thick}_1(\alpha,\cC)$  has  the   structure  of  a
(ramified) covering over a  direct  product $\cH_1(\alpha')\times
B(\epsilon)$ of  the  corresponding  stratum of surfaces obtained
after  degeneration   of   $\cC$   and   a  two-dimensional  disc
$B(\epsilon)$ of radius $\epsilon$. Moreover, the  measure on the
thick  part  $\cH^{\epsilon,thick}_1(\alpha,\cC)$   is  just  the
product  of   corresponding   measures   on   the   two   factors
$\cH_1(\alpha')$ and  $B(\epsilon)$.  The  degree of the covering
(the {\it combinatorial constant}) $M$ corresponds  to the number
of ways to perform the surgeries on $S'$ to obtain $S$. Hence, we
shall obtain the following answer in all cases
$$
\Vol(\cH^{\epsilon}_1(\alpha,\cC)) =
M\cdot \pi\epsilon^2 \cdot \Vol(\cH_1(\alpha')) + o(\epsilon^2).
$$
Applying~\eqref{eq:thin:over:all} we  then  are able to prove the
following result

\begin{equation}
\label{eq:c:answer}
c(\cC)=M\cdot \frac
{\Vol (\cH_1(\alpha'))}
{\Vol(\cH_1(\alpha))}.
\end{equation}

Analogous considerations can be applied to  a configuration $\cC$
of  closed  saddle connections. Since we are  interested  in  the
answers for all connected components of the strata  we shall need
in  certain  situations  to  pay  attention  to  parity  of  spin
structures and hyperellipticity.

One  can  think  of every stratum  as  an  analogue  of a complex
polyhedron;  the  boundary  of  the polyhedron has  ``faces''  of
complex codimension  $1$,  ``edges''  of complex codimension $2$,
etc. A collection of all admissible types $\cC$ of configurations
of homologous  saddle  connections  describes  all  {\it generic}
degenerations of our  surfaces, that is all possible ``faces'' of
the boundary. Thus, though we do not construct a compactification
of the strata  $\cH_1(\alpha)$ in full generality we describe the
{\it principal boundary} of the strata $\cH_1(\alpha)$

\begin{remark}
One  can consider  these  same problems in  the  context of  {\em
quadratic differentials}.  To  avoid  overloading  the  paper  we
prefer to treat them in a subsequent one.
\end{remark}

\begin{remark}
A list of  the geometric types of admissible configurations $\cC$
is shared  by almost all surfaces  in any connected  component of
any  stratum.  However,  this  list  varies  when  we  change the
stratum,  or  a connected component of the  stratum.  Thus,  such
lists can be used as an invariant of a connected component. There
are other orbifolds  in the moduli space of Abelian differentials
invariant under  the  $SL(2,\R{})$-action;  they  have  their own
lists of geometric types of admissible configurations (understood
more  generally).  These  lists  have  proved   to  be  important
invariants of such orbifolds:  in~\cite{Hubert:Schmidt}  they are
used  to  give  an  example  of  different  Veech  surfaces (with
different Teichm\"uller discs) sharing isomorphic Veech groups.
\end{remark}

\subsection{Arithmetic of Siegel---Veech Constants}
\label{ss:Arithmetics:of:Siegel:Veech:Constants}
As remarked earlier the number of (homotopy classes) of primitive
closed geodesics  of length  at most $L$ on a  flat torus of area
one (considered  up to a  symmetry of the torus) is approximately
$\cfrac{3}{\pi^2} \cdot \pi L^2$. The number of geodesic segments
of length  at most $L$ joining a generic  pair of distinct points
$P_1, P_2$ on the same torus is approximately $1\cdot \pi L^2$.

The general  Siegel---Veech  constants  $c(\cC)$  share  the same
arithmetic  property!   Those  Siegel---Veech  constants,   which
correspond to  configurations  of {\it closed} saddle connections
always have the  form  of a  rational  divided by $\pi^2$;  those
Siegel---Veech constants, which correspond  to  configurations of
saddle connections joining distinct points are always rational.

We have no simple geometric  explanation  of  this phenomenon. It
rather follows  directly from ~\eqref{eq:c:answer}, the fact that
the combinatorial factor $M$ is  by  its  meaning always rational
(mostly integer)  and the fact that  the volume of  any connected
component of any stratum $\cH_1(\alpha)$ of Abelian differentials
is a  rational multiple of $\pi^{2g}$, where $g$  is the genus of
the  surface.  This latter fact was conjectured by  M.~Kontsevich
and proved by A.~Eskin and A.~Okounkov, who, actually, calculated
these volumes  (see~\cite{Eskin:Okounkov} for connected  strata).
This  property  of the  volume  is true  not  only for  primitive
strata, but  for strata of surfaces  with marked points  and even
for        strata         of        disconnected        surfaces.
\section{Readers Guide}

In the first half of the paper we solve these problems for saddle
connections joining distinct zeroes. We describe the multiplicity
one  case  first, because it is  the  easiest one and because  it
illustrates    the    computations.    This     is     done    in
Sections~\ref{ss:collapsing}--\ref{ss:comuting:const}.  We   then
describe the principal boundary in the  higher multiplicity case,
describing the  {\it  slit}  and  {\it  gluing} constructions and
again    give     the     constants.     This    is    done    in
Sections~\ref{ss:principal:boundary}--\ref{ss:slit:construction}
and                                    \ref{ss:c:in:higher:mult}.
Sections~\ref{ss:parity:in:slit}--\ref{ss:other:strata}  as  well
as       the       corresponding        preliminaries        from
Section~\ref{ss:parity:of:spin}  are devoted  to  describing  the
constants  in  the cases of strata with  several  components  and
could be skipped  by the reader  interested only in  the  general
thrust of the paper.

In the  second half  of the paper we consider  the case of closed
saddle connections. We describe the principal  boundary using the
{\it  creating  a  pair  of   holes}   and   {\it  figure  eight}
constructions. The reader interested only in the structure of the
principal boundary may consult
Sections~\ref{s:principal:boundary:II}--\ref{s:constructions:zero:to:itself}.
The  computation  of  the  combinatorial factor is  performed  in
Section~\ref{ss:combinatorial:factor:II} and the  computation  of
the        constant         $c$        is        given         in
Section~\ref{ss:computation:of:c:II}.      The      rest       of
Section~\ref{s:zero:to:itself}  is  devoted to  numerous examples
and Section~\ref{s:not:connected:II}  is  devoted  to the special
case of strata with several connected components.

We are illustrating  all  constructions by numerous examples. All
possible strata and  all  possible configurations in genera $g=2$
and  $g=3$  and  all  related computations are presented  as  the
examples in the  main  body  of the text. The  complete  list  of
admissible configurations in genus  $g=4$,  and the values of the
corresponding  constants  are presented in the Appendix. We  also
present in  the Appendix the values of the  constants for the two
distinguished strata $\cH(1,\dots,1)$ and $\cH^{hyp}(g-1,g-1)$ up
to genus $g=30$.


\vspace*{1truecm}
\addcontentsline{toc}{part}
{Part 0. Structure of the Strata}
{\Large\bf Part 0. Structure of the Strata}
\vspace*{0.5truecm}

\section{Connected Components of the Strata}

\subsection{Parity of a Spin Structure}
\label{ss:parity:of:spin}
As noted earlier, the strata $\cH_1(\alpha)$  are not necessarily
connected. To classify the components one  needs  the  notion  of
parity of spin  structure. Consider a smooth simple curve $\beta$
on the flat  surface  $S$ which does not  contain  any zeroes. We
define  the  index  $ind(\beta)\in  \natls$ of the  vector  field
tangent to  $\beta$ to be the  degree of the  corresponding Gauss
map: the total change of the angle between the vector  tangent to
the curve  and the vector  tangent to the horizontal foliation is
$2\pi\cdot ind(\beta)$.

Now  take   a   symplectic   homology   basis   $\{  a_i,  b_i\},
i=1,\ldots,g$,  such   that   the  intersection  matrix  has  the
canonical form: $ a_i\circ a_j = b_i\circ b_j =0$, $ a_i\circ b_j
=  \delta_{ij}$,  $1\le i,j  \le  g$. Though  such  basis is  not
unique,  traditionally  it  is  called a {\it  canonical  basis}.
Consider  a  collection of smooth closed curves representing  the
chosen basis. Denote them by the same symbols $ a_i, b_i$.

When  all   zeroes   of   the   Abelian   differential   $\omega$
corresponding to the flat structure have only even degrees we can
define the {\it parity of the spin structure} $\phi(S)$ as
$$
\phi(S):=\sum_{i=1}^g
(ind( a_i)+1)(ind( b_i)+1)\quad (mod\ 2)
$$

It follows from the results of  D.Johnson~\cite{Johnson} that the
parity of the  spin structure $\phi(S)$  does not depend  on  the
choice of representatives, nor on the  choice  of  the  canonical
homology bases.

Moreover, it  follows  from the results of M.Atiyah~\cite{Atiyah}
that $\phi(S)$ is invariant under continuous  deformations of the
flat      structure      of     the      translation     surface,
(see~\cite{Kontsevich:Zorich} for details).

Since we do not need the notion of a spin structure itself, we do
not discuss it in this paper; the reader may consider  the notion
``parity of  the spin structure'' as a single  term. Note that it
is applicable only to those flat structures $\omega$ which have a
set of zeros $\alpha$ containing even numbers only.

\subsection{Classification of Connected Components of the Strata}
\label{ss:con:comp}
The      components      are      classified      as      follows
(see~\cite{Kontsevich:Zorich}).

Consider  the  general  case,  when  $g  > 3$,  and  $\alpha$  is
different  from  two  exceptional  cases:  $\alpha\neq   (2g-2)$,
$\alpha\neq(g-1,g-1)$.  Under   these  assumptions  the   stratum
$\cH_1(\alpha)$   is   connected  whenever   collection  $\alpha$
contains  an  odd  number.  When  $g  > 3$, $\alpha\neq  (2g-2)$,
$\alpha\neq(g-1,g-1)$, and $\alpha$ has only even elements, there
are  exactly  two connected components corresponding to even  and
odd    spin    structures.    We    shall    denote    them    by
$\cH_1^{even}(\alpha)$ and $\cH_1^{odd}(\alpha)$ correspondingly.

The strata $\cH_1(2g-2)$ and $\cH_1(g-1,g-1)$ have  in addition a
special component ---  the hyperelliptic one. It consists of flat
structures  on   hyperelliptic   surfaces.   In  both  cases  the
hyperelliptic involution $\tau$ sends the flat structure $\omega$
to  $-\omega$;  for $\cH_1(g-1,g-1)$ we also require that  $\tau$
interchanges two zeroes of $\omega$.

For  $g  >  3$  the  stratum  $\cH_1(2g-2)$  has  three connected
components: the  hyperelliptic  one  --- $\cH_1^{hyp}(2g-2)$, and
two  nonhyperelliptic  components  ---  $\cH_1^{even}(2g-2)$  and
$\cH_1^{odd}(2g-2)$   corresponding   to  even   and   odd   spin
structures. For even  $g> 3$ the stratum $\cH_1(g-1,g-1)$ has two
connected   components:   the   hyperelliptic   one,   and    the
nonhyperelliptic one  $\cH_1^{nonhyp}(g-1,g-1)$.  For  odd $g> 3$
the stratum  $\cH_1(g-1,g-1)$ has three connected components: the
hyperelliptic  one,  denoted  by $\cH_1^{hyp}(g-1,g-1)$, and  two
nonhyperelliptic    ones    ---    $\cH_1^{even}(g-1,g-1)$    and
$\cH_1^{odd}(g-1,g-1)$  corresponding  to  even   and   odd  spin
structures.

Genera $2$ and $3$ are special. In genus $2$, all flat structures
are hyperelliptic, so each of  the  $\cH_1(2)$  and  $\cH_1(1,1)$
have one component which is hyperelliptic.  In  genus  $3$,  only
hyperelliptic flat structures may have  even  parity  of the spin
structure.  Thus   $\cH_1(4)$  and  $\cH_1(2,2)$  each  have  two
connected components, the  hyperelliptic one and the one with odd
spin structure. The other strata in genus $3$ are connected.

The  flat  surfaces from hyperelliptic strata have the  following
parity of the  spin-structure (see~\cite{Kontsevich:Zorich}). The
parity of the spin structure determined by a  flat structure from
the hyperelliptic component $S\in \cH_1^{hyp}(2g-2)$ equals

\begin{equation}
\label{eq:spin:2g:minus2}
\phi(S) =
\left[\frac{g+1}{2}\right]\ (mod\ 2)
\end{equation}
where the  square brackets denote the  integer part of  a number.
The parity  of the spin  structure determined by a flat structure
$\omega$       from       the       hyperelliptic       component
$\omega\in\cH_1^{hyp}(g-1,g-1)$, for odd genera $g$ equals
\begin{equation}
\label{eq:spin:hyp}
\phi(S) =
\left(\frac{g+1}{2}\right)\ (mod\ 2) \quad \text{for odd $g$}
\end{equation}

\section{Nonprimitive Strata}

\subsection{Strata of Surfaces with Marked Points}

In this paper  we shall often consider the strata $\cH_1(\alpha)$
of surfaces  $S=(M,\omega)$ where we not  only fix the  zeroes of
the Abelian  differential $\omega$, but  we also mark one or more
regular  points  on the surface. Say, $\cH_1(3,1,0)$ will  denote
the surfaces of  genus $g=3$ endowed with an Abelian differential
with zeroes  of orders $3$ and  $1$ having one  additional marked
point (``zero of order $0$ of the Abelian differential'').

Throughout this  paper we assume that  all zeroes and  all marked
points  are always  numbered.  Thus a stratum  with  one or  more
marked points  has the natural structure  of a fiber  bundle over
the corresponding stratum without marked  points  with  a  direct
product  (minus  diagonals) of several copies of the  translation
surface as a fiber, where the number of copies equals  the number
of  marked  points.  For  example,  the   {\it  universal  curve}
$\cH_1(3,1,0)$ fibers over $\cH_1(3,1)$ with a fiber $S$.

In particular, there is  the  following formula for the dimension
of the  strata with  marked points. Let $g$ be  the genus and let
$card(\alpha)$ denote the number of entries  in $\alpha$. Suppose
that $0\notin\alpha$. Then
\begin{equation}
\label{eq:dim:marked:points}
\dim_\cx\cH(\alpha,\underbrace{0,\dots,0}_{n})=
\dim_\cx\cH(\alpha)+n=2g-1+card(\alpha)+n
\end{equation}
\begin{Convention}
By convention  we {\it  always} mark a point on  a flat torus. We
denote the corresponding stratum $\cH(\torusemptyset)$.

Note   that   $\dim_\cx\cH(\torusemptyset)=2$    which    matches
formula~\eqref{eq:dim:marked:points}.
\end{Convention}

The  natural  measure  on  the   stratum   with   marked   points
disintegrates into a product measure, where the measure along the
fiber  is  just the  Lebesgue  measure  on  $S$  (correspondingly
product of  several copies of $S$) induced by  the flat metric on
$S$, and  the measure  on the base is the  natural measure on the
corresponding stratum taken without marked points.

Recall that from the point  of  view of volumes we have  confined
ourselves to  the  subspaces  $\cH_1(\alpha)$  of  the strata for
which the area of every surface $S$ (measured in the flat metric)
is equal to one. The observation above implies  that {\it volumes
of the strata with marked points coincide with the volumes of the
corresponding strata without marked points}:
\begin{equation}
\label{eq:vol:with:marked:points}
\Vol(\cH_1(\alpha))=\Vol(\cH_1(\alpha,0))=\Vol(\cH_1(\alpha,0,0))=\dots
\end{equation}
   %
\subsection{Strata of Disconnected Surfaces}
\label{ss:Strata:of:Disconnected:Surfaces}
It   will   be  convenient  to  consider  sometimes  the   strata
$\cH(\alpha')=\prod_{i=1}^p\cH(\alpha'_i)$,   of   closed    flat
surfaces $S$  having  $p$  connected  components $S_1\sqcup \dots
\sqcup S_p$ of prescribed types.

\begin{Convention}
\label{conv:alpha:prime}
Using notation $\alpha'=\sqcup_{i=1}^p \alpha'_i$ for the  strata
$\cH(\alpha')$  of  disconnected surfaces we assume that we  keep
track  of   how   $\alpha'$   is   partitioned  into  collections
$\alpha'_i$. \end{Convention}

We shall need the expressions for the volume element and  for the
total volume of such strata.

We  write  $S_i  =  r_i  S_i'$,  where  $area(S_i')  =  1$.  Then
$area(S_i)  =  r_i^2$. Let  $d_i:=\dim_\reals\cH(\alpha'_i)$; let
$d:=\dim_\reals\cH(\alpha')=\sum_{i=1}^p d_i$; let  $d\nu'_i$  be
the volume element on the stratum  $\cH(\alpha'_i)$ (which should
not be confused with the volume element on a ``unit hyperboloid''
$\cH_1(\alpha'_i)$     in      the     same     stratum,      see
Convention~\ref{conv:vol:and:nu}). We have
$$
d\nu(S) =\prod_{i=1}^p  d\nu_i'(S_i)  =
\prod_{i=1}^p \left( r_i^{d_i-1} \, dr_i\right) \,
\prod_{i=1}^p d\vol'(S_i') \
$$
Let $D(1)$ be the unit ball $r_1^2+\dots+r_p^2\le 1$; set
$$
W=\prod_{i=1}^p \Vol(\cH_1(\alpha'_i)).
$$
Then,
$$
\nu(C(\cH_1(\alpha'))  =
W\cdot \int_{D(1)} \prod_{j=1}^p r_j^{d_j - 1}dr_j
$$
We now make the change of  variable $x_i = r_i^2$ to evaluate the
integral. For each $i$, let $b_i=d_i/2-1$, so that $r_i^{d_i-1}\,
dr_i=(1/2)x_i^{b_i} \, dx_i$. Then the above integral becomes
$$
\nu(C(\cH_1(\alpha'))  =
W\cdot
\frac{1}{2^p} \int_{ \sum_i x_i \le 1} x_1^{b_1} \dots x_p^{b_p}
\,
dx_1 \dots dx_p
$$
where  now  we integrate  over  the  standard  simplex.  Repeated
application of the identity
$$
\int_0^u  x^a  (  u  -  x)^b  \, dx  =  \frac{a!\,  b!}{(a+b+1)!}
u^{a+b+1}
$$
yields
$$
\nu(C(\cH_1(\alpha'))  =
W\cdot
\frac{1}{2^p}\frac{ b_1! \dots b_p!}{(b_1 + \dots + b_p + p )!}
$$
Since $b_1+\dots+b_p + p = \sum (d_i/2)=d/2$ we obtain
$$
\nu(C(\cH_1(\alpha') ) = \frac{W}{2^p}\cdot
\frac{(\frac{d_1}{2}-1)!\dots (\frac{d_p}{2}-1)!}
{(\frac{d}{2})!}
$$
Hence,
\begin{multline}
\label{eq:total:volume:of:nonprimitive:stratum}
\Vol(\cH_1(\alpha'))=d\cdot \nu(C(\cH_1(\alpha') )=
\frac{2\cdot\frac{d}{2}\cdot W}{2^p}\cdot
\frac{(\frac{d_1}{2}-1)!\dots (\frac{d_p}{2}-1)!}
{(\frac{d}{2})!} = \\
= \frac{1}{2^{p-1}}\cdot
\frac{(\frac{d_1}{2}-1)!\dots (\frac{d_p}{2}-1)!}
{(\frac{d}{2}-1)!}\cdot
\prod_{i=1}^p \Vol(\cH(\alpha'_i))
\end{multline}

Repeating  literarily   the   same   arguments   we   obtain  the
corresponding formula for the volume elements:
\begin{equation}
\label{eq:volume:element:of:nonprimitive:stratum}
d\vol=\frac{1}{2^{p-1}}\cdot
\frac{(\frac{d_1}{2}-1)!\dots (\frac{d_p}{2}-1)!}
{(\frac{d}{2}-1)!}\cdot d\vol'_1 \cdots d\vol'_p
\end{equation}

\section{Thick---Thin Decomposition, Volume Estimates,
and Computation of the Siegel---Veech Constants}
\label{s:volume:estimates}

In this section we prove Proposition~\ref{prop:thick} and justify
the key  Formula(~\ref{eq:thin:over:all}) for the  Siegel---Veech
constant. We also describe more  precisely  the  structure of the
thick---thin  decomposition   of   neighborhoods   of  the  cusps
$\cH^\epsilon_1(\alpha,\cC)$. Our estimates of the volumes of the
thick and the thin part are based on the following result:

\begin{lemma}
[H.~Masur, J.~Smillie]
\label{lemma:short:saddle:connections}
There is a constant $M$ such that for all $\epsilon,\kappa>0$ the
subset of  $\cH_1(\alpha)$  consisting  of  those  flat surfaces,
which have a saddle connection of length at  most $\epsilon$, has
volume  at most  $M\epsilon^2$.  The volume of  the  set of  flat
surfaces with a saddle connection  of  length  at most $\epsilon$
and  a  nonhomologous  saddle  connection  with  length  at  most
$\kappa$ is at most $M\epsilon^2\kappa^2$.
\end{lemma}
\begin{proof}
The   proof   is  contained  in  the  proof   of   Theorem   10.3
in~\cite{Masur:Smillie}.
\end{proof}

As we  have  seen  in Section~\ref{ss:Siegel:Veech:Formula} it is
convenient  to decompose  the  set  $\cH^\epsilon_1(\alpha,\cC)$,
which  plays  the  role   of   a  neighborhood  of  the  ``cusp''
corresponding to  the  configuration  $\cC$,  into  two  disjoint
subsets: the ``thick'' and the ``thin'' part,
$
\cH_1^\epsilon(\alpha,\cC)=
\cH_1^{\epsilon, thick}(\alpha,\cC)\sqcup
\cH_1^{\epsilon,thin}(\alpha,\cC).
$
Sometimes it will be  convenient  to vary slightly this partition
making the thin part  a bit larger or a bit smaller  depending on
the  consideration.  This  does  not  affect  the  sense  of  the
thick---thin  decomposition,  but  simplifies  the  proofs.  This
variations  can  be described as follows. We  use  the  parameter
$\epsilon$,  $0<\epsilon<1$  to bound the length of the  shortest
saddle      connection.     We      introduce      the      bound
$\kappa:=\lambda\cdot\epsilon^r$ for the length  of  the shortest
saddle  connection  nonhomologous to  the  first  one.  Here  the
parameters $\lambda$  and  $r$  satisfy the following conditions:
$\lambda\geq     1,      0<r\le     1$,     which      guaranties
$\kappa(\lambda,r)\geq\epsilon$ for all $\lambda$ and $r$.

The  subset  $\cH^\epsilon_1(\alpha,\cC)\subset\cH_1(\alpha)$  is
comprised of those surfaces which have at least one collection of
short (shorter  than $\epsilon$) homologous saddle connections of
the  type   $\cC$.   The   {\it   thin   part}  $\cH_1^{\epsilon,
thin}(\alpha,\cC)$ of this subset consists of surfaces $S$ having
at least one additional saddle connection  of  any  type  shorter
than $\kappa$.

The complement  to the thin part in $\cH^\epsilon_1(\alpha,\cC)$,
the  {\it  thick   part}  $\cH_1^{\epsilon,  thick}(\alpha,\cC)$,
consists of surfaces $S$  having  exactly one collection of short
homologous saddle  connections. This collection is necessarily of
the  type  $\cC$;  the  saddle connections from  this  collection
(which are  all of the  same length) are shorter then $\epsilon$;
any other saddle connection on $S$ is longer then $\kappa$.

Lemma~\ref{lemma:short:saddle:connections} implies the  following
immediate corollary.

\begin{corollary}
\label{cor:thick:almost:equals:all}
For any connected component of  any  stratum  $\cH(\alpha)$,  any
configuration $\cC$  and  any  choice  of  parameters $\lambda,r$
defining the thick---thin decomposition we have
$$
\Vol(\cH_1^{\epsilon}(\alpha,\cC))=      \Vol
(\cH_1^{\epsilon,thick}(\alpha,\cC))+o(\epsilon^2).
$$
\end{corollary}

To prove Formula(~\ref{eq:thin:over:all}) it remains to prove the
following Lemma.

\begin{lemma}
\label{lm:can:neglect:thin:part}
Let  $f$  be the  characteristic  function of  a  disc of  radius
$\epsilon$ centered at the  origin  of $\R{2}$. For any connected
component of  any  stratum $\cH(\alpha)$, any configuration $\cC$
and for an appropriate choice of  parameters $\lambda,r$ defining
the thick---thin decomposition the integral of the function $\hat
f_\cC$ over the thin part is negligible:
$$
\int_{\cH^{\epsilon,thin}_1 (\alpha,\cC)}
\hat{f_\cC}(S) \,d\vol(S)=o(\epsilon^2)
$$
\end{lemma}
\begin{proof}
Recall that the nonnegative function $\hat  f_\cC(S)$ counts only
those saddle collections which are arranged  in the configuration
$\cC$; it  can be defined  as the cardinality of the intersection
of the discrete set $V_{\cC}(S)$ with the disc $B(\epsilon)$.
$$
\hat f_\cC(S) := |V_{\cC}(S) \cap B(\epsilon)|
$$
Consider the analogous function
$$
\hat f(S) := |V_{sc}(S) \cap B(\epsilon)|
$$
which counts  all short saddle connections (without multiplicity)
regardless of  which  configuration  they  correspond to. Clearly
$\hat f \ge  \hat  f_\cC$ since $V_{\cC}(S)\subset V_{sc}$. Thus,
it is sufficient to prove the above Lemma for the  function $\hat
f$. We use the following estimate proved in~\cite{Eskin:Masur}.

\begin{Th}[A.~Eskin, H.~Masur]
Let $l(S)$ denote the length of the shortest saddle connection on
$S\in\cH_1(\alpha)$. For  any  connected component of any stratum
$\cH_1(\alpha)$  there  exist constants  $c'$  and  $0<\delta<1$,
depending    only    on   the   stratum,   so   that   for    any
$S\in\cH_1(\alpha)$ for  which  $l(S)$ is sufficiently small, the
following bound is valid:
\begin{equation}
\label{eq:number}
\hat f(S)\leq \cfrac{c'}{\left(l(S)\right)^{1+\delta}}
\end{equation}
\end{Th}

We      now       can       complete       the      proof      of
Lemma~\ref{lm:can:neglect:thin:part}  Choose  any $0<r\leq  1$ so
that
\begin{equation}
\label{eq:thick}
 2r>1+\delta ;
\end{equation}
and   let   $\lambda=1$.  This  choice  of  $\lambda,  r$   gives
$\kappa=1\cdot\epsilon^r$;     consider     the      thick---thin
decomposition corresponding  to  this  choice  of  parameters. To
prove the estimate
\begin{equation}
\label{eq:oepsilon}
\int_{\cH^{\epsilon,thin}_1 (\alpha,\cC)}
\hat{f}(S) \,d\vol(S)=o(\epsilon^2)
\end{equation}
we decompose  the set $\cH^{\epsilon,thin}_1 (\alpha,\cC)$ into a
disjoint union of subsets $U_n$,  such  that  the shortest saddle
connection  for  $S$  in  $U_n$  has   length  $l(S)$  satisfying
$\epsilon/2^{n+1} <  l(S)  \le  \epsilon/2^{n}$,  where  $n$ is a
non-negative integer. Since by definition,  on  each  surface  in
$U_n$   there   is  a  saddle  connection  with  length   between
$\epsilon/2^{n+1}$  and  $\epsilon/2^{n}$   and  a  nonhomologous
saddle connection  with  length  at  most $\kappa=\epsilon^r$, by
Lemma~\ref{lemma:short:saddle:connections}  there  is a  constant
$M$  so   that  the  measure  of  $U_n$  is   at  most  $M  \cdot
2^{-2n}\epsilon^{2+2r}$.

Together  with  (\ref{eq:number}) this implies that for some  new
constant $M'$, the integral of $\hat f$ over $U_n$ is bounded by
$$ 
M'2^{(\delta-1)n}\epsilon^{1+2r-\delta}
$$ 
Summing over $n$, and using~\eqref{eq:thick}  we  find  that  the
estimate~\eqref{eq:oepsilon}                               holds.
Lemma~\ref{lm:can:neglect:thin:part} is proved. \end{proof}
Proposition~\ref{prop:thick}        now        follows       from
Lemma~\ref{lm:can:neglect:thin:part}                          and
Corollary~\ref{cor:thick:almost:equals:all}.

$$
\ast\quad\ast\quad\ast
$$

We complete  this section with  a proof of the statement promised
in the introduction.

\begin{proposition}
\label{pr:no:parallel:nonhomologous}
For  almost  all  flat  surfaces $S$ from any  stratum  $\mathcal
H(\alpha)$ of Abelian differentials one cannot find on $S$ a pair
of parallel saddle connections of different lengths.
\end{proposition}
\begin{proof}
As coordinates  in the stratum me may locally  choose a domain in
the  relative   cohomology   space   $H^1(S,\{P_1,  \dots,  P_k\}
;\cx{})$. Let $c_1,c_2,...,c_n$  be a basis of relative cycles in
$H^1(S,\{P_1, \dots, P_k\} ;\Z)$ The relative periods
$$A_i+\sqrt{-1}\cdot B_i = \int_{c_i} \omega, \qquad i=1,...,n$$
serve as the local coordinates in $\cH(\alpha)$.

Suppose that we have two parallel saddle connections of different
lengths on  a flat surface $S$. They  give us  a pair of  integer
(relative) cycles

$$s_1, s_2 \in H^1(S,\{P_1, \dots, P_k\} ;\Z).$$
(We don't care,  whether  they are both loops,  or  they both are
segments, or  one of  them is a segment, and  another --- a loop,
the argument works for  any  combination.) By assumption they are
not homologous,  so we may assume  that $s_1\neq \pm  s_2$. Since
the cycles $s_1,s_2$ are represented by  simple connected curves,
they are  primitive (i.e., neither  of them can be represented as
an integer multiple of another  integer  cycle).  Hence, they are
not collinear
$$s_1 \neq \lambda \cdot s_2$$
even with a real $\lambda\in \R{}$.

Let $  s_1 =  \sum_{i=1}^n k_i\cdot c_i $, where  $k_i \in \Z$; $
s_2 = \sum_{i=1}^n l_i\cdot c_i $, where $l_i \in \Z$.  Since the
cycles are not collinear, these  two  linear  combinations of the
basic cycles are linearly independent  over  reals.  This  means,
that the  rational  function $f(x_1,\dots,x_n)$ of real variables
$x_1, \dots, x_n$

$$ f(x_1,\dots,x_n)=
\cfrac{\sum_{i=1}^n   k_i\cdot x_i}{ \sum_{i=1}^n   l_i\cdot x_i}
$$
is nonconstant.

The fact, that the two  saddle  connections  are parallel, means,
that  $$\int_{s_1}\omega=  \lambda \cdot  \int_{s_2}\omega, \quad
\text{ with some real }\lambda\neq 0.$$

We have
$$
\int_{s_1}\omega=  \sum_{i=1}^n   k_i\cdot(A_i+\sqrt{-1}\cdot B_i) \qquad
\int_{s_2}\omega=  \sum_{i=1}^n   l_i\cdot(A_i+\sqrt{-1}\cdot B_i)
$$
This implies that $$f(A_1,\dots,A_n)=f(B_1,\dots,B_n),$$ which is
an algebraic condition  on our coordinates $A_1, \dots, A_n, B_1,
\dots, B_n$. (To be absolutely rigorous we have to avoid  the set
of measure zero defined by the algebraic condition $\sum l_i\cdot
A_i=0$ or by $\sum l_i\cdot  B_i=0$.)  Thus,  the set, satisfying
this  condition,  has  measure  zero.  Taking  a  union  over the
countable collection of possible  conditions  (countable, because
we  have  to  consider  all  possible  pairs  of  integer vectors
$(k_1,...,k_n)$, $(l_1,\dots,l_n)$) we still get a set of measure
zero. \end{proof}


\vspace*{1truecm}
\addcontentsline{toc}{part}
{Part 1. Saddle Connections Joining Distinct Zeroes}
{\Large\bf Part 1. Saddle Connections Joining Distinct Zeroes}
\vspace*{0.5truecm}

In this part we describe  the  possible  configurations of saddle
connections joining distinct zeroes and compute  the constants in
the corresponding asymptotics.

To separate the  basic construction and numerous details we start
with the easiest  case.  In  the first section of  this  part  we
assume that the saddle connection joining a pair of zeroes $P_i$,
$P_j$ has {\it multiplicity one}, i.e., there are no other saddle
connections  in  the  same  direction  joining  the same pair  of
zeroes.  In  this  section  we also assume that  the  translation
surface under consideration belongs to a connected stratum.

In the second  section  of this part we  consider  the problem in
full generality  for  the  translation  surfaces  from  connected
strata. In particular we give  the  explicit  general formula for
the surfaces from the principal stratum.

In the last section  of this part we treat the surfaces  from the
strata which are not connected.

\section{Saddle  Connections   of  Multiplicity  One.   Connected
Strata}


\subsection{Breaking up a Zero}
\label{ss:breaking:zero}

Suppose we are given a flat surface $S'$ defined by  $\omega'$, a
zero $w$ of order $m\geq  2$,  a pair of positive integers  $m_1,
m_2$, such that $m=m_1+m_2$, and a  vector $\gamma\in\reals^2$ of
length $2\delta\leq \epsilon$. Further suppose that $S'$ does not
have any saddle connection  or  closed geodesic of length smaller
than  $2\epsilon$.  Let $w, z_1, \dots,  z_l$  be the set of  all
zeroes of the flat structure on  $S'$. If $w$ is the only zero of
$\omega'$  choose  a basis of cycles for  the  relative  homology
group  $H_1(S',\{w\};\Z)$  all of which miss $w$.  If  there  are
other zeroes, we  may  choose a basis of  cycles  in the relative
homology group $H_1(S',\{w,z_1,\dots,z_l\};\Z)$ such that exactly
one curve $\beta_1$ contains $w$ and  $\beta_1$  is  not  closed;
that is, $w$ is an endpoint of $\beta_1$.

\begin{figure}[ht]
%
\includegraphics{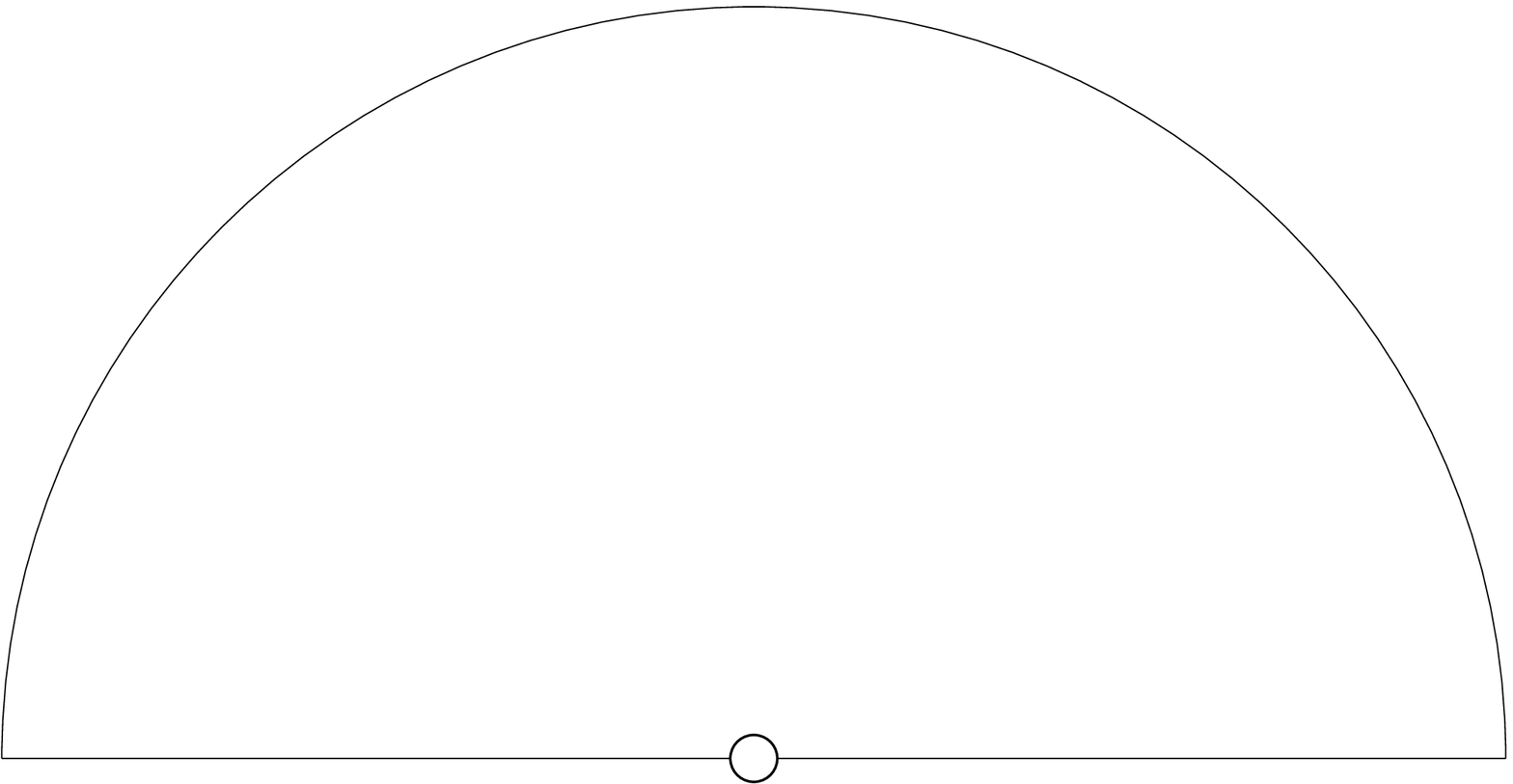}
\includegraphics{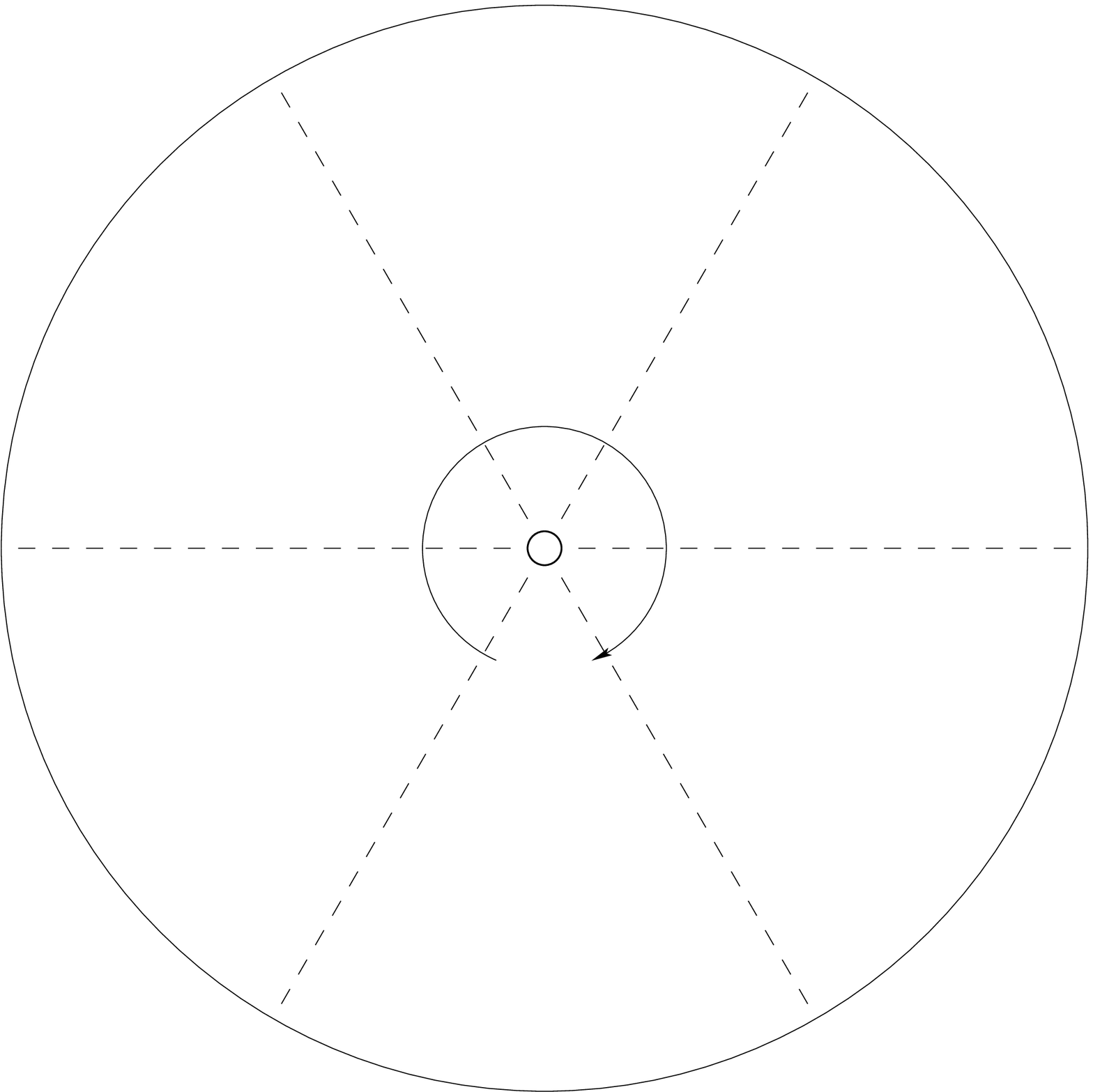}
%
%
%
\includegraphics{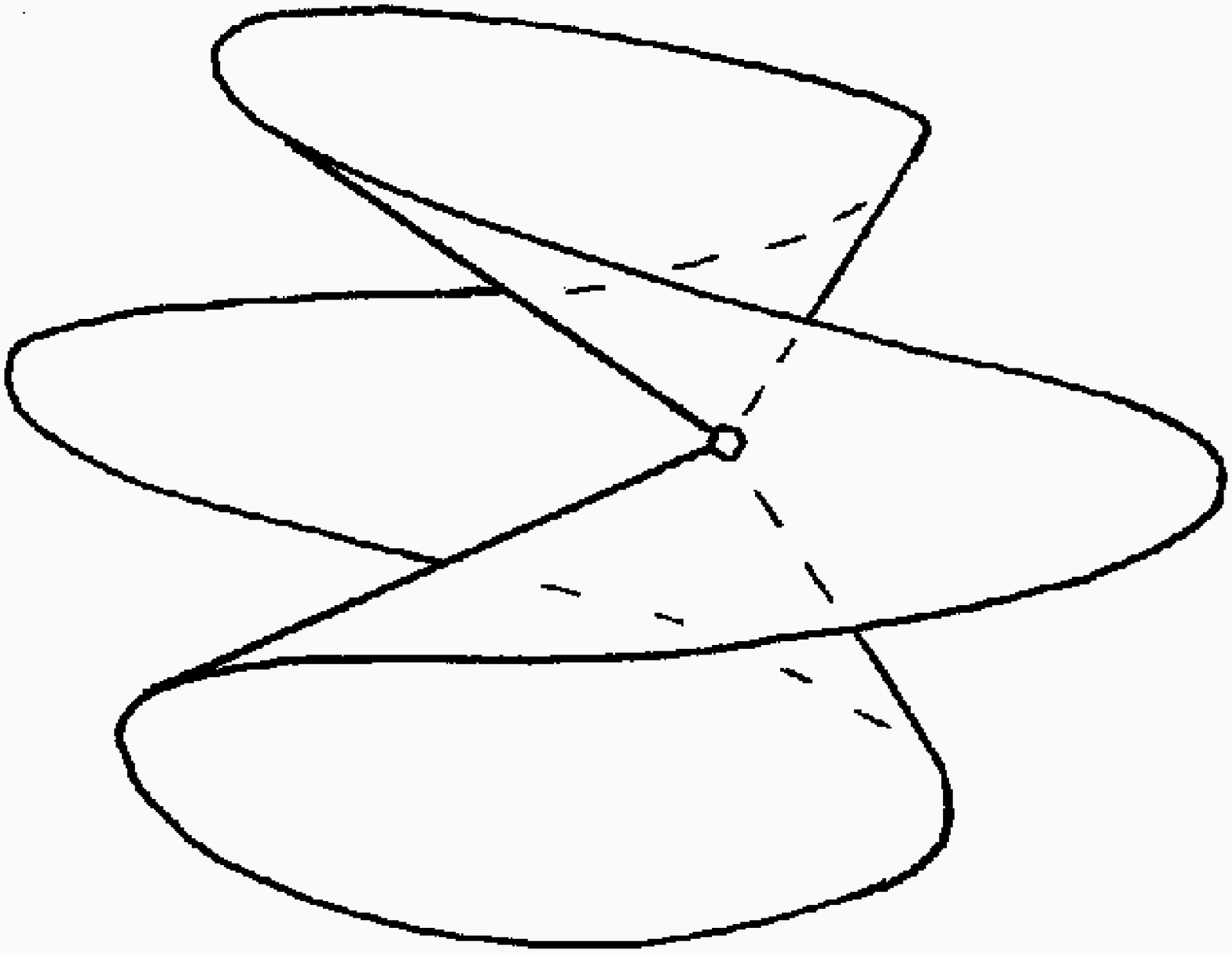}
%
%
\includegraphics{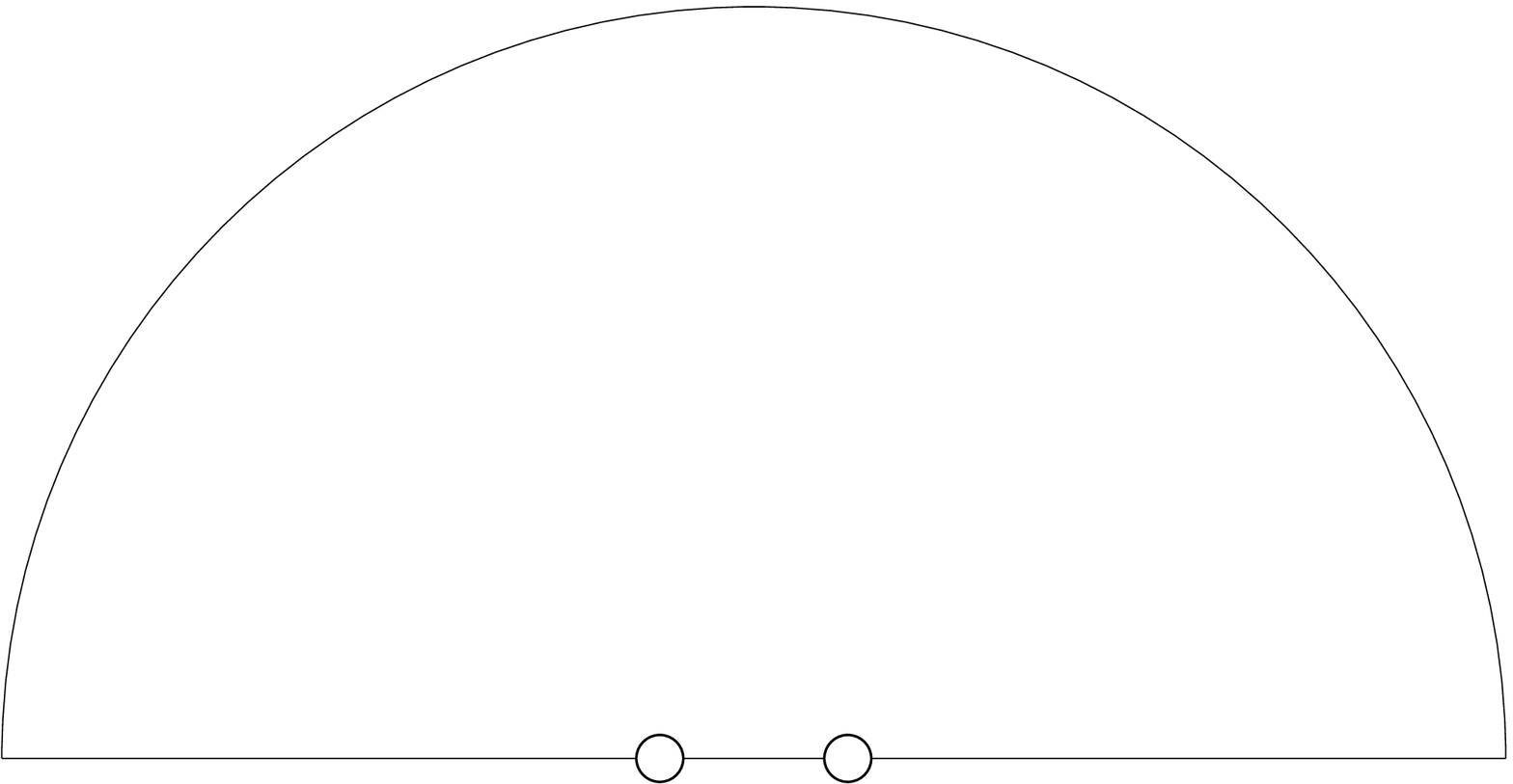}
\includegraphics{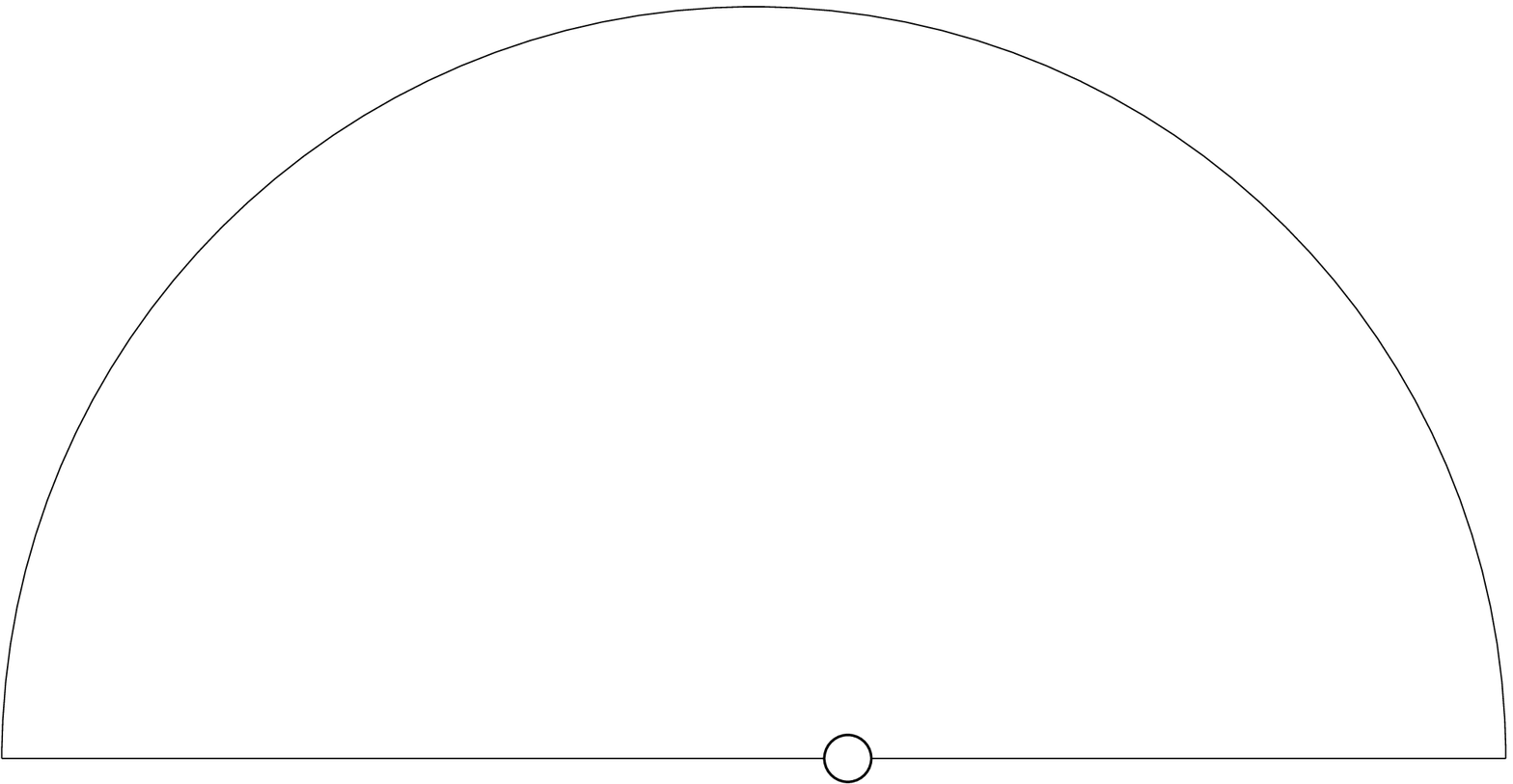}
\includegraphics{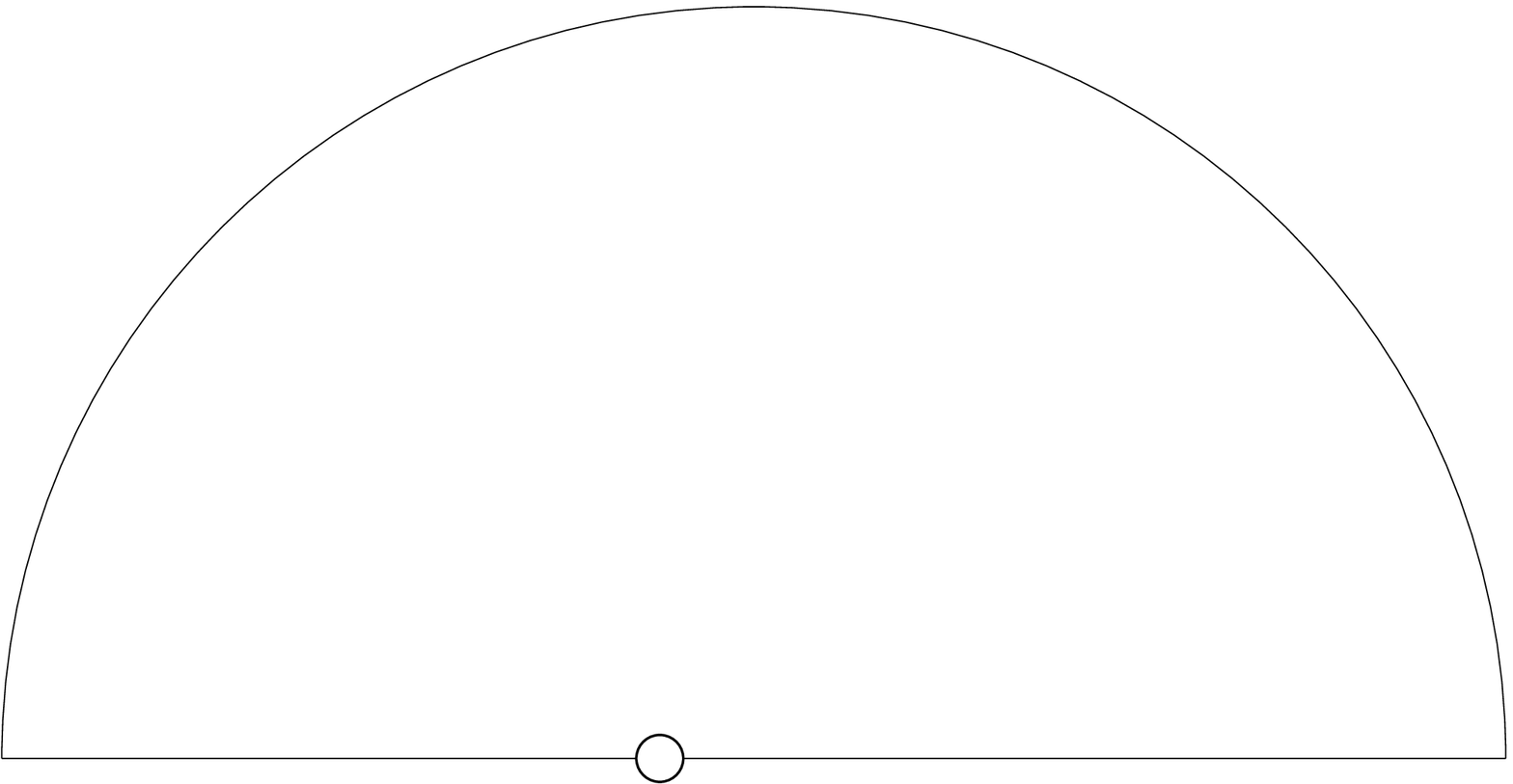}
\includegraphics{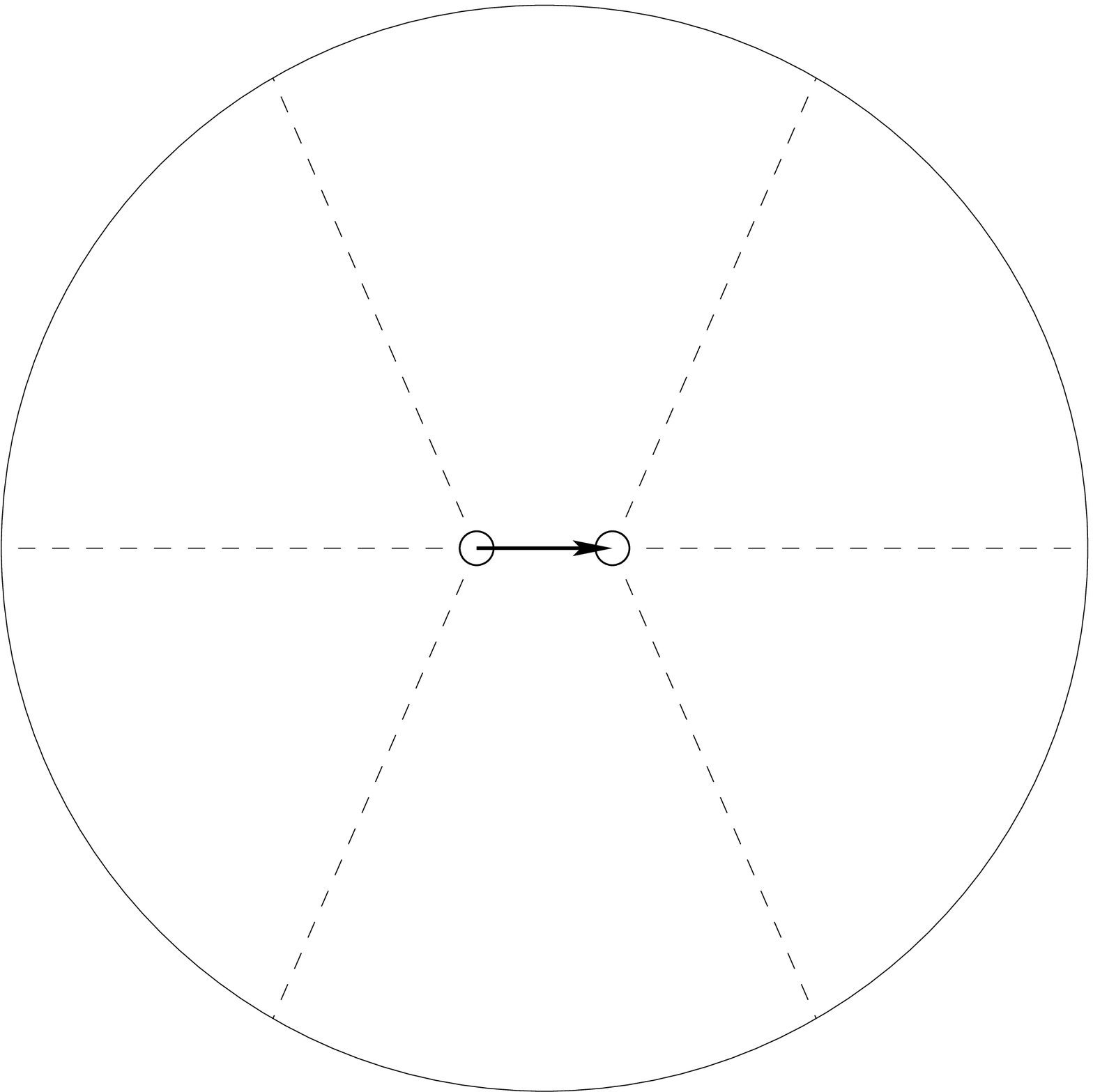}
%
%
\includegraphics{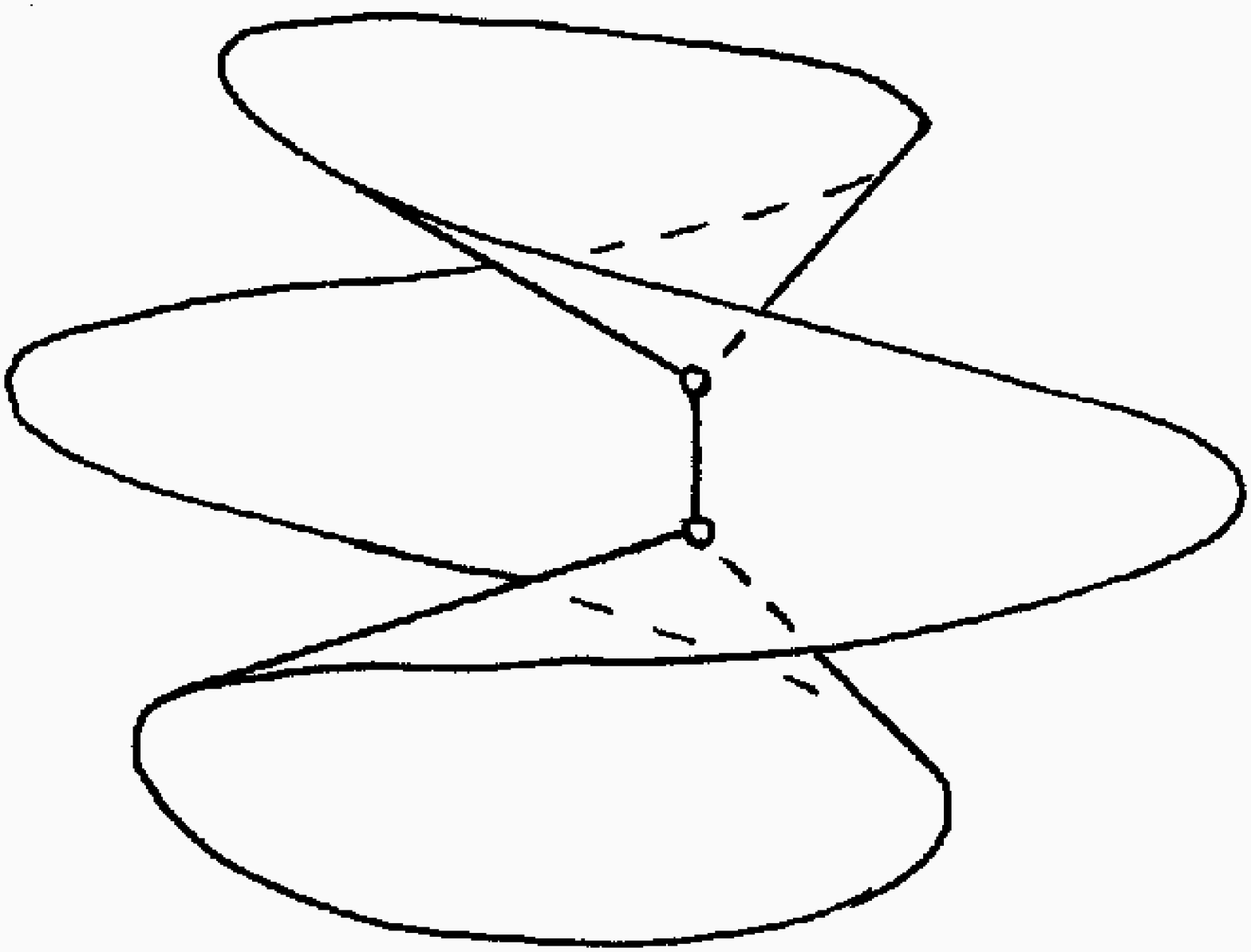}
%
%
\begin{picture}(0,0)(0,0)
\put(10,-10)
 {\begin{picture}(0,0)(0,0)
 \put(-154,-55){$\scriptstyle \varepsilon$}
 \put(-125,-55){$\scriptstyle \varepsilon$}
 \put(-38,-66){$\scriptstyle 6\pi$}
 \end{picture}}
\put(10,-30)
 {\begin{picture}(0,0)(0,0)
 \put(-145,-148){$\scriptstyle 2\delta$}
 \put(-160,-197){$\scriptstyle \varepsilon+\delta$}
 \put(-130,-197){$\scriptstyle \varepsilon-\delta$}
 \put(-160,-247){$\scriptstyle \varepsilon-\delta$}
 \put(-130,-247){$\scriptstyle \varepsilon+\delta$}
 \put(-39,-197){$\scriptstyle 2\delta$}
 \put(-65,-197){$\scriptstyle \varepsilon+\delta$}
 \put(-19,-197){$\scriptstyle \varepsilon+\delta$}
 \put(-64,-173){$\scriptstyle \varepsilon-\delta$}
 \put(-19,-173){$\scriptstyle \varepsilon-\delta$}
 \put(-65,-217){$\scriptstyle \varepsilon-\delta$}
 \put(-17,-217){$\scriptstyle \varepsilon-\delta$}
\end{picture}}
\end{picture}
\vspace{280bp} 
\caption{
\label{pic:degsad}
Breaking up a  zero  of degree $2$ into  two  simple zeroes. Note
that the surgery  is  local:  we do not change  the  flat  metric
outside of the neighborhood of the zero.
}
\end{figure}

We take a  disc  of radius $\epsilon$ about  $w$  that misses all
other zeroes. We may break up the zero $w$ of  order  $m$ on $S'$
into   two   zeroes  $z',z''$   of   orders   $m_1$   and   $m_2$
correspondingly with  a vector $\gamma$ joining them constructing
a flat  surface $S$. We can describe this  breakup as a Whitehead
move  on  the foliation  in  direction $\gamma$.  We  do this  by
forming $2m+2$ half discs of  radius  $\epsilon$.  Along the real
axis of two of them we  mark points at distance $\delta$ from the
origin.  These   two   discs   are   glued   together  along  the
corresponding segment of  length $2\delta$ leaving a pair of free
segments  of  length $\epsilon-\delta$ on each. On  each  of  the
remaining discs we  mark a point  at distance $\delta$  from  the
origin   leaving   segments   of  length  $\epsilon+\delta$   and
$\epsilon-\delta$. We now glue the segments isometrically to each
other in a circular fashion, see Figure~\ref{pic:degsad}.

If $m_1\neq m_2$, the number of ways of effecting the  breakup is
$2m+2$.
However, we  have $\overrightarrow{z'z''}=\gamma$ for only half of
the resulting surfaces, i.e. for  $m+1$  ones;  for another $m+1$
surfaces    we     have    $\overrightarrow{z''z'}=\gamma$.    If
$m_1=m_2=m/2$  the  number  of  the resulting surfaces  is
$m+1$.  For  every such  surface  there are  two  ways to  assign
``names'' $z',z''$ to the newborn zeroes  of  order  $m/2$.  This
doubles the number  of  resulting surfaces with ``named'' zeroes.
However,   we    again    have    only    $m+1$    surfaces   with
$\overrightarrow{z'z''}=\gamma$; for another $m+1$  ones  we have
$\overrightarrow{z''z'}=\gamma$. We again choose only those $m+1$
surfaces for which we have $\overrightarrow{z'z''}=\gamma$.

By convention we let the curve which had the endpoint at  $w$  (if it was
present)  keep  the  corresponding  endpoint at $z'$  during  the
deformation. The fact that this construction was local means that
except for the  curve $\beta_1$ (if  it exists) the  holonomy  is
preserved  along  the  homology  basis of $S'$. The  holonomy  of
$\beta_1$ is  changed  by  $-\gamma/2$.  Furthermore every saddle
connection  other  than $\overrightarrow{z'z''}$  has  length  at
least $\epsilon$.

%

We denote the assignment by
$$ 
(S',\gamma,m)\to (S,m_1,m_2)
$$ 

We  denote  by  $\cH_1(\alpha)$  the stratum which  contains  the
resulting  flat   surface  $S$.  By  construction  the  partition
$\alpha$  is  obtained from the partition $\alpha'$ by  replacing
the entry $m$ by two entries $m_1$ and $m_2$.

Recall that all the zeroes on the surface $S'$ are ``named''. For
all the zeroes  on  the new surface $S$  different  from the zero
which was just broken we keep  the same ``names'' on $S$ as their
initial ``names'' on $S'$.

\subsection{Collapsing a Pair of Zeroes, Principal Boundary}
\label{ss:collapsing}

Now conversely, suppose we have a surface $S\in \cH_1(\alpha)$, a
saddle   connection   of  length   $2\delta\leq\epsilon$  joining
distinct zeroes $w_1,w_2$ of orders $m_1$ and $m_2$ with holonomy
$\gamma$, and no other  saddle  connection of length smaller than
$3\epsilon$. Let $z_1,\ldots,z_l$  be  the other zeroes. Choose a
basis   of    cycles    in    the    relative    homology   group
$H_1(S,\{w_1,w_2,z_1,\dots,z_l\};\Z)$.   This   basis   can   be
represented by a collection of  curves  $\beta_i$  on the surface
$S$. One  of these curves,  $\beta_0$ is the saddle connection joining $w_1$
and $w_2$. If these are  the  only zeroes, then we can choose a basis of cycles in such way
that all other  cycles
are closed  and miss $w_1$ and $w_2$. If  there are other zeroes,
then a  single curve $\beta_1$ intersects  one of the  $w_i$, say
$w_1$. The curve $\beta_1$  is not  closed; it has $w_1$ as  an endpoint.

One can now exactly reverse the breaking up procedure to collapse
the saddle connection of length $2\delta$ to a zero $w$  of order
$m=m_1+m_2$  to  construct a flat surface $S'$.  Namely,  we  can
describe a  neighborhood of the saddle  connection as a  union of
$2m+2$ half discs of radius $\epsilon$ glued along  pieces of their
boundary. Two of  the half discs are  glued along segments  of  length
$2\delta$ to  form the saddle connection. This leaves a pair of
segments of length $\epsilon-\delta$ on  each side.  The  other  discs
have   marked   segments   of   length   $\epsilon-\delta$    and
$\epsilon+\delta$ which  are  glued  isometrically  to  form  the
neighborhood. We  deform $S$ in  the neighborhood by taking $2m+2$
half discs of radius  $\epsilon$  and gluing them cyclically along
segments of length $\epsilon$. The saddle connection $\beta_0$ is
collapsed to  a point  which becomes a zero $w$  of order $m$. We
may perform this deformation while keeping the flat structure in
the complement of the neighborhood fixed, see section 8.1 and Figures 3 and 4.

If $\beta_1$  exists then we can  replace $\beta_1$ with  a curve
joined to $w$. The fact that the deformation is local  means that
the holonomy is fixed along the entire basis other than $\beta_0$
and $\beta_1$. The holonomy of $\beta_1$  is  changed  by  adding
$\gamma/2$. Furthermore,  every  saddle  connection  on  $S'$ has
length at least $2\epsilon$. All named zeroes on $S$ not affected
by the collapse are given the same names on $S'$. We may think of
this surgery as of a  {\em  Whitehead} move  on  the foliation in  direction
$\gamma$.  The  resulting  surface  $S'$  belongs  to  a  stratum
$\cH_1(\alpha')$.  By  construction the  partition  $\alpha'$  is
obtained from the partition $\alpha$ by  replacing  the  pair  of
entries $m_1,m_2$ by the entry $m$.  The stratum $\cH_1(\alpha')$
is in  the closure of  $\cH_1(\alpha)$ inside the moduli space of
all  flat  structures  on  the  surfaces  of genus  $g$.  We  say that
$\cH_1(\alpha')$   is   the   {\em    principal}    boundary   of
$\cH_1(\alpha)$ corresponding to this configuration.

Now choose a  simply connected subset of $\cH_1(\alpha')$ of full
volume  and  remove  the  set  of  flat surfaces  with  a  saddle
connection or closed geodesic of length at most $2\epsilon$. Call
the    resulting    set    $\cF'\subset    \cH_1(\alpha')$.    By
Lemma~\ref{lemma:short:saddle:connections}
$$ 
\Vol(\cH_1(\alpha')\setminus \cF')=O(\epsilon^2)
$$ 
Choose a homology basis that is valid for all $S'\in  \cF'$. Thus
distinct  $S_1',S_2'$  have  different  holonomy  on  some  basis
element.    As    in     Section~\ref{s:volume:estimates}     for
$\kappa=\epsilon$           or           $3\epsilon$          let
$\cH_1^{\epsilon,\kappa}(\alpha,\cC)\subset\cH_1(\alpha)$ be the  set
of  flat surfaces  with  a saddle connection  of  length at  most
$\epsilon$  joining  the   named   zeroes  and  no  other  saddle
connections of length smaller than $\kappa$; the configuration $\cC$
corresponds to a single saddle connection. We have

\begin{lemma}
\label{lm:breaking:zero}
For $\gamma$ a vector in $\reals^2$ of length at most $\epsilon$,
except  for  a  set  of  $S'\in\cF'$  of volume  $0$,  there  are
precisely $m+1$ surfaces  in  $\cH_1^{\epsilon,\epsilon}(\alpha,\cC)$
that  are   the   result   of  the  assignment  $(S',\gamma,m)\to
(S,m_1,m_2)$.       Moreover,       every        surface        in
$\cH_1^{\epsilon,3\epsilon}(\alpha,\cC)$
is the result of such an assignment.

\end{lemma}

\begin{figure}[hb]
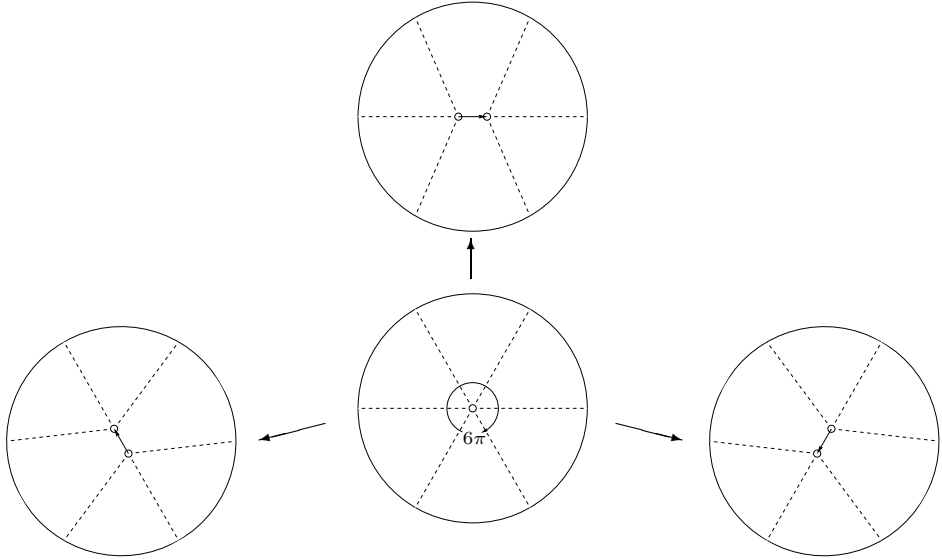

%
%
\includegraphics{b0zero.eps}
\includegraphics{bzero.eps}
%
\includegraphics{bzero.eps}
%
\includegraphics{bzero.eps}
%
%
\begin{picture}(0,0)(0,0)
\put(0,20) 
 {\begin{picture}(0,0)(0,0)
   \put(2,-186){$\scriptstyle 6\pi$}
 \put(60,-178){\vector(4,-1){25}}
 \put(-50,-178){\vector(-4,-1){25}}
 \put(5,-123){\vector(0,1){15}} 
 \end{picture}}
\end{picture}
%
%
\vspace{205bp} 
\caption{
\label{pic:bzero}
The cone angle corresponding to a zero of order $m=2$ is equal to
$(m+1)\cdot 2\pi=6\pi$. Thus we have $(m+1)=3$  different ways of
breaking up a zero of order $2$ in a direction $\vec{\gamma}$. In
this way (generically) we get $m=3$ different flat surfaces
}
\end{figure}

\begin{figure}[htb!]
%
 %
 %
 %
 \includegraphics{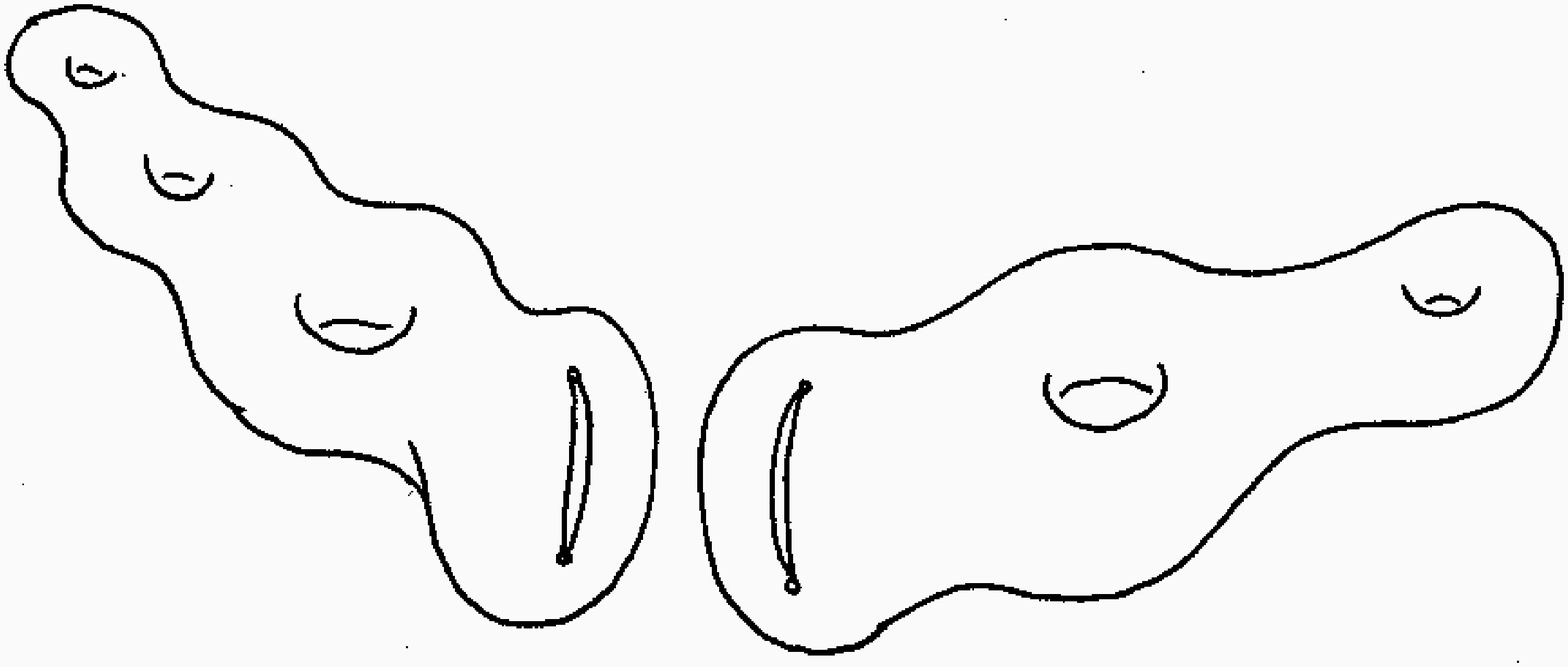}
 %
 %
 %
 \includegraphics{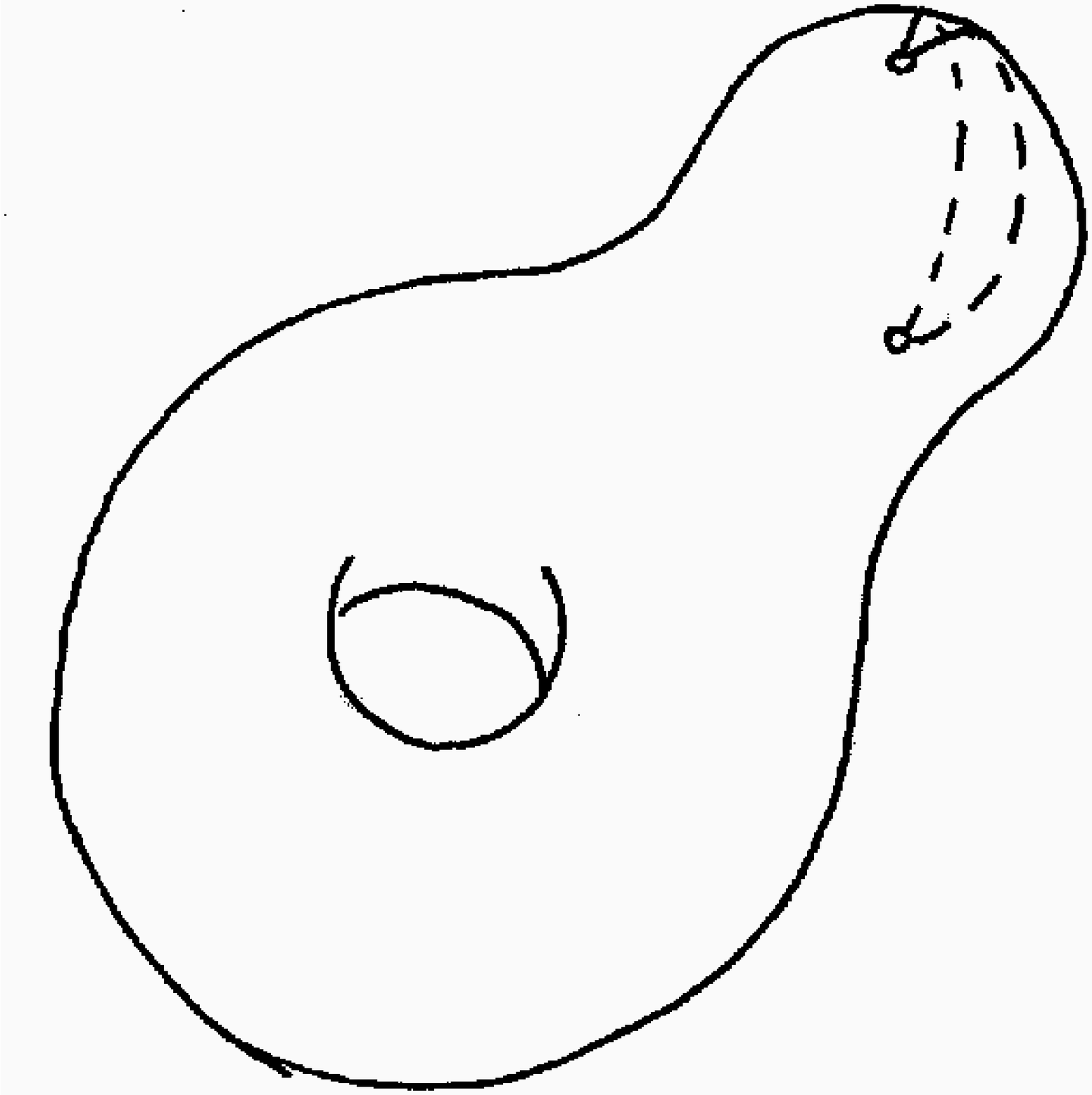}
 %
 %
 %
 \includegraphics{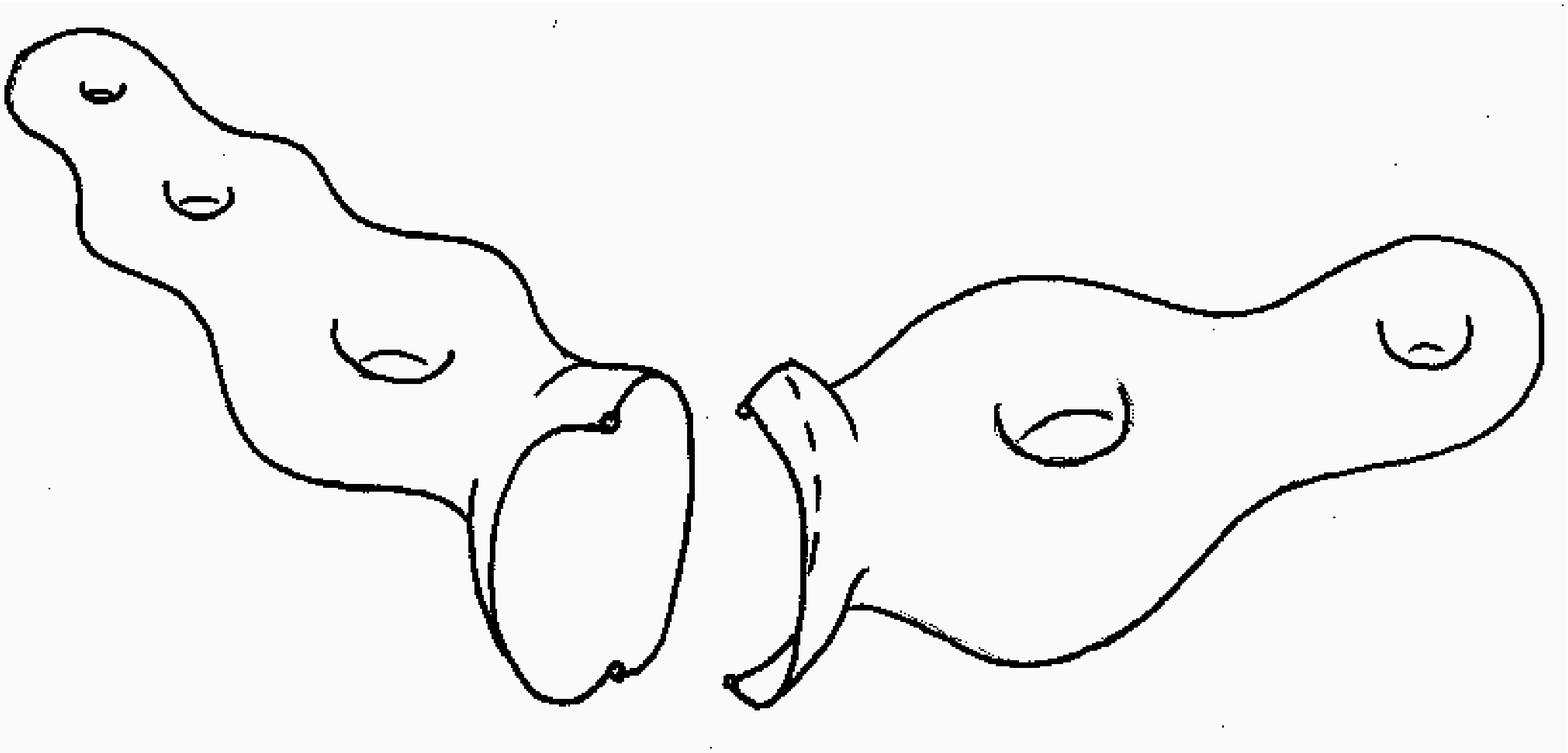}
 %
 %
 \includegraphics{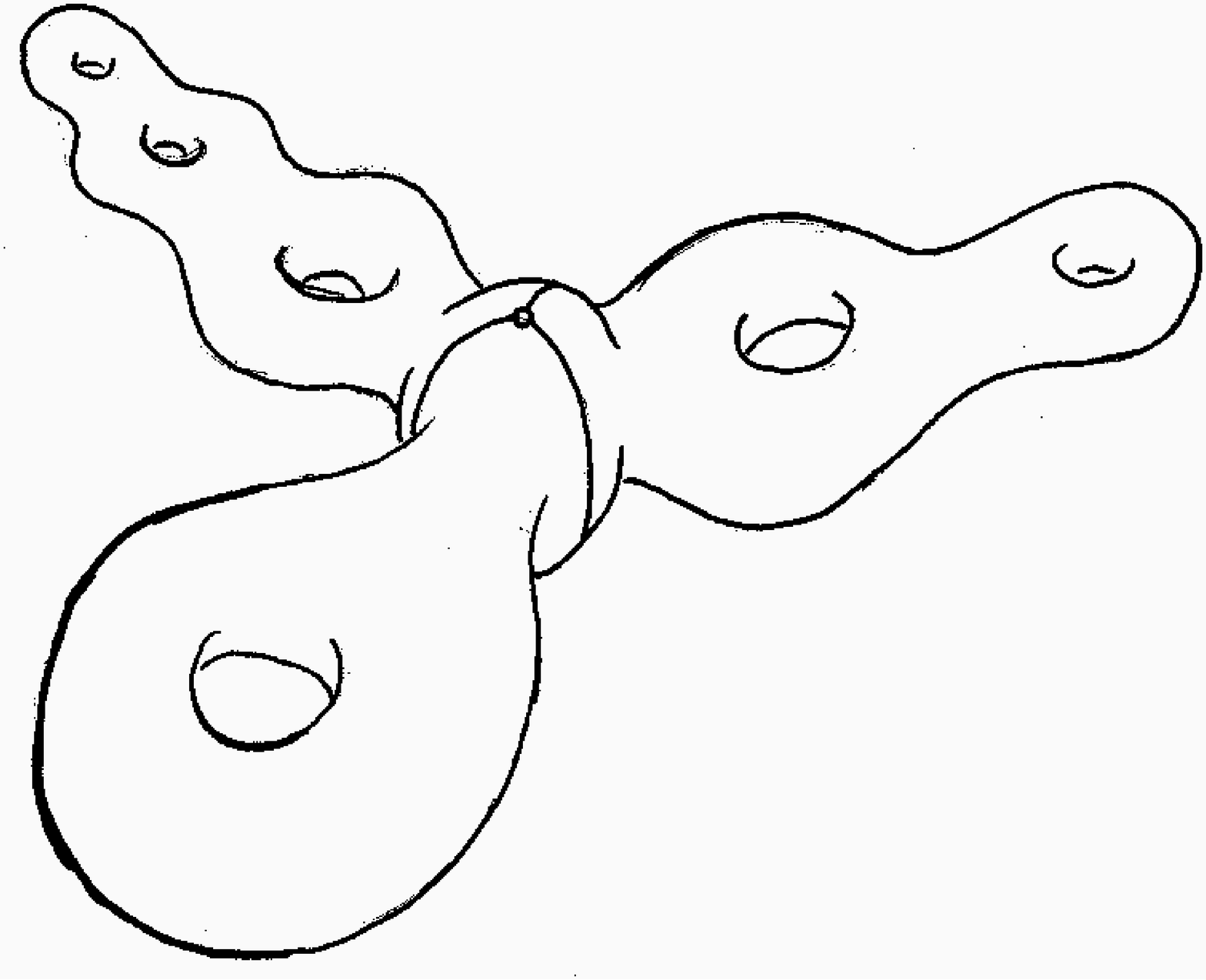}
 %
%
\begin{picture}(0,0)(0,0)
\put(0,0){
 \begin{picture}(0,0)(0,0) 
 \put(-155,-35){$S_3$}
 \put(-35,-45){$S_2$}
 \put(-129,-84){$\scriptscriptstyle \gamma'_3$}
 \put(-118,-84){$\scriptscriptstyle \gamma''_3$}
 \put(-105,-84){$\scriptscriptstyle \gamma'_2$}
 \put(-94,-84){$\scriptscriptstyle \gamma''_2$}
 \end{picture}
}
\put(5,0){
 \begin{picture}(0,0)(0,0) 
 \put(-125,-185){$S_1$}
 \put(-122,-121){$\scriptscriptstyle \gamma''_1$}
 \put(-112,-126){$\scriptscriptstyle \gamma'_1$}
 \end{picture}
}
\put(0,0){
 \begin{picture}(0,0)(0,0) 
 \put(50,-12){$S_3$}
 \put(152,-22){$S_2$}
 \put(90,-40){$z''$}
 \put(84,-81){$z'$}
 \end{picture}
}
\put(0,0){
 \begin{picture}(0,0)(0,0) 
 \put(80,-229){$S_1$}
 \put(50,-117){$S_3$}
 \put(152,-127){$S_2$}
 \put(84,-140){$z_2$}
 \put(85,-162){$\scriptstyle \gamma_2$}
 \put(76,-158){$\scriptstyle \gamma_1$}
 \end{picture}
}
\end{picture}
\vspace{230bp} 
\caption{
\label{pic:dist:sad:mult}
Multiple homologous saddle connections. }
\end{figure}

\begin{proof}
We  have  addressed every issue except the  statement that there  are precisely
$m+1$ surfaces obtained as the result of the assignment. Suppose that for fixed
$S'\in \cF'\subset\cH_1(\alpha')$ two of the  $m+1$ surfaces  $S$  built from $S'$
are  isomorphic. Since by construction each of these surfaces has a single short
saddle connection, the isomorphism sends the newborn saddle connection on one
surface to the newborn saddle connection on the other surface. Hence, it sends the
corresponding ``disc'' on one surface to the corresponding ``disc'' on the other
surface, see Figure~\ref{pic:bzero}. Hence, it is an isomorphism of the complements of the
``discs'', which implies that it induces an automorphism of the surface $S'$. It is
sufficient to note that the set of flat  surfaces, that have  automorphisms, has
measure $0$.
\end{proof}

\begin{remark}
\label{rm:cusp:is:a:ramified:covering}
Note     that      the     cohomology     class      $[\omega]\in
H^1(S,\{w,z_1,\dots,z_l\};\cx)$           together           with
$\overrightarrow{z'z''}=\gamma$  determines  an  element  of  the
relative          cohomology          group          $[\omega]\in
H^1(S,\{z',z'',z_1,\dots,z_l\};\cx)$.
The  lemma  above  claims, actually,      that       the  local    mapping
$\cH_1^{\varepsilon,3\varepsilon}(\alpha,\cC)\to\cF'\times B(\varepsilon)$
is a  ramified covering of  order $m+1$ almost everywhere and that
$d\nu(S)=d\nu(S')\, d\gamma$.
\end{remark}

\subsection{Computing the Siegel---Veech Constants}
\label{ss:comuting:const}

In  this  section we  derive  formulae for  the  constant $c$  in
quadratic asymptotics of the number of saddle connections joining
a zero  of order  $m_1$ to a zero of  order $m_2$ in multiplicity
one.
According to Proposition~\ref{prop:thick}
this means that we have to
find the asymptotics of
$\Vol(\cH_1^{\epsilon,\epsilon}(\alpha, \cC))$,
where the configuration $\cC$ corresponds to a single saddle connection.
Recall  also that  we
are  assuming  that the  zeroes  are  named.  Since  we  may have
numerous zeroes of orders $m_1$ and $m_2$ we,  actually, have two
different counting problems:

{\bf Problem 1.}  Count the constant in the quadratic asymptotics
for the  number of saddle  connections joining a {\it fixed} zero
$z_1$ of order $m_1$ to a  {\it fixed} zero $z_2$ of order $m_2$.
(When $m_1=m_2$ we require that $z_1\neq z_2$.)

{\bf Problem 2.}  Count the constant in the quadratic asymptotics
for the  number of saddle  connections joining {\it some} zero of
order $m_1$ to {\it some other} zero of order $m_2$. (This is, in
fact, equivalent  to  counting  the  saddle  connections  joining
unnumbered zeroes.)

Let us start by assuming that $\cH_1(\alpha)$ is {\em connected}.
The case of strata that are not connected will be postponed until
the end of the section on higher multiplicity.

In the computations to follow we will obscure
the distinction between $\gamma$ as a saddle connection  in the
configuration and its holonomy.  Thus we will use $\gamma$ as a vector
in the disc $B(\epsilon)$ and as a variable of integration.

Again  let  $n =  \dim_\reals  \cH(\alpha)$,  $n'  =  \dim_\reals
\cH(\alpha')$. Then  $n' = n  -2$. Recall that $d\nu'(S')$ is the
natural measure on $\cH(\alpha')$. For  $S'\in  C(\cF')$  we  set
$S'=rS''$  where  $\area(S'')=1$  so  $\area(S')=r^2$.  Then   by
Remark~\ref{rm:cusp:is:a:ramified:covering},
$$ 
d\nu(S) = d\nu'(S') \, d\gamma = r^{n'-1} \,
dr \,  d \vol'(S'') \,
d\gamma
$$ 
By Lemma~\ref{lemma:short:saddle:connections}
$$
\Vol\left(
 \cH_1^{\epsilon,\epsilon}(\alpha,\cC)\setminus
 \cH_1^{\epsilon,3\epsilon}(\alpha,\cC)
\right) =O(\epsilon^4)
$$
so  by Lemma~\ref{lm:breaking:zero}
\begin{align*}
\nu(C(\cH_1^{\epsilon,\epsilon}(\alpha,\cC))) & =  (m_1+m_2+1)\Vol(\cF')
\int_0^1 r^{n'-1}
  \int_{B(\epsilon r)}
  \, d\gamma \, dr +O(\epsilon^4) \\
& = \pi \epsilon^2 (m_1+m_2+1)\Vol(\cF') \int_0^1
r^{n'+1} \,
dr +O(\epsilon^4)
\\
& = \frac{\pi \epsilon^2}{n'+2}
(m_1+m_2+1)\Vol(\cF')+O(\epsilon^4).
\end{align*}
Since  $n = n' + 2$ we have
$$ 
\Vol( \cH_1^{\epsilon,\epsilon}(\alpha,\cC )) = n
\nu(C(\cH_1^{\epsilon,\epsilon}(\alpha))) = (m_1+m_2+1) \pi
\epsilon^2
\Vol(\cF')+O(\epsilon^4)
$$ 
Thus, the constant from the Problem 1 has the following form:
\begin{multline*}
c = \lim_{\epsilon \to 0}
\frac{\Vol(\cH_1^{\epsilon,\epsilon}(\alpha,\cC))}{\pi
  \epsilon^2 \Vol(\cH_1(\alpha))} =\lim_{\epsilon\to 0}
\frac{
(m_1+m_2+1)\Vol(\cF'))}{
\Vol(\cH_1(\alpha))}=\\
=\frac{
(m_1+m_2+1)\Vol(\cH_1(\alpha'))}{
\Vol(\cH_1(\alpha))}.
\end{multline*}

Note  that  having  solved the Problem  1  we  get  a solution of
Problem  2   by  an  elementary  combinatorial  calculation.  Let
$o(m_i)$ be the number of  zeroes  of order $m_i$ in the  stratum
$\cH(\alpha)$. If $m_1\neq  m_2$, then there are $o(m_1)$ ways of
choosing the  zero of order $m_1$  and $o(m_2)$ ways  of choosing
the  zero  of order  $m_2$.  Thus  we  get  an  additional factor
$o(m_1)\cdot o(m_2)$ in comparison with Problem 1.

\begin{formula}
The constant in  any  connected  stratum $\cH_1(\alpha)$ for the
number of  saddle  connections  of  multiplicity  one joining two
zeroes of orders $m_1\neq m_2$
is equal to
\end{formula}
%
\begin{equation}
\label{eq:sad:conn:mult:1}
c = o(m_1)o(m_2)\cdot(m_1+m_2+1)\cdot
\frac{\Vol(\cH_1(\alpha'))}{ \Vol(\cH_1(\alpha))}
\quad \text{when }m_1\neq m_2
\end{equation}

If  $m_1=m_2$,  then  there  are  $o(m_1)(o(m_1)-1)/2$  ways  of
choosing  an  unordered pair of zeroes of  order  $m_1$.  Having
chosen an  unordered pair, we  choose, which of two zeroes would
called be $z_1$, and which one would be $z_2$ in  arbitrary way:
the number of saddle connections  is  obviously  symmetric  with
respect  to  interchange of the names $z_1,  z_2$.  Thus,  when
$m_1=m_2$ the answer  to  Problem 2 is given  as  the answer for
Problem    1    multiplied    by     an     additional    factor
$o(m_1)(o(m_1)-1)/2$.

\begin{formula}
The
constant for the number of  saddle  connections  of  multiplicity
one joining two
distinct   zeroes   of  the  same  orders  $m_1$  is equal to
\end{formula}
%
\begin{equation}
\label{eq:sad:conn:mult:1:m1:m1}
c = \frac{o(m_1)(o(m_1)-1)}{2}\cdot(2m_1+1)\cdot
\frac{\Vol(\cH_1(\alpha'))}{ \Vol(\cH_1(\alpha))}
\qquad \text{ when }m_1= m_2
\end{equation}

Note that we may view the  entire  combinatorial  calculation  in
terms  of  the  naming  of  the zeroes.  We  can  assume  that on
$S'$ we  have broken up  the zero of order $m=m_1+m_2$ at the
first ``named'' zero of order $m$. The remaining  zeroes are then
given  the  same names on $S$. Other  zeroes  of  orders
$m_1,m_2$  already   have   names   on  $S$.  There  are
$o(m_1)o(m_2)$  ways  of  adding  additional names to  zeroes  of
orders $m_1,m_2$  if  $m_1\neq  m_2$  and  $o(m_1)(o(m_1)-1)$  if
$m_1=m_2$.

Note that in  the  consideration above the stratum
$\cH_1(\alpha')$
might be  nonconnected. In this case $\Vol(\cH_1(\alpha'))$ is the
sum of the  volumes of connected  components. For example  if  we
start with the stratum $\cH_1(5,3)$, then  $\cH_1(\alpha')$ has
three connected
components  $\cH_1^{hyp}(8),\cH_1^{even}(8)$  and $\cH_1^{odd}(8)$,
and
the sum of the volumes of these must be in the numerator.

\subsection{Examples: Constants for Connected Strata in Genus $3$}
\label{ss:examples:part1}

As an illustration of the formulae above we now give the explicit
values for these constants in a number of cases.

\begin{table}[ht!]
\small
\caption{Normalized   volumes
$\cfrac{1}{\pi^{2g}}\Vol(\cH_1(\alpha))$ of the strata in small
genera}
\label{tab:volumes}
%
\scriptsize
\vspace*{-12truept}
$$
\begin{array}{|c|c|c|}
\multicolumn{3}{c}{\text{Genera $g=1,2$}}\\
[-\halfbls]\multicolumn{3}{c}{\text{}}\\
\hline &&\\
\cH_1(\torusemptyset) & \cH_1(2)         & \cH_1(1,1) \\
[-\halfbls] &&\\ \hline && \\ [-\halfbls]
\cfrac{1}{3}    & \cfrac{1}{120} & \cfrac{1}{135}\\
&&\\ \hline
\end{array}
$$

\bigskip
%
$$
\begin{array}{|c|c|c|c|c|c|c|}
\multicolumn{7}{c}{\text{Genus $g=3$}}\\
[-\halfbls]\multicolumn{7}{c}{\text{}}\\
\hline &&&&&&\\
\cH_1^{hyp}(4)   & \cH_1^{odd}(4)     & \cH_1(3,1)             &
\cH_1^{hyp}(2,2)       & \cH_1^{odd}(2,2)
& \cH_1(2,1,1)           & \cH_1(1,1,1,1) \\
[-\halfbls] &&&&&&\\ \hline &&&&&& \\ [-\halfbls]
\cfrac{1}{6720}& \cfrac{1}{2430}  & \cfrac{16}{42525}    &
\cfrac{1}{9450}      &
\cfrac{1}{4320}  & \cfrac{1}{3780}      & \cfrac{1}{4860}\\
&&&&&&\\ \hline
\end{array}
$$

\bigskip
%
%
$$
\begin{array}{|c|c|c|c|}
\multicolumn{4}{c}{\text{Genus $g=4$}}\\
[-\halfbls]\multicolumn{4}{c}{\text{}}\\ \hline &&&\\
\cH_1^{hyp}(6)     & \cH_1^{odd}(6)     & \cH_1^{even}(6)        &\cH_1(5,1)
\\ [-\halfbls] &&&\\ \hline &&& \\ [-\halfbls]
\cfrac{1}{580608}& \cfrac{1}{37800} &
\cfrac{32}{1913625}       & \cfrac{1}{36750}
%
\\&&&\\ \hline\hline &&&\\
\cH_1^{odd}(4,2) &   \cH_1^{even}(4,2) &
\cH_1(4,1,1) & \cH_1^{nonhyp}(3,3)
\\ [-\halfbls] &&&\\ \hline &&& \\ [-\halfbls]
\cfrac{1}{79380} & \cfrac{1}{107520} & \cfrac{11}{653184} &\cfrac{1}{51030}
%
\\&&&\\ \hline\hline &&&\\
\cH_1^{hyp}(3,3)  & \cH_1(3,2,1)       & \cH_1(3,1,1,1) & \cH_1^{odd}(2,2,2)
\\ [-\halfbls] &&&\\ \hline &&& \\ [-\halfbls]
\cfrac{1}{992250} & \cfrac{1}{70875} &\cfrac{62}{5740875} & \cfrac{31}{4354560}
%
\\&&&\\ \hline\hline &&&\\
\cH_1^{even}(2,2,2) &
\cH_1(2,2,1,1) & \cH_1(2,1,1,1,1)    &\cH_1(1,\dots,1) \\
\\ [-\halfbls] &&&\\ \hline &&& \\ [-\halfbls]
\cfrac{37}{6804000}  &
\cfrac{131}{13608000} & \cfrac{1}{136080}
&\cfrac{377}{67359600}\\ &&&\\ \hline
\end{array}
$$
\end{table}

In  Table~\ref{tab:volumes}  we present  the  normalized  volumes
$\cfrac{1}{\pi^{2g}} \Vol(\cH_1(\alpha))$ of the strata with {\it
numbered} zeroes, see~\cite{Eskin:Okounkov}, see also some values
in~\cite{Zorich}. These values will be used for computations in all
the examples.


Both  strata  in  genus  $g=2$  are  hyperelliptic; they will  be
treated                later                on                 in
section~\ref{s:distinct:zeroes:nonconnected:strata}.   In   genus
$g=3$ there  are  three connected strata: $\cH_1(3,1),
\cH_1(2,1,1)$,
and $\cH_1(1,1,1,1)$.

\begin{example}
{\bf Stratum} $\boldsymbol{\cH_1(3,1)}.$
After collapsing zeroes  we obtain a flat surface $S'\in
\cH_1(4)$,
where  both  components $\cH_1^{hyp}(4)$  and  $\cH_1^{odd}(4)$
occur.
We  have  $m_1=3,  m_2=1$,  thus  the  combinatorial
factor equals $o(m_1)o(m_2)(m_1+m_2+1)=5$. We get
$$ 
c=5\cdot\frac{\Vol(\cH_1(4))}{\Vol(\cH_1(3,1))}=
5\cdot\frac{\Vol(\cH_1^{hyp}(4))+\Vol(\cH_1^{odd}(4))}{\Vol(\cH_1(3,1))}=
\frac{7625}{1024}\approx 7.45
$$ 
\end{example}

\begin{example}
{\bf Stratum} $\boldsymbol{\cH_1(2,1,1)}.$
There are two cases here in multiplicity $1$. The first case is a
saddle connection joining the zero of degree $2$ with any  of the
two zeroes of degree $1$.  Collapsing  the zeroes we obtain a
surface  $S'\in  \cH_1(3,1)$.  For  $m_1=2$,  $m_2=1$,
$o(m_1)=1$, $o(m_2)=2$.  Thus  the  combinatorial constant equals
$8$. We get
$$ 
c=8\cdot\frac{\Vol(\cH_1(3,1))}{\Vol(\cH_1(2,1,1))}=
\frac{512}{45}\approx 11.4
$$ 

The second case consists of collapsing the pair of simple zeroes
to
a zero of order $2$.  We  obtain  a surface $S'\in
\cH_1(2,2)$.  Any  surface in any component of  $\cH_1(2,2)$  can
be
found by collapsing a pair of simple zeroes in this  manner. Here
$m_1=m_2=1$ and $o(m_1)=2$ so the combinatorial constant
is $3$. We get
$$ 
c=3\cdot\frac{\Vol(\cH_1(2,2))}{\Vol(\cH_1(2,1,1))}=
3\cdot\frac{\Vol(\cH_1^{hyp}(2,2))+\Vol(\cH_1^{odd}(2,2))}{\Vol(\cH_1(2,1,1))}=
\frac{153}{40}\approx 3.83
$$ 
\end{example}


\subsection{Constants for the Principal Stratum  $\cH_1(1,\dots,1)$}

The
surface of  genus $g$ has  $2g-2$
simple   zeroes;    $o(1)=2g-2$.   The   corresponding    surface
$S'\in\cH_1(2,1,\dots,1)$ has $2g-4$  simple zeroes.     The
factor
$o(m_1)\left(o(m_1)-1\right)/2$  gives   $(g-1)(2g-3)$,  and  the
factor        $(2m_1+1)$        gives        $3$.        Applying
formula~\eqref{eq:sad:conn:mult:1:m1:m1}) we get

\begin{formula}
\label{f:dist:zeroes:principal:mult:1}
 The
constant for the number of  saddle  connections  of  multiplicity
one joining two
distinct  zeroes  is equal to
\end{formula}
%
        %
$$ 
c=3(g-1)(2g-3)\cdot\frac
{\Vol(\cH_1(2,\overbrace{1,\dots,1}^{2g-4}))}
{\Vol(\cH_1(\underbrace{1,\dots,1}_{2g-2}))}
$$ 

In the table  below  we  present some  values of  the constant  in
this case.
   %
%

\begin{table}[hbt!]
\caption{
Principal stratum $\cH_1(1,\dots,1)$; values of the
constants for  saddle  connections  of  multiplicity  one joining
distinct zeroes}
\label{tab:c:distinct:princ:mult1}
\scriptsize
$$
\begin{array}{|l|c|c|c|c|c|c|c|}
\hline &&&&&&&\\
& g=2 & g=3 & g=4 & g=5 & g=6 & g=7 & g=8  \\
[-\halfbls] &&&&&&& \\ \hline &&&&&&& \\ [-\halfbls]
c = &   \cfrac{27}{8} & \cfrac{162}{7} & \cfrac{22275}{377} &
\cfrac{2594700}{23357} &
\cfrac{2954056635}{16493303} & \cfrac{13229971542}{50280671} &
\cfrac{14740938123723}{40593663941}
\\  [-\halfbls]
&&&&&&& \\ \hline &&&&&&& \\ [-\halfbls]
c \approx &  3.375 &  23.14 & 59.08 & 111.1 & 179.1 & 263.1 &
363.1
\\ [-\halfbls]
&&&&&&&\\ \hline
\end{array}
$$
\end{table}

\section{Multiple Homologous Saddle Connections. Connected Strata}

\subsection{Principal Boundary}
\label{ss:principal:boundary}
Consider a   surface $S\in\cH_1(\alpha)$ with a fixed pair of
zeroes $z_1$ and $z_2$ of orders  $m_1$ and $m_2$ correspondingly.
Suppose that we have a configuration $\cC$ (see Section 3) of
precisely $p$ homologous saddle connections $\gamma_1,\dots,
\gamma_p$  joining   $z_1$ to $z_2$ of length at most $\epsilon$. Assume there are no other saddle connections shorter than $3\epsilon$.    A pair $\gamma_i$ and
$\gamma_{i+1}$ bounds a surface $S_i$. The surfaces $S_i$ and
$S_{i+1}$  share   the   saddle connection $\gamma_{i+1}$.  By
convention the cyclic order  of $\gamma_i$ at $z_1$ is clockwise
in the orientation defined by the flat structure.  The angle
between   $\gamma_i$   and $\gamma_{i+1}$ at  $z_1$  is
$2\pi(a'_i+1)$;  the   angle between $\gamma_i$ and $\gamma_{i+1}$
at  $z_2$  be $2\pi(a''_i+1)$.

Cut the surface  $S$ along all $\gamma_i$, $i=1,\dots,p$. Gluing
together the
two sides $\gamma_i$  and $\gamma_{i+1}$ of the boundary of $S_i$
we get a  flat surface with two distinguished zeroes $z'_i,z''_i$
of orders $a'_i$ and  $a''_i$  correspondingly joined by a saddle
connection,  $\overrightarrow{z'_i z''_i}=\gamma$.  When  one  of
$a'_i, a''_i$ is equal to zero the corresponding point is  just a
marked point.  If  $a'_i=a''_i=0$,  then  both  points are marked
points.

By construction any resulting flat surface $S_i$ has a
saddle connection $\gamma_i$ shorter than $\epsilon$   joining the pair of zeroes and does
not  have  any other saddle connections homologous  to  it.  Every other saddle connection is longer than $3\epsilon$. This
implies  that  contracting the saddle connection $\gamma_i$ to  a
point, see section~\ref{ss:collapsing},
we    get    a   flat   surface
$S'_i$ with a distinguished zero (marked point) $w_i$
of order $a_i=a'_i+a''_i$ and no saddle connection shorter than $2\epsilon$.

Let $S'_i\in \cH(\alpha'_i)$. We say that $\sqcup_{i=1}^p\cH(\alpha_i')$ is the
principal boundary of this configuration.
Note that if $a_i=0$ we mark corresponding point on $S'_i$,
so collection $\alpha'_i$ contains $0$ in this case. We use the notation
$\alpha'=\sqcup \alpha'_i$, see Convention~\ref{conv:alpha:prime}.

\subsection{Slit Construction}
\label{ss:slit:construction}

We would like to reverse the above degeneration and build surfaces
with
multiple homologous  saddle
connections out of surfaces in $\cH(\alpha_i')$. We need the
 following  construction.  Let  $S'$ be a
surface with a zero of  order  $a\geq 0$.  Suppose there are no
saddle connections of length shorter than $2\epsilon$. Let
$a',a''\geq 0$  such
that $a'+a''=a$. Let $\gamma$ be a  vector of length at most $\epsilon$.
If $a',a''\neq 0$ we may
break up the  zero  of order $a$ into  zeroes  of orders $a',a''$
with a saddle connection determining  $\gamma$  joining them,
see section~\ref{ss:breaking:zero}. If $a'=0$ take a point $z'$
of the  form $w-\gamma$ (this means along a  geodesic from $w$ in
direction  $-\gamma$  and  distance  $|\gamma|$) and join  it  to
$z''=w$  on  $S'$. If $a''=0$ we  take  a segment from $z'=w$  to
$z''=w+ \gamma$ on $S'$.

In  either  case we  may  slit the  resulting  surface along  the
saddle connection joining the two zeroes $z',z''$. In  this way we
build a surface with one boundary  component  consisting  of  two
arcs, denoted $\gamma'$ and $\gamma''$, joining  the endpoints of
the  slit.  The angles between $\gamma'$ and  $\gamma''$  at  the
points  $z'$   and   $z''$  are  $2\pi(a'+1)$  and  $2\pi(a''+1)$
correspondingly.

Note  that  the flat structure on  $S'$  fixes the choice of  the
orientation.  By  convention we give the ``names'' $\gamma'$  and
$\gamma''$ to the arcs in such  a way that turning around $z'$ in
a  clockwise  direction from $\gamma''$ to $\gamma'$  we  do  not
leave the surface.

Conversely,  consider  a flat  surface  with  a  single  boundary
component consisting of two arcs joining a pair  of vertices. The
convention  on  the   choice   of  the  ``names''  $\gamma'$  and
$\gamma''$ means that as soon as we know which of two vertices is
$z'$, and  which is  $z''$ we can determine which  of two arcs is
$\gamma'$, and which is $\gamma''$.

We describe now how to build  surfaces  with  multiple
homologous saddle connections using the slit construction.

\subsection{Building Surfaces With Multiple Homologous Saddle
Connections}
\label{ss:building:surfaces}

Suppose     $S'=S_1'\sqcup\ldots     \sqcup
S_p'$ is  a  disconnected  flat  surface  and on each
$S_i'$ there  is a zero or  a regular point  which we
think  of  as a {\em marked  point}  $w_i$ of order $a_i$.  Assume no surface contains a saddle connection shorter than $2\epsilon$. Given
pairs $a_i', a_i''$ with $a_i=a_i'+a_i''$, and a vector $\gamma$ of length smaller than $\epsilon$,
we  perform  the  slit  construction on each surface.  We  obtain
surfaces with one boundary component  each  of  which consists of
two arcs $\gamma_i'$ and $\gamma_i''$.  We  glue  $\gamma_i'$  to
$\gamma_{i+1}''$, calling this curve $\gamma_{i+1}$, (and $\gamma_p'$
to $\gamma_1''$). This gives a  closed  surface  $S$,  a
pair of  zeroes $z_1$ and $z_2$ of orders  $m_1$ and $m_2$, where
$\sum  a_i'=m_1+1-p$  and $\sum  a_i''=m_2+1-p$,  and  a  set  of
homologous  curves  $\gamma_i$, $i=1,\dots,p$  joining  $z_1$  to
$z_2$. The angle  between  $\gamma_{i}$ and $\gamma_{i+1}$ at $z_1$
is  $2\pi(a_i'+1)$  and  the  angle  between  them  at  $z_2$  is
$2\pi(a_i''+1)$. We denote this assignment by
$$ 
(S',\gamma,a_i',a_i'')\to   (S,m_1,m_2,a_i',a_i'')
$$ 

The resulting surface $(S,m_1,m_2,a_i',a_i'')$ has no saddle connections shorter than $\epsilon$ other than the $\gamma_i'$.
Our convention on  the choice of $\gamma_i', \gamma_i''$ (see the
section  above)   implies   that   the   cyclic  order  $\dots\to
\gamma_{i-1}\to \gamma_i\to\gamma_{i+1}\to \dots$ is clockwise at
$z_1$  with  respect  to  the  orientation  defined  by  the flat
structure.

Fix a  basis for the relative homology on  each surface $S_i'$. A
relative homology basis for $S$ is given by:
\begin{itemize}
\item A relative homology basis of curves $\beta_i$  for $S'$.
\item A curve connecting the two zeroes. The integral of $\omega$
along this curve is $\gamma$.
\end{itemize}

We  introduce  $\bar{\alpha}'_i$ as $\alpha'$
from which we remove $a_i$.

Having a collection $\alpha$ of integers  we denote by $|\alpha|$
their sum, say, $|\alpha|=2g-2$.

The  construction  above  implies  the  following  conditions  on
partitions $\alpha'_i$ and the numbers $a_i', a_i''$.
\begin{lemma}
\label{lm:admissible:split}
An         assignment          $(S',\gamma,a_i',a_i'')\to
(S,m_1,m_2,a_i',a_i'')$ satisfies  the following necessary
conditions on the collection $(\alpha'_i,a_i',a''_i)$:
\begin{gather*}
a'_1+\dots+a'_p   = m_1 + 1  -p \\
a''_1+\dots+a''_p = m_2 + 1  -p \\
\bar{\alpha}'_1 \sqcup \dots \sqcup  \bar{\alpha}'_p
    \sqcup \{m_1\}\sqcup\{m_2\} = \alpha \\
|\bar{\alpha}'_i| = a'_i+a''_i\ (mod\ 2) \text{ for } i=1,\dots,p
\end{gather*}
Moreover, when  the  stratum  $\cH_1(\alpha)$  is  connected,
these
conditions       are        sufficient:       every       surface
$S\in\cH_1(\alpha)$ with a configuration of $p$ homologous saddle
connections
joining  the pair  of  zeroes of length at most $\epsilon$ and no other saddle connection with length smaller than $3\epsilon$ can be  obtained  by an  assignment
$(S',\gamma,a_i',a_i'')\to (S,m_1,m_2,a_i',a_i'')$
with appropriate $(S',\gamma,a_i',a_i'')$.
\end{lemma}

We shall introduce   notation to describe the degeneration pattern
for the surfaces $S_i$.

{\bf Notation.} The integer represented as a sum of two integers
corresponds to $a_i$ represented as $a_i'+a_i''$.  Here the order
of appearances of the
summands $a'_i+a''_i$ is significant. The  cyclic  order  of  the
surfaces is represented by
$S_1\succ S_2\succ S_3\succ \ldots \succ S_p\succ
S_1$.

\begin{example}
\label{ex:degeneration:patterns}
Let $S\in  \cH_1(4,3,2,1)$;  $m_1=3,  m_2=4$; $p=3$. There
are  $15$   possible   pictures   for   the three  homologous
saddle
connections joining the zero of  order  $3$ to the zero of  order
$4$  depending  on  the  returning angles (Table 3).
Thus for example,  the first line in table 3
refers to $3$ surfaces in cyclic order. The first is a torus (no
zeroes) where we have broken up a marked point into  $(0+0)$; the
second surface has a single zero of order $2$ which we have broken
into two simple zeroes indicated by $(1+1)$, and the third surface
has zeroes of $2,1,1$, where we have broken one of the simple
zeroes into a simple zero and a zero of order $0$.

\end{example}

\begin{table}
\begin{displaymath}
\begin{array}{c|ccc}
\text{Degeneration pattern} & \alpha'_1  &  \alpha'_2  &  \alpha'_3  \\
\hline                      &             &             &
   \\
(0+0)\succ(1+1)\succ(0+1,2,1)\succ         & (\torusemptyset) & (2)
& (2,1,1)

    \\
(1+1)\succ(0+0)\succ(0+1,2,1)\succ         & (2)         &
(\torusemptyset) & (2,1,1)

   \\
(0+0)\succ(0+0)\succ(1+2,2,1)\succ         & (\torusemptyset) &
(\torusemptyset) & (3,2,1)

   \\
(0+2)\succ(0+0)\succ(1+0,2,1)\succ         & (2)         &
(\torusemptyset) & (2,1,1)

   \\
(0+0)\succ(0+2)\succ(1+0,2,1)\succ         & (\torusemptyset) & (2)
& (2,1,1)

   \\
(1+1)\succ(0+1,1)\succ(0+0,2)\succ         & (2)         & (1,1)
& (2)
    \\
(0+0)\succ(1+0,1)\succ(0+2,2)\succ         & (\torusemptyset) & (1,1)
& (2,2)
    \\
(0+2)\succ(1+0,1)\succ(0+0,2)\succ         & (2)         & (1,1)
& (2)
    \\
(0+0)\succ(1+2,1)\succ(0+0,2)\succ         & (\torusemptyset) & (1,3)
& (2)
    \\
(0+0)\succ(0+1,1)\succ(1+1,2)\succ         & (\torusemptyset) & (1,1)
& (2,2)
    \\
(1+1)\succ(0+0,2)\succ(0+1,1)\succ         & (2)         & (2)
& (1,1)
    \\
(0+0)\succ(1+1,2)\succ(0+1,1)\succ         & (\torusemptyset) & (2,2)
& (1,1)
    \\
(0+0)\succ(0+0,2)\succ(1+2,1)\succ         & (\torusemptyset) & (2)
& (1,3)
    \\
(0+2)\succ(0+0,2)\succ(1+0,1)\succ         & (2)         & (2)
& (1,1)
    \\
(0+0)\succ(0+2,2)\succ(1+0,1)\succ         & (\torusemptyset) & (2,2)
& (1,1)
\end{array}
\end{displaymath}
\medskip

\caption{
\label{4321:mult:3}
Possible degenerations for
$S\in \cH_1(4,3,2,1)$ with $m_1=3,  m_2=4$;  $p=3$}
\end{table}

\subsection{Stratum Interchange and $\gamma\to-\gamma$ Symmetry}
\label{ss:stratum:interchange}
In  this  section  we  discuss  the  possible symmetries of  the
assignment
$$ 
(S',\gamma,a'_i,a''_i)\to(S,m_1,m_2,a_i',a_i'')
$$ 
Let us first specify our problem.

{\bf Problem 1.}
We  assume  that  the  zeroes  $z_1,\dots,z_l$  of  the  surface
$S$  are  numbered. We  fix  the  zeroes  $z_1,  z_2\in
S$,  $z_1\neq   z_2$,   of   orders   $m_1$  and  $m_2$
correspondingly.

We start  with the  case when we fix also  the following data. We
assume that the zeroes $z_3, \dots, z_{l_1}$ belong to $S_1$, the
zeroes $z_{l_1+1},\dots,  z_{l_1+l_2}$  belong to $S_2$, ..., the
zeroes $z_{l_1+\dots+l_{p-1}+1}, \dots, z_{l_1+\dots+l_p}$ belong
to $S_p$.  Here $l_1-1=card(\alpha'_1)$ if $a_1\neq 0$ and $l_1=card(\alpha'_1)$ when $a_1=0$ and for $i>1$,   $l_i+1=card(\alpha'_i)$,  when  $a_i\neq 0$, and
$l_i=card(\alpha'_i)$, when  $a_i=0$. We have  $l = l_1 + \dots +
l_p = card(\alpha)$.

In this setting the only thing that is not determined is the case
when $S$ has exactly  two  zeroes: $z_1$ and $z_2$. Then
every $S'_i$ either  has a single zero of order $a_i$
or it is a torus with  a marked (regular) point. Note that we fix
only the cyclic order  of  the collection $(a'_i,a''_i)$, but not
the  numbering of  each  of these pair  of  numbers. The  natural
action  of  the cyclic group of  order  $p$ on the collection  of
ordered pairs $(a'_i,a''_i)$ organized in  a  cyclic  order,  may
have nontrivial stabilizer, which  we  denote by $\Gamma$. We get
symmetry of  order $|\Gamma|$. This  type of symmetry we call the
{\it stratum interchange}.

{\bf Problem 2.}  Now consider the problem with fewer constrains.
We count the number  of  occurrences of the following phenomenon:
{\it some} zero of $S$ of order $m_1$ is joined  to {\it
some}  zero  of order $m_2$ by precisely  $p$  homologous  saddle
connections $\gamma_1,\dots,\gamma_p$ as above.
   %
       %
We   fix   the    type    $\alpha'_i$   of   the   flat   surface
$S'_i$  obtained  as described above, but we make  no
assumptions on the distribution of the  zeroes $z_1,\dots,z_l$ by
the surfaces $S'_i$. As before we fix only the cyclic
order of appearances of the surfaces $S_i$.

\begin{remark}
Problem 2 may be considered as generalization of Problem 1 to the
case  when  the zeroes of $S$ and of  $S'_i$
are {\it not} numbered.
\end{remark}

\begin{remark}
In both settings the $\cH_1(\alpha_i')$   might  be nonconnected.
We could   also    specify    the   connected component   of   the
$\cH_1(\alpha'_i)$  when  this occurs. We prefer the setting where
this data is not specified.
\end{remark}

Consider the  natural action of the cyclic group  of order $p$ on
the collection
$$
(\alpha'_1,a'_1,a''_1) \succ \dots \succ (\alpha'_p,a'_p,a''_p)
\succ
$$
organized in a cyclic order. (Recall the notation $\succ$ from the previous section).  If it has a nontrivial stabilizer we
denote it by the same symbol $\Gamma$. We get a symmetry of order
$|\Gamma|$ which we also call the {\it stratum interchange}.

\begin{example}
Let $\alpha=(27,15,1,1)$, $m_1=27$, $m_2=15$.    Consider     the
degeneration of $S$  into  $6$ surfaces of the following
types $\alpha'_1 \succ \dots \succ \alpha'_6 \succ$
$$ 
(2)\succ (7,1) \succ (12) \succ (2)\succ (7,1) \succ (12) \succ
$$ 
The possible symmetry in this  case  is a symmetry of order  $2$,
where subgroup $\Gamma$ is a shift by $3$.

The degeneration
$$ 
(0+2)\succ(4+3,1)\succ(7+5)\succ(0+2)\succ(4+3,1)\succ(7+5)\succ
$$ 
possesses this symmetry, while the degeneration
$$ 
(0+2)\succ(3+4,1)\succ(7+5)\succ(2+0)\succ(4+3,1)\succ(7+5)\succ
$$ 
does not.
\end{example}

When $m_1=m_2$  we may have  an additional symmetry in Problem 2.
Let $P, Q$ be two zeroes of $\omega$ of order $m_1=m_2$ joined by
precisely $p$  homologous  saddle  connections.  Let the homology
class  of  the  saddle  connection  be  represented  by  a vector
$\gamma$.  Assigning  the ``names'' $z_1:=P$ and $z_2:=Q$ to  $P$
and        $Q$        we        get        a        decomposition
$(S'_i,\gamma,\alpha'_i,a'_i,a''_i)$.

Consider now the very same surface  $S$  with  the  same
configuration of  $p$  homologous  saddle connections joining the
same  two  zeroes  $P$  and  $Q$. We  may  declare  now  that the
correspondence of zeroes  is inverse with respect to the previous
assignment: $z_1:=Q$  and  $z_2:=P$. Since our saddle connections
are oriented (from $z_1$ to $z_2$) the homology class of the same
saddle connection (from  $P$ to $Q$)  is represented now  by  the
vector  $-\gamma$.  Since  the  cyclic  order  of  ``pieces''  is
determined  by  the cyclic  order  at  the  point  $z_1$  the new
identification  reverses  the  cyclic  order  in  the  collection
$(\alpha'_i)$, as well as the order in the pairs $(a'_i,a''_i)$.

We say,  that we have  a $\gamma\to-\gamma$ {\it symmetry} if and
only  if   the   assignment   $z_1:=Q$   and   $z_2:=P$  gives  a
decomposition  with  the same  (up  to  a  cyclic  reenumeration)
collection $(\alpha'_i,a'_i,a''_i)$ as before.

It is  easy to see the  following combinatorial criterion  of the
$\gamma\to-\gamma$  symmetry.  Interchange   simultaneously   all
$a'_i\leftrightarrow a''_i$, and  change  the cyclic order of the
resulting  collection $(\alpha'_i,a''_i,a'_i)$  to  the  opposite
one. We possess a $\gamma\to-\gamma$ symmetry if and  only if the
result  of  this  operation  gives  us   the  initial  collection
$(\alpha'_i,a'_i,a''_i)$ up  to a cyclic reenumeration. We denote
$\gamma\to-\gamma$ symmetry by  $\Gamma_-$, where $|\Gamma_-|=2$.
We  let   $|\Gamma_-|=1$  when  there  is  no  $\gamma\to-\gamma$
symmetry.

The condition $m_1=m_2$ is an  obvious  necessary  condition  for
$\gamma\to-\gamma$ symmetry. It will be  always  assumed  in  the
discussions of $\gamma\to-\gamma$ symmetry.

Under this  condition  every  degeneration  of  multiplicity  $1$
possesses   a   $\gamma\to-\gamma$   symmetry.    There    is   a
$\gamma\to-\gamma$ symmetry in multiplicity $2$ if and only if at
least one the following conditions is valid: either\newline
--- $a'_1=a''_1$ (and then $a'_2=a''_2$ automatically);
or\newline
---   $\alpha'_1=\alpha'_2$    and   $a'_1=a''_2$   (and    then
$a''_1=a'_2$ automatically).

\begin{example}
Let    $\alpha=(11,11)$.     Consider degeneration  of
$S$  into   $4$   surfaces   of   the   following  types
$\alpha'_1\succ \alpha'_2\succ \alpha'_3\succ \alpha'_4$:

$$ 
(2)\succ (6) \succ (2)\succ (6) \succ
$$ 
The degeneration
$$ 
(0+2)\succ(3+3)\succ(2+0)\succ(3+3)\succ
$$ 
possesses the $\gamma\to-\gamma$ symmetry, but does  not have the
stratum interchange, so $|\Gamma_-|=2$, $|\Gamma|=1$.

The degeneration
$$ 
(0+2)\succ(4+2)\succ(0+2)\succ(4+2)\succ
$$ 
does  not  have  a  $\gamma\to-\gamma$ symmetry, but  allows  the
stratum interchange, so $|\Gamma_-|=1$, $|\Gamma|=2$.

The degeneration
$$ 
(1+1)\succ(3+3)\succ(1+1)\succ(3+3)\succ
$$ 
has both symmetries: $|\Gamma_-|=2$, $|\Gamma|=2$.
\end{example}

\begin{remark}
Describing                    the                     assignments
$(S',\gamma,a'_i,a''_i)\to(S,m_1,m_2,a_i',a_i'')$
one can consider {\it all}  possible  assignments,  and then take
into account corresponding symmetries,  or  one can deal with the
classes, eliminating the symmetry whenever it is possible.

We have chosen the second way. For example the assignment
$$
(2+0)\succ(0+0)\succ(0+1,2,1)
$$
is  not  presented  in  the Table~\ref{4321:mult:3} since  it  is
symmetric to the assignment
$$
(0+0)\succ(0+2)\succ(1+0,2,1)
$$
(which  is  in  the  list) by composition  of  $\gamma\to-\gamma$
symmetry with the stratum interchange.
\end{remark}
%
$$
\ast\quad\ast\quad\ast
$$
\medskip

Let us return to the degeneration construction described in
Section~\ref{ss:principal:boundary}. To every surface $S$ in
the thick part $\cH_1^{\epsilon,thick}(\alpha,\cC)$ we
associate a surface $S'=\sqcup_{i=1}^p S'_i$. We
denote the corresponding stratum by $\cH(\alpha')$,
where $\alpha'=\sqcup \alpha'_i$, see Convention~\ref{conv:alpha:prime}.

However, our consideration of the symmetries shows, that the surface
$S'$ is defined up to a finite symmetry. In other words,
we may have a nontrivial symmetry group $\Gamma_\pm$
acting on the stratum $\cH(\alpha')$. This symmetry group
is generated by the subgroup $\Gamma$ of
stratum interchange symmetries and subgroup $\Gamma_-$ of
$\gamma\to-\gamma$ symmetries. (Any of the subgroups, or both might be trivial.)
The order $|\Gamma_\pm|$ of the symmetry group $\Gamma_\pm$ is equal to the product of the orders
$|\Gamma|$ and $|\Gamma_-|$.

Note that the symmetry group $\Gamma_\pm$ preserves the natural measure
on $\cH(\alpha')$. It obviously preserves the area of the surface,
hence $\Gamma_\pm$ preserves the ``unit hyperboloid'' $\cH_1(\alpha')$.
Thus we get a natural volume element on the quotient
$\cH_1(\alpha')/\Gamma_\pm$.

The degeneration construction can be considered
as a map
$$
\cH_1^{\epsilon,thick}(\alpha,\cC)\to
(\cH_1(\alpha')/\Gamma_\pm)\times B(\epsilon)
$$
In the next section we shall study this map. This will allow us
to compute the volume of the thick part and thus compute the corresponding
Siegel---Veech constant.

\subsection{Siegel---Veech Constants in Higher Multiplicity, Connected Strata}
\label{ss:c:in:higher:mult}

Consider the set $\cH_1^{\epsilon,thick}(\alpha,\cC)$.
It consists  of  surfaces  $S$  having   exactly one
collection of short homologous saddle connections.
This collection is necessarily of the type $\cC$;
the saddle connections from this collection
(which are all of the same length) are shorter
then $\epsilon$; any other saddle connection on $S$
is longer then $\kappa$.

Actually, we shall need two variants of this set: the one,
which we denote $\cH_1^{\epsilon,\epsilon}(\alpha,\cC)$
corresponds to the value $\kappa=\epsilon$
of parameter $\kappa$; the other, denoted
$\cH_1^{\epsilon,3\epsilon}(\alpha,\cC)$,
corresponds to $\kappa=3\epsilon$.

Clearly,
\begin{equation}
\label{eq:epsilon:3epsilon}
\cH_1^{\epsilon,3\epsilon}(\alpha,\cC)\subset
\cH_1^{\epsilon,\epsilon}(\alpha,\cC).
\end{equation}

\begin{lemma}
\label{lemma:1234}
Consider the map
$$
p:\cH_1^{\epsilon,\epsilon}(\alpha,\cC) \to
\cH_1(\alpha')/\Gamma_\pm \times B(\epsilon)
$$
corresponding to the assignment
$(S',\gamma,a_i',a_i'')\to
(S,m_1,m_2,a_i',a_i'')$.
Let
$$
\cF'_{2\epsilon}:=\left(\cH_1(\alpha') - \cH_1^{2\epsilon}(\alpha')\right)/\Gamma_\pm.
$$

1. The image of the map $p$
contains the subset
$\cF'_{2\epsilon} \times B(\epsilon)$.

2. Let $\cU:=p^{-1}\left(\cF'_{2\epsilon} \times B(\epsilon)\right)$.
The restriction $p|_{\cU}$ is a (ramified) covering of degree
\begin{equation}
\label{eq:prod:a}
M_0:=\prod_{i=1}^p (a_i+1)
\end{equation}
where $a_i=a_i'+a_i''$ is the degree of the distinguished zero
on the surface $S'_i$.

3. The map $p|_{\cU}$ is volume preserving:
the volume element on $\cU$ induced from
$d\vol'\!\times\, d\gamma$ on
$\cF'_{2\epsilon} \times B(\epsilon)$ coincides with
the volume element
$d\vol$ on $\cU\subset\cH_1(\alpha)$.

4. The set $\cU$ contains the set
$\cH_1^{\epsilon,3\epsilon}(\alpha,\cC)$.
\end{lemma}
\begin{proof}
In fact, the statement of the Lemma is almost
tautological reformulation
in geometric terms
of the properties which
we already know.

Part 1 follows immediately from the construction
in the beginning of Section~\ref{ss:building:surfaces},
where having a surface
$S'\in \cH_1(\alpha') - \cH_1^{2\epsilon}(\alpha')$
and a vector $\gamma\in B(\epsilon)$ we constructed a surface
$S\in\cH_1^{\epsilon,\epsilon}(\alpha,\cC)$. We have $p(S)=S'$.

The fact that $p$ is nondegenerate almost everywhere follows
immediately from consideration of $p$ in cohomological coordinates,
so in this sense it is a (ramified) covering. To prove statement 2
it remains to compute the degree of this covering, i.e. to prove
that for almost all pairs $(S',\gamma)$ we can construct
exactly $M_0$ pairwise nonisometric surfaces.
For each surface $S'_i$ there  are  $a_i+1$ ways to
perform  a  slit
construction  at   a   zero   of   order   $a_i=a'_i+a''_i$ of $S'_i$
providing $a_i+1$ surfaces $S_i$. According to Lemma~\ref{lm:breaking:zero}
for almost all $S'\in \cF'_{2\epsilon}$ the resulting $a_i+1$ surfaces $S_i$
are pairwise nonisometric.
For almost all $S'$
the surfaces $S_i$ and $S_j$ corresponding
to different indices $i\ne j$ are nonisometric as well.
Thus generically there are no isomorphisms between the resulting surfaces $S$
which map the collection $\{\gamma_i\}$ on one surface to the collection $\{\gamma_i\}$
on the other surface. Since by the choice
of $(S',\gamma)$ such a collection is unique on each of the
$M_0$ surfaces $S$, they are generically
pairwise nonisometric.

Statement 3 follows immediately from the fact that as a basis
of cycles in $H_1(S,\{P_1, \dots, P_k\};\Z)$ we can choose
a union of basic cycles for $S'_i$ and one of the homologous
saddle connections $\gamma_i$, where by definition
$hol(\gamma_i)=\gamma$. Thus,
this statement is completely analogous to the
corresponding statement in multiplicity one.

Statement 4 follows from the fact that the surgery, which associates
a surface $S'$ to a surface $S$ either does not change the holonomy of
saddle connections different from $\gamma_i$ or changes it by holonomy
of $\gamma_i$. Since $|\gamma_i|<\epsilon$ it means that
a saddle connection, which is longer than $3\epsilon$ would stay
longer than $2\epsilon$ after the surgery.
\end{proof}

\begin{corollary}
\label{cr:volume:thick:in:problem1}
Suppose that the stratum $\cH(\alpha)$ is connected;
let $\cC$ be an admissible configuration of saddle connections
joining a pair of distinct zeroes; let
$(S',\gamma,a_i',a_i'')\to
(S,m_1,m_2,a_i',a_i'')$ be the corresponding
assignment. Then
$$
\Vol(\cH_1^{\epsilon}(\alpha,\cC)) =
M\cdot\Vol(\cH_1(\alpha'))\cdot \pi\epsilon^2 + o(\epsilon^2),
$$
where the integer constant $M$ is defined as
\begin{equation}
\label{eq:factor:M}
M=\cfrac{\prod_{i=1}^p (a_i+1)}{|\Gamma_{-}|\cdot|\Gamma|}\
\end{equation}
Here $a_i=a_i'+a_i''$; $|\Gamma_{-}|=2$ if the assignment has
$\gamma\to-\gamma$ symmetry and $|\Gamma_{-}|=1$ otherwise;
$|\Gamma|$ is the order of the ``stratum interchange'' symmetry.
\end{corollary}
\begin{proof}
Combining statement 4 of the Lemma above with the obvious
inclusions~\eqref{eq:epsilon:3epsilon} we obtain the following
inclusions
$$
\cH_1^{\epsilon,thick}(\alpha,\cC) =
\cH_1^{\epsilon,3\epsilon}(\alpha,\cC) \subset
\cU \subset
\cH_1^{\epsilon}(\alpha,\cC)
$$
Hence, by Corollary~\ref{cor:thick:almost:equals:all}
$$
\Vol\left(\cH_1^{\epsilon}(\alpha,\cC)\right)=
\Vol(\cU)+o(\epsilon^2).
$$
On the other hand, by statements 2 and 3 of the Lemma above,
$$
\Vol(\cU)=M_0\cdot\Vol(\cF'_{2\epsilon})\cdot \Vol(B(\epsilon))
$$
The volume of the disc $B(\epsilon)$ equals $\pi \epsilon^2$. The volume
of the quotient over a finite group of isometries equals to
the volume of the total space divided by the order of the group.
Since $|\Gamma_\pm|=|\Gamma_-|\cdot|\Gamma|$, and
the set $\cF'_{2\epsilon}$ is defined as
$$
\cF'_{2\epsilon}:=\left(\cH_1(\alpha') - \cH_1^{2\epsilon}(\alpha')\right)/\Gamma_\pm
$$
we obtain
$$
\Vol(\cF'_{2\epsilon})=
\left(\Vol(\cH_1(\alpha'))-\Vol(\cH_1^{2\epsilon}(\alpha',\cC))\right)/
(|\Gamma_-|\cdot|\Gamma|)
$$
It follows from Lemma~\ref{lemma:short:saddle:connections}
that $\Vol(\cH_1^{2\epsilon}(\alpha'))=O(\epsilon^2)$. Hence,
$$
\Vol(\cF'_{2\epsilon})=
\Vol(\cH_1(\alpha'))/
(|\Gamma_-|\cdot|\Gamma|) + O(\epsilon^2)
$$
Summarizing the above arguments we obtain the statement of the
corollary.
\end{proof}

Now everything is ready to compute the
Siegel---Veech constant $c(\cC)$. Note that we work
in the setting where
we  fix  the  partition  of  numbered  zeroes of $S$  by
components     $S'_i$     (see    Problem     1    in
section~\ref{ss:stratum:interchange}).  We  assume   that
$\cH_1(\alpha)$ is {\em connected}.
Applying Proposition~\ref{prop:thick} we get the following
value of the $c(\cC)$:
$$
c(\cC)=M\cdot\cfrac{ \Vol(\cH_1(\alpha')) }{ \Vol(\cH_1(\alpha)) }
$$
Using expression~\eqref{eq:factor:M} for $M$ and
Formula~\eqref{eq:total:volume:of:nonprimitive:stratum}
from Section~\ref{ss:Strata:of:Disconnected:Surfaces}
for the volume of a nonprimitive stratum $\cH_1(\alpha')$
we finally obtain the following
formula:
$$
c = \frac{1}{|\Gamma|}\frac{1}{|\Gamma_-|} \cdot
\prod_{j=1}^p (a_j+1)  \cdot
\frac{1}{2^{p-1}}\cdot
\frac{ \prod_{i=1}^p (\tfrac{d_i}{2}-1)! }
     { (\frac{d}{2}-2)! } \cdot
\frac{ \prod_{i=1}^p \Vol( \cH_1(\alpha_i)) }
     { \Vol(\cH_1(\alpha)) }
$$
where
  %
   %
\begin{equation}
\label{eq:di}
d_i=
\dim_\reals\cH_1(\alpha'_i)
\end{equation}
Note, that if $a_i=0$, the corresponding surface $S'_i$ has a marked
point, so $\alpha'_i$ contains $0$. Note also, that by convention
when $S'_i$ is a torus, we use the point which is already marked,
so $\alpha'_i = {0}$.

$$
\ast\quad\ast\quad\ast
$$
\medskip

Let us study now the constant  $c$ in the setting where we do not
specify which zeroes  $S$ are to be joined. (see Problem
2 in Section~\ref{ss:stratum:interchange}). This counting problem
can be reduced  to  the previous  one  by a purely  combinatorial
computation. Consider  $m\neq  m_1$,  $m\neq  m_2$,  $m\neq a_i$,
$i=1,\dots,p$. Then these zeroes are  all  ``inherited''  by  the
surfaces $S_i$.  Since all the  zeroes $z_1, z_2, \dots, z_l$ are
named, the  number of ways to  distribute $o(m)$ zeroes  of order
$m$ into groups of $o_1(m), \dots, o_p(m)$ zeroes equals
$$ 
\frac{o(m)!}{\prod_{i=1}^p o_i(m)!}
$$ 
If $m\neq m_1$ and $m\neq m_2$  but $m=a_j$ for some $j$ then one
of the $o_j(m)$  zeroes  of  order $m$ which lives  on  $S_j$  is
newborn; that is,  arising from the degeneration, while the other
$o_j(m)-1$ zeroes of order $m$ come from the corresponding zeroes
of  order  $m$ on  $S$.  Thus  the  corresponding  factor  in the
denominator  becomes  $(o_j(m)-1)!$.  Multiplying  numerator  and
denominator by $o_j(m)$ we get the following factor:
$$ 
\frac{o(m)!}{\prod_{j=1}^p o_j(m)!} \cdot \prod_{j\ |\ a_j=m}
o_j(m)
$$ 

(Equivalently,  we  may  say  that  we  have  performed  the slit
construction at the first ``named'' zero of order $a_j$ on $S_j$.
The rest of the zeroes of  that order inherit their names on $S$,
but then  we count  the number of ways of  merging the names from
different $S_i$ onto $S$.)
Consider now the case when $m=m_1$ or $m=m_2$. When $m_1\neq m_2$
and $m=m_1$ (resp. $m=m_2$) we have the arrangements of the other
$o(m_1)-1$ (resp. $o(m_2)-1$) zeroes of order $m$. When $m_1=m_2$
we have  an arrangement of the  other $o(m_1)-2$ zeroes  of order
$m$.  This changes  the  numerator in the  formula  above by  the
corresponding factorial. However, in these cases we have to count
the number  of ways to choose a zero of degree $m_1$ and a zero of
degree $m_2$ on the initial surface $S$. This number is  equal to
$o(m_1)\cdot o(m_2)$  if  $m_1\neq m_2$ and $o(m_1)(o(m_1)-1)$ if
$m_1=m_2$. Thus  for zeroes of these  degrees we have  this extra
factor  $o(m_1)\cdot   o(m_2)$   which   combined  with  previous
considerations gives $o(m)!$ in this  case  as  well. Taking into
account        the         possible        symmetries        (see
section~\ref{ss:stratum:interchange}),   we   finally   get   the
following expression for the constant  $c$  in  the setting where
the zeroes of  $S$ are not  numbered (see Problem  2  in
section~\ref{ss:stratum:interchange}).

\begin{formula}
   %
        %
\begin{multline}
\label{eq:sad:conn:const}
c= \frac{1}{|\Gamma_-|} \cdot \frac{1}{|\Gamma|} \cdot
\prod_{m\in\alpha}
\left(\frac{o(m)!}{\prod_{j=1}^p
o_j(m)!}\right) \cdot \prod_{\substack{j=1\\a_j\neq 0}}^p o_j(a_j) \cdot \\
\cdot \prod_{j=1}^p (a_j+1)  \cdot
\frac{1}{2^{p-1}}\cdot
\frac{ \prod_{i=1}^p (\tfrac{d_i}{2}-1)! }
     { (\frac{d}{2}-2)! } \cdot
\frac{ \prod_{i=1}^p \Vol( \cH_1(\alpha_i)) }
     { \Vol(\cH_1(\alpha)) }
\end{multline}
\end{formula}

We give some examples of  this  formula in the case of  connected
strata $\cH_1(\alpha)$.

\begin{example}
{\bf Stratum} $\boldsymbol{\cH_1(3,1)}$;
   %
The multiplicity  is    $p=2$.
After degeneration we get a surface $S'_1$ of genus
$2$ with a single zero of order $a_1=2$, where now  $a'_1=2,
a''_1=0$, and we get a
torus $S'_2$, with $a_2=a'_2=a''_2=0$.

The saddle  connections  $\gamma_1$  and $\gamma_2$ partition the
cone angle $8\pi$ at $z_1$ in two sectors with   angles $2\pi$
and $6\pi$; they partition the cone angle $4\pi$ at $z_2$  in two
sectors with angles $2\pi$ and $2\pi$.
There  is  no  stratum  interchange,  so  $|\Gamma|=1$.  We get
$$ 
c=3 \cdot \frac{1}{2} \cdot \frac{3!\cdot 1!}{5!}
\frac{\Vol(\cH_1(2))\cdot\Vol(\cH_1(\torusemptyset))}{\Vol(\cH_1(3,1))}
=
\frac{567}{1024} \approx 0.554
$$ 
\end{example}

\begin{example}
{\bf    Stratum}    $\boldsymbol{\cH_1(2,1,1)}$;    genus
$g=3$,
multiplicity $p=2$.
First  suppose  that  we   have   a  pair  of  homologous  saddle
connections joining  the zero of degree  $m_1=2$ with one  of two
zeroes of degree  $m_2=1$. After degeneration we get two surfaces
$S'_1\in\cH_1(a_1,1)$,                              and
$S'_2\in\cH_1(a_2)$,  where  $a_1+a_2= m_1+m_2+2-2p=1$.
Thus $a_1=1$, $a_2=0$.
There
is no  stratum  interchange,  so  $|\Gamma|=1$. Altogether we get
$$ 
c=4 \cdot \frac{1}{2} \cdot \frac{4!\cdot 1!}{6!}\cdot
\frac{\Vol(\cH_1(1,1))\Vol(\cH_1(\torusemptyset))}{\Vol(\cH_1(2,1,1))}
=
\frac{28}{45} \approx 0.622
$$ 

Suppose now that we have a pair of homologous saddle connections
joining  two  simple  zeroes.  After  degeneration  we  get  two
surfaces          $S'_1\in\cH(a_1,2)$,           and
$S'_2\in\cH(a_2)$, where $a_1+a_2=  m_1+m_1+2-2p=0$.
Thus $S'$ is genus
 $2$ with  a marked  point and with  a zero  of
degree $2$, while $S'_2$ is  a  torus  with a marked
point. Here we have $a'_i=a''_i=0$, $i=1,2$.

There  is no stratum interchange, so $|\Gamma|=1$. We now get
$$ 
c=\frac{1}{2} \cdot 1 \cdot \frac{4!\cdot 1!}{6!}\cdot
\frac{\Vol(\cH_1(2))\Vol(\cH_1(\torusemptyset))}{\Vol(\cH_1(2,1,1))} =
\frac{7}{40} = 0.175
$$ 

\end{example}

\subsection{Principal Stratum}

The cone angle  at  any simple zero equals  $4\pi$.  Thus for the
principal  stratum  we  may  have at most two  homologous  saddle
connections  joining  a  pair  of distinct zeroes. Hence  in  the
higher        multiplicity        case,        the        surface
$S,\in\cH_1(1,\dots,1)$ of  genus $g$ degenerates to a pair of
surfaces $S'_1, S'_2$ of  genera  $g_1+g_2=g$.  Each of $S'_1,
S'_2$ has only simple zeroes and one additional  marked point.
We    always    have    the
$\gamma\to-\gamma$ symmetry, so $|\Gamma_-|=2$. When $g_1=g_2$ we
have a  symmetry  $|\Gamma|=2$  due  to  the stratum interchange;
there is no other  symmetry  otherwise.

   %
   %
   %

Applying    formula~\eqref{eq:sad:conn:const} and  ~\eqref{eq:di}
we   finally  get  the
following answer.

\begin{formula}
 The constant for the number of   saddle  connections  of
multiplicity  two  joining  two  distinct  zeroes in the principal
stratum  is
\end{formula}
\begin{equation*}
\begin{split}
c & =\frac{1}{2|\Gamma|}\cdot
\frac{(2g-2)!}{(2g_1-2)! (2g_2-2)!}
\cdot 1 \cdot 1 \cdot \frac{1}{2^{2-1}}\cdot
\frac{(4g_1-3)!(4g_2-3)!}{(4g-5)!} \cdot \\
& \qquad\qquad\qquad
\cdot \frac{\Vol(\cH_1(\overbrace{1,\dots,1}^{2g_1-2}))
      \Vol(\cH_1(\overbrace{1,\dots,1}^{2g_2-2}))}
{\Vol(\cH_1(\underbrace{1,\dots,1}_{2g-2}))} =\\
& = \frac{1}{4|\Gamma|} \cdot
\frac{(2g-2)!\,(4g_1-3)!\,(4g_2-3)!}{(2g_1-2)!\,(2g_2-2)!\,(4g-5)!}
\cdot\frac{\Vol(\cH_1(\overbrace{1,\dots,1}^{2g_1-2}))
      \Vol(\cH_1(\overbrace{1,\dots,1}^{2g_2-2}))}
{\Vol(\cH_1(\underbrace{1,\dots,1}_{2g-2}))}
\end{split}
\end{equation*}
$$
\text{where}\quad g_1+g_2=g,\ g_1,g_2\ge 1,\qquad |\Gamma|=
\begin{cases}
  1&\text{when }g_1\neq g_2\\
  2&\text{when } g_1=g_2
\end{cases}
$$

\begin{figure}[hb!]
%
%
%
%
\includegraphics{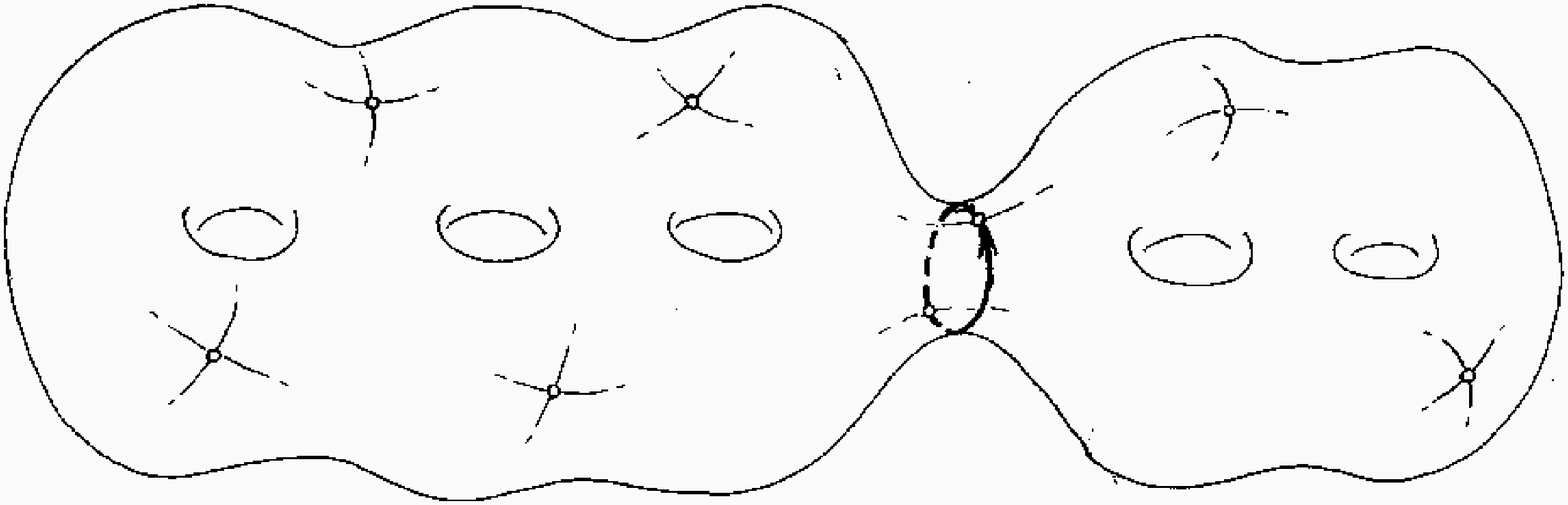}
\vspace{90bp} 
\caption{
\label{pic:h1111:mult2} Saddle connection of multiplicity 2 on a
surfaces from the principal stratum $\cH_1(1,\dots,1)$. }
\end{figure}

Table~\ref{tab:c:distinct:principal:mult2:exact}   presents   the
numerical values for  these  constants. More values are presented
in Table~\ref{tab:c:distinct:principal:mult2:approx} in the
Appendix. (We   present   only   approximate   values   of
constants   in Table~\ref{tab:c:distinct:principal:mult2:approx}
just  to  save the space.)

\begin{table}[htb!]
\small
\caption{Principal stratum;  values  of  the constants for saddle
connections  of  multiplicity two  joining  a  pair  of  distinct
zeroes}
\label{tab:c:distinct:principal:mult2:exact}

\bigskip
\scriptsize

$$ 
\begin{array}{|l|c|c|c|c|c|}
\hline &&&&&\\
& g_2=1 & g_2=2 & g_2=3 & g_2=4 & g_2=5 \\
[-\halfbls] &&&&& \\ \hline &&&&& \\ [-\halfbls]
 g_1=1 & \cfrac{5}{8}  & \cfrac{6}{7}  & \cfrac{315}{377}  & \cfrac{19604}{23357}  & \cfrac{13897415}{16493303} \\
 [-\halfbls] &&&&& \\ \hline &&&&& \\ [-\halfbls]
 g_1=2 & \cfrac{6}{7}  & \cfrac{30}{377}  & \cfrac{1680}{23357}  & \cfrac{686140}{16493303}  & \cfrac{9529656}{351964697} \\
 [-\halfbls] &&&&& \\ \hline &&&&& \\ [-\halfbls]
 g_1=3 & \cfrac{315}{377}  & \cfrac{1680}{23357}  & \cfrac{154350}{16493303}  & \cfrac{352872}{50280671}  & \cfrac{389127620}{121780991823} \\
 [-\halfbls] &&&&& \\ \hline &&&&& \\ [-\halfbls]
 g_1=4 & \cfrac{19604}{23357}  & \cfrac{686140}{16493303}  & \cfrac{352872}{50280671}  & \cfrac{336277214}{365342975469}  & \cfrac{435911877856}{704782787198207} \\
 [-\halfbls] &&&&& \\ \hline &&&&& \\ [-\halfbls]
 g_1=5 & \cfrac{13897415}{16493303}  & \cfrac{9529656}{351964697}  & \cfrac{389127620}{121780991823}  & \cfrac{435911877856}{704782787198207}  & \cfrac{19865637635886}{249020069093788675} \\ [-\halfbls]
&&&&&\\ \hline
\end{array}
$$ 

\end{table}

\section{Strata Having Several Connected Components}
\label{s:distinct:zeroes:nonconnected:strata}

The  analysis  of  admissible  assignments,  structure  of  local
surgeries,  and  most of the elements of  our  constructions  are
basically  the  same for  the  strata  having  several  connected
components. However, we need to establish the correspondence
between connected components and parts of the boundary of the
stratum
they are adjacent to. In other words we need to describe the
admissible constructions, which lead to the flat surfaces in
the prescribe connected component.

\subsection{Parity of the Spin Structure in the Slit Construction}
\label{ss:parity:in:slit}

We start with outlining the relation  between  the
parities  of  spin structures of the resulting  surface  and  its
components in the slit construction.

\begin{lemma}
\label{lm:spin:in:slit:constr}
Suppose  $S$  constructed above has only even zeroes  so
has   a    spin    structure.   Then   each   component   surface
$S_i'$ has  only even zeroes and  so also has  a spin
structure. The parity of the  spin  structure  $\phi(S)$  is
equal to  the sum  of the parities of the  spin structures of the
components:
\begin{equation}
\label{eq:parities:in:slits}
\phi(S)=\sum_{i=1}^p \phi(S_i')\ (mod\ 2)
\end{equation}
\end{lemma}

\begin{proof}
For each $S'_i$, all the zeroes except the distinguished one come
from the corresponding zeroes on  $S$  unchanged;  in  particular
they have the same degrees as the initial zeroes of $S$. Thus for
each $S'_i$ all the zeroes  but  at most one, have even  degrees.
Note that the total sum  of  the degrees  of the zeroes of $S_i'$
is
even: it equals $2g_i+2$, where $g_i$ is the genus of the surface
$S'_i$. Hence the  remaining zero has  even degree as  well.  (By
usual convention we consider a marked point as a ``zero of degree
$0$''.)

For each surface $S'_i$ consider  a  collection  of $2g_i$ smooth
simple closed curves  on  $S'_i$ representing the canonical basis
of cycles. Breaking up the zero  $z_i$ we make a local surgery of
the flat structure,  i.e., we do  not change the  flat  structure
outside of a  small domain $U'_i\subset S'_i$ containing the zero
$z_i$. Thus for each $i$  we  can deform if necessary the  curves
from the  corresponding collection in such  a way that  they stay
outside of  a neighborhood of the zero $z_i$.  Thus we may assume
that for  each $i=1,\dots,p$ none  of the chosen curves on $S'_i$
intersect $U'_i$. This means that  all  the  curves survive under
the surgery, moreover, the union of collections for all $i$ gives
us a canonical  basis  on $S$. By construction  the  index of any
resulting curve in the flat structure $\omega$ is the same as the
index   of   the  corresponding  curve  in  the  flat   structure
$\omega'_i$.       Hence       we      get       the      desired
relation~\eqref{eq:parities:in:slits}.
        %
        %
The lemma is proved.
\end{proof}

\begin{lemma}
\label{lm:even:and:odd:assignments}
Let  $S$  belong to one of the    components
$\cH_1^{odd}(\alpha)$  or   $\cH_1^{even}(\alpha)$. Suppose that
$\alpha\neq
(g-1,g-1)$.  Then the  surface $S$  can  be  obtained  by  an
assignment                     $(S',\gamma,a_i',a_i'')\to
(S,m_1,m_2,a_i',a_i'')$, where  $m_1,m_2\in\alpha$, if and
only if the  collection  $(\alpha'_i, a'_i, a''_i)$ satisfies the
conditions  of  Lemma~\ref{lm:admissible:split}, all  $\alpha'_i$
are  even,  and  the  parities  of   the   spin   structures   of
$S'_i$      satisfy      the       condition       of
Lemma~\ref{lm:spin:in:slit:constr}.
\end{lemma}

\subsection{Admissible Assignments for Surfaces from Hyperelliptic Components}

Note  that   the   surfaces   $S'_i$  from  the  {\it
hyperelliptic}  component  $\cH_1^{hyp}(2g-2)$ and  from  the
{\it
hyperelliptic} component $\cH_1^{hyp}(g-1,g-1)$, where $g$ is
odd,
also    have     parity    of    the    spin    structure    (see
formulae~\eqref{eq:spin:2g:minus2},   \eqref{eq:spin:hyp}).    In
general these  surfaces  are  also  involved  in assignments from
Lemma~\ref{lm:even:and:odd:assignments}      producing       {\it
nonhyperelliptic}     components     $\cH_1^{even}(\alpha)$
and
$\cH_1^{odd}(\alpha)$.  However,   in   the   exceptional  case
of
$\alpha=(g-1,g-1)$ with odd $g$ some assignments lead to surfaces
$S$   from   the   hyperelliptic   connected   component
$\cH_1^{hyp}(g-1,g-1)$.  These  assignments  are classified in
the
lemma below.

\begin{lemma}
\label{lm:slit:hyp}
Let a flat surface $S$ from  the hyperelliptic connected
component $\cH_1^{hyp}(g-1,g-1)$ be obtained by an assignment
$$
(S',\gamma,a_i',a_i'')\to (S,g-1,g-1,a_i',a_i'').
$$
Then for almost every $S$ the assignment has  one of the
following two types:
\begin{itemize}
\item[---]  The   multiplicity   is   $1$;   the   flat  surface
$S'$   belongs   to    the   hyperelliptic   component
$\cH_1^{hyp}(2g-2)$, and $a'=a''=g-1$.
\item[---]  The   multiplicity   is   $2$;   the  flat  surfaces
$S'_i$,   $i=1,2$  belong   to   the   hyperelliptic
components   $\cH^{hyp}(2g_i-2)$,   where    $g_1+g_2=g$,    and
$a'_i=a''_i=g_i-1$, $i=1,2$.
\end{itemize}
The
assignment is $|\Gamma|$ to $1$ and is onto the corresponding set
of  hyperelliptic  surfaces.
\end{lemma}
\begin{proof}
The set  of  hyperelliptic  structures  such  that the involution
interchanges some components has positive codimension,  so we may
disregard these and  assume  the involution fixes each component.
This implies that all flat surfaces  $S'_i$ belong to
components   $\cH^{hyp}(2g_i-2)$.   Since    the    hyperelliptic
involution sends a vector $v$  to  $-v$, but does not change  the
orientation of the surface,  if  the multiplicity is greater than
one, it cannot  map  any of  the  $\gamma_i'$ or $\gamma_i''$  to
themselves.  It  must  therefore   interchange   $\gamma_i'$  and
$\gamma_i''$ for each $i$. This implies that $p=2$.
\end{proof}

\begin{corollary}
\label{cr:no:high:mult:in:hyp}
Almost any  flat  surface  $S$  from  the  hyperelliptic
connected   component   $\cH_1^{hyp}(g-1,g-1)$    has   no
saddle
connections of multiplicity greater then two.
\end{corollary}

\subsection{Constants for the Components of the Stratum $\cH_1(g-1,g-1)$}

We start with  the  stratum $\cH_1(g-1,g-1)$  which  is
exceptional:
%
it has the hyperelliptic connected component.

\subsubsection{Hyperelliptic Connected Component $\cH_1^{hyp}(g-1,g-1)$}

The  admissible   assignments   for   the   surfaces   from   the
hyperelliptic   connected   component  $\cH_1^{hyp}(g-1,g-1)$
are
described by Lemma~\ref{lm:slit:hyp}. The possible multiplicities
are      one      and     two.      In      multiplicity      one
modifying~\eqref{eq:sad:conn:mult:1:m1:m1}  in  accordance   with
Lemma~\ref{lm:slit:hyp} we get the following formula.
\begin{formula}
\label{eq:dist:sad:hyp:mult:1}
   %
The constant for the  number  of
saddle connections  of  multiplicity  one joining distinct zeroes
in  $\cH_1^{hyp}(g-1,g-1)$ is given by
\end{formula}
$$ 
c = (2g-1)\cdot
\frac{\Vol(\cH_1^{hyp}(2g-2))}{\Vol(\cH_1^{hyp}(g-1,g-1))}
$$ 

In  other  words, collapsing a saddle connection of  multiplicity
one on  a  surface  $S\in\cH_1^{hyp}(g-1,g-1)$  we  get  a
surface $S'$ from the hyperelliptic connected component
$\cH_1^{hyp}(2g-2)$.  Breaking   up   the   zero   of   a
surface
$S'\in\cH_1^{hyp}(2g-2)$  we necessarily  get  a  surface
$S$      from      the      hyperelliptic      component
$\cH_1^{hyp}(g-1,g-1)$.

\begin{table}[htb!]
\small
\caption{\small  Hyperelliptic  component
$\cH_1^{hyp}(g-1,g-1)$; values  of  the constants for saddle
connections of  multiplicity one joining distinct zeroes}
\label{tab:c:distinct:hyp:mult1}
   %
\scriptsize
$$ 
\begin{array}{|l|c|c|c|c|c|c|c|c|}
\hline &&&&&&&&\\ & g=2 & g=3 & g=4 & g=5 & g=6 & g=7 & g=8
& g=9  \\ [-\halfbls] &&&&&&&& \\ \hline &&&&&&&& \\
[-\halfbls]
 c = & \cfrac{27}{8}  & \cfrac{225}{32}  & \cfrac{6125}{512}  & \cfrac{297675}{16384}  & \cfrac{3361743}{131072}  & \cfrac{9018009}{262144}  & \cfrac{372683025}{8388608}  & \cfrac{59836330125}{1073741824}
 \\  [-\halfbls] &&&&&&&& \\ \hline &&&&&&&& \\ [-\halfbls]
 c \approx & 3.375 & 7.03125 & 11.9629 & 18.1686 & 25.6481 & 34.401 & 44.4273 & 55.7269 \\
 [-\halfbls]&&&&&&&&\\ \hline
\end{array}
$$ 
\end{table}

\begin{example}{\bf Component} $\boldsymbol{\cH_1^{hyp}(1,1)}$.
In genus $g=2$ the component $\cH_1^{hyp}(1,1)$  coincides with
the
whole principal stratum  $\cH_1(1,1)$.  The formula above gives
the
same constant as formula~\ref{f:dist:zeroes:principal:mult:1} for
the principal stratum in genus $g=2$:
$$ 
c = 3\cdot
\frac{\Vol(\cH_1^{hyp}(2))}{\Vol(\cH_1^{hyp}(1,1))}=
3\cdot
\frac{\Vol(\cH_1(2))}{\Vol(\cH_1(1,1))} = \cfrac{27}{8}=3.375
$$ 
\end{example}

\begin{example}
{\bf Component} $\boldsymbol{\cH_1^{hyp}(2,2)}$.
In  genus  $3$ we have the connected component
$\cH_1^{hyp}(2,2)$.
After  collapsing   zeroes  we  obtain  a  hyperelliptic  surface
$S'\in \cH_1^{hyp}(4)$.
   %
%
%
We get
$$ 
c=5\cdot\frac{\Vol(\cH_1^{hyp}(4))}{\Vol(\cH_1^{hyp}(2,2))}=
\frac{225}{32}\approx 7.03
$$ 
\end{example}

We  present  more values  in
table~\ref{tab:c:distinct:hyp:mult1}.

Consider  now  multiple  saddle  connections for a  flat
surface $S$ from  the connected component
$\cH_1^{hyp}(g-1,g-1)$. By  Lemma~\ref{lm:slit:hyp}  we  only
have $p=2$ in  the  higher multiplicity,     and
$S'_i\in\cH^{hyp}(2g_i-2)$, $i=1,2$,   where
$g_1+g_2=g$.   Here   we   always   have   the
$\gamma\to-\gamma$ symmetry, so $  |\Gamma_-|=2$.  When
$g_1=g_2$ we also have  the  stratum interchange symmetry, so
$|\Gamma|=2$, and $|\Gamma|=1$  otherwise.
 Thus,   after
appropriate  modification   of
formula~\eqref{eq:sad:conn:const} we get the following
constant
\begin{formula}
\label{f:hyp:dist:mult2}
The constant in multiplicity $2$ for saddle connections joining
distinct zeroes   in  $\cH_1^{hyp}(g-1,g-1)$ is given by
   %
   %
   %
\end{formula}
%
        %
$$ 
c= \frac{(2g_1-1)(2g_2-1)}{2|\Gamma|} \cdot
\frac{(2g_1-1)!(2g_2-1)!}{(2g-1)!} \cdot
\frac{ \Vol(\cH_1^{hyp}(2g_1-2))\! \Vol(\cH_1^{hyp}(2g_2-2))
}
     { \Vol(\cH_1^{hyp}(g-1,g-1)) }
$$ 
     %

\begin{example}
{\bf Component} $\boldsymbol{\cH_1^{hyp}(1,1)}$.
Saddle connections of multiplicity two
joining distinct zeroes on
a surface from   $\cH_1^{hyp}(1,1)=\cH_1(1,1)$
were actually, already considered in the section treating the
principal
stratum. Here after decomposition
we      get     a pair of tori $S'_1,S'_2$  each with  a
marked point.
We get the
following value for the corresponding constant:
$$ 
c=\frac{1 \cdot 1}{2\cdot 2} \cdot \frac{1!\cdot 1!}{3!}\cdot
\frac{\Vol(\cH_1(\torusemptyset))\Vol(\cH_1(\torusemptyset))}{\Vol(\cH_1(1,1))}
=
\frac{5}{8} = 0.625
$$ 
\end{example}

\begin{table}[htb]
\small
\caption{\small  Hyperelliptic  component
$\cH_1^{hyp}(g-1,g-1)$;
values  of  the constants for saddle connections of  multiplicity
two joining distinct zeroes}
\label{tab:c:distinct:hyp:mult2}
\scriptsize
$$ 
\begin{array}{|l|c|c|c|c|c|c|}
\hline &&&&&&\\
& g_2=1 & g_2=2 & g_2=3 & g_2=4 & g_2=5 & g_2=6 \\
[-\halfbls] &&&&&& \\ \hline &&&&&& \\ [-\halfbls]
 g_1=1 & \cfrac{5}{8}  & \cfrac{63}{32}  & \cfrac{375}{128}  & \cfrac{67375}{16384}  & \cfrac{361179}{65536}  & \cfrac{1867635}{262144} \\
 [-\halfbls] &&&&&& \\ \hline &&&&&& \\ [-\halfbls]
 g_1=2 & \cfrac{63}{32}  & \cfrac{567}{512}  & \cfrac{22275}{8192}  & \cfrac{441441}{131072}  & \cfrac{1083537}{262144}  & \cfrac{10496871}{2097152} \\
 [-\halfbls] &&&&&& \\ \hline &&&&&& \\ [-\halfbls]
 g_1=3 & \cfrac{375}{128}  & \cfrac{22275}{8192}  & \cfrac{96525}{65536}  & \cfrac{875875}{262144}  & \cfrac{4027725}{1048576}  & \cfrac{1187146125}{268435456} \\
 [-\halfbls] &&&&&& \\ \hline &&&&&& \\ [-\halfbls]
 g_1=4 & \cfrac{67375}{16384}  & \cfrac{441441}{131072}  & \cfrac{875875}{262144}  & \cfrac{14889875}{8388608}  & \cfrac{2083228875}{536870912}  & \cfrac{18432519875}{4294967296} \\
 [-\halfbls] &&&&&& \\ \hline &&&&&& \\ [-\halfbls]
 g_1=5 & \cfrac{361179}{65536}  & \cfrac{1083537}{262144}  & \cfrac{4027725}{1048576}  & \cfrac{2083228875}{536870912}  & \cfrac{8749561275}{4294967296}  & \cfrac{74923166241}{17179869184} \\
 [-\halfbls] &&&&&& \\ \hline &&&&&& \\ [-\halfbls]
 g_1=6 & \cfrac{1867635}{262144}  & \cfrac{10496871}{2097152}  & \cfrac{1187146125}{268435456}  & \cfrac{18432519875}{4294967296}  & \cfrac{74923166241}{17179869184}  & \cfrac{624359718675}{274877906944} \\ [-\halfbls]
&&&&&&\\ \hline
\end{array}
$$ 
\end{table}

\begin{example}
{\bf Component} $\boldsymbol{\cH_1^{hyp}(2,2)}$.
Consider saddle  connections of multiplicity two joining distinct
zeroes on a  surface from $\cH_1^{hyp}(2,2)$.
After decomposition we get a  surface $S'_1\in\cH_1(2)$ of
genus  $2$  with  a  single  zero  of  order  $2$  and   a  torus
$S'_2\in\cH(\torusemptyset)$ with  a  marked  point.  The  cone  angle
between $\gamma_1$ and $\gamma_2$ is the same at $z_1$ and $z_2$;
it  equals  $2\pi$.   The   cone  angle  between  $\gamma_2$  and
$\gamma_1$ is also the same at $z_1$ and $z_2$, but now it equals
$4\pi$. The constant is given by
$$ 
c=\frac{3 \cdot 1}{2\cdot 1} \cdot \frac{3!\cdot 1!}{5!}\cdot
\frac{\Vol(\cH_1(2))\Vol(\cH_1(\torusemptyset))}{\Vol(\cH_1^{hyp}(2,2))}
=
\frac{63}{32} \approx 1.97
$$ 
\end{example}

We present more values in table~\ref{tab:c:distinct:hyp:mult2} and
in table~\ref{tab:c:dist:hyp:mult2} in the Appendix.

\subsubsection{Connected Component $\cH_1^{nonhyp}(g-1,g-1)$; Even Genus}

For $g\ge 4$, $g$ even,  the  stratum  $\cH_1(g-1,g-1)$ has
exactly
two     connected     components:    $\cH_1^{hyp}(g-1,g-1)$
and
$\cH_1^{nonhyp}(g-1,g-1)$. Thus,  dealing with the connected
component
$\cH_1^{nonhyp}(g-1,g-1)$ we  have  to exclude all the assignments
producing      surfaces   $S$   in   the   stratum
$\cH_1(g-1,g-1)$ those that actually   belong  to  the
component
$\cH_1^{hyp}(g-1,g-1)$. The  latter problem has  already been
solved. Thus
we get the following answers.
\begin{formula}
\label{eq:dist:sad:nonhyp:mult:1}
For any  even genus $g\ge  4$ any flat surfaces $S$ from the
connected component $\cH_1^{nonhyp}(g-1,g-1)$  having  a saddle
connection of multiplicity one joining  distinct  zeroes  may be
obtained from  the  corresponding surface $S'$ from one of the
components $\cH_1^{even}(2g-2)$, $\cH_1^{odd}(2g-2)$.

The corresponding constant is equal to
\end{formula}
\begin{multline*}
c = (2g-1)\cdot
\frac{\Vol(\cH_1^{even}(2g-2))+\Vol(\cH_1^{odd}(2g-2))}{\Vol(\cH_1^{nonhyp}(g-1,g-1))}=
\\
= (2g-1)\ \cdot\ \frac{\Vol(\cH_1(2g-2))-\Vol(\cH_1^{hyp}(2g-2))}
{\Vol(\cH_1(g-1,g-1))-\Vol(\cH_1^{hyp}(g-1,g-1))}
\end{multline*}

We proceed with  multiplicity  two. Here  we get the  formula
analogous  to formula~\ref{f:hyp:dist:mult2}  for  the
component $\cH_1^{hyp}(g-1,g-1)$.
\begin{formula}
\label{f:nonhyp:dist:mult2}
%
The constant for multiplicity $2$ joining two zeroes in the
stratum $\cH^{nonhyp}(g-1,g-1)$
is given by
   %
   %
   %
\end{formula}
%

  %
\begin{multline*}
c= \frac{(2g_1-1)(2g_2-1)}{2|\Gamma|}\cdot\frac{(2g_1-1)!(2g_2-1)!}{(2g-1)!} \cdot \\
\cdot \frac{\Vol(\cH_1(2g_1-2)) \Vol(\cH_1(2g_2-2))
- \Vol(\cH_1^{hyp}(2g_1-2))
 \Vol(\cH_1^{hyp}(2g_2-2))}
     { \Vol(\cH_1^{nonhyp}(g-1,g-1)) },
\end{multline*}
  %

For multiplicity $p\geq 3$ it is  easy to account  for the
impact      of      the     hyperelliptic      component     (see
Corollary~\ref{cr:no:high:mult:in:hyp}).  We  get  the  following
particular case of formula(~\ref{eq:sad:conn:const}) in the case of
 $\cH^{nonhyp}(g-1,g-1)$
\begin{formula}
   %
%
  %
$$ 
c= \frac{1}{2|\Gamma|} \cdot \prod_{i=1}^p (2g_i-1) \cdot
\frac{\prod_{i=1}^p (2g_i-1)!}{(2g-1)!} \cdot
\frac{ \prod_{i=1}^p \Vol( \cH_1(2g_i-2)) }
     { \Vol(\cH^{nonhyp}_1(g-1,g-1)) }
$$ 
\end{formula}


\subsubsection{Components $\cH^{even}(g-1,g-1)$ and $\cH^{odd}(g-1,g-1)$}

We  consider  the  multiplicity  one  case   first.  Assume  that
$S$ belongs  to  one  of the nonhyperelliptic components
$\cH^{even}(g-1,g-1)$  or  $\cH^{odd}(g-1,g-1)$. Contracting the
 saddle connection  joining the two zeroes, we  merge the two zeroes into
one  and  we get a surface  $S'$  in
$\cH^{even}(2g-2)$  or  $\cH^{odd}(2g-2)$  correspondingly.
Conversely
 breaking up a single zero of a surface $S'$ from
a nonhyperelliptic  connected  component  $\cH^{even}(2g-2)$   or
$\cH^{odd}(2g-2)$  into  two  zeroes  of degrees $g-1$ we  get  a
surface     from      the component     $\cH^{even}(g-1,g-1)$
or
$\cH^{odd}(g-1,g-1)$  correspondingly.  Thus we  have  to  modify
formula~\eqref{eq:sad:conn:mult:1:m1:m1} in the following way.

\begin{formula}
   %

The constant is equal to
\end{formula}
$$ 
c = (2g-1)\cdot
\frac{\Vol(\cH^{even}_1(2g-2))}{\Vol(\cH^{even}_1(g-1,g-1))}
$$ 
$$ 
c = (2g-1)\cdot
\frac{\Vol(\cH^{odd}_1(2g-2))}{\Vol(\cH^{odd}_1(g-1,g-1))}
$$ 

\begin{example}
{\bf Component }  $\boldsymbol{\cH^{odd}(2,2).}$
   %
After collapsing  zeroes  we  obtain  a  surface $S'\in
\cH^{odd}(4)$. We get
$$ 
c=5\cdot\frac{\Vol(\cH_1^{odd}(4))}{\Vol(\cH_1^{odd}(2,2))}=
\frac{80}{9}\approx 8.89
$$ 
\end{example}

\begin{remark}
Note that  the similar nonhyperelliptic connected components with
{\it even} parity of the  spins  structure do not exist in  genus
$g=3$.
\end{remark}


Now we consider higher multiplicity.  We  again  have to consider
multiplicity   two   separately   because  some  assignments   in
multiplicity  two   produce   surfaces   from  the  hyperelliptic
connected        component        $\cH^{hyp}(g-1,g-1)$,       see
Lemma~\ref{lm:spin:in:slit:constr}.  We  again  get  the  formula
analogous  to formula~\ref{f:hyp:dist:mult2}  for  the  component
$\cH^{hyp}(g-1,g-1)$ and to formula~\ref{f:nonhyp:dist:mult2} for
the component $\cH^{nonhyp}(g-1,g-1)$.

It will be convenient  to  introduce the following notation. We
introduce a function $\delta(\alpha,\phi)$ which has values $0$ or
$1$. It is equal to $1$ when all of the following three conditions
are satisfied:
all the zeroes of $\alpha$ (if any) have even degrees;
$\cH(\alpha)$
contains the hyperelliptic component $\cH^{hyp}(\alpha)$;
the parity of the spin structure of $\cH^{hyp}(\alpha)$ coincides
with $\phi$. Otherwise $\delta(\alpha,\phi)=0$.

Formula~\eqref{eq:spin:2g:minus2}  for the parity of the spin
structure of $\cH^{hyp}(\alpha')$ shows that
\begin{equation}
\label{eq:delta:hyp}
 \delta(\alpha,\phi)=
\begin{cases}
1&\text{when } \alpha=(2g-2)\text{ and }\left[\cfrac{g+1}{2}\right] = \phi\ (\mod\ 2)\\
1&\text{when } \alpha=(g-1,g-1),\ g \text{ is odd},\\
 & \text{and }\left(\cfrac{g+1}{2}\right) = \phi\ (\mod\ 2)\\
0&\text{otherwise}
\end{cases}
\end{equation}

We  shall also  use  the following convention:  the  volume
of  a nonexistent component is equal to zero. In this
notation we get the following six ``dummy'' volumes for the
strata in small genera:
\begin{alignat}4
\label{eq:dummy:volumes}
&\Vol(\cH_1^{odd}(\torusemptyset)) =0\qquad  &\Vol(\cH_1^{odd}(2))
=0\qquad&\Vol(\cH_1^{even}(4))   &=0
\notag\\
&\Vol(\cH_1^{even}(\torusemptyset)) =0\qquad
&\Vol(\cH_1^{even}(2))=0\qquad  &\Vol(\cH_1^{even}(2,2)) &=0
\end{alignat}
and also we have
$$
\Vol(\cH^{hyp}(\alpha)) = 0\quad \text{if } \alpha\neq (2g-2),(g-1,g-1)
$$

\begin{formula}
\label{f:dist:even:odd:mult2}
A surface $S$ from one of the nonhyperelliptic connected
components $\cH^{even}(g-1,g-1)$ or $\cH^{odd}(g-1,g-1)$ may have
a saddle connection  of multiplicity two joining the two distinct
zeroes if and only if it can be obtained from surfaces $S_1'$ and
$S_2'$ of genus $g_1,g_2$ where $g_1+g_2=g$
   %
   %
and  where  the
following additional  requirements  are satisfied.  At  least
one  of  the surfaces $S'_1,S'_2$
does not belong  to the hyperelliptic  component
$\cH^{hyp}(2g_i-2)$,  $i=1,2$;  the parities  of  the   spin
structures  of  the  surfaces  satisfy
relation~\eqref{eq:parities:in:slits}:
$\phi(S'_1)+\phi(S'_2)=\phi(S)$.  For $S\in \cH^{odd}(g-1,g-1)$
the spin structures of $S_1',S_2'$ should have opposite parities.
This means in particular that there is no longer a stratum
interchange symmetry. For $g_1\neq g_2$
\end{formula}
\begin{equation*}
\begin{split}
c &=
\frac{(2g_1-1)(2g_2-1)}{2}\cdot\frac{(2g_1-1)!(2g_2-1)!}{(2g-1)!}\cdot
\frac{1}{ \Vol(\cH_1^{odd}(g-1,g-1)) } \cdot
   \\&
\Bigg(
\Big(\big(\Vol(\cH^{even}_1(2g_1-2)) +
\delta((2g_1-2),even)\Vol(\cH^{hyp}_1(2g_1-2))\big)\cdot
   \\&\phantom{\Bigg(}
\cdot\big(\Vol(\cH^{odd}_1(2g_2-2)) +
\delta((2g_2-2),odd)\Vol(\cH^{hyp}_1(2g_2-2))\big)\Big) +
   \\& \quad\
      +\Big(\big(\Vol(\cH^{odd}_1(2g_1-2)) + \delta((2g_1-2),odd)\Vol(\cH^{hyp}_1(2g_1-2))\big)\cdot
   \\&\ \quad\phantom{ + }
\cdot\big(\Vol(\cH^{even}_1(2g_2-2)) +
\delta((2g_2-2),even)\Vol(\cH^{hyp}_1(2g_2-2))\big)\Big) -
   \\&\quad
-\delta((2g_1-2),even)\delta((2g_2-2),odd)
\Vol(\cH^{hyp}_1(2g_1-2))\Vol(\cH^{hyp}_1(2g_2-2)) -
   \\&\quad
-\delta((2g_1-2),odd)\delta((2g_2-2),even)\Vol(\cH^{hyp}_1(2g_1-2))\Vol(\cH^{hyp}_1(2g_2-2))
\Bigg)
\end{split}
\end{equation*}
For $g_1=g_2=g/2$, where $g$ is even, we get
\begin{multline*}
c =
\frac{(2g_1-1)^2}{2}\cdot\frac{\big((2g_1-1)!\big)^2}{(2g-1)!}\cdot
\frac{1}{ \Vol(\cH_1^{odd}(g-1,g-1)) } \cdot
   \\
\cdot
\big(\Vol(\cH^{even}_1(2g_1-2)) +
\delta((2g_1-2),even)\Vol(\cH^{hyp}_1(2g_1-2))\big)\cdot
   \\
\cdot\big(\Vol(\cH^{odd}_1(2g_2-2)) +
\delta((2g_2-2),odd)\Vol(\cH^{hyp}_1(2g_2-2))\big)
\end{multline*}

\begin{example}
{\bf Component} $\boldsymbol{\cH^{odd}(2,2)}$.
   %
   %
From    Formula~\ref{f:dist:even:odd:mult2}    we    know
that the degeneration of multiplicity $2$ of  a  surface  from
$\cH(2,2)$  has   the  form
$\big(\overline{0}+\overline{0}\big) \succ
\big(\overline{1}+\overline{1}\big)  \succ$. Since  the
strata $\cH(\torusemptyset)$ and $\cH(2)$ are (hyper)elliptic
the resulting surface in $\cH(2,2)$ is also hyperelliptic, see
Lemma~\ref{lm:slit:hyp}. Hence degeneration of multiplicity
$2$ is not realizable for flat surfaces from
$\cH^{odd}(2,2)$.
\end{example}

The  expression  for  the  constant   for   the   surfaces   from
$\cH^{even}(g-1,g-1)$ is analogous to the one above with the only
exception  that  now we  get  the  stratum  interchange  symmetry
$|\Gamma|=2$ for $g_1=g_2$:
\begin{equation*}
\begin{split}
c &=
\frac{(2g_1-1)(2g_2-1)}{2|\Gamma|}\cdot\frac{(2g_1-1)!(2g_2-1)!}{(2g-1)!}\cdot
\frac{1}{ \Vol(\cH_1^{even}(g-1,g-1)) } \cdot
   \\&
\Bigg(
\Big(\big(\Vol(\cH^{even}_1(2g_1-2)) +
\delta((2g_1-2),even)\Vol(\cH^{hyp}_1(2g_1-2))\big)\cdot
   \\&\phantom{\Bigg(}
\cdot\big(\Vol(\cH^{even}_1(2g_2-2)) +
\delta((2g_2-2),even)\Vol(\cH^{hyp}_1(2g_2-2))\big)\Big) +
   \\&\quad
      +\Big(\big(\Vol(\cH^{odd}_1(2g_1-2)) + \delta((2g_1-2),odd)\Vol(\cH^{hyp}_1(2g_1-2))\big)\cdot
   \\& \quad\phantom{ + }
\cdot\big(\Vol(\cH^{odd}_1(2g_2-2)) +
\delta((2g_2-2),odd)\Vol(\cH^{hyp}_1(2g_2-2))\big)\Big) -
   \\&
-\delta((2g_1-2),even)\delta((2g_2-2),even)
\Vol(\cH^{hyp}_1(2g_1-2))\Vol(\cH^{hyp}_1(2g_2-2)) -
   \\&
-\delta((2g_1-2),odd)\delta((2g_2-2),odd)\Vol(\cH^{hyp}_1(2g_1-2))\Vol(\cH^{hyp}_1(2g_2-2))
\Bigg)
\end{split}
\end{equation*}
   %


Finally, for from multiplicity  at least three  it is easy to
account
for   the    impact   of   the   hyperelliptic   component   (see
Corollary~\ref{cr:no:high:mult:in:hyp}).  We  get  the  following
particular case of formula(~\ref{eq:sad:conn:const}):
\begin{formula}
Almost   every   surface   in    $\cH^{even}(g-1,g-1)$    or   in
$\cH^{odd}(g-1,g-1)$  having  saddle connections  of multiplicity
$p\ge 3$ joining two zeroes can be obtained by assignment
$$
\big(\overline{(g_1-1)})+\overline{(g_1-1)}\big) \succ \dots
\big(\overline{(g_p-1)}+\overline{(g_p-1)}\big) \succ
$$
where  $g_i\ge  1$,  $i=1,\dots,p$, and $g_1+\dots+g_p=g$.
The  corresponding  constant equals
\begin{multline*}
c=  \frac{1}{|\Gamma|} \cdot \prod_{i=1}^p (2g_i-1) \cdot
\frac{1}{2^{p-1}}\cdot
\frac{\prod_{i=1}^p (2g_i-1)!}{(2g-1)!} \cdot
\frac{1}{ \Vol(\cH^{\phi}_1(g-1,g-1)) } \cdot \\
\cdot
\sum_{\substack{\phi_1,...,\phi_p\in\{even, odd\}\\
\phi_1+\dots+\phi_p=\phi\, (mod\, 2)}}
\prod_{i=1}^p \Big(\Vol(\cH^{\phi_i}_1(2g_i-2)) +
\delta((2g_i-2),\phi_i)\Vol(\cH^{hyp}_1(2g_i-2))\Big)
\end{multline*}
\end{formula}
%

Note that here we always have the $\gamma\to-\gamma$ symmetry,
$|\Gamma_-|=2$,
which is cancelled by $o(m)!=o(g-1)!$
(see~\eqref{eq:sad:conn:const}). Note that
in this particular case we prefer to use summation over {\it all}
combinations of parities $\phi_i$ producing the proper parity of
the sum.
Thus the stratum interchange
symmetry $\Gamma$ depends only on the collection of $g_i$, and on
their
cyclic order.

\begin{example}
{\bf Component} $\boldsymbol{\cH^{odd}(2,2).}$
Consider a saddle connection joining the two zeroes  of a surface
$S\in\cH^{odd}(2,2)$,  and   suppose  that  the   saddle
connection has multiplicity $3$. Each of $S'_i$, for $i=1,2,3$ is
a torus  with a marked point.  The angle between  any consecutive
saddle connections  at $z_1$ and  $z_2$ equals $2\pi$. There is a
symmetry  of  order  $3$  coming  from  stratum  interchange,  so
$|\Gamma|=3$. Taking into consideration $1/4$ from $1/2^{p-1}$ we
get
   %
$$
c=\frac{1}{3} \cdot \frac{1}{4} \cdot \frac{1!\cdot 1!\cdot
1!}{5!}\cdot
\frac{\bigl(\Vol(\cH_1(\torusemptyset))\bigr)^3}{\Vol(\cH_1^{odd}(2,2))}
=
\frac{1}{9} \approx 0.111
$$
\end{example}

\begin{remark}
Let $g\ge  4$ be  odd. We see that the  maximal multiplicity of a
saddle connection joining the  distinct  zeroes of a surface from
the stratum  $\cH(g-1,g-1)$  is different for different connected
components:   it   equals  $2$  for  almost  all  surfaces   from
$\cH^{hyp}(g-1,g-1)$, it equals $g-1$ for almost  all surfaces in
$\cH^{even}(g-1,g-1)$; it equals  $g$  for almost all surfaces in
$\cH^{odd}(g-1,g-1)$.
\end{remark}

\subsection{Other Strata}
\label{ss:other:strata}
The description for all other strata $\cH(\alpha)$, where all the
entries  of  $\alpha$  are  even,  is  similar to  the  one  just
presented.
\begin{formula}
   %
The constant for the  nonhyperelliptic connected component
$\cH^{\phi}(\alpha)$,    $\phi\in\{even,    odd\}$,   $\alpha\neq
(g-1,g-1)$, is given by
   %
\end{formula}
\begin{multline*}
c=
\frac{1}{|\Gamma_-|} \cdot
\frac{1}{|\Gamma|} \cdot
\prod_{m\in\alpha} \left(\frac{o(m)!}{\prod_{j=1}^p
o_j(m)!}\right)
\cdot \prod_{\substack{j=1\\a_j\neq 0}}^p o_j(a_j)
\cdot \prod_{j=1}^p (a_j+1)  \cdot \frac{1}{2^{p-1}}
\cdot \frac{ \prod_{i=1}^p (\tfrac{d_i}{2}-1)! }
     { (\frac{d}{2}-2)! } \cdot
     \\
\cdot \frac{1}{ \Vol(\cH^{\phi}_1(\alpha)) } \cdot
\sum_{\substack{\phi'_1,...,\phi'_p\in\{even,\, odd\}\\
\phi'_1+\dots+\phi'_p=\phi\, (mod\, 2)}}
\prod_{i=1}^p \Big(\Vol(\cH^{\phi_i}_1(\alpha'_i)) +
\delta(\alpha'_i,\phi'_i)\Vol(\cH^{hyp}_1(\alpha'_i))\Big)
\end{multline*}

\begin{example} {\bf Stratum} $\boldsymbol{\cH^{even}(4,2)}.$
 In multiplicity one, after collapsing  zeroes  we  obtain  a  surface
$S'\in \cH(6)$,  where either   of  the two components
$\cH^{hyp}(6)$  and $\cH^{even}(6)$ having  even  parity  of
the  spin structure are possible. We get
  %
$$ 
c =
7\cdot\frac{\bigl(\Vol(\cH_1^{hyp}(6))+\Vol(\cH_1^{even}(6))\bigr)}
               {\Vol(\cH_1^{even}(4,2))}=
\frac{253001}{18225}\approx 13.88
$$ 

In multiplicity two we  get  two surfaces $S_1'\in\cH(\alpha_1')$,
$S'_2\in\cH(\alpha_2')$,  where $\alpha_1'+\alpha_2'=
m_1+m_2-p=4$.
   %
   %

There are two possibilities. In the first case
$\alpha_1'=\alpha_2'$.      Note
that        any $S'\in\cH(2)$ has  odd parity of
the  spin structure, see~\eqref{eq:spin:2g:minus2}. Thus the
parity of the resulting flat      structure $S$ is
even, see Lemma~\ref{lm:spin:in:slit:constr}.

The saddle  connections  $\gamma_1$  and $\gamma_2$ partition the
cone angle $10\pi$ at $z_1$ in two sectors with angles $6\pi$
and $4\pi$; they partition the cone angle $6\pi$ at $z_2$  in two
sectors with the angles $4\pi$ and $2\pi$. Here $a_1'=2,a_1''=0$
while
$a_2'=a_2''=1$
so that while
$\cH(\alpha'_1)=\cH(\alpha'_2)=\cH(2)$,   there   is  no   stratum
interchange since  $a'_1\neq  a'_2$.  Thus $|\Gamma|=1$.  We get
   %
$$ 
c  = 9 \cdot \frac{1}{2} \cdot \frac{3!\cdot 3!}{7!}
\cdot\frac{\bigl(\Vol(\cH_1(2))\bigr)^2}{\Vol(\cH_1^{even}(4,2))}
= \frac{6}{25} = 0.24
$$ 

In the second case after  degeneration we get
$S'_1      \in\cH_1(\torusemptyset)$,
$S'_2\in\cH_1(4)$
Since
$\cH(\torusemptyset)$  has   odd   parity   of   the  spin  structure,
$S'_2\in \cH_1^{odd}(4)$ in
order  to  result  in a surface $S\in\cH^{even}(4,2)$  (see
Lemma~\ref{lm:spin:in:slit:constr}).

The saddle  connections  $\gamma_1$  and $\gamma_2$ partition the
cone angle $10\pi$ at $z_1$ in two sectors with the angles $2\pi$
and $8\pi$; they partition the cone angle $6\pi$ at $z_2$  in two
sectors with the angles $2\pi$ and $4\pi$.
We get
   %
$$ 
c=5 \cdot \frac{1}{2} \cdot \frac{1!\cdot 5!}{7!}
\cdot\frac{\Vol(\cH_1(\torusemptyset))\Vol(\cH^{odd}_1(4))}
          {\Vol(\cH_1^{even}(4,2))} =
\frac{640}{729} \approx 0.8779
$$ 

Multiplicity   $3$   is   not   realizable   in   the   component
$\cH^{even}(4,2)$ (see the next subsection). \end{example}

\begin{example} {\bf Stratum} $\boldsymbol{\cH^{odd}(4,2)}.$
This component is similar in many aspects to the previous one, so
we skip  those details of  calculations which are common for both
components.

The first case is $\cH^{odd}(4,2)$ with multiplicity $p=1$. After
collapsing zeroes  we  obtain a surface $S'\in \cH(6)$,
the only component that results is $\cH^{odd}(6)$.  We get
$$ 
c = 7\cdot\frac{\Vol(\cH_1^{odd}(6))}
               {\Vol(\cH_1^{odd}(4,2))}=
\frac{147}{10}= 14.7
$$ 

Consider  multiplicity
$2$.   The    case    of    degeneration    to   two   surfaces
$S'_1, S'_2\in\cH(2)$  does not take place here,
since  the  surface  $S\in\cH(4,2)$  obtained  from  two
surfaces  from  $\cH(2)$  has even  parity  of spin
structure.

Thus,     the     only    possibility     is    $S'_1
\in\cH_1(\torusemptyset)$;
$S'_2\in\cH_1(4)$,    $a'_2=3$,    $a''_2=1$.   Since
$\cH(\torusemptyset)$  has   odd   parity   of   the  spin  structure,
$S'_2$ must be  in components of $\cH(4)$ having
even  parity  of spin  structure  in order  to  result a  surface
$S\in\cH^{odd}(4,2)$                                (see
Lemma~\ref{lm:spin:in:slit:constr}).  There  is   only  one  such
component: $\cH^{hyp}(4)$.

Thus we get
$$ 
c= 5 \cdot \frac{1}{2} \cdot \frac{1!\cdot 5!}{7!}
\cdot\frac{\Vol(\cH_1(\torusemptyset))\Vol(\cH^{hyp}_1(4))}
          {\Vol(\cH_1^{odd}(4,2))} =
\frac{15}{64} \approx 0.2344
$$ 
(see analogous calculation for the component
$\cH^{even}(4,2)$.

Consider  now   multiplicity   $3$.
It is not hard to check that $S'_1\in\cH(2)$,  $S'_2,
S'_3\in\cH(\torusemptyset)$.    Note     that     both
$\cH(\torusemptyset)$  and  $\cH(2)$  have  odd  parity  of  the  spin
structure. Thus  the  resulting surface $S\in\cH(4,2)$ always has
odd      parity      of      the     spin     structure,      see
Lemma~\ref{lm:spin:in:slit:constr}.

The saddle connections $\gamma_1,  \gamma_2,  \gamma_3$ partition
the cone angle $10\pi$ at $z_1$ in three sectors with  the angles
$6\pi$, $2\pi$, and  $2\pi$; they partition the cone angle $6\pi$
at $z_2$ in three sectors with the angles $2\pi$.
We find

   %
$$ 
c=3 \cdot \frac{1}{4} \cdot \frac{3!\cdot 1!\cdot 1!}{7!}
\cdot\frac{\Vol(\cH(2))\bigl(\Vol(\cH_1(\torusemptyset))\bigr)^2}
          {\Vol(\cH_1^{odd}(4,2))} =
\frac{21}{320} \approx 0.06563
$$ 
\end{example}


\vspace*{1truecm}
\addcontentsline{toc}{part}
{Part 2. Saddle Connections Joining a Zero to Itself}
{\Large\bf Part 2. Saddle Connections Joining a Zero to Itself}
\vspace*{0.5truecm}

   %
In this part  we  consider the second problem  of  this paper: we
count closed saddle connections joining  a  zero to itself.

\section{Approaching the Principal Boundary by Shrinking Closed Geodesics}
\label{s:principal:boundary:II}
\subsection{Configurations of Closed Saddle Connections and
Corresponding Surface Decompositions}
\label{ss:configurations:of:closed:sad:connections}
\begin{figure}[htb!]
%
 %
 %
  %
 %
 \includegraphics{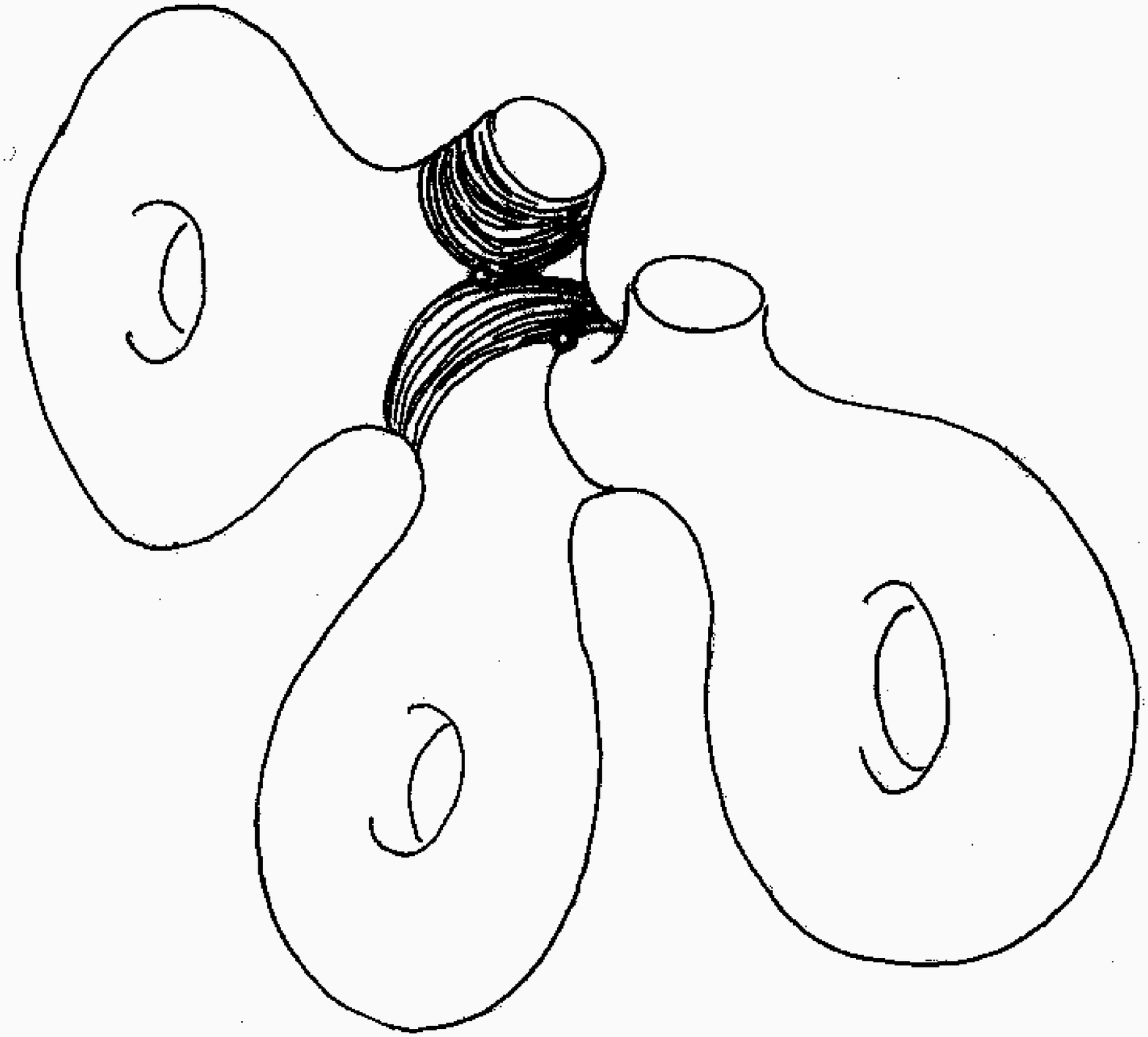}
\includegraphics{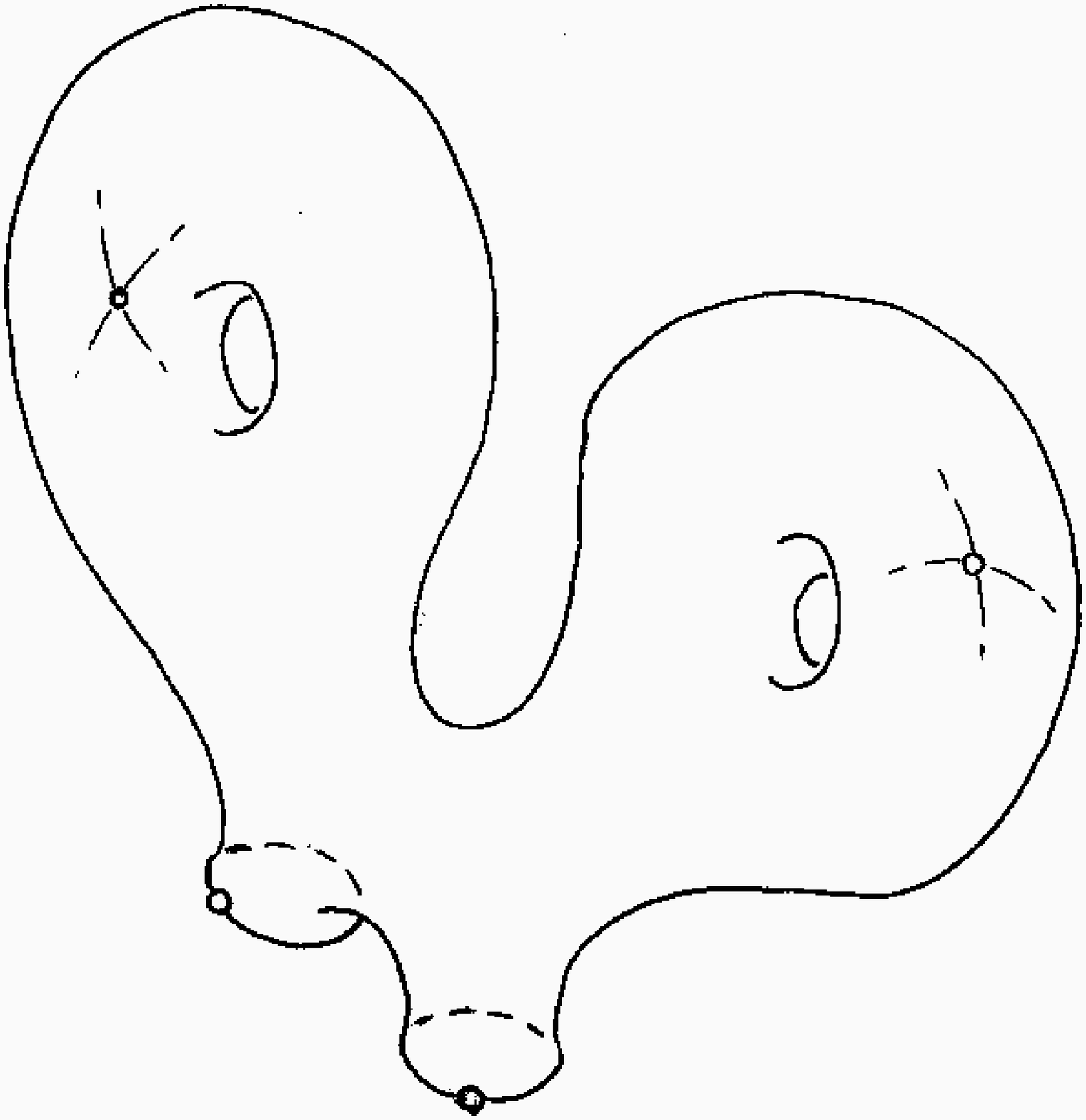}
\begin{picture}(0,0)(0,0)
  \begin{picture}(0,0)(0,0)
     \put(-60,-150){$S_1$}
     \put(90,-115){$S_2$}
     \put(10,-30){$S_3$}
     \put(-110,-25){$S_4$}
  \end{picture}
\end{picture}
\vspace{160bp} 
\caption{
\label{pic:dance2}
Homologous  closed  geodesics  of high multiplicity.  Topological
picture.}
\end{figure}

We have already indicated
in Section~\ref{sec:configurations} that together
with a closed saddle connection $\gamma_1$ joining a
zero  to  itself some   other   saddle   connections
$\gamma_2,\dots,\gamma_p$ homologous  to $\gamma_1$ might be present
on a surface $S$,
see  Figure~\ref{pic:dance2}.  There  might  be  also  some  metric
cylinders filled  with  regular  closed  geodesics  homologous to
$\gamma_1$; such cylinders are bounded  by  the  singular  closed
geodesics  from  the  collection  $\gamma_1,  \dots,   \gamma_p$.
Suppose  that  the  curves  $\gamma_j$ with indices
from a set $J=\{i_1,\ldots,i_l\}\subset \{1,\ldots,p\}$ bound $q$  cylinders.  The
complement  of  the  union  of  the  curves  $\gamma_j$  and  the
cylinders splits into $p=m-q$ disjoint surfaces $S_1,\ldots,S_p$.
For example, the surface $S$ presented in Figure~\ref{pic:dance2}
is composed of $p=4$ surfaces $S_j$ and of $q=2$ cylinders, where
one cylinder  joins $S_2$ to $S_3$  and the other  cylinder joins
$S_3$ to  $S_4$. When $p\ge 2$ the boundary  of the closure $\bar
S_j$ of any  surface  $S_j\subset S$ is a  union  of two singular
closed geodesics $\gamma_{j-1} \cup \gamma_j$.

The sum of the  genera of the $S_j$ is $g-1$. The  surfaces $S_j$
are attached  in some cyclic order. Thus, for  any $j$ the ``next
surface'' $S_{j+1}$  and  the ``previous surface'' $S_{j-1}$ make
sense mod $p$.  We use the  following convention for  the  cyclic
order: let $\vec{n}$ be the vector orthogonal to $\vec{\gamma}_j$
and pointing {\it from} $S_j$ {\it to} $S_{j+1}$.  Then the frame
$(\vec{n},\vec{\gamma}_j)$   corresponds   to    the    canonical
orientation of the surface induced by the complex structure.

We distinguish two types of surfaces $S_j$.

If the  boundary of the closure $\bar S_k$  of $S_k\subset S$ has
two connected components, we get a surface of the first  type. In
this  case the  boundary  $\partial \bar S_k$  is  formed by  two
disjoint  singular  closed geodesics  $\gamma''=\gamma_{k-1}$ and
$\gamma'=\gamma_k$.    The surfaces     $S_1$    and     $S_2$    in
Figure~\ref{pic:dance2} are  of  that  type. The closed geodesics
$\gamma''$ contains a single conical point $z_k$  the closed
geodesics $\gamma'$ contains a  single  conical point $w_k$.  We
denote by $(2b_k'+3)\pi$ the cone angle  at  $w_k$  {\it  inside}
$\bar S_k$,  and by $(2b_k''+3)\pi$  the cone angle at $z_k$ {\it
inside} $\bar S_k$.

Suppose  now  that the  boundary  of the  closure  $\bar S_k$  of
$S_k\subset S$ has a  single  connected component produced by two
singular closed geodesics sharing the  same  conical  point  $P$.
Cut $\bar S_k$ at this  conical point.
If  the  boundary of the resulting surface
now has two connected components, that  means  that  the  surface
$\bar S_k$ is obtained from a surface described  in the paragraph
above by  identification of the  conical points $w\sim z$, and we
do not  distinguish a surface of  this type from  those discussed
above.

We also put in the first group the surface $S_1$ if $\bar S_1=S$;
that  is, if  we  have a singular  closed  geodesic $\gamma_1$  of
multiplicity  one  which  does  not bound a cylinder  of  regular
closed geodesics. Cutting $S$ by $\gamma_1$  we  again  obtain  a
surface with  two disjoint boundary  components $\gamma'\neq \gamma''$ as
above. The initial unique conical point $P\in\gamma_1$ produces
two distinct conical points; one on each of the two boundary components.

Consider now the remaining case. Suppose  that  the  boundary  of
$\bar  S_i$  has  a  single connected component produced  by  two
singular    closed    geodesics    $\gamma_{i-1}=\gamma''$    and
$\gamma_{i}=\gamma'$ sharing the  same  conical point $P$. If the
boundary of  the surface $\bar S_i $  cut at $P$ still
has a single connected component, we say that we get a surface of
the  second  type.  For  example,  the surfaces  $S_3$  and  $S_4$ at
Figure~\ref{pic:dance2} are  of  that  type. The curves $\gamma',
\gamma''$ bound  angles  $2(a_i'+1)\pi$  and $2(a_i''+1)\pi$ {\it
inside} $S_i$  at the conical  point $P$. The two boundary curves
$\gamma'$ and $\gamma''$  are joined together at the zero forming
a ``figure eight''.

Therefore we  get two  types of surfaces $S_k, S_i$;  a boundary of each
surface  is  produced  by  a  pair  of singular closed  geodesics
$\gamma', \gamma''$.
Altogether the configuration consists of the data
\begin{displaymath}
(J,a_i',a_i'',b_k',b_k'')
\end{displaymath}

We will show how to {\it metrically} shrink the boundary components
of each $S_i$ to produce a closed surface.  For surfaces of the first
type the boundary curves will shrink to zeroes of
orders $b_k',b_k''$ and for surfaces of the second type to a single zero
of order $a_i'+a_i''$.  Each such closed surface will then lie
in a stratum $\cH(\alpha_i')$ and the union of these strata comprises
the {\it principal boundary}.

The goal of the rest of this  section is to describe how, for each $S_i$,
one  can naturally  ``shrink''  the boundary components $\gamma',
\gamma''$ to  obtain a  closed regular flat surface. We start with
surfaces of the second type.

\subsection{Shrinking a Pair of Adjacent Holes}
\label{ss:shrinking:adjacent:holes}
Suppose that  the boundary of $\bar  S_i$ has a  single connected
component   produced   by   two    singular    closed   geodesics
$\gamma',\gamma''$ sharing  the  same  conical  point  $P$; the
boundary of such surface has  a form of  figure eight. By assumption the surface
$S_i$ is of the second type,  so cutting $\bar S_i$ at $P$ we get
a surface with a connected boundary. The surgery corresponding to cutting
$\bar S_i$ at $P$ is illustrated at Figure~\ref{pic:figure:8:constr}.
We start with the surface  presented   on   the   right   of
Figure~\ref{pic:figure:8:constr} (there is an extra flat cylinder drawn in the picture,
this flat cylinder will be discussed later)
and we get the surface presented in  the
middle of Figure~\ref{pic:figure:8:constr}.  By construction, the
boundary of the surface obtained after the surgery is made  up of
two parallel geodesic segments of  the  same  length (coming from
$\gamma'$ and $\gamma''$) joining two singular  points $P_1, P_2$
(obtained by breaking $P$ into two  distinct  points).  The  cone
angles  {\it  inside}  $S_i$  at  these  points  are  denoted  by
$2\pi(a_i'+1)$ and  $2\pi(a_i''+1)$.  Now  we  identify  the  two
geodesic  segments  which  form  the boundary; see  the  left
picture  at  Figure~\ref{pic:figure:8:constr}.  We get a  closed
surface with a saddle connection $\gamma$ joining the two zeroes.
By  construction   this   saddle  connection  has    direction
$\gamma_1$ and   length $|\gamma_1|$; moreover, by construction
the resulting surface  does not have any other saddle connections
homologous to $\gamma$.  Thus, if the initial  surface  $S$  belongs  to
$\cH^{\epsilon,\,thick}_1(\alpha)$  then  the
surface obtained from $S_i$ by  our surgery does  not have any  other  short
saddle  connections,  except $\gamma$. Hence,
following the construction of section~\ref{ss:collapsing} we  can
collapse the saddle connection $\gamma$ into  a  single  zero  of
order    $a_i=a_i'+a_i''$.  We have proved that  for any $S$   from
$\cH^{\epsilon,\,thick}_1(\alpha)$ there is a canonical way to  associate
to every component $S_i$  of the second  type a closed flat surface $S'_i$.

In other words, we have proved that for a surface $S_i$
of the second type (one having
boundary formed by two small adjacent
holes) there is a canonical way to shrink metrically
the holes,  obtaining as a result a closed flat surface $S'_i$.

\subsection{Transporting a Small Hole Along a Flat Surface}
\label{ss:transporting:holes}
Intuitively one can interpret the surgery described in the previous section as
annihilation of a ``$+$'' and a ``$-$'' hole which are joined  together. In this
section we show how one can continuously deform the location of an isolated small
hole on a flat surface. In the next section we consider a flat surface with  a pair
of small holes located at different places on the surface. Using the hole-transport
construction we move ``$+$'' and ``$-$'' holes in such way that they become
adjacent  to each other and then we
make  them annihilate  as was described in the previous section.

In this section we consider the case of multiplicity one. Moreover, we assume that
the short closed saddle connection $\gamma$ of multiplicity one does not bound a
metric cylinder. Choose some orientation of $\gamma$; let $\vec{v}=hol(\gamma)$.
Let $\Delta$ be
the length of the shortest saddle connection
different from $\gamma$. Assume that $\Delta\gg\epsilon=|\gamma|$. Now cut the surface $S$
along $\gamma$. Since $\gamma$ is not homologous to zero, we get a connected
surface $S_1$ with boundary having two connected components. Each component is
formed by a single geodesic segment containing exactly one conical point. By
construction the size $\epsilon$ of the holes is small with respect to the scale of
the surface $S_1$, and the holes are located relatively far away one from the
other.

\begin{figure}[ht]
%
\includegraphics{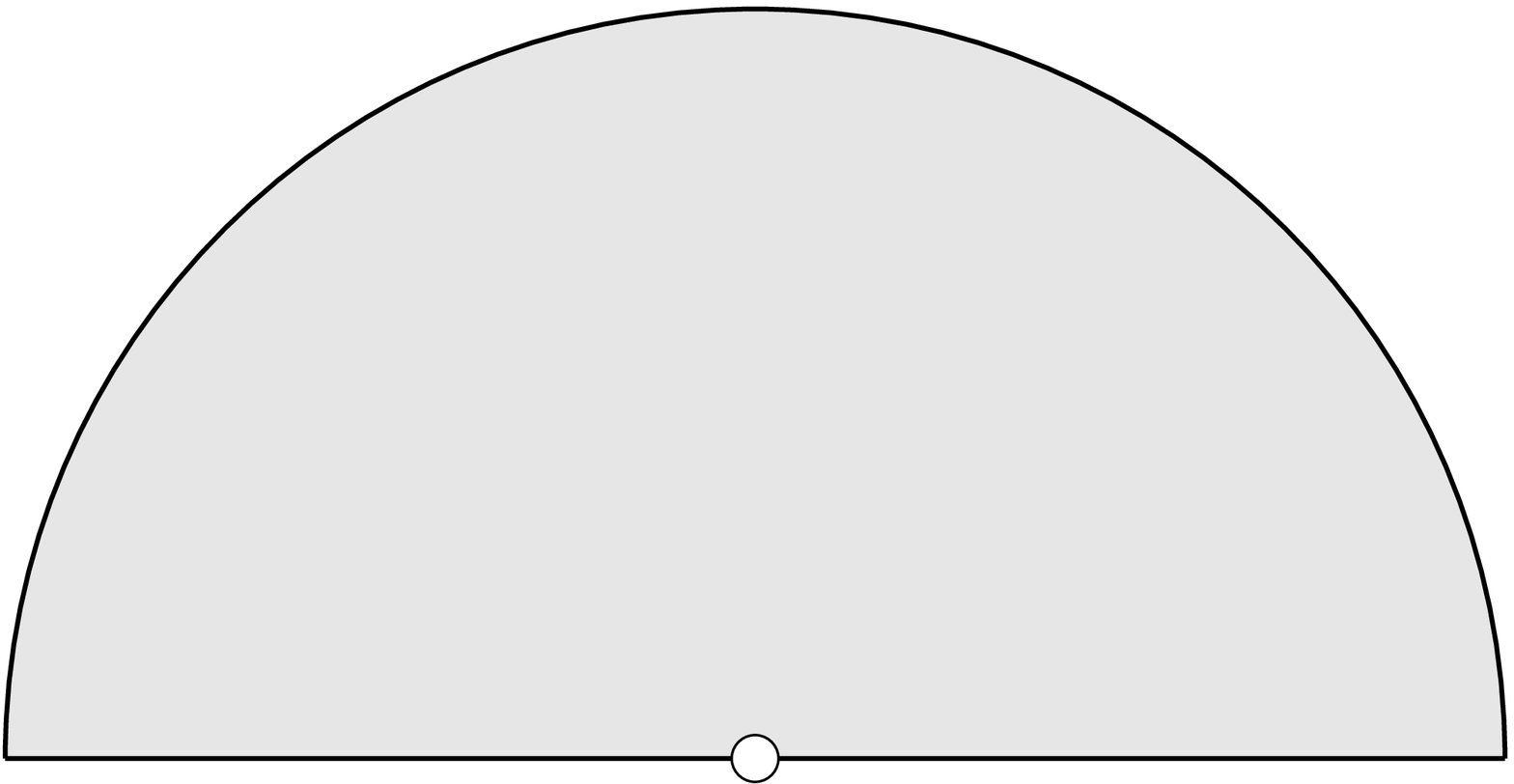}
\includegraphics{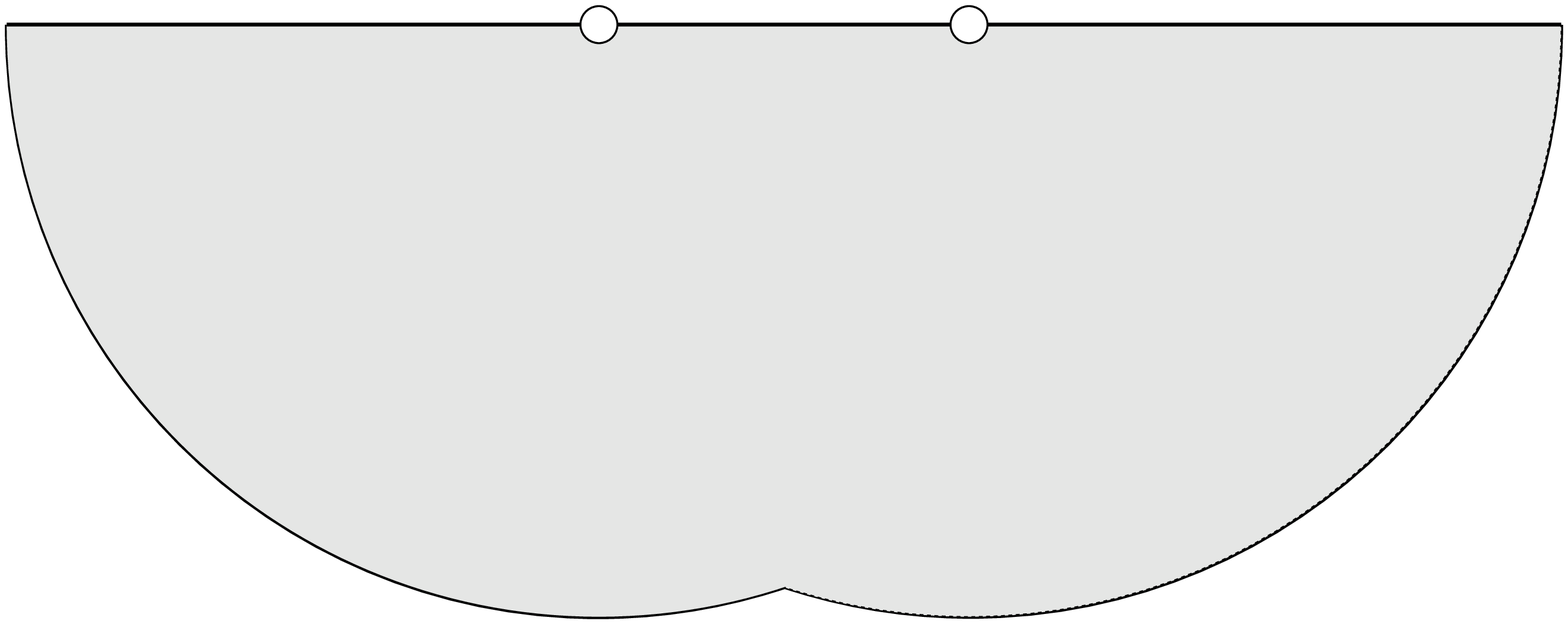}
\includegraphics{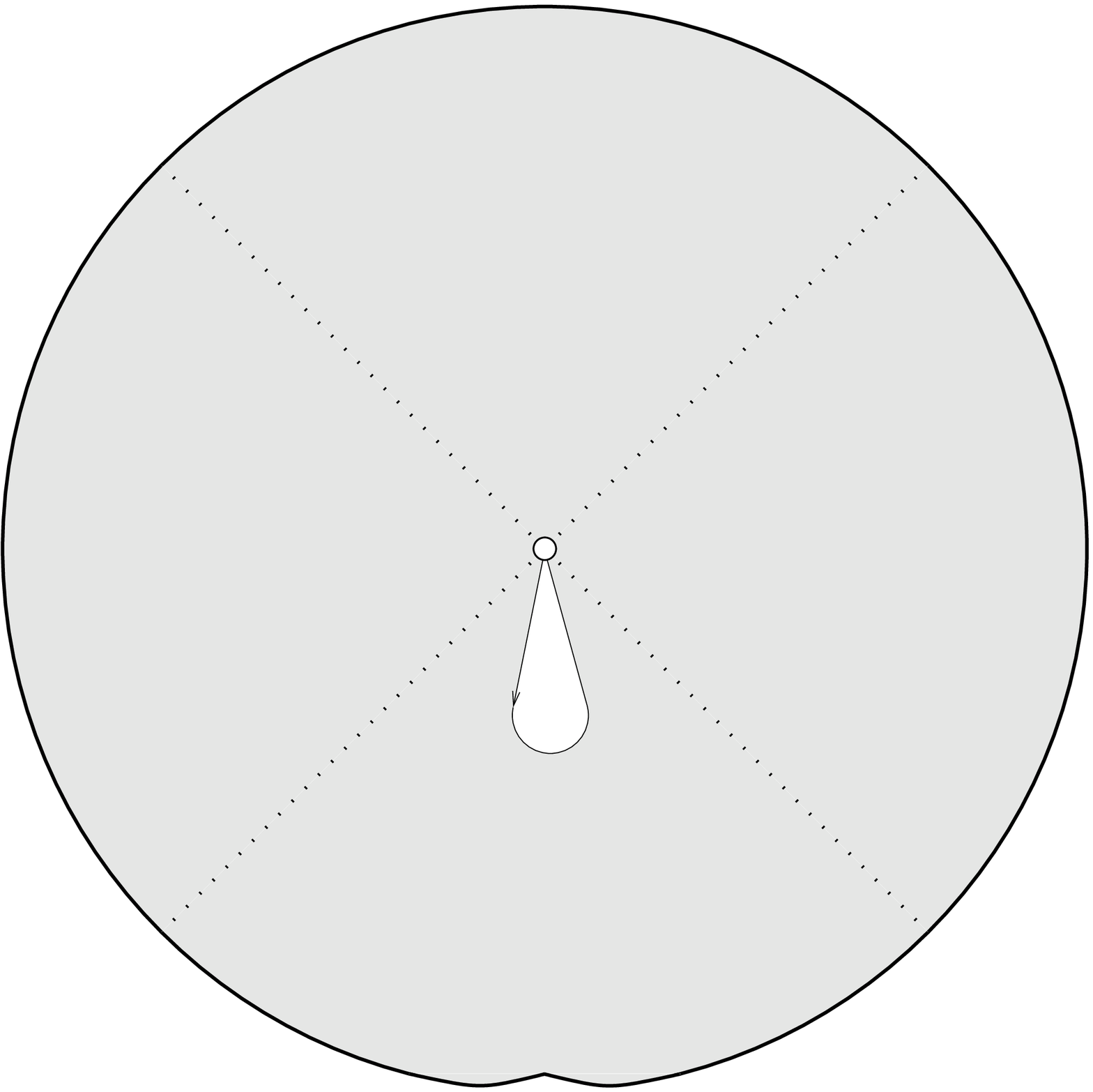}
%
%
\begin{picture}(0,0)(0,0)
\put(60,5) 
 {\begin{picture}(0,0)(0,0)
 \put(-155,-52){$\scriptstyle \Delta$}
 \put(-127,-52){$\scriptstyle \Delta$}
 \put(-155,-82){$\scriptstyle \Delta$} 
 \put(-132,-82){$\scriptstyle \epsilon$} 
 \put(-112,-82){$\scriptstyle \Delta$} 
 \put(-141,-61){$\scriptstyle P_{}$}
 \put(-141,-74){$\scriptstyle P_{}$} 
 \put(-123,-74){$\scriptstyle P_{}$} 
 \put(-169,-20){$\Omega_i,\ {\scriptstyle i=1,\dots,2b+1}$}
 \put(-139,-115){$\Omega_0$}
 \put(-17,-95){$\Omega_0$}
 \put(10,-65){$\Omega_1$}
 \put(-17,-45){$\Omega_2$}
 \put(-48,-65){$\Omega_{2b+1}$}
 \end{picture}}
\end{picture}
\vspace{110bp} 
\caption{
\label{fig:a:hole}
Neighborhood $D_\Delta(P)$ of a boundary singularity $P$ is composed of
a special sector $\Omega_0$ and of $2b+1$ regular sectors $\Omega_1, \dots
\Omega_{2b+1}$. The cone angle at $P$ is $(2b+3)\pi$.}
\end{figure}

Choose one of the holes, say, the {\it positive} one for which the
orientation induced by the natural orientation of the boundary
coincides with the chosen orientation of $\gamma$. We denote this hole
as $\gamma_+$. Let $P$ be the unique
singularity point located on $\gamma_+$. Consider the set
$D_{\Delta}(P)\subset S_1$ of
points in $S_1$ located at  distance less than or equal to
$\Delta$ from $P$.
Denote the cone angle at $P$ by $(2b+3)\pi$;
it is easy to see that the cone angle is always an odd multiple of $\pi$,
with $b$ being some nonnegative integer.
The set $D_{\Delta}(P)$ is composed of $2b+1$ regular sectors,
isometric to metric half-discs of radius $\Delta$
and of one irregular sector.
The irregular sector can be
constructed metrically as follows: superpose two copies of
metric half-discs of radius $\Delta$; then shift
one with respect to the other by $\epsilon$ in direction of
the common diameter; finally, take the union of the resulting
metric half-discs;  see Figure~\ref{fig:a:hole}.
The distinguished
sector contains the hole, see Figure~\ref{fig:a:hole},
so its impact to the cone angle at $P$ is $2\pi$ instead of $\pi$.
We choose  the boundaries
of the sectors to have directions $\pm\vec{v}$.
We denote the distinguished sector by $\Omega_0$; we enumerate
the other sectors $\Omega_1, \dots, \Omega_{2b+1}$ in a cyclic order.

\begin{figure}[htb] 
%
\includegraphics{break0_sh.eps}
\includegraphics{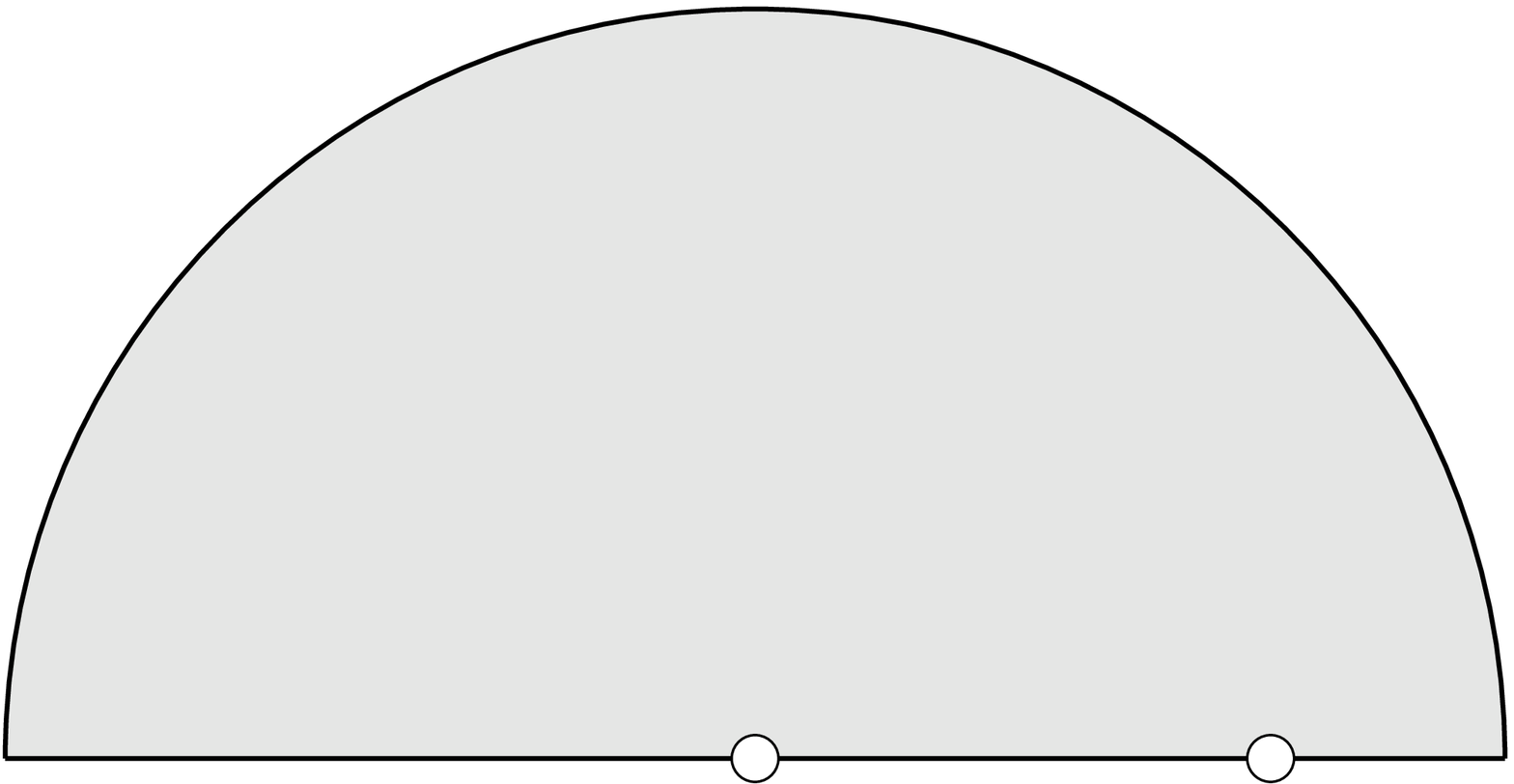}
\includegraphics{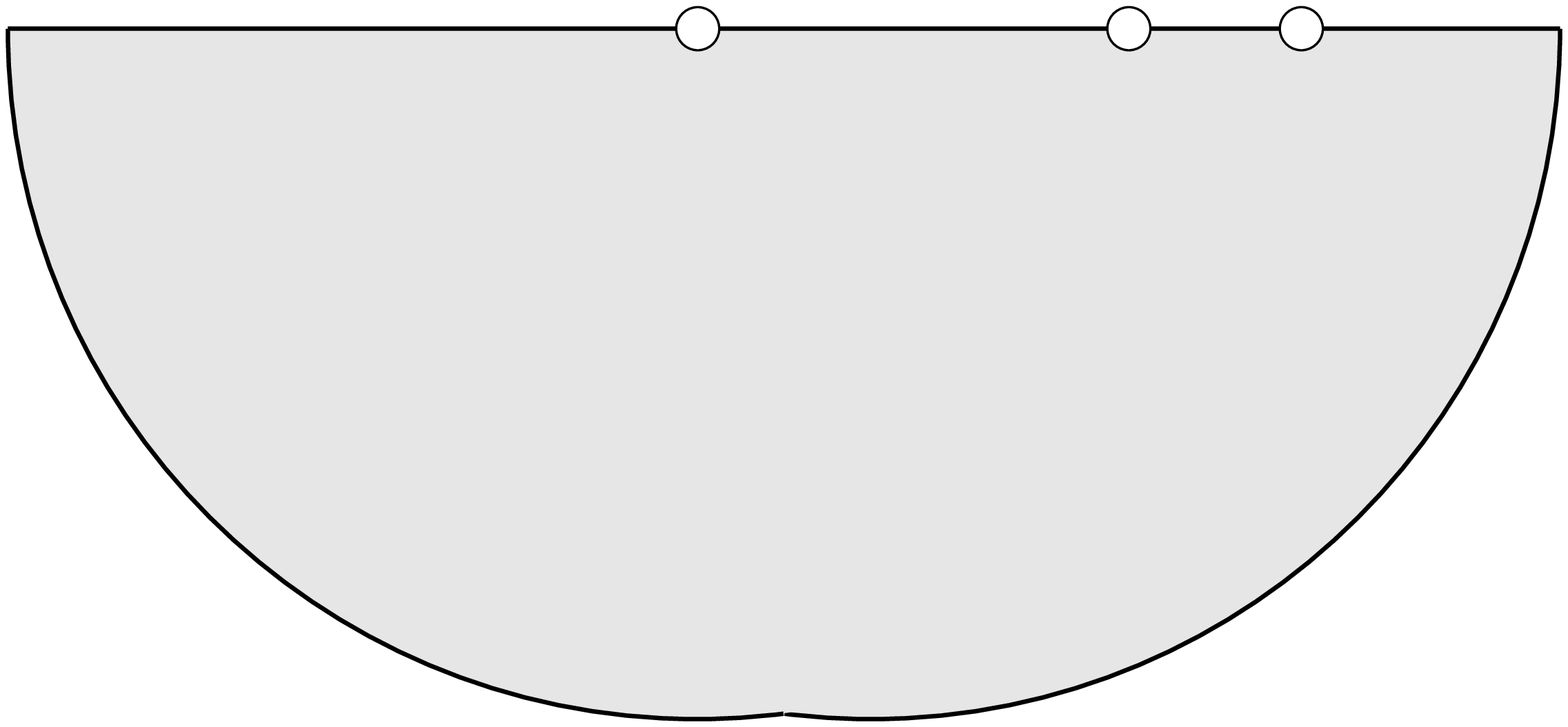}
\includegraphics{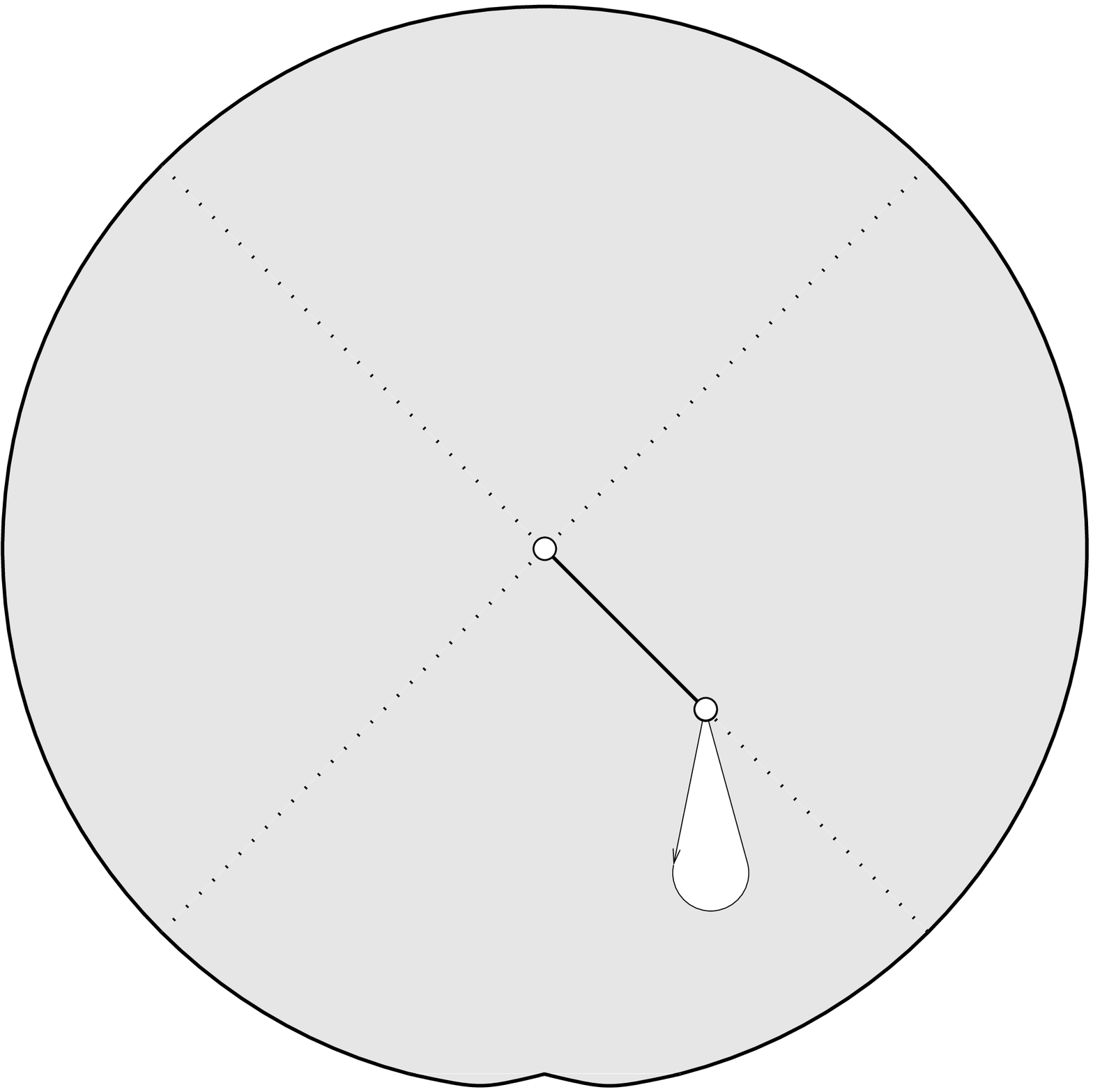}
%
\includegraphics{break0_sh.eps}
\includegraphics{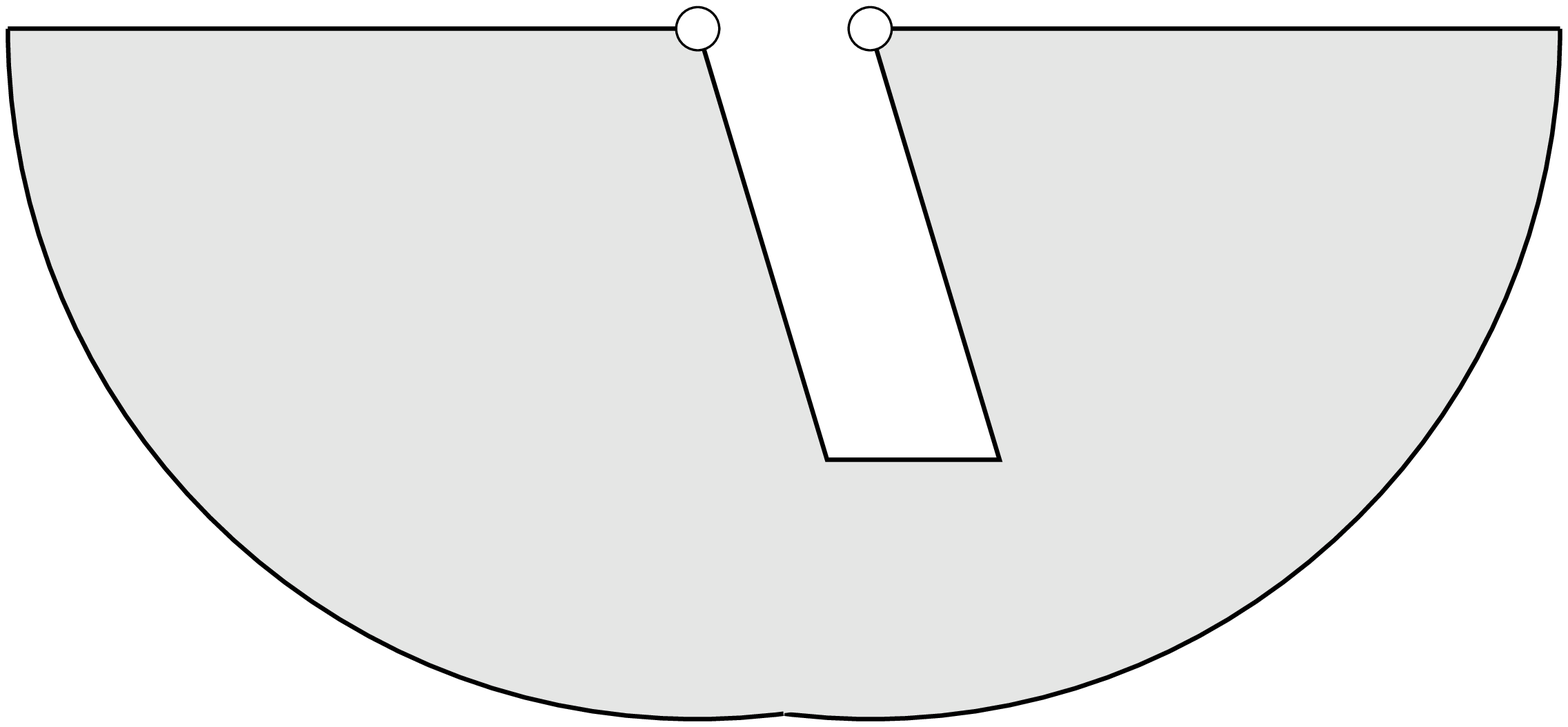}
\includegraphics{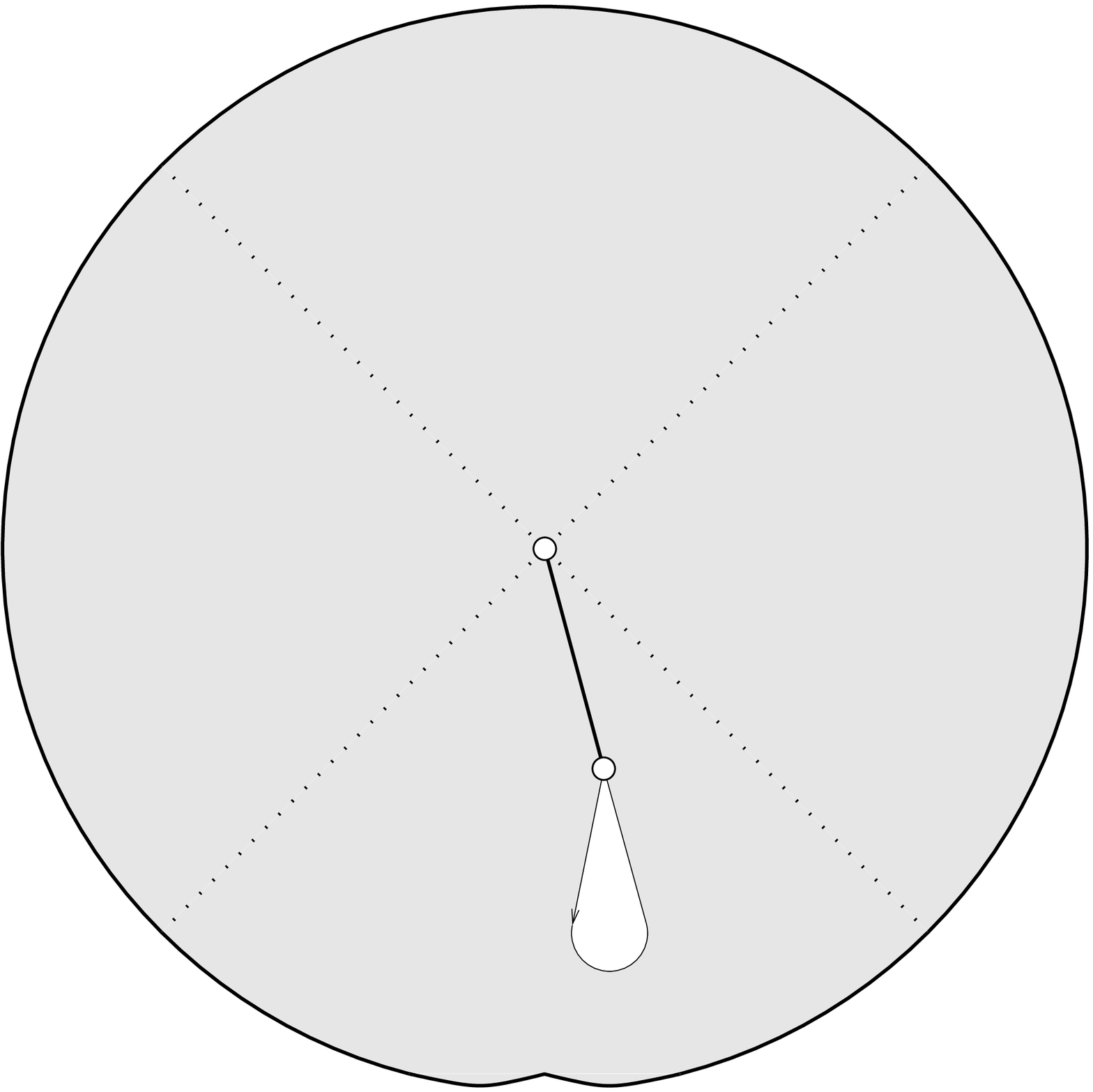}
%
\includegraphics{break0_sh.eps}
\includegraphics{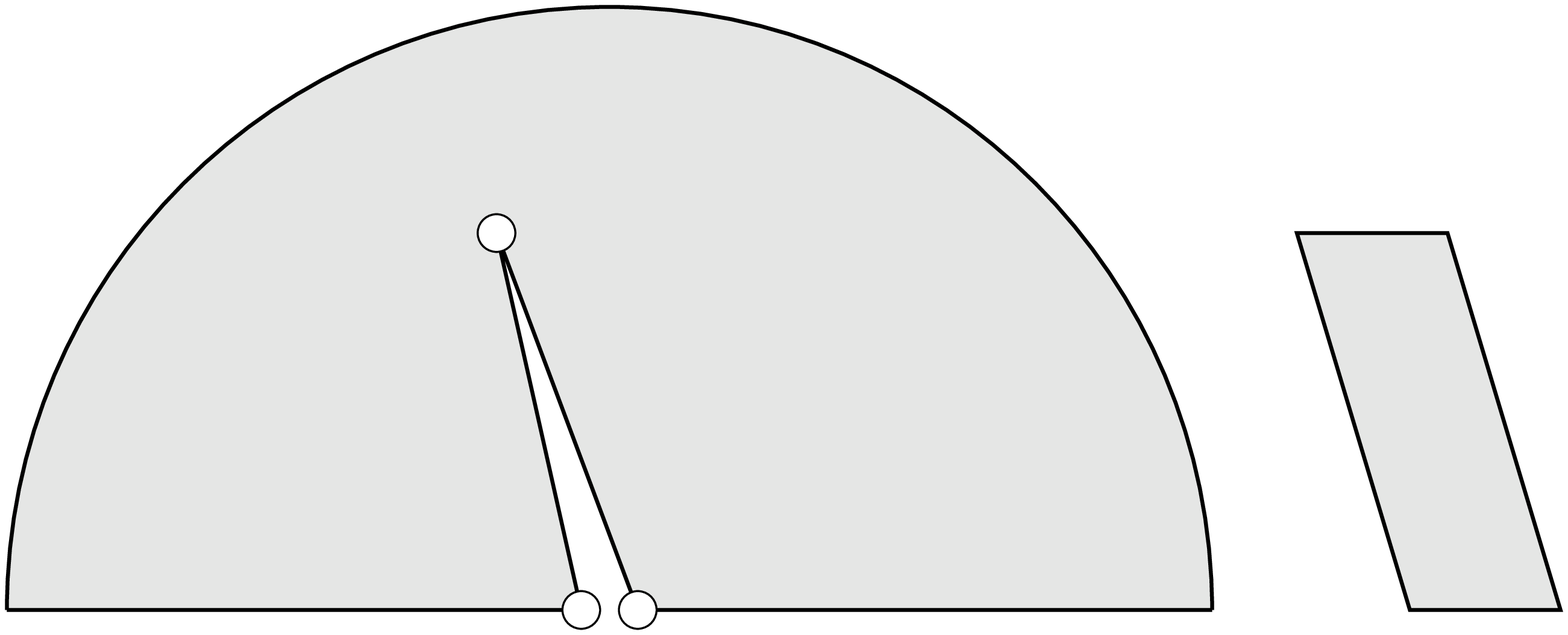}
\includegraphics{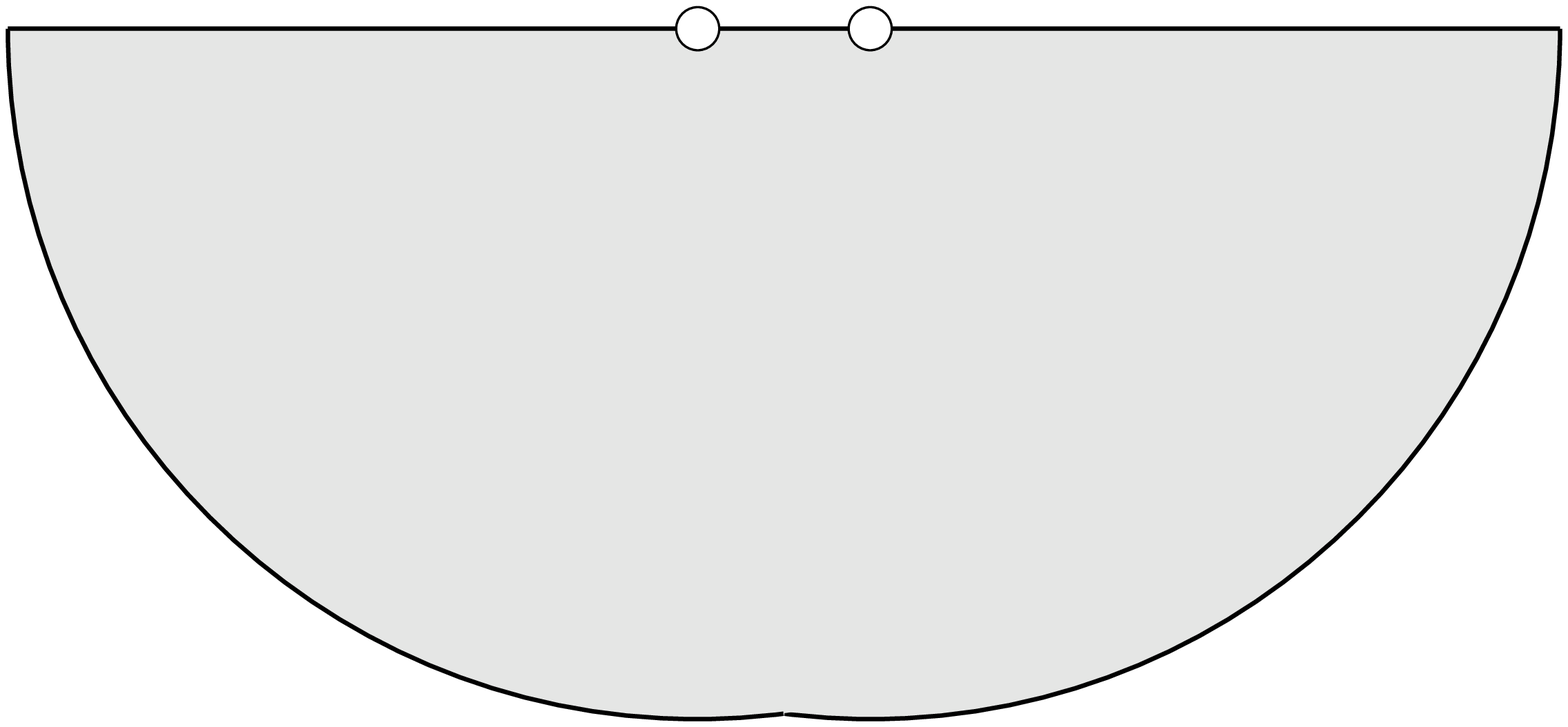}
\includegraphics{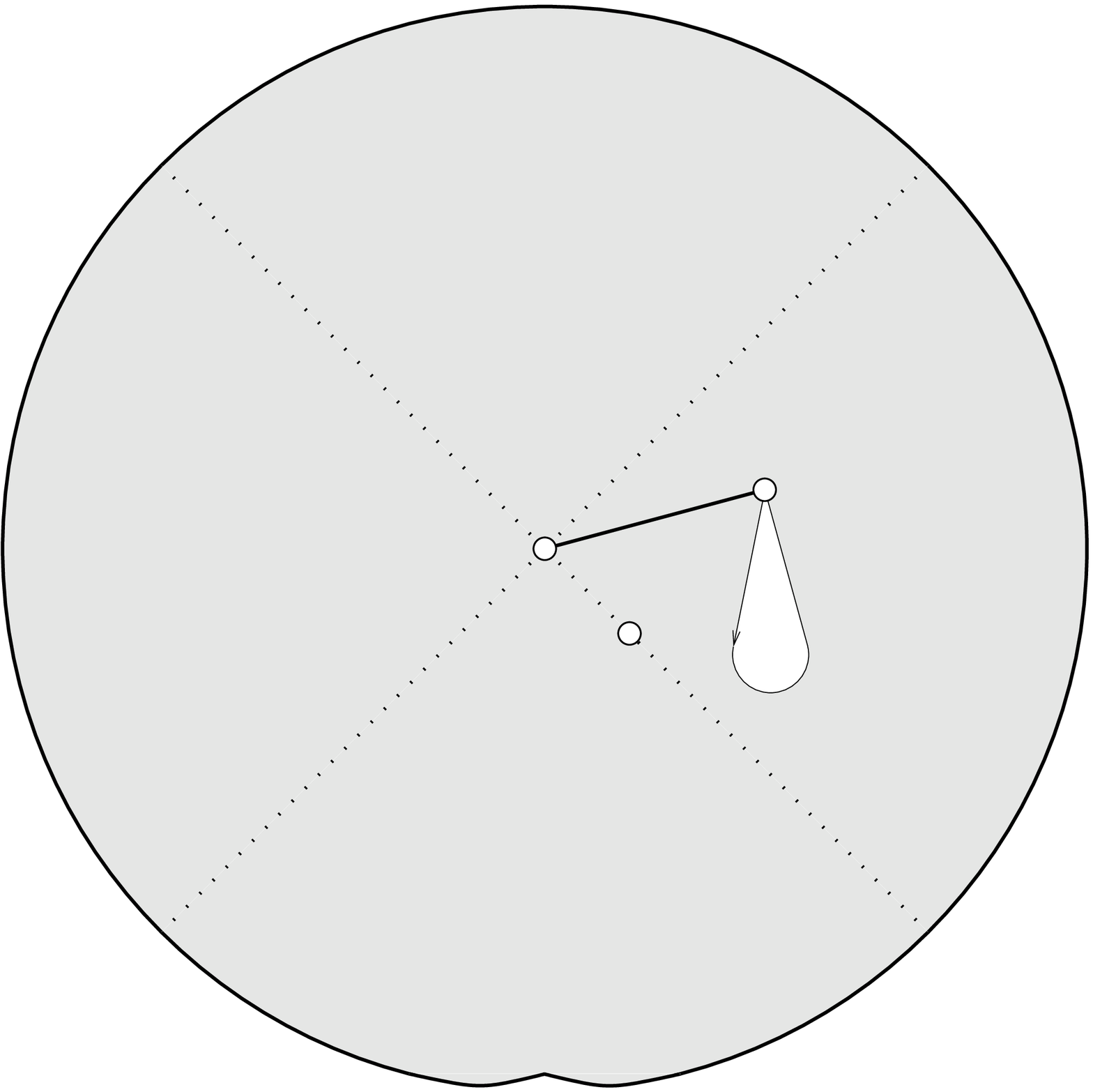}
%
\includegraphics{break0_sh.eps}
\includegraphics{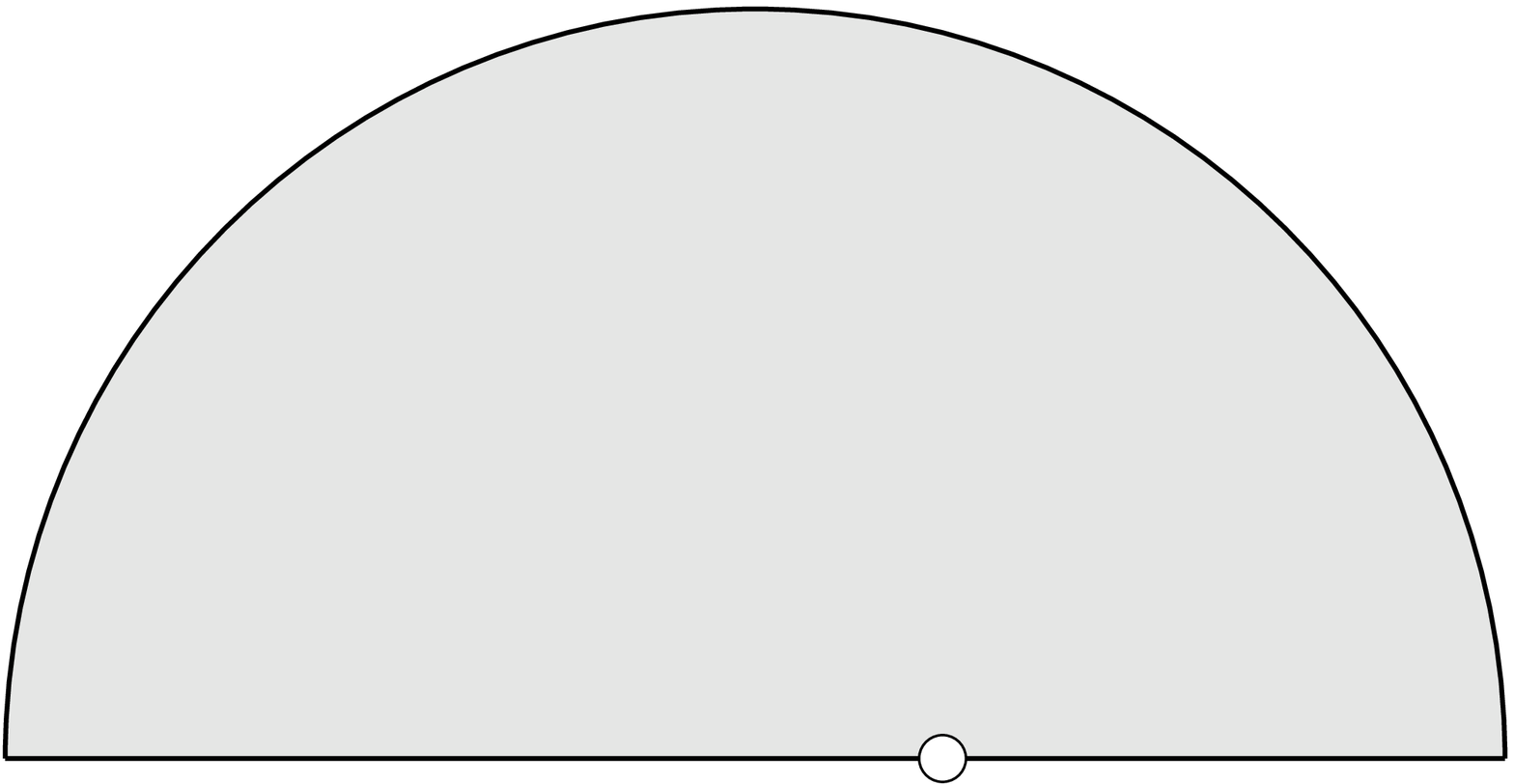}
\includegraphics{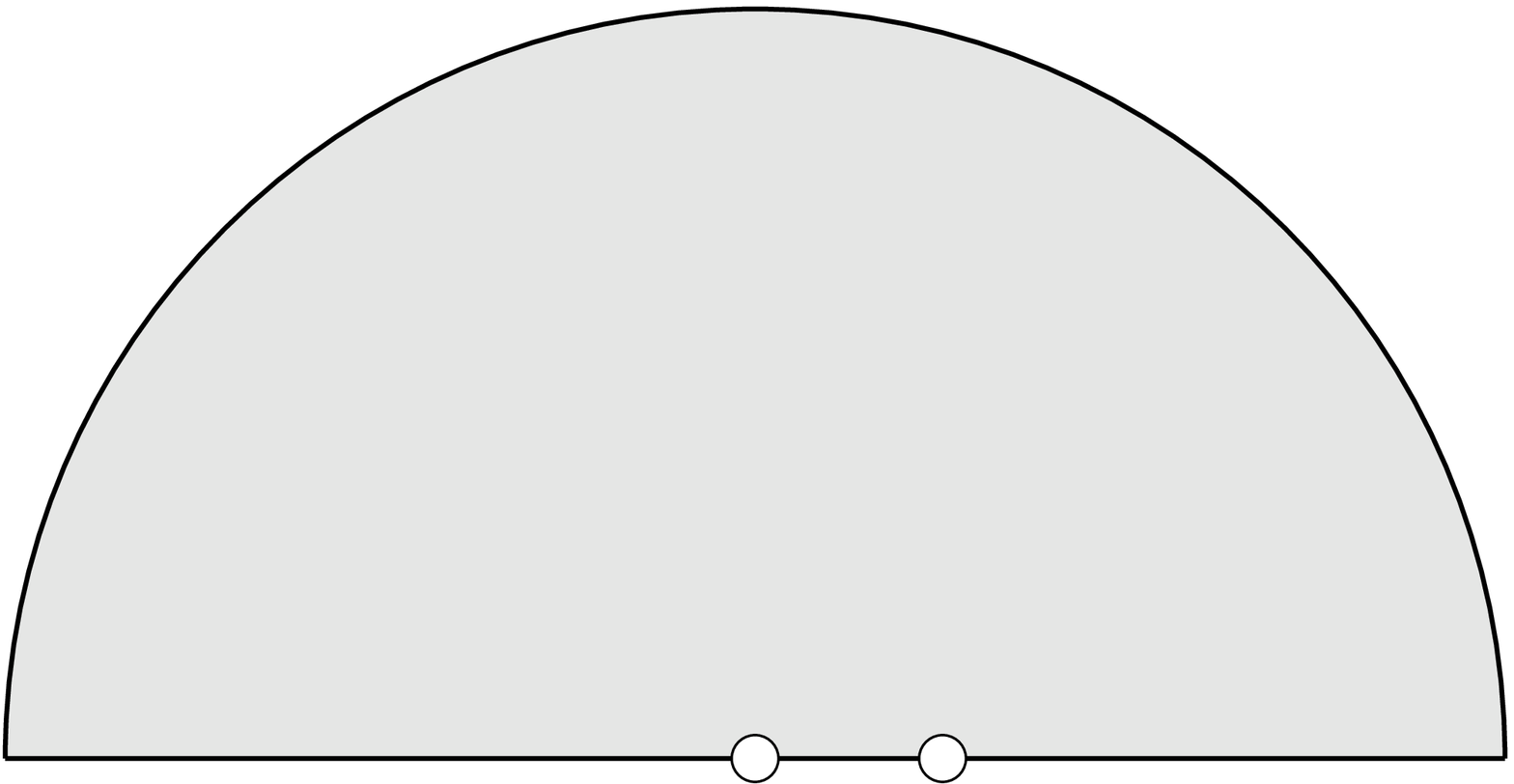}
\includegraphics{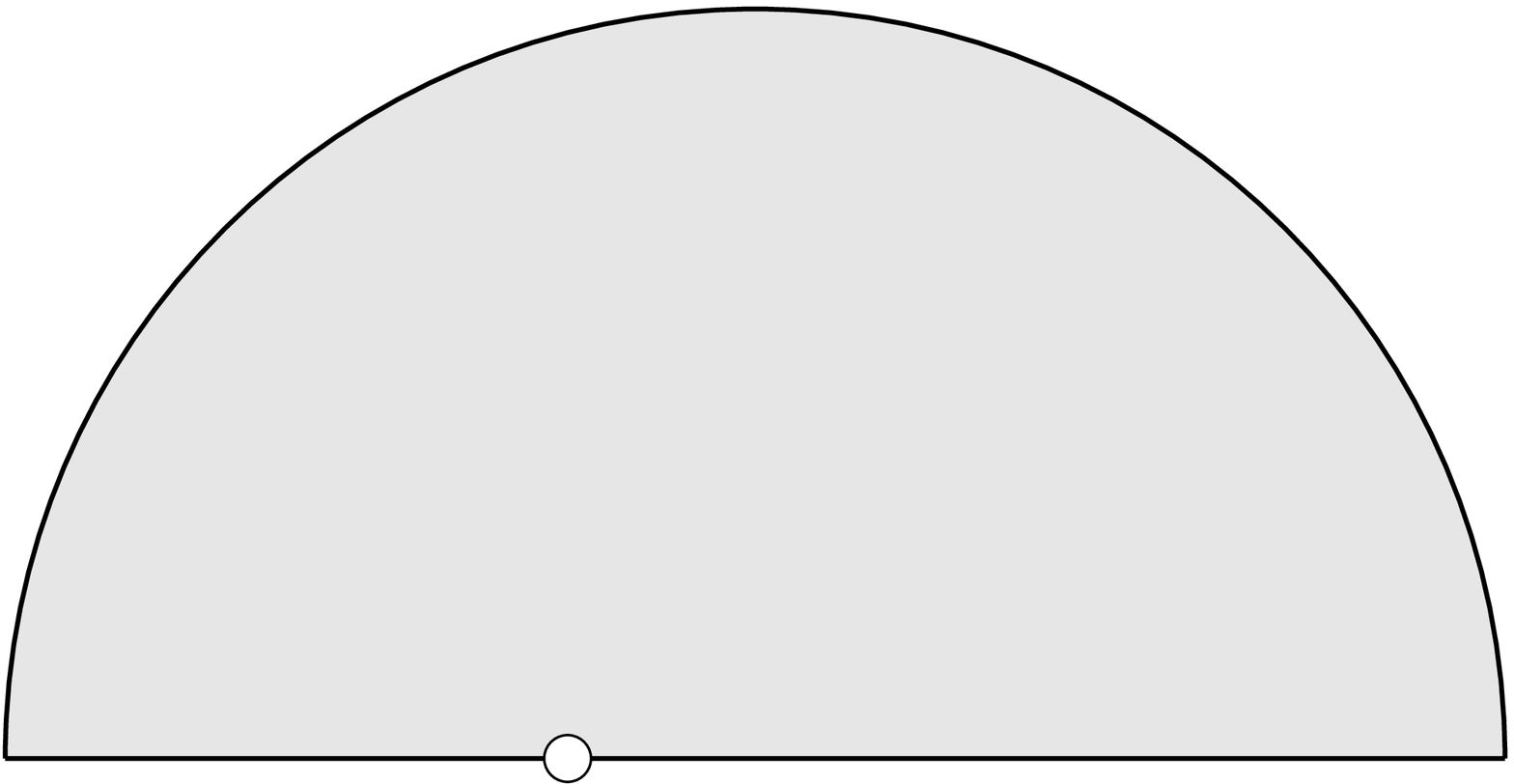}
\includegraphics{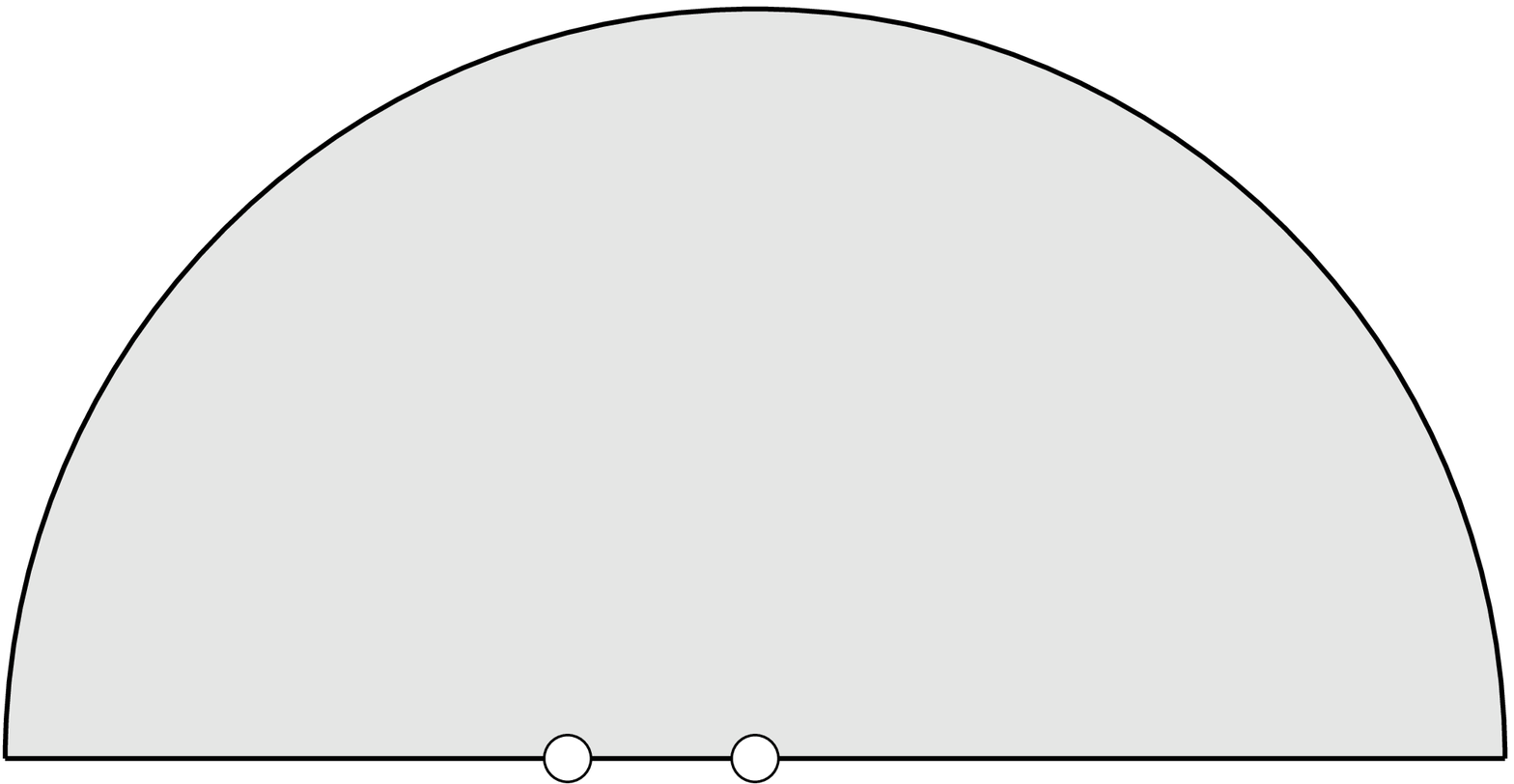}
\includegraphics{omega0.eps}
\includegraphics{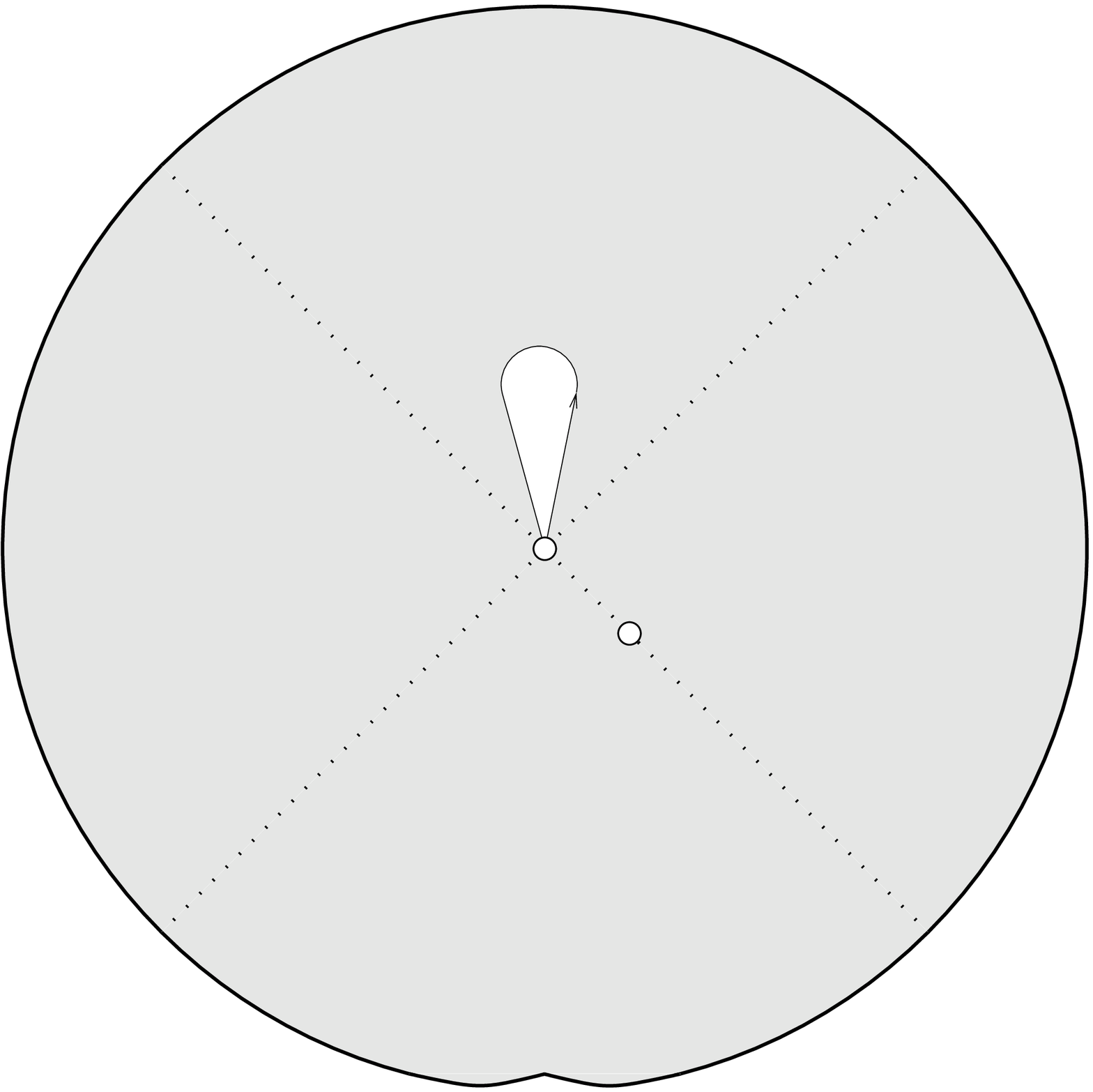}
%
%
\begin{picture}(0,0)(0,0)
\put(10,-10) 
 {\begin{picture}(0,0)(0,0)
 \put(-146,-38){$\scriptstyle \Delta$} %
 \put(-119,-38){$\scriptstyle \Delta$} %
 \put(-133,-46){$\scriptstyle P_{}$} 
 \put(-150,-6){$\Omega_i,\ {\scriptstyle i=2,\dots,2b+1}$} %
 \put(-65,-38){$\scriptstyle \Delta$} %
 \put(-41,-38){$\scriptstyle \delta$} %
 \put(-53,-46){$\scriptstyle P_{}$}  %
 \put(-33,-46){$\scriptstyle P'$} %
 \put(-53,-6){$\Omega_1$} %
 \put(-105,-68){$\scriptstyle \Delta$} 
 \put(-82,-68){$\scriptstyle \delta$}
 \put(-70,-68){$\scriptstyle \epsilon$}
 \put(-91,-59){$\scriptstyle P_{}$}  
 \put(-74,-59){$\scriptstyle P'$} %
 \put(-64,-59){$\scriptstyle P'$} %
 \put(-89,-100){$\Omega_0$} %
 \put(63,-46){$\scriptstyle P_{}$}  %
 \put(77,-59){$\scriptstyle P'$} %
 \put(-180,-6){\bf I.}
 \end{picture}}
\put(10,-130) 
 {\begin{picture}(0,0)(0,0)
 \put(-104,-38){$\scriptstyle \Delta$} 
 \put(-77,-38){$\scriptstyle \Delta$} 
 \put(-93,-46){$\scriptstyle P_{}$} 
 \put(-102,-6){$\Omega_i,\ {\scriptstyle i=1,\dots,2b+1}$} 
 \put(-105,-68){$\scriptstyle \Delta$} 
 \put(-70,-68){$\scriptstyle \Delta$}
 \put(-82,-84){$\scriptstyle \epsilon$}
 \put(-93,-59){$\scriptstyle P_{}$}  
 \put(-82,-59){$\scriptstyle P_{}$} %
 \put(-93,-81){$\scriptstyle P'$}  
 \put(-76,-81){$\scriptstyle P'$} %
 \put(-89,-100){$\Omega'_0$} %
 \put(63,-46){$\scriptstyle P_{}$}  %
 \put(71,-64){$\scriptstyle P'$} %
 \put(-180,-6){\bf II.}
 \end{picture}}
\put(10,-251) 
 {\begin{picture}(0,0)(0,0)
 \put(-146,-38){$\scriptstyle \Delta$} %
 \put(-119,-38){$\scriptstyle \Delta$} %
 \put(-133,-46){$\scriptstyle P_{}$} 
 \put(-150,-6){$\Omega_i,\ {\scriptstyle i=2,\dots,2b+1}$} %
 \put(-75,-38){$\scriptstyle \Delta$} 
 \put(-48,-38){$\scriptstyle \Delta$} 
 \put(-67,-46){$\scriptstyle P_{}$}  
 \put(-58,-46){$\scriptstyle P_1$} 
 \put(-63,-19){$\scriptstyle P'$} 
 \put(-63,-6){$\Omega_1$} 
 \put(-19,-38){$\scriptstyle \epsilon$} %
 \put(-23,-25){$\scriptstyle \epsilon$} %
 \put(-25,-46){$\scriptstyle P_{}$}  %
 \put(-16,-46){$\scriptstyle P_1$} %
 \put(-28,-19){$\scriptstyle P'$} %
 \put(-19,-19){$\scriptstyle P'$} %
 \put(-26,-6){$\Pi$} %
 \put(-105,-63){$\scriptstyle \Delta$} %
 \put(-87,-63){$\scriptstyle \epsilon$}
 \put(-70,-63){$\scriptstyle \Delta$}
 \put(-93,-53){$\scriptstyle P_{}$}  %
 \put(-82,-53){$\scriptstyle P_1$} %
 \put(-89,-95){$\Omega_0$} %
 \put(63,-46){$\scriptstyle P_{}$}  %
 \put(82,-42){$\scriptstyle P'$} %
 \put(69,-63){$\scriptstyle P_1$} %
 \put(-180,-6){\bf III.}
 \end{picture}}
\put(10,-385)                      
 {\begin{picture}(0,0)(0,0)
 \put(-146,-54){$\scriptstyle \Delta$} 
 \put(-119,-54){$\scriptstyle \Delta$} 
 \put(-135,-54){$\scriptstyle P_{}$}   
 \put(-133,-40){$\Omega_k$} 
 \put(-170,-40){${\scriptstyle k>j}$} 
 \put(-75,-14){$\scriptstyle \Delta-\epsilon$} 
 \put(-45,-14){$\scriptstyle \Delta+\epsilon$}          
 \put(-61,-22){$\scriptstyle P_{}$}            
 \put(-53,  0){$\Omega_i$} 
 \put(-20,  0){${\scriptstyle i=3,5,\dots}$} 
 \put(-15,  -10){${\scriptstyle i<j}$} 
 \put(-70,-54){$\scriptstyle \Delta+\epsilon$} 
 \put(-39,-54){$\scriptstyle \Delta-\epsilon$}          
 \put(-48,-62){$\scriptstyle P$}             
 \put(-53,-40){$\Omega_i$} 
 \put(-20,-40){${\scriptstyle i=2,4,\dots}$} 
 \put(-15,  -50){${\scriptstyle i<j}$} 
 \put(-75,-94){$\scriptstyle \Delta-\epsilon$} 
 \put(-54,-94){$\scriptstyle \epsilon$}        
 \put(-39,-94){$\scriptstyle \Delta$}          
 \put(-61,-102){$\scriptstyle P_{}$}            
 \put(-50,-102){$\scriptstyle P_1$}             
 \put(-53,-80){$\Omega_1$} 
 \put(-145,-14){$\scriptstyle \Delta$} 
 \put(-127,-14){$\scriptstyle \epsilon$} 
 \put(-119,-14){$\scriptstyle \Delta-\epsilon$} 
 \put(-134,-22){$\scriptstyle P_{}$}  
 \put(-123,-22){$\scriptstyle P_1$} 
 \put(-133, 0){$\Omega_j$} 
 \put(-145,-73){$\scriptstyle \Delta$} 
 \put(-127,-73){$\scriptstyle \epsilon$}
 \put(-110,-73){$\scriptstyle \Delta$}  
 \put(-134,-65){$\scriptstyle P_{}$}    
 \put(-123,-65){$\scriptstyle P_1$}     
 \put(-129,-90){$\Omega_0$} 
 \put(58,-65){$\scriptstyle P_{}$}  %
 \put(74,-66){$\scriptstyle P_1$} %
 \put(45,-75){$\scriptstyle \Delta$}  %
 \put(69,-64){$\scriptstyle \epsilon$} %
 \put(83,-77){$\scriptstyle \Delta$}
 \put(45,-50){$\scriptstyle \Delta$}  %
 \put(80,-50){$\scriptstyle \Delta-\epsilon$} %
 \put(-180,18){\bf IV.}
 \end{picture}}
\end{picture}
\vspace{480bp} 
\caption{
\label{fig:neighborhood:of:a:hole} Transporting a hole from $P_{}$ to $P'$ along
$\tau$ and tunnelling it to a nonadjacent sector.\newline
\hspace*{1truecm} I. $\tau\subset\Omega_0$ and $\tau$ has direction $\pm \vec{v}$;
\newline
\hspace*{1truecm} II. $\tau\subset\Omega_0$ and $\tau$ has direction different from
$\pm\vec{v}$;
\newline
\hspace*{1truecm} III. $\tau\subset\Omega_{\pm 1}$ and $\tau$ has direction
different from $\pm\vec{v}$;
\newline
\hspace*{1truecm} IV. $\tau\subset\Omega_i$ and $i\neq {-1,0,1}$
}
\end{figure}

Consider a geodesic segment $\tau\subset D_\Delta(P)$ having the conical
point $P$ as one of the endpoints; let $P'$ be the other endpoint of $\tau$;
let $\delta$ be the length of $\tau$.
We assume that
$\delta< \Delta-\epsilon$, so
$P'$ is located in the interior of $D_\Delta(P)$.
Our goal is to ``move the hole'' to the point $P'$. We treat separately
the four cases listed below which describe all possible locations of $\tau$.
Each time we modify some sectors $\Omega_i$ or the identification
of the boundaries of the sectors in such way that after regluing
of the modified sectors
a neighborhood
of the boundary of
the modified domain $D'_\Delta(P)$ stays isometric to a neighborhood
of the boundary of
the initial domain $D_\Delta(P)$. In this way we can cut $D_\Delta(P)$
out of the flat surface; make a surgery on $D_\Delta(P)$ and then
{\it metrically}
paste back
the modified domain $D'_\Delta(P)$ into the surface.
Here is a complete list of possible locations of $\tau$
(see also Figure~\ref{fig:neighborhood:of:a:hole});
as usual indices are taken modulo the largest one,
so $\Omega_{-1}=\Omega_{2b+1}$.
   %

   %
\begin{itemize}
\item[I.]
$\tau\subset\Omega_0$ and $\tau$ has direction $\pm \vec{v}$;
\item[II.]
$\tau\subset\Omega_0$ and $\tau$ has direction different from $\pm\vec{v}$;
\item[III.]
$\tau\subset\Omega_{\pm 1}$ and $\tau$ has direction
different from $\pm\vec{v}$;
\item[IV.]
$\tau\subset\Omega_i$ where $i\neq {-1,0,1}$
\end{itemize}
\medskip

{\bf I. Moving the hole along $\pmb{\vec{v}}$}.\newline
\indent
To move the hole in direction $\vec{v}$
it is sufficient to modify the identification
of the boundary of $\Omega_0$ with the boundary of $\Omega_1$.
Similarly, to move the hole in direction $-\vec{v}$
it is sufficient to modify the identification
of the boundary of $\Omega_0$ with the boundary of $\Omega_{-1}$.
The first case is illustrated at
Figure~\ref{fig:neighborhood:of:a:hole}, case I. We identify
the segments $[P_{},P']$ on $\Omega_0$ and on
$\Omega_1$ (this is exactly
our geodesic segment $\tau$ of length $\delta$). Then we leave
a segment of length $\epsilon$ on $\Omega_0$
without identification, but we identify its endpoints
getting a single point $P'$.
This is the new hole. Then we identify the remaining segments
of length $\Delta-\delta$ (the ones to the right of the point
$P'$ at Figure~\ref{fig:neighborhood:of:a:hole}, case I)
on $\Omega_0$ and on $\Omega_1$.

In other words, our surgery consists of
marking a superfluous segment of
length $\epsilon$ on the boundary of $\Omega_0$
at the distance $\delta$ to the right (correspondingly to the
left) of the initial location.
\smallskip

{\bf II. Moving the hole inside the distinguished sector}.\newline
\indent
Consider a parallelogram $\Pi$ in $\Omega_0$ having $\gamma$ and $\tau$
as a pair of sides (see Figure~\ref{fig:neighborhood:of:a:hole}, II);
let $\gamma'$ and $\tau'$ be the other two sides of $\Pi$ parallel to
$\gamma$ and $\tau$ correspondingly. We can think of $\tau'$ as of a
parallel translation  of $\tau$ by $\gamma$.  Cut $\Pi$ out of $\Omega_0$ and identify
the opposite sides $\tau$ and $\tau'$ in $\Omega'_0=\Omega_0-\Pi$.
Perform the initial identification of the boundaries of the sectors using $\Omega'_0$ instead
of $\Omega_0$. As a result of this surgery we get a new domain
$D'_\Delta(P)$ with a hole located at the point $P'$. Note that
performing this surgery we have reduced the area of the surface
by the area of the parallelogram $\Pi$ which was cut out.
\smallskip

{\bf III. Moving the hole inside one of the two adjacent sectors
$\pmb{\Omega_{-1}, \Omega_{1}}$}.\newline
\indent
If the segment $\tau$ transversal to the direction $\vec{v}$
is located in one of the sectors $\Omega_{-1}$, $\Omega_1$ we can
perform an operation inverse to the one presented above.
Suppose that, say, $\tau\subset \Omega_1$.
Consider a parallelogram $\Pi$ in the Euclidean plane
with a pair of sides having the same lengths and
directions as $\gamma$ and $\tau$ correspondingly.
Slit $\Omega_1$ along $\tau$ from $P$ to $P'$;
let $\tau$ and $\tau'$ be the sides of the resulting slit. Identify
the corresponding sides of $\Pi$ with $\tau$ and $\tau'$. We have
modified the sector $\Omega_1$ in the following way: the
modified sector $\Omega'_1=\Omega_1+\Pi$ has an extra segment $\gamma$ at the base,
where $hol(\gamma)=\vec{v}$, and a new hole $\gamma'$ at the point $P'$,
where $hol(\gamma')=\vec{v}$.
Perform the initial identification of the boundaries of the sectors using $\Omega'_1$ instead
of $\Omega_1$; add a new one identifying the segment $\gamma$ of length $\epsilon$
at the base of
$\Omega'_1$ with the segment $\gamma$ of length $\epsilon$ at the base of $\Omega_0$.
As a result of this surgery we get a new domain
$D'_\Delta(P)$ with a hole located at the point $P'$. Note that
performing this surgery we have increased the area of the surface
by the area of the parallelogram $\Pi$ which was pasted in.
\smallskip

{\bf IV. Tunnelling the hole to a nonadjacent sector}.\newline
\indent
If the segment $\tau$
is located in a sector $\Omega_i$ where $i\neq -1,0,1$ we first ``push
the hole'' either to $\Omega_i$ or to its neighbor. This reduces
the situation to one of the cases I--III.
We push the hole to the sector $\Omega_{j}$,
where $j=i$ if $i$ is even and $j=i-1$ if $i$ is odd. To do this
we change the identifications of the boundary segments of the
sectors $\Omega_0, \dots, \Omega_j$ as indicated at
Figure~\ref{fig:neighborhood:of:a:hole}, case IV.
The operation of
tunnelling the hole to a nonadjacent sector does not change the
area of the surface.
\medskip

\subsection{Shrinking a Pair of Nonadjacent Holes}
\label{ss:shrinking:nonadjacent:holes}

To simplify the construction above we assumed that the path
$\tau$, along which we have transported the hole, is a geodesic
segment which is shorter
than the shortest saddle connection. It is quite
clear how to generalize the  construction in several ways.

A geodesic segment $\tau$ transversal to $\vec{v}$ can be replaced by a
path $\tau$ transversal to $\vec{v}$; the result of the surgery depends
only on the homotopy class of $\tau$, where the homotopy is performed
with the fixed endpoints; the path remains  transversal
to $\vec{v}$ during the homotopy; the homotopy is performed inside
the surface with punctured singularities.

We can invert our construction: having a hole located close to a zero
we can move it to the zero.

The segment (path) $\tau$ need not be necessarily short: an appropriate
surgery can be performed using a finite covering of $\tau$ by
domains $D_{\Delta}(P_0), D_{\Delta_1}(P_1), \dots,\! D_{\Delta_s}(P_s)$.
However, we require that $\tau$ is either a segment parallel
to $\vec{v}$ or a segment (path) transversal to $\vec{v}$. We also
 require that $\tau$ does not have self-intersections, and, moreover,
that it does not make ``almost loops'' when a pair of points located
far one from another with respect to parameterization of $\tau$ by
length occur to be close to each other on the surface.

Finally, we can consider a situation when $\tau$ is a broken
line of geodesic segments joining singularities
$P_0, P_1, \dots, P_s$. We can
consecutively move the hole along the first segment $\tau_{01}$
from singularity
$P_0$ to singularity $P_1$; perform if necessary
tunnelling of the hole from one sector at $P_1$ to another;
move the hole along the second segment $\tau_{12}$
from singularity
$P_1$ to singularity $P_2$; etc, and finally bring the hole
to the appropriate sector of $P_s$.

We can now return to our original problem.  We have a surface $S_1$ with two small holes which we want to shrink {\it metrically}. We move one hole
to another and then apply the construction of shrinking a pair of adjacent holes.

Note that the operation of moving a hole from $P_0$ to $P_1$  typically
changes the area of the surface. However, the change of the area
is at most the product of the length of circumference of the hole $\gamma$
by the length of the path $\tau$.

\section{Constructing Surfaces with Homologous Closed Saddle
Connections}
\label{s:constructions:zero:to:itself}

Our objective now is to  reverse the shrinking operation to show how to  build
homologous saddle
connections with associated  configuration
$(J,a_i',a_i'',b_k',b_k'')$ out
of simpler  surfaces. We shall start  with a collection  of closed
surfaces  and   via   a   ``figure  eight''  and  ``creating   a pair of holes
construction'' build a  closed  surface with  a  curve or set  of
homologous curves joining a zero to itself.

Fix $0<r<1/4$ and let
$\cF'$ denote a compact simply connected subset of
$\cH_1(\sqcup\alpha_i')$ such that for  $S'\in\cF'$ every saddle
connection has
 length at least  $2\epsilon^r$.
By Lemma~\ref{lemma:short:saddle:connections} its complement has volume   $O(\epsilon^{2r})$.
Choose a homology basis $\beta_i$ that is valid for all $S'\in \cF'$.

\begin{figure}[hbt!]
%
%
\includegraphics{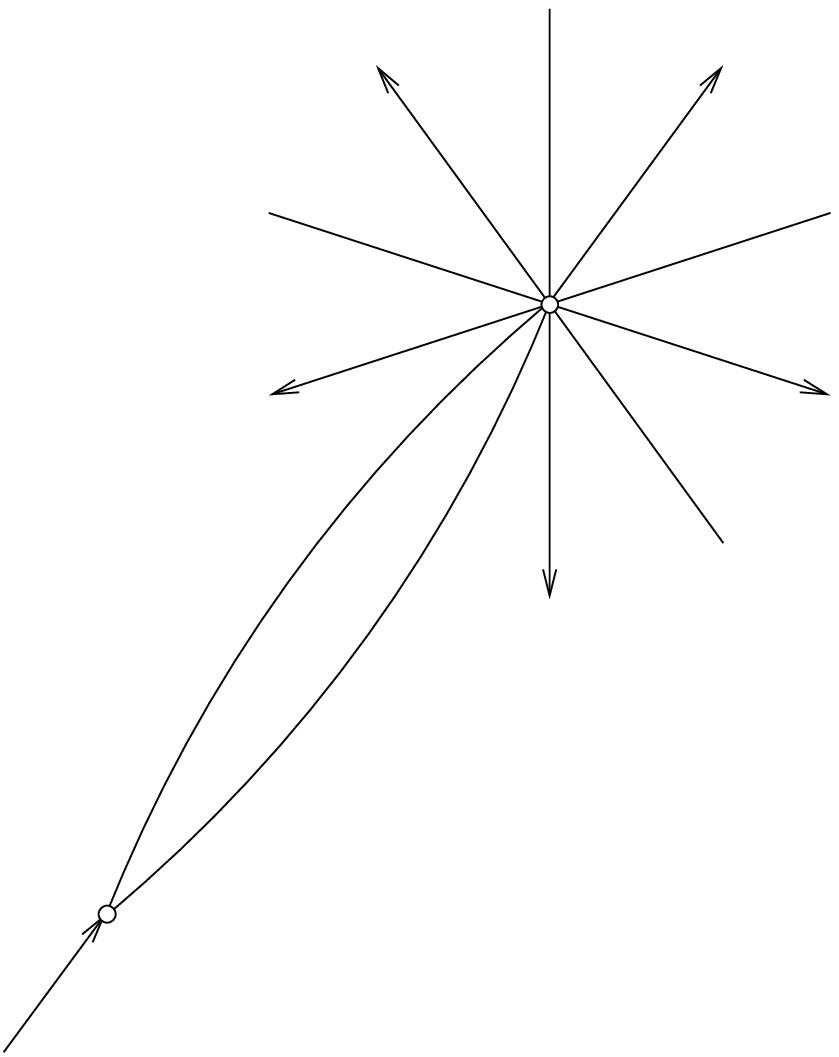}
\includegraphics{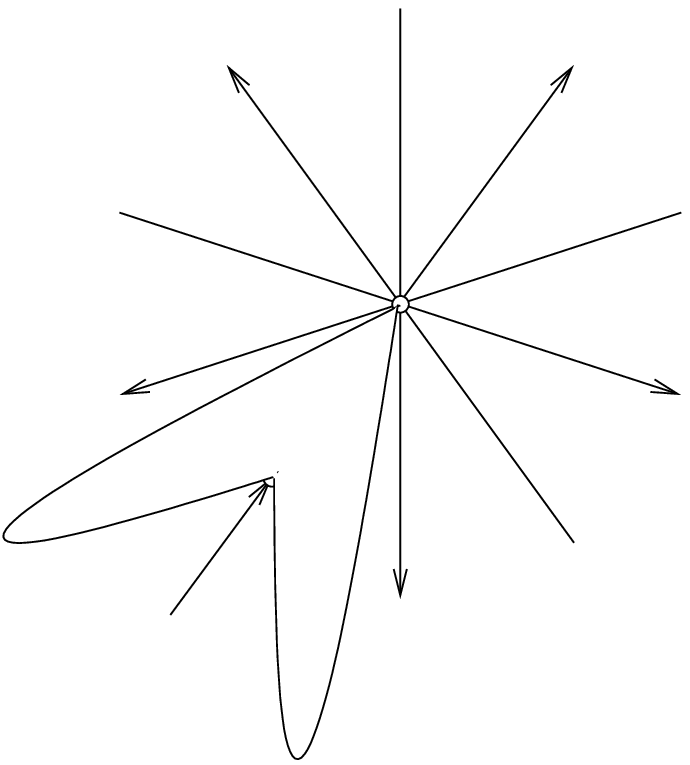}
\includegraphics{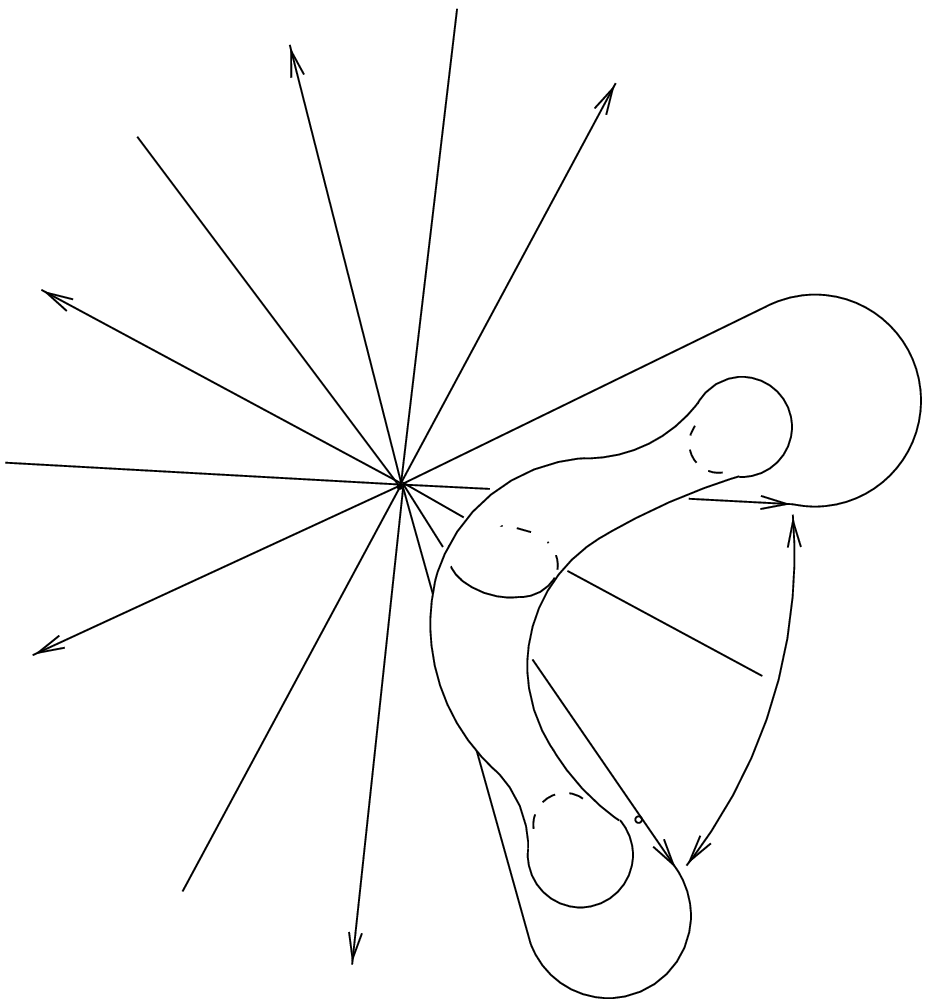} \vspace{115bp}
\begin{picture}(0,0)(0,0)
\put(0,0){\begin{picture}(0,0)(0,0) \put(70,14){$2\pi(a'+1)$}
\end{picture}}
\end{picture}
\caption{
\label{pic:figure:8:constr} Here we applied the figure eight
construction with $a'=0$ and $a''=4$ to a zero of order
$a=a'+a''=4$, and then we glued in a cylinder.
The degree of the resulting zero is $6$.
}
\end{figure}

\subsection{Figure Eight Construction}
\label{ss:figure:8}

Given a  surface $S_i'$,  a zero (or a marked point)
$z_i$ of  order $a_i$, and  a pair of numbers $a'_i+a''_i=a_i$ we
may break  up the zero and perform the  slit construction. We may
then identify the two endpoints of the slit. This gives a surface
with two  circular  boundary  components $\gamma'_i, \gamma''_i$,
which are joined at a point, see Figure~\ref{pic:figure:8:constr}.
The two circular boundary components
separate the total  angle at the point into angles $2\pi(a_i'+1)$
and $2\pi(a_i''+1)$. By  convention  on the choice of $\gamma'_i,
\gamma''_i$   (see   Section~\ref{ss:slit:construction})  turning
around  $z_i$  in a  clockwise  direction  from  $\gamma''_i$  to
$\gamma'_i$  inside  the  surface  $S'_i$  we  turn by the  angle
$2\pi(a'_i+1)$.

Recall that  we may  break the zero of order  $a_i$ up in $a_i+1$
ways. Thus  the  figure  eight  construction is not uniquely defined;
the construction gives $a_i+1$ surfaces.

\subsection{Creating a Pair of Holes}

In this construction we
 are  given a pair of distinct points and a
vector.   We shall ``open'' the surface at the two points  using the
vector  to  build  a  surface with two boundary  components.
The
construction depends on choosing a path between the two points, and so  is not quite canonical.
Suppose then that $z$ and $w$ are zeroes of orders $b',b''\geq 0$ on a
surface  $S'$,   $\gamma$   is   a   vector  of  length
$\epsilon$.  Let $\rho$
be the shortest path joining $z$ and $w$. It consists of saddle connections. On a set
of surfaces of full measure none of the saddle connections are in the direction $\gamma$.
The endpoints of the segments are located
at conical points.  We need first to estimate the length of $\rho$
and establish a few simple properties of the path. We use
Corollary 5.7 from paper~\cite{Masur:Smillie} (in which
we adjusted the notations).

\begin{StarProposition}[H.~Masur, J.~Smillie]
Let $S$ be a closed flat surface with diameter $d$. Let $\gamma$ be any geodesic
segment joining singularities whose length is minimal among all geodesic
segments that join singularities. Assume that $\gamma$ is not a simple
closed curve which bounds a metric cylinder. Let $\Delta$ be the length
of the shortest geodesic element $\beta$ joining singularities not equal
to $\gamma$. Then there is a constant $C$ depending
only on the genus of the surface $M$ such that $d\le C/\Delta$.
\end{StarProposition}

From the  assumption on $\cF'$, and the above Proposition
$|\rho|\leq C\epsilon^{-r}$ for a new constant $C$.
Since  $\rho$ is a shortest path it does not have any self-intersections, so a conical point
$P_i$ may have at most two segments of $\rho$ adjacent to it.
If $\rho$ comes within $\epsilon+\epsilon^{2r}$ of a singularity other than
at its endpoints, we may
modify $\rho$ choosing a broken line
going through that singularity, increasing the
length of $\rho$ by at most $2(\epsilon+\epsilon^{2r})$.
Since the shortest loop on the surface has length at least
$\epsilon^r$, and there are at most $2g-2$ conical points,
for small $\epsilon$, after at most $2g-2$ replacements we can
assume that the new path $\rho$ does not pass within $\epsilon+\epsilon^{2r}$
of any singularity except at its endpoints. The  total increase in length is
at most $2g(\epsilon+\epsilon^{2r})$ and so its  length
is bounded by $C\epsilon^{-r}$ for a slightly different $C$.
(The path might now be not shortest.)

\begin{proposition}
\label{prop:trans2}
\begin{enumerate}

\item Any line segment of length $\epsilon+\epsilon^{2r}$ in direction $\vec{v}$ through any point on $\rho$ does not meet $\rho$ again.
\item For small $\epsilon$ and any saddle connection $\beta$, the algebraic intersection
$|\rho\cdot \beta|\leq |\beta|\epsilon^{-2r}$.
\end{enumerate}
\end{proposition}

\begin{proof}
To prove (1) notice that if the intersections points $Q_1,Q_2$ of the line segment
with $\rho$ were further than $(2g+1)(\epsilon+\epsilon^{2r})$ apart on $\rho$, we could shorten $\rho$ by more than $2g(\epsilon+\epsilon^{2r})$ following the line segment between $Q_1$ and $Q_2$ instead of $\rho$, which is impossible.  If $Q_1,Q_2$ are within $2g(\epsilon+\epsilon^{2r})$ of each other on $\rho$ we have a loop of length at most $(2g+1)(\epsilon+2\epsilon^{2r})
<2\epsilon^r$ for small $\epsilon$, which is also a contradiction.

The proof of (2) is similar.
\end{proof}

Now returning to the problem of creating a pair of holes, suppose
we are given a pair of  separatrices  leaving  $z$ and $w$ in direction
$\gamma$.  The    notion of $z'$ as the point  $z+\gamma$ makes sense for small $\gamma$.
Perform the figure eight construction using
the points $z$ and $z'$ creating a surface with a pair of holes joined at a point $P$.

We can
now transport a hole from $P$ to $w$ along $\rho$ as  was
described in the previous section. Conclusion 1) of the Proposition says that the hole transport is well-defined.
We get a surface with two holes.
We may then tunnel  the hole at $w$ so that it is adjacent to the given separatrix.
Note that the operation of creating a pair of holes  typically
changes the area of the surface. However, the change of the area
is at most the product of the length of circumference of the hole $\gamma$
by the length of the path $\rho$.
Our estimates show that  the change of the
area is bounded by $C\epsilon^{1-r}$.

We note that the resulting  flat  structure  depends  on   the   choice   of
the separatrices leaving each zero, on the homotopy class of the path $\rho$
as well as on the
other given data.
From (2) we see that for any saddle connection $\beta$, the
 length of $\beta$ on the new  surface is at least
$|\beta|-\epsilon|\beta\cdot\rho|\geq |\beta|(1-\epsilon^{1-2r})$ and since
$|\beta|\geq \epsilon^{2r}$ on $S'$, on the new surface $\beta$ must have length at
least $\epsilon^r$.

\subsection{Admissible Constructions}
\label{ss:admissible:constructions}

We formalize the above construction as follows.
Consider the following data:

\begin{itemize}
\item
For  $1 \le  i  \le p$: surfaces  $S_i'  \in \cH(\alpha_1')$  with
$\area S_i' < 1$ and no saddle connections shorter than $2\epsilon^r$
\item
A vector  $\gamma  \in  \reals^2$,  with  $|\gamma| < \epsilon$.
\item
For $1  \le  i \le p$: a pair of points $z_i, w_i \in S_i'$ which
are zeroes of orders $b'_i,b''_i$  or  marked  points. The points
may  coincide.  If  they  are  distinct,  we  are  also  given  a
broken line  $\rho_i$ joining $z_i$ and $w_i$. For  any      $\beta$ including those in a fixed homology basis,
\begin{equation}
\label{eq:hol:bound}
|\beta\cdot\rho_i|\leq |\beta|\epsilon^{-2r}
\end{equation}
where $\beta\cdot \rho_i$
is the algebraic intersection number.
\item $|\rho_i|\leq C\epsilon^{-r}$, where $C$ is a constant.
\item
If  $w_i=z_i$  is a  zero  of order  $a_i$, then  a pair  of
numbers $a'_i,a_i''$ with $a_i=a'_i+a_i''$.
\item
A subset $J=\{i_1,\dots,i_q\}\subset\{1,\dots,p\}$,  where  $0\le
q\le p$. By convention, when $J=\emptyset$ we let $q=0$.
\item
For $1  \le  j  \le  q$:  a number $h_{i_j} > 0$ such that $ \sum
h_{i_j} |\gamma | + \area S_{i_j} < 1$.
\item
For  $1  \le  j \le q$: a number $t_{i_j}$ with $0 \le t_{i_j} <
|\gamma|$.
\end{itemize}

Given  the  above data  we  construct  a  flat  surface  $S$ with
homologous saddle connections   of  multiplicity  $m$  with  associated  data
$(J,a_i',a_i'',b_k',b_k'')$.
Namely,  for
each $i$,  if $z_i=w_i$ we  use the figure eight construction and
if $z_i\neq w_i$ we use the construction of creating a pair of holes giving
two
boundary    curves    $\gamma_i'$    and
$\gamma_i''$.  For  each  $i$  either   glue
$\gamma_i'$     directly     to    $\gamma_{i+1}''$     or,    if
$i=i_j\in\{i_1,\dots,i_q\}=J$  glue  in  cylinders  $C_{i_j}$  of
``width''   $|\gamma|$,  ``height''   $h_{i_j}$   and   ``twist''
$t_{i_j}$   to    connect   $\gamma_i'$    to
$\gamma_{i+1}''$ ($i+1$ considered mod $p$).

The resulting surface $S$  also  depends on the data
$h_{i_j},t_{i_j}$ and
$\rho_i$.
Denote  by $S'$  the union  of  the pieces  $S_i'$,  $h$
  the vector of  $h_{i_j}$,  $t$  the vector of
$t_{i_j}$ and $\rho$ the collection of  $\rho_i$.
Then we have an
{\em admissible construction}

$$ 
(S',J,a_i',a_i'',b_i',b_i'',\rho,\gamma,h,t)\to S
$$ 

A relative homology  basis for $S$ is given by:
\begin{itemize}
\item  A  relative  homology  basis  of  curves $\beta$ for  each
component $S_i'$ of $S'$.
\item If one of the $z_i,w_i$ is a marked point, the curve $\rho_i$.
If both are marked points,
then in addition to $\rho$ a curve from some singularity on $S_i'$ to
$z_i$ (or to $w_i$) unless $S_i'$ is a torus.
\item  A    closed
curve in the homology class of $\gamma_i$.
The holonomy along this curve is $\gamma$.
\item For each cylinder a curve joining its boundary singularities.
The holonomy   of
this curve is $h+it$.
\end{itemize}

We summarize the discussion as follows.
\begin{lemma}
\label{lm:admissible:to:itself}
Suppose that we\newline
--- either use at least one  construction of creating a pair of holes;
or\newline
--- if we use only  figure  eight constructions, we use at  least
one cylinder.
\newline

Then the resulting  surface  $S$ is a nondegenerate flat
surface with no saddle connections other than the $\gamma_i$ shorter than $\epsilon^r$
and thus it lies in $\cH^{\epsilon,\epsilon^r}(\alpha,\cC)$.

The holonomy vectors  $v_S(\beta)$ on $S$  of those  curves
$\beta$, which correspond to the original homology basis on $S_i'$ are given by
\begin{equation}
\label{eq:hol}
v_S(\beta)=v_{S'}(\beta)+(\beta\cdot\rho_i)\gamma
\end{equation}

\end{lemma}

\begin{figure}[ht]
%
\includegraphics{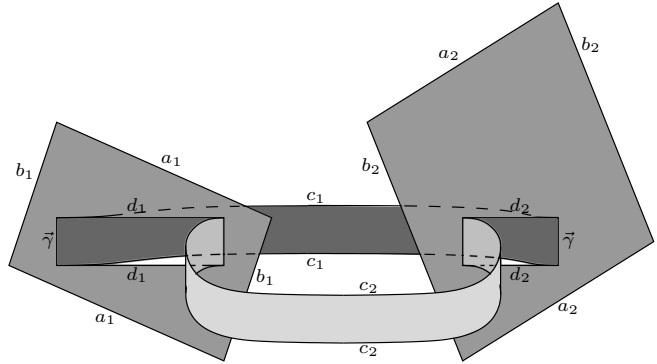}
%
%
\begin{picture}(0,0)(0,0)
\put(0,0)
 {\begin{picture}(0,0)(0,0)
 \put(-65,-60){$\scriptstyle a_1$}
 \put(-90,-120){$\scriptstyle a_1$}
 \put(-120,-65){$\scriptstyle b_1$}
 \put(-29,-105){$\scriptstyle b_1$}
 \put(40,-20){$\scriptstyle a_2$}
 \put(85,-116){$\scriptstyle a_2$}
 \put(11,-63){$\scriptstyle b_2$}
 \put(94,-17){$\scriptstyle b_2$}
 \put(-10,-74){$\scriptstyle c_1$}
 \put(-10,-99){$\scriptstyle c_1$}
 \put(10,-108){$\scriptstyle c_2$}
 \put(10,-132){$\scriptstyle c_2$}
 \put(-78,-104){$\scriptstyle d_1$}
 \put(-78,-78){$\scriptstyle d_1$}
 \put(67,-78){$\scriptstyle d_2$}
 \put(67,-104){$\scriptstyle d_2$}
 \put(-110,-91){$\scriptstyle \vec{\gamma}$}
 \put(87,-91){$\scriptstyle \vec{\gamma}$}
\end{picture}}
\end{picture}
\vspace{140bp}
\caption{
\label{pic:parallel:constr}
In this example we apply the {\it construction of creating  a pair of holes}
simultaneously to two tori (represented by large parallelograms
with identified opposite sides) and then glue in a pair of
cylinders (represented by the narrow stripes with identified
opposite sides).
}
\end{figure}

\begin{example}
In the example presented at Figure~\ref{pic:parallel:constr} we
take two flat tori. In the picture they are represented by the
large parallelograms with identified opposite sides $a_1$ to
$a_1$, $b_1$ to $b_1$, and $a_2$ to $a_2$, $b_2$ to $b_2$. In
fact, the zeroes $z_1$, $w_1$, and $z_2$, $w_2$ are fake here: all
these points are just  marked points of the tori. In this
example we choose the vector $\vec{\gamma}$ to be vertical.

We cut a rectangular domain out of the interior of each of the two
tori in such way that the vertical side of each of the two removed
rectangles is equal to $\vec{\gamma}$. Then we identify the
horizontal sides $d_1$ to $d_1$ and $d_2$ to $d_2$. We get a pair
of tori, each provided with a pair of holes.

Finally, we join the hole $\overline{w_1 w'_1}$ on the first torus
to the hole $\overline{z_2 z'_2}$ on the second torus by a
cylinder (in the picture the cylinder is represented by a narrow
stripe with identified opposite sides $c_1$ to $c_1$); we also
join the hole $\overline{w_2 w'_2}$ on the second torus to the
hole $\overline{z_1 z'_1}$ on the first torus by a cylinder (in
the picture this second cylinder is represented by a narrow stripe
with identified opposite sides $c_2$ to $c_2$).

We get a flat surface of genus $g=3$ with four simple zeroes. By
construction the surface contains two cylinders filled with closed
regular homologous geodesics with holonomy  $\vec{\gamma}$ (these cylinders
are represented by narrow stripes with horizontal sides $c_1$ and
$c_2$). In other words the saddle connection represented by
$\gamma$ has multiplicity $2$. Note that the two cylinders would
survive under any sufficiently small deformation of the surface,
which means that this configuration is generic.
\end{example}


Now we return to  the discussion at the beginning of part II. We are
given a  set of homologous saddle connections
$\gamma_i$ and a surface $S\in\cH^{\epsilon,3\epsilon}(\alpha,\cC)$.  We
may cut $S$ along the $\gamma_i$ to form a collection of surfaces $S_i$, each with a pair of boundary components. We may then choose a path $\tau_i$ in each $S_i$ joining the cone points and shrink the boundary circles giving  surfaces $S_i'$ in  $\cH(\alpha_i')$.  An analysis similar to that given in the section on creating a pair of holes  shows that we may choose the $\tau_i$ so that the resulting surfaces lie in $\cF'$.

Our main objective is  now to show that the degeneration of a
surface to a point on the principal boundary  and the admissible
construction of a surface from a point on the principal boundary
are essentially inverse operations.  This will allow us to compute
the constants for the counting problem. The complicating
issue  is  that the degeneration  depends on the paths $\tau_i$,
while the admissible construction depends on the choice of $\rho_i$
and these curves might be not the same.  In particular, since lengths
of curves may change by $\epsilon$ during the construction, the shortest
curve may change.
However, changing from one path to
another
amounts to a change of basis and can be described.

Recall that we have fixed a homology basis $\beta_i$ valid for all $S'\in \cF$.  The holonomy along the basis allows one to define a metric $d(\cdot,\cdot)$ on $\cF'$. The length of the $\beta_i$ satisfy
$|\beta_i|\leq C\epsilon^{-r}$ for a fixed constant $C$. Fix a
pair of zeroes of a given $S_i'$ and a path  $\rho_i$ joining them which does
not pass within $\epsilon+\epsilon^{2r}$ of a singularity or of itself.
We may find a ball of radius $\epsilon^{2r}$ about $S_i$ in the metric $d$ such that for all
other surfaces in the ball we may perform the admissible construction   using the same paths  $\rho_i$ since for any surface in the ball the corresponding paths  do not come within      $\epsilon$ of a singularity or themselves. Thus we can find a finite collection of balls
$B(x_j,r_j)$ centered at such points $x_j$ of radius $\epsilon^{2r}$ that  cover $\cF'$
such that  for  each   $B(x_j,r_j)$,  for
each component $S_i'$ of $S'$ for which $z_i\neq w_i$ there is a fixed
curve  $\rho_{i,j}$ joining $z_i$ to $w_i$ and so that we may form the admissible construction for each $x\in x_j$.

We may use the cover  to define a
partition of $\cF'$ by a finite collection of sets $U_j\subset
B(x_j,r_j)$.
Now we  can give a well defined  assignment
$$ 
\sigma:\cF'\to\cH_1^{\epsilon,\epsilon^r}(\alpha,\cC)
$$ 

Namely each $S'$ belongs to a unique   $U_j$.
Then use  $\rho_{i,j}$ to assign to $S'$ a closed surface
$S=\sigma(S')$ by the admissible construction.

As  in  the  previous  section,  we  have the  possibility  of  a
construction for  $-\gamma$  coinciding  with a construction with
$\gamma$ and also  the possibility of stratum interchange.
However, even without considering these symmetries, it is possible
because of the change of transversals from one set $U_j$ to another,
 that the construction is
not $1-1$.
We prove

\begin{proposition}   \label{proposition:injective:onto}  There  is   $E'\subset \cF'$ such  that
\begin{itemize}
\item $\Vol(\cF'\setminus E')=O(\epsilon^{1-3r})$,
\item On a subset of $E'$ of full volume,  $\sigma$ is $1-1$ except for the $\gamma\to -\gamma$ and stratum  interchange symmetries.
\item $\Vol(\cH_1^{\epsilon,\epsilon^r}(\alpha,\cC)\setminus \sigma(E'))=o(\epsilon^2)$
\end{itemize}
\end{proposition}
\begin{proof}
   Define  $E'$ to be the set of points which are at least $C\epsilon^{1-3r}$, from the boundary of any
ball.
This   means that if $S'\in E'\cap U_j$ then the $C\epsilon^{1-3r}$
neighborhood of $S'$ is contained in $U_j$. It also implies
$$ 
\Vol(\cF'\setminus E')= O(\epsilon^{1-3r})
$$ 

We prove $\sigma$ is $1-1$ on almost all of $E'$ except for  possible stratum interchange and $\gamma\to-\gamma$ symmetries.
  If $S_1'\neq S_2'$ belong to the same $U_j$
then the
collection of surfaces $\sigma(S_1')$ and $\sigma(S_2')$  are
constructed with the same broken curve and yet differ on holonomy for a set of basis curves. The surfaces can be isomorphic only if the isomorphism does not
preserve the basis. If the set of configuration curves are preserved then
 we have the
stratum
interchange or $\gamma\to -\gamma$ symmetry. If they are not preserved, then
  $S'$ contains a curve of length at most $\epsilon$ and this possibility has been ruled out by the definition of $\cF'$.    The remaining  possibility is if
for a fixed $S'$ two or more of the remaining surfaces in $\sigma(S')$ are
isomorphic by an isomorphism that preserves the set of curves in the configuration. But such an isomorphism
must arise from an automorphism of $S'$ and such $S'$ have volume $0$.

  If $S_1',S_2'\in E'$ belong to
different $U_j$, then their distance apart in the metric $d$ is at
least  $2C\epsilon^{1-3r}$ so their holonomy on some basis curve $\beta_i$ differs by at least $2C\epsilon^{1-3r}$.
However,  by (\ref{eq:hol}) the change in holonomy in passing to $\sigma(S')$ from  $S'$ and $\sigma(S'')$ from $S''$ is at most $\epsilon|\beta_i\cdot \rho_{i,j}|\leq \epsilon|\beta_i|\epsilon^{-2r}\leq C\epsilon^{1-3r}$ and so again  $\sigma(S_1')$ and $\sigma(S_2')$
again differ  on $\beta_i$ and we have the same conclusion as before.

We now   show the third statement.
Given   the
transversals $\tau_i$,  let $E$ be the
subset of  $\cH_1^{\epsilon,3\epsilon^r}(\alpha,\cC)$ consisting of those $S$ whose degeneration yields surfaces $S'$ which are at least $3C\epsilon^{1-3r}$
from the boundary of any ball.  We have
$$ 
\Vol(\cH_1^{\epsilon,3\epsilon^r}(\alpha,\cC)\setminus E)=o(\epsilon^2)
$$ 
For some $j$, the  resulting surface $S'\in U_j\cap B(x_j,r_j-3\epsilon^{1-3r})$.
However, since it may not be the case that
$\tau_i=\rho_{i,j}$,  it  may not  be  true
 that
$ \sigma(S')=S$.
Now we claim     that
   %
   %
there      is some     $S''\in  E'$      such      that
 $S=\sigma(S'')$.
To
prove  the claim, notice that since  $B(x_j,r_j)$  is
a ball    and  since
$|\beta\cdot\tau_i|\leq C\epsilon^{-3r}$,  $|\beta\cdot\rho_{i,j}|\leq
C\epsilon^{-3r}$, and $|\gamma|\leq\epsilon$,
  the  point  in $\reals^{n'}$ defined  by  the
collection of vectors
$$ 
v_{S'}(\beta)+ (\sum_i(\beta\cdot \tau_i)
-(\beta\cdot\rho_{i,j})) \gamma
$$ 
as $\beta$ runs over the homology basis, corresponds to the
holonomy for a
unique    $S''\in B(x_j,r_j-C\epsilon^{1-3r})$.
so that  $S''\in E'$. But  now
by (\ref{eq:hol}),
$\sigma(S'')$
has  the same
holonomy as $S$  along  the  basis and so is  isomorphic  to  it, proving the claim. The third statement now follows from the estimate
$\Vol(\cH_1^{\epsilon,\epsilon^r}(\alpha,\cC)
\setminus \cH_1^{\epsilon,3\epsilon^r}(\alpha,\cC))=o(\epsilon^2)$, which follows from Lemma~\ref{lemma:short:saddle:connections}.
\end{proof}

\section{Computing the Siegel---Veech Constants for Connected Strata}
\label{s:zero:to:itself}

Recall that $q$ is the  number  of cylinders we shall attach.  Let
$n_i =  \dim_{\reals}  \cH(\alpha_i)$.  We  also  define $d_i$ to
be equal to $n_i$ except in the following situations. If $S_i$ has
a  single  marked  point which is not a zero and $S_i$ is not the
torus, we have $d_i=n_i+2$.
If there is a pair of  marked points,
and the surface is not the torus $d_i=n_i+4$. If the surface is a
torus then $d_i=n_i+2$. Then,
$$ 
\dim_{\reals} \cH(\alpha) = 2q+ 2+ \sum_{i=1}^p d_i,
$$ 
where  $\alpha$  is   the   stratum  such  that  $S  \in
\cH(\alpha)$.

\subsection{Computation of the Constants}
\label{ss:computation:of:c:II}
Recall that $d\nu(S)$ is the measure in  $\cH(\alpha)$ and
$d\nu(S_i')$ is  the  measure   in   $\cH(\alpha_i)$.
As in Part I we will obscure the distinction between $\gamma$
a saddle connection and its holonomy so that we will use $\gamma$ to refer to a vector as a variable of integration in a disc of radius $\epsilon$.
We let $M$ be the number which computes of fixed data
$(S',J,a_i',a_i'',b_k',b_k'',\gamma,h,t)$ the number
of $S=\sigma(S)$ that can be built. This number will be computed in the next sections. We  first
assume  that $\cH(\alpha)$ is connected. If we choose a homology
basis $\beta_i$ for the $S_i'$, and let $hol_S(\beta_i)$ be the
holonomy of $\beta_i$ on the surface $S$, then the measure $d\nu$
restricted to the image  $\sigma(E')$ under the admissible
construction is the product measure of $Md\gamma\ dt_idh_i$ and
$\prod d(hol_S(\beta_i))$. By~\eqref{eq:hol}
\begin{displaymath}
\prod d(hol_S(\beta_i))d\gamma=d\nu'd\gamma
\end{displaymath}
and again since $|\gamma|\leq\epsilon$ we can write
$$ 
d\nu(S) = d\gamma \, \prod_{i=1}^p d\nu'(S_i') \,  
\prod_{i=1}^q dh_idt_i
$$ 
except for  the following situations.  If $z_i$ is a marked point
and $w_i$ is a zero then we have the additional  factor $dz(S_i)$
in  the  term  $d\nu_i$. If $z_i$ is a zero and $w_i$ is a marked
point we have $dw(S_i)$. If $w_i=z_i$ is a marked point,  we have
$dz(S_i)$ and if $z_i\neq w_i$ are both marked points we have the
factor $dz(S_i)\ dw(S_i)$.

Note that we could also  consider  strata with one or two  marked
points.  The stratum  would  be a fiber  space  over the  stratum
$\cH_1(\alpha)$. The fiber over a point would be the flat surface
if  there  is one marked point  or  product of flat surfaces  if
there  are  two,  the  flat  surface  representing that point  in
$\cH_1(\alpha)$. The measure on  the  fiber space is then locally
the product measure of  the  base with the area  (or  product of area
measures) on the fiber.

Now for each $S_i$ we either perform the figure eight or creating a pair of holes  construction. In the latter case we remove or add a parallelogram $P_i$.
Thus  $\area  S = \sum_i  (\area  S_i) +  h_i
|\gamma|\pm\area P_i$, (where $\area P_i=0$ if we perform the figure eight construction).
If   $S   \in
C(\cH_1^{\epsilon,3\epsilon^r}(\alpha,\cC))$  then the normalization says that
$|\gamma| \le \epsilon \sqrt{(\area S)}$. Hence,
$$ 
|\gamma|^2 \le \epsilon^2 (\sum_i (\area S_i + h_i |\gamma|\pm\area P_i)),
$$ 
i.e.
$$ 
\sum_i h_i \ge \frac{|\gamma|}{\epsilon^2} - \frac{\sum_i(\area
  S_i\pm\area P_i)}{|\gamma|}.
$$ 
The expression on the right is positive when  $|\gamma| > \epsilon
\sqrt{ (\sum_i \area S_i\pm\area P_i)}$.

Let
$$ 
W = \prod_i \Vol(\cH_1(\alpha_i))
$$ 
For $s  > 0$, let $D(s) = \{ (r_1, \dots, r_p) \st \sum_i r_i^2 <
s \}$. Then using the fact that $\area P_i=O(\epsilon^{1-r})$,
Proposition~\ref{proposition:injective:onto},
the fact that
$\Vol(\prod_i \cH_1(\alpha_i)\setminus \cF')=O(\epsilon^{2r})$ and $|\gamma|\leq \epsilon$,
we have

\begin{multline*}
\nu (C(\cH_1^{\epsilon,\epsilon^r}(\alpha)))=WM\cdot\\
\Bigg(
\int_{D(1)} \biggl(\prod_{i=1}^p r_i^{d_i-1} \, dr_i
\biggr)
\int_{ \gamma \in B(\epsilon\sqrt{\sum_i r_i^2})} \int_{0 \le
\sum_i h_i \le
  \frac{1 - \sum_i r_i^2}{|\gamma|}} \prod_{i=1}^q \biggl(
\int_0^{|\gamma|}
    dt_i \, dh_i \, \biggr) d\gamma +  \\
 +  \int_{D(1)}
\prod_{i=1}^p r_i^{d_i-1} \, dr_i
\int_{ \gamma \in B(\epsilon) \setminus B(\epsilon\sqrt{\sum_i
    r_i^2})} \int_{\frac{|\gamma|}{\epsilon^2} - \frac{\sum_i
    r_i^2}{|\gamma|} \le \sum_i h_i \le
  \frac{1 - \sum_i r_i^2}{|\gamma|}} |\gamma|^q \prod_{i=1}^q dh_i \, d\gamma
\Bigg)+\\
+o(\epsilon^2)  =
\end{multline*}
\begin{multline*}
= WM\cdot\Bigg(
 \int_{D(1)}
 \prod_{i=1}^p r_i^{d_i-1} \, dr_i
\int_{ \gamma \in B(\epsilon \sqrt{\sum_i r_i^2})} \frac{1}{q!}
\bigl(
  1 - \sum_i r_i^2 \bigr)^q \,
d\gamma + \\
+ \int_{D(1)}\hspace*{-1.5truept}
 \prod_{i=1}^p r_i^{d_i-1} \, dr_i
\int_{ \gamma \in B(\epsilon) \setminus B(\epsilon\sqrt{\sum_i
    r_i^2})} \frac{1}{q!} \bigg[\bigl( 1 - \sum_i r_i^2 \bigr)^q
-
\biggl( \frac{|\gamma|^2}{\epsilon^2} - \sum_i r_i^2 \biggr)^q
\biggr]
d\gamma
\Bigg)+\\
+o(\epsilon^2)
\end{multline*}

Performing the integral over $\gamma$ we
obtain

\begin{align*}
\nu(C(\cH_1^{\epsilon,\epsilon^r}(\alpha,\cC)))  = \hspace*{-57pt} & \\
& = \quad \frac{\pi \epsilon^2 WM }{q!} \int_{D(1)} \left(\prod_{i=1}^p
  r_i^{d_i-1} \right) \left( 1 - \sum_i r_i^2 \right)^q
\left(\sum_i
r_i^2 \right) \prod_i dr_i \\
& \quad + \frac{\pi \epsilon^2 WM }{q!} \int_{D(1)} \left(\prod_{i=1}^p
  r_i^{d_i-1} \right) \left( 1 - \sum_i r_i^2 \right)^q \left( 1
  - \sum_i r_i^2 \right) \prod_i dr_i \\
& \quad - \frac{\pi WM }{q!} \int_{D(1)} \left(\prod_{i=1}^p
  r_i^{d_i-1} \right) \int_{ \epsilon \sqrt{ \sum_i
    r_i^2}}^\epsilon \left( \frac{s^2}{\epsilon^2} - \sum_i
  r_i^2 \right)^q 2 s \, ds \prod_i dr_i +o(\epsilon^2)\\
& = \quad \frac{\pi \epsilon^2 WM }{q!} \int_{D(1)} \left(\prod_{i=1}^p
  r_i^{d_i-1} \right) \left( 1 - \sum_i r_i^2 \right)^q \prod_i
dr_i \\
& \quad - \frac{ \pi \epsilon^2 WM}{(q+1)!} \int_{D(1)}
\left(\prod_{i=1}^p
  r_i^{d_i-1} \right) \left( 1 - \sum_i r_i^2 \right)^{q+1}
\prod_i dr_i +o(\epsilon^2)\\
\end{align*}
where we have used the change of variable $u =  (s^2/\epsilon^2 -
\sum r_i^2)$ in the last line. We now make the change of variable
$x_i =  r_i^2$.  Let $b_i = d_i/2 - 1 $, so that $r_i^{d_i -1} \,
dr_i = (1/2) x_i^{b_i} \, dx_i$. Then,
\begin{multline*}
\nu(C(\cH_1^{\epsilon,\epsilon^r}(\alpha,\cC))) = \frac{\pi \epsilon^2
WM}{2^p} \int_{
\sum_i
  x_i \le 1} x_1^{b_1} \dots x_p^{b_p} \Bigg[ \frac{1}{q!}
\left( 1  -
  \sum_i x_i\right)^q \\
   - \frac{1}{(q+1)!} \left( 1 - \sum_i x_i
  \right)^{q+1}\Bigg] \, dx_1 \dots dx_p+o(\epsilon^2)
\end{multline*}
Evaluating  the  integral  from  the  inside   out  via  repeated
application of the identity
$$ 
\int_0^u  x^a  (  u  -  x)^b  \, dx  =  \frac{a!\,  b!}{(a+b+1)!}
u^{a+b+1}
$$ 
yields
\begin{align*}
\nu(C&(\cH_1^{\epsilon,\epsilon^r}(\alpha, \cC)))  = \\
& =
\frac{\pi \epsilon^2
WM}{2^p} \left[
\frac{ b_1!
    \dots b_p!}{(b_1 + \dots + b_p + p + q )!} - \frac{b_1! \dots
    b_p!}{(b_1 + \dots + b_p + p+ q+1)!} \right]+o(\epsilon^2) = \\
& = \frac{\pi \epsilon^2 WM}{2^p}  \frac{b_1! \dots
    b_p!}{(b_1 +  \dots + b_p + p+ q - 1 )!} \cdot
\frac{1}{b_1+\dots +
  b_p + p+ q + 1}+o(\epsilon^2) \\
\end{align*}
Since $b_1+\dots+b_p + p + q + 1 = n/2$, we have
$$ 
\nu(C(\cH_1^{\epsilon,\epsilon^r}(\alpha,\cC))) = \frac{\pi \epsilon^2
WM}{2^{p-1}}
\frac{b_1!
  \dots b_p!}{ n(\frac{n}{2} - 2)!}+o(\epsilon^2)
$$ 
Now taking the limit as $\epsilon\to 0$, we get
\begin{multline*}
c = \lim_{\epsilon\to 0 }
\frac{\nu_1(C(\cH_1^{\epsilon,\epsilon^r}(\alpha,\cC)))}{\pi \epsilon^2
\nu_1(\cH_1(\alpha))} =
\\
= M \cdot \frac{1}{2^{p-1}} \cdot
\frac{\prod_{i=1}^p (\tfrac{d_i}{2}-1)!}{(\frac{n}{2}-2)!} \cdot
\frac{\prod_{i=1}^p \Vol(\cH_1(\alpha'_i) )}{\Vol(\cH_1(\alpha))}
\end{multline*}

 In order to  compute the combinatorial constant $M$ let us
study in more detail how, given  the admissible construction
$(S',J,a'_i,a''_i,b'_k,b''_k,\rho,\gamma,h,t)\to S$,  the resulting flat surface $S$ is made  up of
subsurfaces.
The resulting surface  $S$ has two groups of zeroes: the
ones inherited from the components $S_i'$ without any
changes, and  the newborn ones. Every  newborn zero has  at least
one curve $\gamma_i$ passing through it; the curves $\gamma_i$ do
not pass  through the zeroes of the first  group. Actually, it is
easy to see from our construction, that the  collection of curves
$\gamma_{i_j}$ passing through a zero  of  $S$  has  the
form   $\gamma_i,  \gamma_{i+1},   \dots,   \gamma_{i+s}$,   with
consecutive indices (where as usual $p+1$ is identified with $1$,
so $p$ is followed by $1$). It is easy to  give  a description of
the nature of the corresponding set $\{i,i+1,\dots,i+s\}$.

Any newborn zero is formed by one of the following three constructions.
The reader may refer to figure 10.

$\bullet$ {\bf Type I.} There is a chain  of consecutive surfaces
$S_{j_r+1} \to S_{j_r+2}\to \dots\to  S_{j_r+s}$  glued directly.
To  each  surface  of  the  chain  we  apply   the  figure  eight
construction. The first surface of the chain $S_{j_r+1}$ is glued
to the preceding surface $S_{j_r}$  by  a  cylinder and similarly
for the last surface of the chain $S_{j_r+s}$  and its successor.
In  our  notation  we  have  $j_r+s=j_{r+1}$,  where  both  $j_r,
j_{r+1}\in J$, and $s\ge 0$.

In this case the newborn zero has order %
\begin{equation}
\label{eq:newborn:zero:I}
\sum_{k=1}^{s} (a_{j_r+k} +2),
\end{equation}
where  $a_{j_r+k}=a_{j_r+k}'+a_{j_r+k}''$.  By  convention  $s=0$
stands  for  the case  when  multiplicity $p$  is  equal to  $1$,
$J=\{1\}$, and  we have a single surface $S_1$  to which we apply
figure eight construction and glue in a cylinder.

The left drawing in figure 10 illustrates the situation when $s=2$. The figure eight construction is applied to each of the two surfaces and they are glued to each other directly. Each is glued to another surface by a cylinder.
\begin{figure}[htb!]
%
 %
 %
 \includegraphics{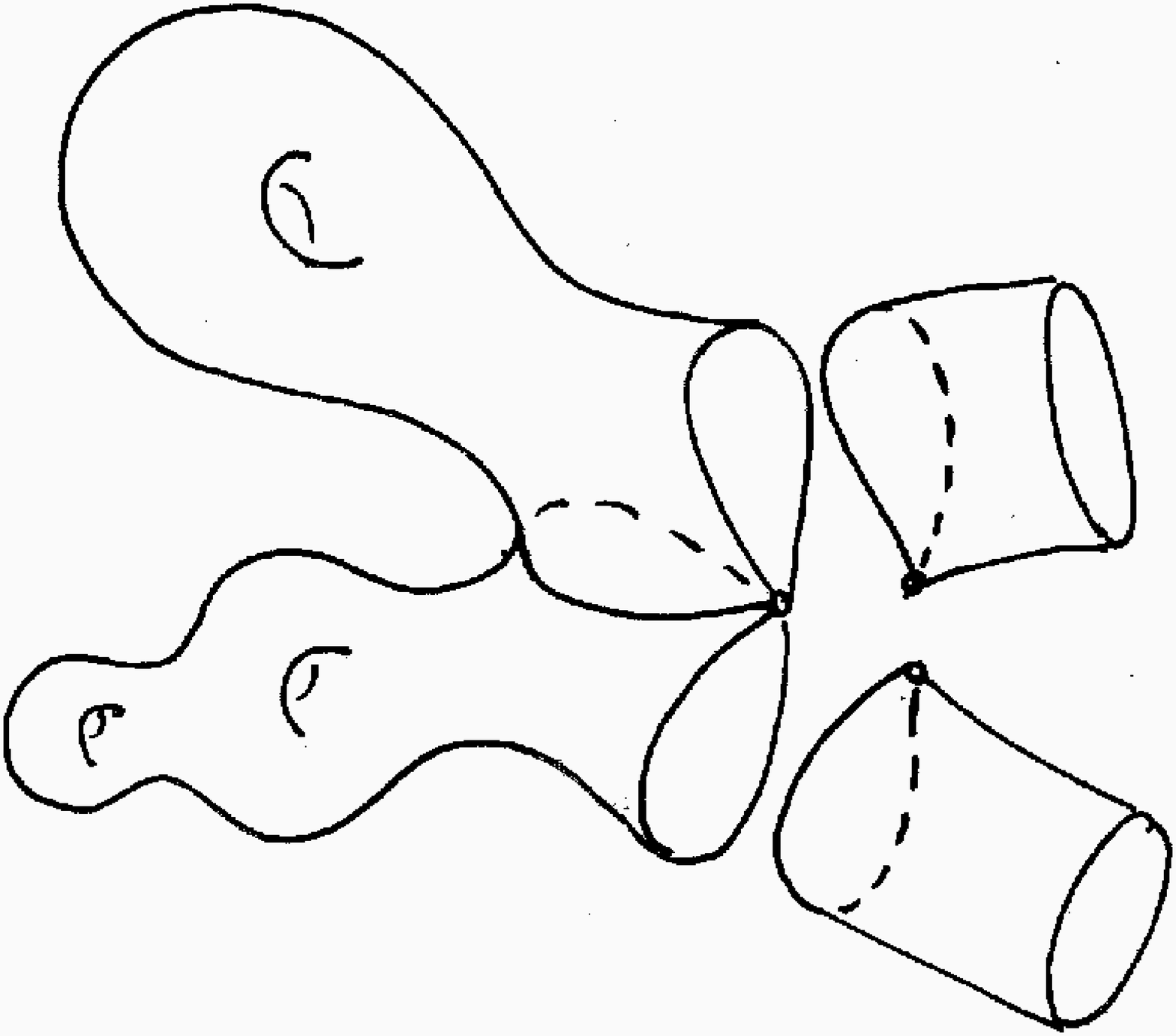}
 \includegraphics{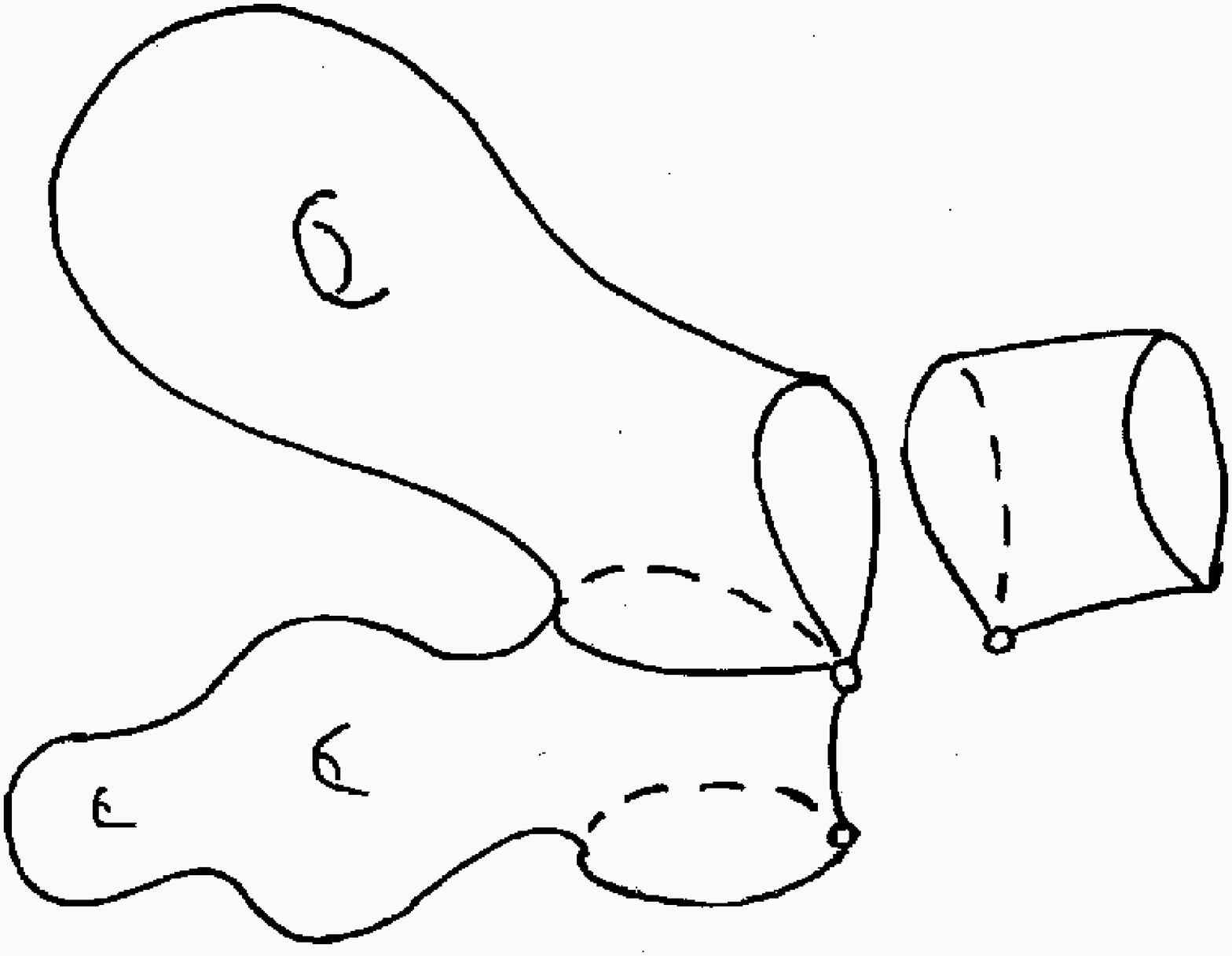}
 \includegraphics{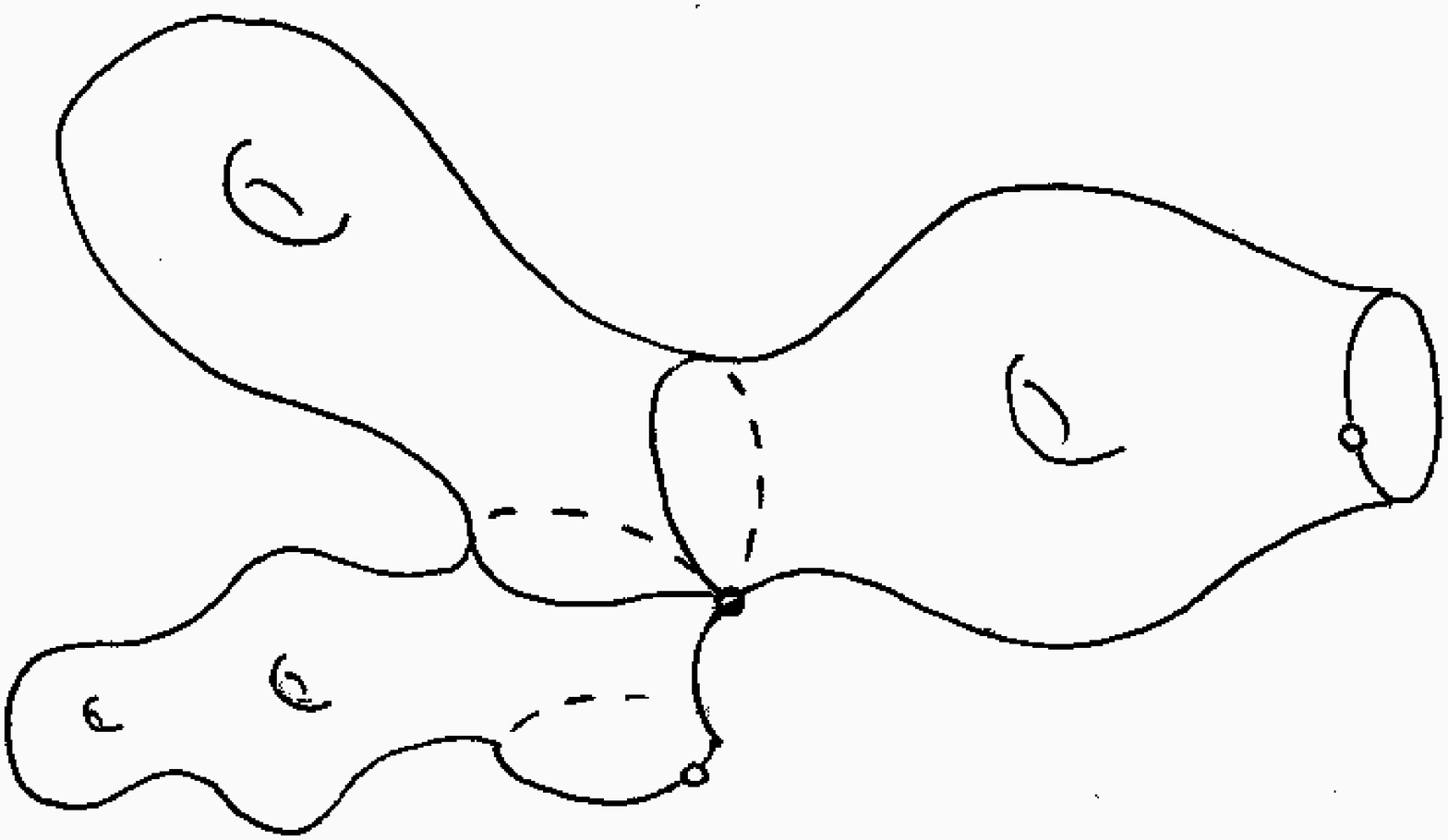}
\vspace{150bp}
\caption{
\label{pic:threetypes}
Newborn zeroes of three different types.
}
\end{figure}

$\bullet$  {\bf Type  II.}  Now we have  the  following chain  of
surfaces   $S_{j_r+1},   S_{j_r+2},  \dots,   S_{j_r+s+1}$  glued
directly. We apply  the figure eight construction to all surfaces
except  the   last,   to   which   we   apply  the
construction of creating a pair of holes.  The first surface $S_{j_r+1}$ of the chain is glued
to $S_{j_r}$  by a cylinder. The value $s=0$  is allowed: in this
case we have  a  single surface $S_{j_r+1}$ to  which  we apply a
construction of creating a pair of holes,
and then attach  a  cylinder  to  the
zero  $z_{j_r+1}$  of order  $b'_{j_r+1}+1$.  In  this  case  the
newborn zero has order
\begin{equation}
\label{eq:newborn:zero:II}
\sum_{k=1}^{s} (a_{j_r+k} +2) \ +\  (b''_{j_r+s+1} +1)
\end{equation}
We may have  the  symmetric picture,  when  the surface $S_j$  to
which we apply the construction of creating a pair of holes is at the beginning
of the chain, then we glue directly the surfaces $S_{j+1}, \dots,
S_{j+s}$ to which we apply figure eight construction,  and to the
last surface $S_{j+s}$ (where $j+s\in  J$  equals  some $j_r$) we
glue a cylinder. We have a symmetric formula for the order of the
zero:
$$ 
(b'_j+1) + \sum_{k=1}^{s} (a_{j+k} +2)
$$ 

The middle drawing in Figure 10 illustrates this with $s=2$. A figure eight
construction is applied to the genus $1$ surface and a construction of creating a pair of holes to the genus $2$ surface. The first surface is glued to a cylinder.

$\bullet$  {\bf  Type   III.}   Finally,  we  may  have  surfaces
$S_k,\ldots, S_{k+s+1}$  to  which  we  apply  the
construction of creating a pair of holes boarding both ends of the chain; we apply the figure
eight  construction  to  all  interior  surfaces  $S_{k+1}  \dots
S_{k+s}$;  we  glue  all  surfaces directly. The value  $s=0$  is
admissible: by convention it means that we have only two surfaces
$S_k$ and  $S_{k+1}$ in the chain  to which we  apply
the construction of creating a pair of holes   and   then   glue   directly   $\gamma''_k$    to
$\gamma'_{k+1}$. In this  case  we get  the  newborn zero of  the
following order
\begin{equation}
\label{eq:newborn:zero:III}
(b'_k+1) +\ \sum_{i=1}^{s} (a_{k+i} +2)\ +(b''_{k+s+1}+1)
\end{equation}

The last drawing in figure 10 illustrates this situation with $s=2$. A
construction of creating a pair of holes is applied to one torus and the surface of genus $2$.  The figure eight construction is applied to the middle torus. One boundary component of it is glued to
 a boundary component of the genus $2$ surface; the other  boundary component glued to a boundary component of the other torus.

We now adopt the following notation to describe how the surfaces
are put together.

{\bf Notation} The symbol $\to$ corresponds to  direct  gluing
of the surfaces;
the  symbol  $\Rightarrow$ corresponds to a cylinder joining  two
consecutive surfaces. If a number is written as a
sum  of  two  numbers  $a_i'+a_i''$ with a bar over the sum,  we
apply the figure  eight
construction to a zero  of  order $a_i=a_i'+a_i''$ breaking it up
into zeroes of order $a_i'$ and $a_i''$;
if there  are two multiplicity
numbers, $\bar  b'_i,  \bar  b''_i$,  we  apply the construction of creating a pair of holes.
By convention $b'_i$ is  the  left  one, $b''_i$ is
the right one. Say,  in $(4,\bar 3, 2, \bar 1)$ we  have $b'_1=3,
b''_1=1$, and in $(\overline{0+0})$ we have $a_2=a_2'=a_2''=0$.

\begin{example}
Consider the following collection of flat  surfaces organized in
a cyclic order
   %
   %
\begin{alignat*}{19}
\to&  (4,\!\bar {3},\!2,\!\bar{1})
&\nsp\to&
(\overline{0\!\!+\!\!0})
&\nsp\to&
(\overline{4\!\!+\!\!2},\!2)
&\nsp\to&
(\bar{9},\!8,\!\bar{7})
&\nsp\to&
(\bar{5},\bar{3})
&\Rightarrow&
(\overline{2\!\!+\!\!2},\!2)
&\nsp\to&
(5,\!\overline{1\!\!+\!\!2})
&\Rightarrow&
(\overline{1\!\!+\!\!0},\!1)
&\nsp\to&
(2,\!\overline{0\!\!+\!\!0})
&\nsp\to
\\
\nsp\to&\quad \ \ S_1
&\nsp\to&
\quad S_2
&\nsp\to& \quad
S_3 &\nsp\to&  \quad S_4            &\nsp\to&  \ \ S_5 &\Rightarrow&
\quad S_6      &\nsp\to&  \quad  S_7     &\Rightarrow&  \quad S_8
&\nsp\to& \quad S_9      &\nsp\to
\end{alignat*}
In the example above we glue  $S_5$ to $S_6$ with a cylinder, and
$S_7$ to $S_8$ with a cylinder; all the other gluings are direct.
Thus in  this example we have  $p=9$, and $J=\{5,7\}$.  The chain
$S_1\to S_2 \to S_3 \to S_4$ gives a newborn zero of type III. It
has order $(1+1)+\bigl( (0+2)+(6+2) \bigr) + (9+1)=22$; there are
$3$ separatrix loops $\gamma_1,\gamma_2,\gamma_3$ passing through
this newborn zero.  The loops $\gamma_1$ and $\gamma_2$ bound the
surface $S_2$ and bound angles $2\pi$ with each  other. The loops
$\gamma_2,\gamma_3$ bound  $S_3$  with angles $10\pi$ and $6\pi$.
The next chain $S_4\to S_5$ gives  a newborn zero of type III; it
has order $(7+1)+(5+1)=14$; there is  only  one  separatrix  loop
$\gamma_4$  passing  through   this   zero.  The  next  chain  is
$S_5\Rightarrow$. It is of type  II;  it gives a newborn zero  of
order $3+1=4$ with  a  separatrix loop $\gamma_5$ passing through
it. The next chain $\Rightarrow S_6 \to S_7  \Rightarrow$ gives a
newborn  zero  of  order  $(4+2)+(3+2)=11$ of type I,  and  three
separatrix loops $\gamma_6,\gamma_7,\gamma_8$ passing through it.
The surface $S_6$ is  bounded  by $\gamma_6$ and $\gamma_7$ which
make angles $8\pi$ and $4\pi$  at  the zero of order $11$,  while
$S_7$ is bounded by  $\gamma_7$  and $\gamma_8$ which make angles
$6\pi$  and  $  4\pi$  at  the  zero.  Finally,  the  last  chain
$\Rightarrow S_8\to S_9\to S_1$  produces  a newborn zero of type
II.  Its   order   is   $(1+2)+(0+2)+(3+1)=9$;  there  are  three
separatrix   loops   $\gamma_9,\gamma_{10},\gamma_{11}$   passing
through it.  The  loops  $\gamma_9$  and  $\gamma_{10}$ bound the
surface $S_8$ making angles $4\pi$ and  $2\pi$  with  each  other
while $\gamma_{10}$  and  $\gamma_{11}$ bound $S_9$ making angles
$2\pi$ and $2\pi$ with each other.

The   resulting   surface   inherits   the   zeroes   of   orders
$4,2,2,8,2,5,1,2$ which  are coming unchanged from the components
$S_i$ (the integers which are not barred.) It also has
the newborn zeroes  of  orders $22,14,4,11,9$. Thus the resulting
surface     $S$      belongs     to     the      stratum
$\cH(22,14,11,9,8,5,4,4,2,2,2,2,1)$.
\end{example}

\subsection{Stratum Interchange and $\gamma\to-\gamma$ Symmetry}
\label{ss:symmetries:in:zero:to:itself}
In  this  section  we  discuss  the  possible symmetries of  the
admissible construction
$$ 
(S,a_i',a_i'',b'_k,b''_k,\rho,\gamma,h,t)\to S
$$ 
In many aspects  these symmetries are analogous to the symmetries
discussed  in  Section~\ref{ss:stratum:interchange}.   We   again
consider two settings:  the first one with more restrictions, the
second one with fewer restrictions.

{\bf Problem 1.}
We  assume  that  all  zeroes  $z_1,\dots,z_l$   of  the  surface
$S$ are numbered. We fix the zeroes $z_1, \dots, z_r \in
S$, $z_i\neq z_j$, for $i\neq j$, of orders $m_1, \dots,
m_r$  correspondingly.   We   assume   that   the   flat  surface
$S$ possesses  exactly  $s$  loops of saddle connections
$\gamma_1,\dots,  \gamma_s$   homologous   to   some  fixed  loop
$\gamma$. We assume that every loop $\gamma_i$ starts and ends at
one of the zeroes $z_1, \dots, z_r$.

Cutting along the loops $\gamma_i$, $i=1,\dots,r$  we perform the
surgery as  described  above  decomposing  $S$  into the
collection $S'_i$.  The  collection is organized in a
cyclic order as described above.

We fix the types $\alpha'_i$ of  the surfaces $S'_i$,
the distinguished zeroes (marked points) $w_i,  z_i\in S'_i$, and
the numbers $(a'_i,  a''_i)$,  $(b'_k,b''_k)$ where both types of
pairs are ordered.  We  fix, whether $S'_j$ is joined
to $S'_j$ by a cylinder, or directly.

We start  with the  case when we fix also  the following data. We
assume that $z_1, \dots, z_r$ are  the  newborn  zeroes  numbered
with respect to  this  cyclic order.  We  assume that the  zeroes
$z_{r+1},   \dots,    z_{l_1}$   lie   in   $S_1$,   the   zeroes
$z_{l_1+1},\dots, z_{l_1+l_2}$  lie  in  $S_2$,  ...,  the zeroes
$z_{l_1+\dots+l_{p-1}+1}, \dots, z_{l_1+\dots+l_p}$ lie in $S_p$.
Here $l = l_1 + \dots + l_p = card(\alpha)$.

{\bf Problem 2.}
Now  consider  the problem with fewer constraints.  We  count  the
number  of  occurrences of  configurations  of  loops  of  saddle
connections  described  in  the  previous problem with  the  only
difference that now we assume all the zeroes to be ``anonymous''.
In  other   words,   we   keep   all   geometric  information  on
configuration (number $r$ of zeroes involved; number $s$ of loops
of  homologous  saddle  connections;  types  $\alpha'_i$  of  the
surfaces  $S'_i$,  and  the  cyclic  order  of  their
appearances; numbers $(a'_i, a''_i)$, $(b'_j,b''_j)$, where  both
types of  pairs are  ordered; the way (by means  of a cylinder or
directly)       $S'_i$       is       joined       to
$S'_{i+1}$).

However, if we have two zeroes $P_i$, $P_j$ of the same degree on
the surface $S$, and, neither of them is ``newborn'', in
the setting of Problem 2 we do not distinguish cases,  when $P_i$
lies in  $S'_q$ and $P_j$ gets to $S'_t$,
and the opposite case.

\begin{remark}
Problem 2  may be considered as  generalization of Problem  1 to
the   case   when   the   zeroes   of    $S$   and   of
$S'_i$ are {\it not} numbered.
\end{remark}

\begin{remark}
In both settings the $\cH(\alpha_i')$  may  be  disconnected.  We
could   also    specify    the   connected   component   of   the
$\cH(\alpha'_i)$  when  this occurs. We prefer the setting  where
this data is not specified.
\end{remark}

Consider the  natural action of the cyclic group  of order $p$ on
the elements of the assignment organized in a cyclic order. If it
has  a nontrivial  stabilizer  we denote it  by  the same  symbol
$\Gamma$ as  before. We get  a symmetry of order $|\Gamma|$ which
we  again  call  the  {\it stratum interchange}. Note  that  this
symmetry might be  different  in the  settings  of Problem 1  and
Problem 2.

\begin{example}
\label{ex:symmetries}
The assignment
$$
(\bar{1},\bar{1},2)\to(\bar{1},\bar{1})\Rightarrow
(\bar{1},\bar{1},2)\to(\bar{1},\bar{1})\Rightarrow
$$
does not possess  the stratum interchange symmetry in the setting
of Problem  1 since the unchanged zeroes of  order $2$ are named,
and  so  they  identify  the  components  of the type  $(2,1,1)$.
However, in  the setting of  Problem 2 this configuration has the
stratum interchange symmetry of order 2.

The assignment
$$
(\bar{1},\bar{1})\to(\bar{1},\bar{1})\Rightarrow
(\bar{1},\bar{1})\to(\bar{1},\bar{1})\Rightarrow
$$
possesses the stratum interchange  symmetry  of order $2$ in both
settings.
\end{example}

Recall that for saddle connections  joining  distinct  zeroes  we
also  might  also have the $\gamma\to-\gamma$ symmetry. It  could
appear only when $m_1=m_2$, for in this case we have no invariant
way of choosing the orientation of the saddle connection. When we
have only  saddle connections joining  a zero to itself, there is
no   geometric   way  of  choosing  the  orientation  of   saddle
connections. The two different choices of the orientation give two
decompositions of $S$. Thus we get  the  action  of  the
group $\integers/2\integers$ on the collection of assignments. We
call this  action  the  $\gamma\to-\gamma$  {\it  symmetry}. This
action is easily described in combinatorial  terms:  one  has  to
change the cyclic order in the assignment to the opposite one, as
well as the order in each pair $(a'_i,a''_i)$ and $(b'_j,b''_j)$.
Some assignments may stay invariant under this operation. In this
case we  say  that  they  have  $\gamma\to-\gamma$ symmetry. This
property also depends  whether we consider the assignments in the
setting of Problem 1 or of Problem 2.

\begin{example}
The first assignment in Example~\ref{ex:symmetries} does not have the
$\gamma\to-\gamma$ symmetry; the second one does.

An assignment  of  multiplicity 1 has $\gamma\to-\gamma$ symmetry
if  and  only  if  $a'=a''$  when  we  apply   the  figure  eight
construction  and  if and  only  if $b'=b''$  when  we the  apply
the creating a hole construction.

The assignment
$$
(\bar{2},\bar{2},2)\Rightarrow(\bar{2},\bar{2},2)\Rightarrow
$$
has the $\gamma\to-\gamma$  symmetry  in both settings, while the
assignment
$$
(\bar{2},\bar{2},2)\Rightarrow(\bar{2},\bar{2},2)\Rightarrow(\bar{2},\bar{2},2)\Rightarrow
$$
has  the  $\gamma\to-\gamma$ symmetry  only  in  the  setting  of
Problem 2.
\end{example}

\begin{remark}
In             describing             the       construction
$(S',J,a'_i,a''_i,b_i',b_i'',\tau,\gamma,h,t)\to S$
one can consider {\it all}  possible  assignments,  and then take
into account corresponding symmetries,  or  one can deal with the
classes, eliminating the symmetry whenever it is possible.

We  have  chosen  the  second  way.  For example, the  assignment
$(\bar{1},\bar{0},1)\to$
(see  Section~\ref{subsec:toitself:genus3})  is symmetric  to the
assignment
$
(\bar{0},\bar{1},1)\to
$
by the $\gamma\to-\gamma$ symmetry.
\end{remark}


\subsection{Combinatorial Factor, Connected Strata}
\label{ss:combinatorial:factor:II}
We now  compute  the   combinatorial factor $M$, which computes
for fixed data  $(S',J,a_i',a_i'',b'_k,b''_k,\gamma,h,t)$ the
number of $S=\sigma(S')$ that can be built.  We start with the
setting of Problem            1 in
section~\ref{ss:symmetries:in:zero:to:itself}.  We  first suppose
that the stratum $\cH(\alpha)$ is connected.

Consider  a  surface $S'_i$ such that $z_i=w_i$ is  a
zero  of  order  $a_i=a'_i+a''_i$.  We  apply  the  figure  eight
construction by breaking  $z_i=w_i$  into zeroes of orders $a'_i$
and  $a''_i$.  There  are   $a_i+1$   ways  of  doing  this,  see
section~\ref{ss:figure:8}. The  corresponding boundary components
are denoted $\gamma'_i, \gamma''_i$.  The component $\gamma'_i$ is
attached  to  $S_{i+1}$  (directly  or  by  means  of  a cylinder
depending on whether or not $i\in J$); the component $\gamma''_i$
is attached to  $S_{i-1}$  (directly or  by  means of a  cylinder
depending on whether or not $i-1\in J$).

If $z_k\neq w_k$  we  perform the construction  of creating a pair of holes.
 The
orders of the zeroes  are  $b''_k$ and $b'_k$ correspondingly, so
there are $b''_k+1$ ways of choosing the vector $\gamma$ at $z_k$
and $b'_k+1$ ways  of choosing the  vector $\gamma$ at  $w_i$  so
this gives $(b'_k+1)(b''_k+1)$ choices.

We  may  have  the  stratum   interchange   symmetry   or   (and)
$\gamma\to-\gamma$                 symmetry,                  see
section~\ref{ss:symmetries:in:zero:to:itself}. Thus the resulting
combinatorial   factor   in   the   setting  of  Problem   1   in
section~\ref{ss:symmetries:in:zero:to:itself} equals
$$ 
c= \frac{1}{|\Gamma_-|}\cdot\frac{1}{|\Gamma|}
\cdot \prod_{\substack{1\le i \le p\\z_i=w_i}} (a_i+1)  \cdot
\prod_{\substack{1\le k \le p\\z_k\neq w_k}} (b'_k+1)(b''_k+1)
$$ 

We now compute the combinatorial factor in the setting of Problem
2  in   Section~\ref{ss:symmetries:in:zero:to:itself}.  For  each
integer $m$, let $o(m)$ denote  its  multiplicity.  For each $i$,
$i=1,\dots,p$, let $o_i(m)$ denote the  multiplicity  of  $m$  in
$\alpha'_i$.

Consider an  integer $m$ different  from the order of any newborn
zero  on  $S$   (see   formulae~\ref{eq:newborn:zero:I},
\ref{eq:newborn:zero:II},     \ref{eq:newborn:zero:III}),     and
different from  any $a_i$ or  $b'_k, b''_k$. Since all the zeroes
$z_1, z_2, \dots, z_l$  of  $S$ are numbered, the number
of  ways  to arrange $o(m)$ zeroes  of  order $m$ into groups  of
$o_1(m), \dots, o_p(m)$ zeroes equals
$$ 
\frac{o(m)!}{\prod_{i=1}^p o_i(m)!}
$$ 
where we let $0!=1$ by convention.

Let  now  $s\ge  1$  be  the  number  of  the newborn  zeroes  in
$S$ of order $m$, and suppose that $m$ is different from
any  $a_i,  b'_k,  b''_k$. Now we  have  to  choose  $s$ of $o(m)$
numbered zeroes to  be  the  newborn ones, and then  we  have  to
arrange the  remaining  $o(m)-s$  ordered  zeroes  into groups of
$o_1(m), \dots, o_p(m)$ zeroes. For such $m$ we get the factor
$$ 
o(m)(o(m)-1)\cdot\dots\cdot(o(m)-s+1)\ \cdot\
\frac{(o(m)-s)!}{\prod_{i=1}^p o_i(m)!} =
\frac{o(m)!}{\prod_{i=1}^p o_i(m)!}
$$ 
which coincides  with the previous one. Hence, we  do not need to
distinguish this case from the previous one.

If $m$ is different from any  $b'_k, b''_k$, but is equal to some
$a_i$, then there  are  only $o_i(m)-1$  zeroes  of order $m$  on
$S'_i$   inherited  from   $S$.   Thus   the
corresponding factor  in  the  denominator  in  the formula above
equals $(o_i(m)-1)!$.  Multiplying  numerator  and denominator by
$o_i(m)$ we get the following factor for such $m$:
$$ 
\frac{o(m)!}{\prod_{i=1}^p o_i(m)!} \cdot
\prod_{\substack{1\leq i\le p\\ \vspace*{-5pt} \\a_i=m\\z_i=w_i}}
o_i(m)
$$ 
We have  to make a similar  correction for those $m$ which
coincide with some $b'_k$ or  $b''_k$. In the case when $m$
equals  one of the  $b'_k\neq  b''_k$ the correction is
completely analogous  to the  previous  one.  In  the  case
when  for  some $k$  we  have $m=b'_k=b''_k$,  we  get  only
$o(m)-2$ zeroes of order  $m$  on $S'_k,$
inherited  from   $S$.   Thus   the corresponding
factor  in  the  denominator  in  the formula above equals
$(o_k(m)-2)!$.  Multiplying  numerator  and denominator by
$o_k(m)(o_k(m)-1)$ and  collecting  all the correction terms
we finally get the following combinatorial factor for the
setting of Problem 2 in
section~\ref{ss:symmetries:in:zero:to:itself}.
\begin{multline*}
M=\frac{1}{|\Gamma_-|}\cdot\frac{1}{|\Gamma|} \cdot
\prod_{m\in \alpha} \left(\frac{o(m)!}
{\prod_{i=1}^p o_i(m)!}\right) \cdot \\
\cdot \prod_{\substack{1\le i \le p\\z_i=w_i\\a_i\neq 0}} o_i(a_i)
\cdot
\prod_{\substack{1\le k \le p\\z_k\neq w_k\\b'_k\neq
b''_k\\b'_k\neq 0}}
o_k(b'_k) \cdot  \prod_{\substack{1\le k \le p\\z_k\neq
w_k\\b'_k\neq b''_k\\b''_k\neq 0}}o_k(b_k'')\cdot
\prod_{\substack{1\le k \le p\\z_k\neq w_k\\b'_k=b''_k\neq 0}}
o_k(b'_k) (o_k(b'_k)-1) \cdot \\
\cdot \prod_{\substack{1\le i \le p\\z_i=w_i}} (a_i+1)  \cdot
\prod_{\substack{1\le k \le p\\z_k\neq w_k}} (b'_k+1)(b''_k+1)
\end{multline*}
Note that in  general the symmetry groups $\Gamma$ and $\Gamma_-$
corresponding to  possible  stratum  interchange  symmetry and to
$\gamma\to-\gamma$  symmetry  are different in Problems 1 and  2,
see Section~\ref{ss:symmetries:in:zero:to:itself}.

There  will  also  be  additional  factors  coming from the  spin
structures  and  the hyperelliptic strata. We shall discuss  these
separately.

We finally  get the following  expression for the constant $c$ in
the setting where the zeroes of  $S$  are  not  numbered
(see Problem 2 in Section~\ref{ss:symmetries:in:zero:to:itself}):

\begin{formula}
\label{f:toitself:connected}
The  list   of   possible  configurations  of  homologous  saddle
connections joining zeroes to  themselves  is the same for almost
all surfaces in any connected stratum $\cH(\alpha)$. The possible
configurations            are            described             in
section~\ref{ss:admissible:constructions}.  The  constant in  the
quadratic  asymptotics  for  the  number  of  saddle  connections
with associated data
$(J,a_i',a_i'',b'_k,b''_k)$   is
presented by the following formula: \end{formula}
   %
        %
\begin{multline}
\label{eq:sad:conn:to:itself:const}
c=\frac{1}{|\Gamma_-|}\cdot\frac{1}{|\Gamma|} \cdot
\prod_{m\in\alpha} \left(\frac{o(m)!}{\prod_{i=1}^p
o_i(m)!}\right) \cdot \\
\cdot \prod_{\substack{1\le i \le p\\z_i=w_i\\a_i\in\alpha}}
o_i(a_i)
\cdot
\prod_{\substack{1\le k \le p\\z_k\neq w_k\\b'_k\neq
b''_k\\b'_k\in\alpha}}
o_k(b'_k) \cdot  \prod_{\substack{1\le k \le p\\z_k\neq
w_k\\b'_k\neq b''_k\\b''_k\in\alpha}}o_k(b_k'')\cdot
\prod_{\substack{1\le k \le p\\z_k\neq w_k\\b'_k=b''_k\in\alpha}}
o_k(b'_k) (o_k(b'_k)-1) \cdot \\
\cdot \prod_{\substack{1\le i \le p\\z_i=w_i}} (a_i+1)  \cdot
\prod_{\substack{1\le k \le p\\z_k\neq w_k}} (b'_k+1)(b''_k+1)
\cdot\\
\cdot \frac{1}{2^{p-1}}\cdot
\frac{ \prod_{i=1}^p (\tfrac{d_i}{2}-1)! }
     { (\frac{d}{2}-2)! } \cdot
\frac{ \prod_{i=1}^p \Vol( \cH_1(\alpha'_i)) }
     { \Vol(\cH_1(\alpha)) }
\qquad
\end{multline}
\bigskip


\subsection{Examples: Constants for Strata in Genus 3}
\label{subsec:toitself:genus3}



\bigskip
\noindent{\large \bf{Stratum} $\boldsymbol{\cH(3,1)}$}
\smallskip


\medskip
\noindent  $\boldsymbol{ \to (\overline{1},\overline{0},1) \to
}$\\ \noindent
The first possibility is multiplicity $1$ with no cylinder. After
degeneration  we  get a genus $2$ surface  with  a  distinguished
simple zero, a  marked  point and  another  simple zero. We  have
$d_1=12$.  We   employ  the  parallelogram  construction  to  the
distinguished simple zero and marked point  producing two circles
which we glue directly. Since $1=b'_1\neq b''_1=0$ we do not have
the $\gamma\to-\gamma$ symmetry. Thus $M=(b'_1+1)(b''_1+1)=2$. We
have
$$ 
c=2\cdot\frac{(5+1-1)!}{(7-2)!} \cdot
\frac{\Vol(\cH_1(1,1))}{\Vol(\cH_1(3,1))} =
\frac{105}{16}\cdot\frac{1}{\zeta(2)}\approx 3.990
$$ 

\medskip
\noindent  $\boldsymbol{ \Rightarrow (\overline{1+0},1)
\Rightarrow }$\\ \noindent
The next possibility is multiplicity  $1$  with  a cylinder where
the other side  also  returns  to the zero of  order  $3$.  After
degeneration we  get a surface  of genus $2$ with a distinguished
simple zero on which we  perform  the  figure eight construction,
gluing in  a cylinder. Since  $1=a'_1\neq a''_1=0$ we do not have
the $\gamma\to-\gamma$ symmetry. Thus $M=(a_1+1)=2$. Now $d_1=10$
so
$$ 
c=2\cdot\frac{(5-1)!}{(7-2)!} \cdot
\frac{\Vol(\cH_1(1,1))}{\Vol(\cH_1(3,1))}=
\frac{21}{16}\cdot\frac{1}{\zeta(2)}\approx 0.7979
$$ 

\medskip
\noindent  $\boldsymbol{ \Rightarrow (\overline{2},\overline{0})
\Rightarrow }$\\ \noindent
If  the  other  side  returns  to  the simple  zero,  then  after
degeneration, we get a genus $2$ surface with a double zero and a
marked point on which we perform  the creating a hole construction.
Again $d_1=10$. Since $2=b'_1\neq b''_1=0$ we  do  not  have  the
$\gamma\to-\gamma$  symmetry.  Thus  $M=(b'_1+1)(b''_1+1)=3$.  We
have
$$ 
c=3\cdot\frac{(4+1-1)!}{(7-2)!} \cdot
\frac{\Vol(\cH_1(2))}{\Vol(\cH_1(3,1))}=
\frac{567}{256}\cdot\frac{1}{\zeta(2)}\approx 1.346
$$ 

\medskip
\noindent  $\boldsymbol{ \Rightarrow (\overline{0+0}) \to
(\overline{0},\overline{0})
\Rightarrow }$\\ \noindent
Now we consider multiplicity $2$ in this stratum. Then we  have a
cylinder  whose  boundary components are curves returning to  the
zero of order  $3$  and zero of order  one.  The homologous curve
returns to  the zero of  order $3$. The degenerating surfaces are
tori, one  with a single marked point, the  other with $2$ marked
points. On the torus with one point, we perform the  figure eight
construction and on the other the  creating a hole construction. We
glue in one cylinder and one pair of circles is glued directly.

There is no $\gamma\to-\gamma$ symmetry in  this  case.


   %
In this case we get
$$ 
c=1\cdot\frac{1}{2^{2-1}}\cdot
\frac{(3-1)!\cdot(2-1)!}{(7-2)!}\cdot
\frac{\Vol(\cH_1(\torusemptyset)))^2}{\Vol(\cH_1(3,1))}=
\frac{105}{256}\cdot\frac{1}{\zeta(2)}\approx 0.2493
$$ 


\bigskip
\noindent{\large \bf{Stratum} $\boldsymbol{\cH(2,1,1)}$}
\smallskip


\medskip
\noindent $\boldsymbol{ \to (1,1,\overline{0},\overline{0}) \to
}$\\ \noindent
The multiplicity $1$ case with no cylinder degenerates to a genus
$2$ surface with $2$ marked  points.  We  apply the parallelogram
construction to the marked points, gluing the circles directly to
each other. There is the $\gamma\to-\gamma$ symmetry.

   %
We get
$$ 
c=\frac{1}{2}\cdot\frac{(5+2-1)!}{(8-2)!}\cdot
\frac{\Vol(\cH_1(1,1))}{\Vol(\cH_1(2,1,1))}=
\frac{7}{3}\cdot\frac{1}{\zeta(2)}\approx 1.418
$$ 

\medskip
\noindent $\boldsymbol{ \Rightarrow (1,1,\overline{0+0})
\Rightarrow }$\\ \noindent
For the case of a cylinder with multiplicity $1$ where  the other
side  returns  to the double zero the  degenerating  surface  has
genus $2$ with a single marked point.
We apply the figure
 eight   construction   at   the   marked  point.  We   have
the
$\gamma\to-\gamma$ symmetry.
   %
$$ 
c=\frac{1}{2} \cdot \frac{(5+1-1)!}{(8-2)!} \cdot
\frac{\Vol(\cH_1(1,1))}{\Vol(\cH_1(2,1,1))} =
\frac{7}{18} \cdot \frac{1}{\zeta(2)} \approx 0.2364
$$ 

\medskip
\noindent $\boldsymbol{ \Rightarrow (1,\overline{1},\overline{0})
\Rightarrow }$\\ \noindent
If the other side of the cylinder returns to a simple zero, then
after degeneration  we get a surface of genus  $2$ with a marked
point and a simple  zero (as well as another simple zero  on the
surface). We apply the creating a pair of holes construction  at the marked
point and a  zero  gluing in a cylinder.  Now  we have $b'_1=0$,
$b''_1=1$.   Since   $b'_1\neq   b''_2$   we  do  not   have   a
$\gamma\to-\gamma$ symmetry.
   %
$$ 
c=4\cdot\frac{(5+1-1)!}{(8-2)!}\cdot
\frac{\Vol(\cH_1(1,1))}{\Vol(\cH_1(2,1,1))}=
\frac{28}{9}\cdot\frac{1}{\zeta(2)}\approx 1.891
$$ 

\medskip
\noindent $\boldsymbol{ \Rightarrow (2,\overline{0},\overline{0})
\Rightarrow }$\\ \noindent
If the cylinder  joins two simple zeroes, then after degeneration
we get  a  surface  of genus $2$ with a zero of order $2$ and two
marked points. We apply the  parallelogram  construction  at  the
marked   points.   We   have   $b'_1=b''_1=0$.   Here   we   have
$\gamma\to-\gamma$   symmetry.
   %
$$ 
c=1\cdot\frac{(4+2-1)!}{(8-2)!}\cdot
\frac{\Vol(\cH_1(2))}{\Vol(\cH_1(2,1,1))}=
\frac{7}{8}\cdot\frac{1}{\zeta(2)}\approx 0.5319
$$ 

\medskip
\noindent $\boldsymbol{ \Rightarrow (\overline{0},\overline{0})
\to (\overline{0},\overline{0})
\Rightarrow }$\\ \noindent
There are two multiplicity $p=2$ cases. The first is if  there is
a cylinder with boundary curves  returning  to  the simple zeroes
and a  homologous curve returning to the double  zero. We get $2$
tori each with $2$ marked points so $p=2,  d_1=d_2=6$. We perform
the  parallelogram  construction  on  each and then glue  in  one
cylinder and glue two circles  directly.  We  have the $\gamma\to
-\gamma$  symmetry,  and we  do  not  have   a  stratum
interchange.
   %
$$ 
c=1\cdot\frac{1}{2^{2-1}}\cdot\frac{(3-1)!\cdot(3-1)!}{(8-2)!}
\cdot
\frac{\Vol(\cH_1(\torusemptyset))^2}{\Vol(\cH_1(2,1,1))}=
\frac{7}{36}\cdot\frac{1}{\zeta(2)}\approx 0.1182
$$ 

\medskip
\noindent $\boldsymbol{ \Rightarrow (\overline{0},\overline{0})
\Rightarrow (\overline{0+0})
\Rightarrow }$\\ \noindent
The  second  case is if there  are  $2$ cylinders. In that  case,
after degeneration,  we get a torus with $2$  marked points and a
torus  with  one  marked   point.   We  apply  the  figure  eight
construction  to  the   torus  with  one  marked  point  and  the
parallelogram  construction  to  the  other.  We   then  glue  in
cylinders to  each pair of  circles. Again we have the $\gamma\to
-\gamma$  symmetry,  and we  do  not  have   a  stratum
interchange.
   %
$$ 
c=1\cdot\frac{1}{2^{2-1}}\cdot\frac{(3-1)!\cdot(2-1)!}{(8-2)!}\cdot
\frac{\Vol(\cH_1(\torusemptyset))^2}{\Vol(\cH_1(2,1,1))}=
\frac{7}{72}\cdot\frac{1}{\zeta(2)}\approx 0.05910
$$ 

\subsection{Principal Stratum}
\label{s:distinct:principal}

Since any figure  eight  construction produces a ``newborn'' zero
of order at  least  two,  here we cannot have  any  figure  eight
constructions. All  the surfaces $S'_i,$ belong to the
principal  strata in  lower  genera $g_i$; there  are  a pair  of
marked points $z_i\neq w_i$ on each surface, and  we always apply
the  parallelogram  construction  with  $b'_i=b''_i=0$.  All  the
curves $\gamma_i$ bound cylinders.

\subsubsection{Saddle  Connections  of  Multiplicity $1$.}

In  multiplicity   one   there   is   a   single  surface  $S'\in
\cH(\alpha')$  of  genus  $g-1$  and $S'$ has two  marked  points
$z\neq w$. We recover $S$ from $S'$ by applying the parallelogram
construction followed by gluing  in  a cylinder; $\alpha'$ is the
partition of $2(g-1)-2$ into ones.

We always have the $\gamma\to-\gamma$ symmetry, so $|\Gamma_-|=2$
in   all   cases.  We   use   the  setting   of   Problem  2   of
section~\ref{ss:symmetries:in:zero:to:itself}. The combinatorial
constant is
$$ 
M=\frac{1}{2} \cdot \frac{(2g-2)!}{(2g-4)!} \cdot (0+1)(0+1)
$$ 
Applying~\eqref{eq:sad:conn:to:itself:const} we get the following
value for the constant $c$:
$$ 
c=\frac{(2g-2)(2g-3)}{2} \cdot
\frac{((4g-3-2)-1)!}{((4g-3)-2)!} \cdot
\frac{\Vol(\cH_1(\alpha'))}{\Vol(\cH_1(\alpha))}
$$ 
Thus, finally,
$$ 
c=\frac{(g-1)(2g-3)}{4g-5} \cdot
\frac{\Vol(\cH_1(\overbrace{1,\dots,1}^{2g-4}))}
{\Vol(\cH_1(\underbrace{1,\dots,1}_{2g-2}))}
$$ 
In  the  table below we present  the  values of the constant  for
saddle connections of  multiplicity one joining a zero to itself,
for flat surfaces living in the principal strata in small genera.

\begin{table}[hbt!]
\caption{
Principal stratum $\cH(1,\dots,1)$; values of the
constants in the quadratic  asymptotics  for the number of closed
geodesics of multiplicity one.}
\label{tab:c:toitself:princ:mult1}
\scriptsize
$$
\begin{array}{|l|c|c|c|c|c|c|c|}
\hline &&&&&&&\\
& g=2 & g=3 & g=4 & g=5 & g=6 & g=7 & g=8 \\
[-\halfbls] &&&&&&& \\ \hline &&&&&&& \\ [-\halfbls]
c\cdot\zeta(2) = &   \cfrac{5}{2} & \cfrac{36}{7} &
\cfrac{3150}{377} & \cfrac{274456}{23357} &
\cfrac{250153470}{16493303} & \cfrac{6531347988}{351964697} &
\cfrac{8007196856750}{365342975469}
%
%
\\  [-\halfbls]
&&&&&&& \\ \hline &&&&&&& \\ [-\halfbls]
c \approx &  1.51982& 3.12648& 5.07950& 7.14344& 9.22041& 11.2812&
13.3239
%
%
\\ [-\halfbls]
&&&&&&&\\ \hline
\end{array}
$$
\end{table}


\subsubsection{Higher Multiplicity in the Principle Stratum}
The surface $S$ breaks up into  a union of $p$ surfaces $S_i$ of
genera $g_i$, $i=1,\dots,p$. Each  $S_i$ has only simple zeroes
and a pair of marked
points $z_i$ and $w_i$.
We reconstruct $S$  by the creating a pair of holes construction
at the marked points $z_i$ and $w_i$,
and then by
gluing in cylinders $C_i$ to connect the boundary circles.


For multiplicity two assignments
there is always  a  factor of $1/2$ due  to  the  symmetry
$\gamma\to-\gamma$. For higher multiplicities an assignment has
the $\gamma\to-\gamma$ symmetry if and only if
the corresponding cyclic ordering
$g_1\to\dots\to g_p\to$ is invariant under reversing arrows
(up to a cyclic shift of the entries).
The assignment has the stratum interchange symmetry if
and only if
the corresponding cyclic ordering
$g_1\to\dots\to g_p\to$ has translational symmetry.
The resulting combinatorial factor is
$$ 
M=\frac{1}{|\Gamma_{-}|\cdot|\Gamma|}\frac{(2g-2)!}{\prod_{i=1}^p (2g_i-2)!}
$$ 

Let $\alpha_i$ denote
the partition of $2g_i-2$ into ones.
Applying our formula with $M$ as above, we get
$$ 
 c =\frac{M}{2^{p-1}(4g-5)!
\Vol(\cH_1(\alpha))}
\prod_{i=1}^p (4
 g_i - 2)! \Vol( \cH_1(\alpha_i))
$$ 

To  complete  this  section  we   present   some   examples   of
these computations.

\begin{example}{\bf Stratum} $\boldsymbol{\cH(1,1,1,1).}$
In  the  stratum  $\cH(1,1,1,1)$  in genus $3$, the  only  higher
multiplicity is $p=2$ so $S'_1, S'_2$ have genus $1$. We have the
stratum  interchange  symmetry factor of $|\Gamma|=2$; as in  all
these  cases   we  also  have  the  $\gamma\to-\gamma$  symmetry,
$|\Gamma_-|=2$. Hence we have the combinatorial factor
$$ 
M=\frac{1}{2}\cdot\frac{1}{2}\cdot\frac{4!}{0!\cdot 0!}=6
$$ 
Thus
$$ 
c=6\cdot\frac{1}{2^{2-1}}\cdot\frac{2!\cdot 2!}{7!}\cdot
\frac{\bigl(\Vol(\cH_1(\torusemptyset))\bigr)^2}{\Vol(\cH_1(1,1,1,1))}
=
\frac{3}{14}\cdot\frac{1}{\zeta(2)} \approx 0.1303
$$ 
\end{example}

\bigskip

\begin{example}{\bf Stratum} $\boldsymbol{\cH(1,1,1,1,1,1).}$
Now $g=4$, and we have two  different  nontrivial  partitions  of
$g-1=3$ arranged in a cyclic  order:  in  multiplicity $2$,
a  surface of  genus $2$ and a surface  of genus $1$ and
in multiplicity $3$ yielding three  tori.  In  the  first
case  we  do  not  have   the   stratum   interchange   symmetry,
$|\Gamma|=1$. Hence we get the following combinatorial factor
$$ 
M=\frac{1}{2} \cdot\frac{6!}{0!\cdot 2!}
$$ 
Thus
$$ 
c=\frac{6!}{4}\cdot\frac{1}{2^{2-1}}\cdot\frac{2!\cdot
6!}{11!}\cdot
\frac{\Vol(\cH_1(\torusemptyset))\Vol(\cH_1(1,1))}{\Vol(\cH_1(1,1,1,1,1,1))}
=
\frac{90}{377}\cdot\frac{1}{\zeta(2)} \approx 0.1451
$$ 

In  the  second case we have the  stratum  interchange  symmetry,
$|\Gamma|=3$. Hence we get the following combinatorial factor
$$ 
M=\frac{1}{2}\cdot\frac{1}{3}\cdot\frac{6!}{0!\cdot 0!\cdot 0!} =
5!
$$ 
Thus
$$ 
c=5!\cdot\frac{1}{2^{3-1}}\cdot\frac{2!\cdot 2!\cdot 2!}{11!}\cdot
\frac{\bigl(\Vol(\cH_1(\torusemptyset))\bigr)^3}{\Vol(\cH_1(1,1,1,1,1,1))}
=
\frac{5}{754}\cdot\frac{1}{\zeta(2)} \approx 0.004031
$$ 
\end{example}

\begin{example}{\bf Stratum} $\boldsymbol{\cH(1,1,1,1,1,1,1,1).}$
Proceeding  to  genus   $5$   we  get  the  following  nontrivial
partitions  of  $g-1=4$  arranged  in  a  cyclic order:  in
multiplicity $2$  giving a surface of genus $1$  and one of genus
$3$; in multiplicity $2$ yielding  two  surfaces  of  genus
$2$, in multiplicity $3$ yielding  two tori and  a genus
$2$ surface and in multiplicity $4$ yielding $4$ tori.

In  the  first  case  we  do  not have  any  stratum  interchange
symmetry, $|\Gamma|=1$.  Hence we get the following combinatorial
factor
$$ 
M=\frac{1}{2} \cdot\frac{8!}{0!\cdot 4!}
$$ 
Thus
$$ 
c=\frac{8!}{2\cdot 4!}\cdot\frac{1}{2^{2-1}}\cdot\frac{2!\cdot
10!}{15!}\cdot
\frac{\Vol(\cH_1(\torusemptyset))\Vol(\cH_1(1,1,1,1))}{\Vol(\cH_1(1,1,1,1,1,1,1,1))}
=
\frac{4200}{23357}\cdot\frac{1}{\zeta(2)} \approx 0.1093
$$ 

For   the    other   multiplicity   $2$
case $\alpha'_1=\alpha'_2=(1,1)$.  We  have  the  stratum
interchange
symmetry, $|\Gamma|=2$.  Hence we get the following combinatorial
factor
$$ 
M=\frac{1}{2}\cdot\frac{1}{2}\cdot\frac{8!}{2!\cdot 2!}
$$ 
Thus
$$ 
c=\frac{8!}{16}\cdot\frac{1}{2^{2-1}}\cdot\frac{6!\cdot
6!}{15!}\cdot
\frac{\bigl(\Vol(\cH_1(1,1))\bigr)^2}{\Vol(\cH_1(1,1,1,1,1,1,1,1))}
=
\frac{720}{23357}\cdot\frac{1}{\zeta(2)} \approx 0.01874
$$ 

For  the  third case  we  do  not  have  any  stratum interchange
symmetry, $|\Gamma|=1$.  Hence we get the following combinatorial
factor
$$ 
M=\frac{1}{2} \cdot\frac{8!}{0!\cdot 0!\cdot 2!}
$$ 
Thus
$$ 
c=\frac{8!}{4}\cdot\frac{1}{2^{3-1}}\cdot\frac{2!\cdot 2!\cdot
6!}{15!}\cdot
\frac{\bigl(\Vol(\cH_1(\torusemptyset))\bigr)^2\Vol(\cH_1(1,1))}{\Vol(\cH_1(1,1,1,1,1,1,1,1))}
=
\frac{120}{23357}\cdot\frac{1}{\zeta(2)} \approx 0.003123
$$ 

In  the  last  case  we  have  stratum  a  interchange  symmetry,
$|\Gamma|=4$. Hence we get the following combinatorial factor
$$ 
M=\frac{1}{2}\cdot\frac{1}{4}\cdot\frac{8!}{0!\cdot 0!\cdot
0!\cdot 0!} = 7!
$$ 
Thus
$$ 
c=7!\cdot\frac{1}{2^{4-1}}\cdot\frac{2!\cdot 2!\cdot 2!\cdot
2!}{15!}\cdot
\frac{\bigl(\Vol(\cH_1(\torusemptyset))\bigr)^4}{\Vol(\cH_1(1,1,1,1,1,1,1,1))}
=
\frac{5}{46714}\cdot\frac{1}{\zeta(2)} \approx 0.00006507
$$ 
\end{example}

\section{Strata that are not Connected}
\label{s:not:connected:II}

\subsection{Parity of the Spin Structure of the Compound Surface}
\label{ss:parities}

Consider a flat surface $S$ constructed  in the previous
section. Throughout this section we  suppose  that  all zeroes of
$\omega$  are  of  even  orders,  so  the resulting flat  surface
$S$ has a spin structure. In this section we discuss the
conditions on components  $S'_i$  , and on the gluing
rules  between  them  which  allow   us   to   obtain  a  surface
$S$ with even $\alpha$. We  also  compute  the parity of
the spin  structure of $S$. This  will allow us  to find
constants in the case of components with spin structures.

\begin{lemma}
\label{lm:newborn:zeroes}
Let a nondegenerate  flat surface $S$ be obtained from a
collection $S'_i$  by  applying  the figure eight and
parallelogram   constructions.   The   resulting   flat   surface
$S$  has  zeroes of  even  orders  if  and  only  if the
following conditions are valid:
\begin{itemize}
\item[---]
All zeroes of those surfaces $S'_i$ to which we apply
the figure eight construction have even orders;
\item[---]
If  there  is  at  least  one  surface  to  which  we  apply  the
parallelogram   construction,   then   there  is  the   following
additional alternative. Either
\begin{itemize}
\item[---]
all zeroes of all surfaces  are  of even order and every  newborn
zero of $S$ is of type III;
\item[or]
\item[---]
the orders  $b'_k$, $b''_k$ of  all zeroes $z_k\neq w_k$ to which
we apply the creating a pair of holes construction are  odd  while  all  the
other zeroes of any $S_i$ are even.
\end{itemize}
\end{itemize}
\end{lemma}
\begin{proof}  First   note   that   all   zeroes   inherited  by
$S$ from  $S'_i$ without changes are of even
orders. Thus, if  we apply to $S'_i$ the figure eight
construction, than all the zeroes of $S'_i$ different
from $z_i$  must have even orders. Since the  total sum of orders
equals $2g_i+2$,  it means that $a_i$ --- the  order of $z_i$ ---
is even as well.

Suppose now that there  is  a surface $S'_k$ to which
we  apply  the creating a pair of holes construction with at  least  one  of
$b'_k, b''_k$ even. Without loss of generality we may assume that
$b'_k$  is  even. By the same  arguments  as above all zeroes  of
$S'_k$ different  from  $z_k$  and  $w_k$  have  even
orders. Since  the total sum $2g_k+2$ of all  orders is even, the
remaining zero has even order $b''_k$ as well.

Consider the newborn zero of $S$ induced from $w_k'$. It
is of one of the types II or III.  If it  were of  type II,  then
by~\eqref{eq:newborn:zero:II} it  would have odd order. Thus this
newborn     zero      is     of     the     type     III.     Let
$S'_{k+s+1}$ be  the surface, to which we apply
creating a pair of holes construction, at the opposite end of the chain from
$S'_k$.  By~\eqref{eq:newborn:zero:III}  $b'_{k+s+1}$
is even.  Repeating our arguments we show that  all zeroes of all
$S'_i$ are even, and  that  all newborn zeroes are of
the type III.

The remaining part of the alternative now becomes obvious.
\end{proof}
Let  us  calculate  the  parity  of  the spin  structure  of  the
resulting surface  $S$  in  all possible cases described
above.
\begin{lemma}
\label{lm:spin:in:8:constr}
Suppose that we  use only figure  eight constructions to  glue  a
nondegenerate surface $S$. Suppose  that  all $\alpha_i$
are even.
The  parity  of  the  spin  structure  of  the  resulting surface
$S$ is defined by the parities of the spin structures of
the components, and by the parities of the numbers $a'_i$  in the
following way:
$$ 
\phi(S)=1+ \sum_{i=1}^{p} \phi(S'_i) +
\sum (a'_i+1) \quad   (mod\ 2)
$$ 
\end{lemma}
\begin{proof} For each surface $S_i$  consider  a  collection  of
$2g_i$  smooth  simple  closed  curves on $S'_i$  representing  a
canonical basis of cycles. As usual we may assume that the curves
stay outside  a neighborhood of  the zeroes and the marked points
of $S'_i$.

We may  construct now the  following canonical basis of cycles on
the resulting surface $S$. Let $g$ be the  genus of $S$.
First take all  the  cycles represented  by  the basic curves  on
components. They give us  $g-1$  pairs of cycles. By construction
they form incomplete canonical basis of cycles on $S$.

Since we are using only figure eight constructions we have to use
at  least  one   cylinder   to  obtain  a  nondegenerate  surface
$S$. As a  cycle  $A_g$ we  may  use the cycle  $\gamma$
represented by the  waist curve of this cylinder. By construction
it  is  independent from  the  previous  ones,  and  it  does not
intersect them. To complete the construction of the basis we have
to choose a  cycle  $B_g$ dual  to  $A_g$. Consider the  following
curve representing $B_g$: \newline --- On any cylinder it follows
the direction transversal to $\gamma$.\newline --- On the surface
$S'_i$ the  curve  starts  at  the boundary component
$\gamma''_i$, then follows in the clockwise direction an arc in a
small neighborhood of  $z_i$, and finally arrives at the boundary
component $\gamma'_i$ in direction transversal to $\gamma$.

The canonical basis  of cycles of $S$ is constructed. By
construction  of  this   basis  the  index  of  any  basic  curve
$A_i,B_i$,  $i=1,\dots,g-1$  in  the flat structure  $S$
does not depend  on  the way in which  we  applied the prescribed
combination of figure eight constructions. Thus we have
$$ 
\sum_{k=1}^{g-1} \bigl(ind(A_k)+1\bigr)\bigl(ind(B_k)+1\bigr) =
\sum_{i=1}^p
\sum_{j_i} \bigl(ind(A_{j_i})+1\bigr)\bigl(ind(B_{j_i})+1\bigr) =
\sum_{i=1}^p \phi(S_i)
$$ 
Recall that the cycle $A_g$ is represented by the waist  curve of
a cylinder. The tangent vector to the waist curve of the cylinder
is  constant  in the flat structure $\omega$. Thus  $ind(A_g)=0$,
which implies that
\begin{multline*}
\phi(S)=
\sum_{k=1}^{g} \bigl(ind(A_k)+1\bigr)\bigl(ind(B_k)+1\bigr) = \\
= 1\cdot (ind(B_g) + 1) +
\sum_{k=1}^{g-1} \bigl(ind(A_k)+1\bigr)\bigl(ind(B_k)+1\bigr)=\\
= ind(B_g) + 1 + \sum_{i=1}^p \phi(S_i) \ (mod\ 2)
\end{multline*}

It remains to compute  the  index $ind(B_g)$. By the construction
of  $B_g$,  we  may  assume  that  the  tangent   vector  to  the
corresponding curve does  not turn in the flat structure $\omega$
while  it  crosses the cylinders. It  makes  a turn by the  angle
$2\pi (a'_i+1)$ while it follows a small arc joining two holes of
a ``figure eight'' at the point $z_i=w_i$. This  shows that $ind(
B_g)= \sum (a'_i+1)\ (mod\ 2)$.
\end{proof}

\begin{lemma}
\label{lm:spin:in:even:par:constr}
Suppose $\alpha$ is even so $S$ has a spin structure. If there is
some $S'_k$  to  which  we  apply  the  parallelogram
construction and at least one of  the points $z_k, w_k$ is a zero
(marked point) of even order, then all $\alpha'_i$  are even, and
the parity of the spin structure of the flat surface $S$
is equal to the sum of the parities of the spin structures of the
components
$$ 
\phi(S)=\sum_{i=1}^p \phi(S_i)
$$ 
\end{lemma}
\begin{proof} We use the canonical basis  on $S$ similar
to the one used in the  Lemma above. To construct the first $g-1$
pairs of  curves we again use the curves  living on the component
surfaces $S'_i$. In addition to previous  assumptions
we require that  when we apply creating a pair of holes construction to some
$S_i$,  the   basic   curves   do   not  approach  the
transversal $\rho_i$ chosen to join $z_i$ to $w_i$. We extend the
description  of  the curve representing the cycle  $B_g$  to  the
surfaces to  which  we  apply  the  creating a pair of holes construction as
follows: on such surfaces $S'_i$ the  curve starts at
the boundary component $\gamma''_i$ in a direction transversal to
$\gamma$ and arrives at the  boundary  component  $\gamma'_i$  in
direction transversal to $\gamma$. We  choose this  curve in such
way   that   it   does   not  intersect  any  basic   curves   on
$S'_i$.

It follows from Lemma~\ref{lm:newborn:zeroes} that all the zeroes
are of type III; we do not have any glued  in  cylinders: all the
gluings are  direct. As a curve representing the  cycle $ A_g$ we
now use  a  small smooth deformation $\tilde{\gamma}\subset S'_k$
of the  boundary  curve  $\gamma'_i$,  such that $\tilde{\gamma}$
does  not  pass   through  any  zeroes  of  $\omega'_i$.  We  get
$ind(A_g)=ind(\tilde{\gamma})= (b'_i+1)$. Since the number $b'_i$
is even, we see that $ind(A_g)+1  = 0\ (mod\ 2)$. Thus the impact
of the last pair of cycles to the sum
$$ 
\phi(S)=
\sum_{i=1}^{g} \bigl(ind(A_i)+1\bigr)\bigl(ind(B_i)+1\bigr)
$$ 
equals to zero, and
\begin{multline*}
\phi(S)=
\sum_{i=1}^{g-1} \bigl(ind(A_i)+1\bigr)\bigl(ind(B_i)+1\bigr) =
\\
= \sum_{i=1}^p
\sum_{j_i} \bigl(ind(A_{j_i})+1\bigr)\bigl(ind(B_{j_i})+1\bigr) =
\sum_{i=1}^p \phi(S_i)
\end{multline*}
Lemma~\ref{lm:spin:in:even:par:constr} is proved.
\end{proof}

Note that there is a discrete freedom left in the figure eight or
in the creating a pair of holes construction. When we perform a figure eight
construction at a zero of  order  $a=a'+a''$  with fixed $a',a''$
there  are  actually   $a+1$  ways  to  perform  a  figure  eight
construction, (see Section~\ref{ss:figure:8}). Similarly, when we
perform  a  parallelogram  construction  at a pair of  zeroes  of
orders  $b_i',  b_i''$  we, actually obtain  $(b'_i+1)\cdot(b''_i+1)$
different surfaces.  Thus, even when we  fix all the  elements of
the collection
$$ 
(S',J,a'_j,a''_j,b'_k,b''_k)
$$ 
we  usually  get numerous surfaces $S$ defined by  such
assignment though some freedom may be cancelled by the symmetry.
The collection $\alpha$ of degrees of  zeros  of  the  resulting
flat structure $S$ is, of course, invariant, but a spin
structure may vary.

\begin{lemma}
\label{lm:spin:in:odd:par:constr}
Let    $N$     be     the     total     number    of    surfaces
$S\in\cH(\alpha)$ obtained by the assignment
$$ 
(S',J,a'_j,a''_j,b'_k,b''_k,\gamma,\rho,h,t)\to S
$$ 
with    fixed    $(S',J,a'_j,a''_j,b'_k,b''_k,\gamma,\rho,h,t)$.
Suppose $\alpha$  is even so  every such $S$ has a spin
structure. If there is some $S'_k$ to which we apply
the creating a pair of holes construction, and one of the points $z_k, w_k$
is  an  odd order  zero,  then exactly  $N/2$  of the  resulting
surfaces  have  even  spin  structure  and  $N/2$ have odd  spin
structure.
\end{lemma}

\begin{proof} Note that  if  one of  $b'_k,  b''_k$ is odd,  the
other one is odd as well.

Consider a canonical basis  of  cycles similar to one constructed
in Lemma~\ref{lm:spin:in:even:par:constr}. The only difference is
that  now  we   may  have  some   cylinders  glued  in.   As   in
Lemma~\ref{lm:spin:in:8:constr}   we  assume   that   the   curve
representing the cycle  $B_g$  stays transversal to the direction
of $\gamma$ while passing through the cylinders.

By the choice of this canonical  basis of cycles the index of any
basic curve  representing  $A_i,B_i$, $i=1,\dots,g-1$ in the flat
structure $S$ does  not depend on  the way in  which  we
applied  the   prescribed   combination   of   figure  eight  and
creating a pair of holes constructions. Thus the number
$$ 
\phi_0 := \sum_{i=1}^{g-1} (ind(A_i)+1)(ind(B_i)+1)
$$ 
is invariant of the discrete freedom in the assignment
$$ 
(S',J,a'_j,a''_j,b'_k,b''_k,\gamma,\rho,h,t)\to S
$$ 

We  construct  a  representative   of   the  cycle  $A_g$  as  in
Lemma~\ref{lm:spin:in:even:par:constr}  using  a   small   smooth
deformation   $\tilde{\gamma}\subset    S'_k$   of   the    curve
$\gamma'_k$.  We  get  $ind(A_g)=ind(\tilde{\gamma})=  (b'_k+1)$.
Since    the    number    $b'_k$    is   odd,   we    see    that
$ind(\tilde{\gamma})+1 = 1\ (mod\ 2)$. Thus  the  parity  of  the
spin  structure   of  $S$  computed  in  the  constructed
canonical basis is represented as follows:
\begin{equation}
\label{eq:tmp}
\phi(S) = \phi_0 + (ind(A_g)+1)(ind(B_g)+1) =
\phi_0 + ind(B_g) +1 \ (mod\ 2)
\end{equation}

Let $z_i\neq w_i$ be the odd  zero  of $\omega'_i$ to which we  apply
creating a pair of holes construction. Since  the  conical angle at $z_i$ is
$2(b''_i+1)\pi$, there are an even number $r=b''_i+1$ of distinct
points  $P_1,\dots,P_r$  of the form $z_i+\gamma$. We number  the
points in the natural cyclic order. Fix all the other elements of
the construction  leaving the only freedom  in the choice  of the
point $P_j$ to  perform the slit  $z_i P_j$. We  get  $r=b''_i+1$
surfaces  $S$.  We claim that $r/2$ of  them  have  even
parity of the spin structure, while another $r/2$ have odd parity
of the spin structure.

To see  this compare the surfaces  obtained from the  slits along
$z_i P_j$ and  along $z_i P_{j+1}$. The two surfaces $S$
share the same collection of curves representing the cycles $A_i,
B_i$,  $i=1,\dots,  g-1$, and $A_g$. The curves representing  the
cycle $B_g$ differ only near the point $z_i$:  the tangent vector
to one  of the curves makes an  extra turn  by the angle  $2\pi$.
Thus by  equation~\eqref{eq:tmp}  these  two  flat  surfaces have
opposite parities of the spin structures.
\end{proof}

\subsection{Decomposition of Surfaces from Hyperelliptic Components}
\label{ss:zero:to:itself:hyp}

In  this  section  we  consider possible decompositions  of  flat
surfaces $S$ from hyperelliptic components corresponding
to  possible  configurations  of  homologous  saddle  connections
joining  a  zero to itself. In  the  next section we compute  the
constants  in  the  quadratic  asymptotics  for   the  number  of
configurations of each type.

We start with the stratum $\cH_1^{hyp}(2g-2)$.

\begin{lemma}
\label{lm:hyp:zero:to:itself:2g:minus2}
Surfaces  $S$  in  the
hyperelliptic component $\cH_1^{2g-2}(\alpha),g\geq 3$ are obtained
by
an assignment
$$ 
(S',J,a'_j,a''_j,b'_k,b''_k,\gamma,\rho,h,t)\to S
$$ 
of  one of the following three types:
\begin{itemize}
\item[i)]  $(\overline{g-2},\ \overline{g-2})\to$

The multiplicity is $1$; the flat  surface $S'$ belongs
to the hyperelliptic component $\cH_1^{hyp}(g-2,g-2)$; we apply
the
creating a pair of holes construction  gluing  the  boundaries directly; the
points $z+\gamma$, $w+\gamma$ in  the  creating a pair of holes construction
are  chosen  to be symmetric with respect  to  the  hyperelliptic
involution on  $S'$.  For  fixed  data  there precisely
$g-1$ surfaces $S$ that can be built.
\item[ii)]  $(\overline{(g-2)+(g-2)})\Rightarrow$

The multiplicity is $1$; the flat  surface $S'$ belongs
to the  hyperelliptic  component  $\cH_1^{hyp}(2g-4)$; we apply
the
figure eight  construction  gluing  in  a cylinder; $a'=a''=g-2$.
There are $2g-3$ surfaces $S$ that can be built.
\item[iii)] $(\overline{(g_1-1)+(g_1-1)})\to(\overline{g_2-1},\ \overline{g_2-1})\to$

The  multiplicity  is  $2$;  the  flat  surfaces  belong  to  the
hyperelliptic       components       $S'_1        \in
\cH^{hyp}(2g_1-2)$,             $S'_2             \in
\cH^{hyp}(g_2-1,g_2-1)$, where $g_1+g_2=g-1$. We apply the figure
eight construction  with  $a'_1=a''_1=g_1-1$ to the first surface
and the creating a pair of holes construction with $b'_2=b''_2=g_2-1$ to the
second one gluing both  pairs  of components directly. The points
$z_2+\gamma$,  $w_2+\gamma  \in S'_2$ are chosen to be  symmetric
with    respect    to    the    hyperelliptic    involution    on
$S'_2$. There  are $g_2(2g_1-1)$ surfaces that can be
built.
\end{itemize}
The surfaces can be related by a stratum  interchange or by
a $\gamma\to-\gamma$ symmetry.
\end{lemma}

The proof is analogous to the proof of Lemma~\ref{lm:slit:hyp}.

Decompositions  of  the  flat  surfaces  from  the  hyperelliptic
connected  component  $\cH_1^{hyp}(g-1,g-1)$  are described by
the
following Lemma.

\begin{lemma}
\label{lm:hyp:zero:to:itself:gm1:gm1}
Surfaces   $S$  in  the
hyperelliptic component  $\cH_1^{hyp}(g-1,g-1)$  are obtained by
an
assignment
$$ 
(S',J,a'_j,a''_j,b'_k,b''_k,\gamma,\rho,h,t)\to S
$$ 
of one of the following two types.
\begin{itemize}
\item[i)]  $(\overline{g-2}, \overline{g-2})\Rightarrow$

The multiplicity is $1$;  $S'\in\cH_1^{hyp}(g-2,g-2)$. We apply
the
parallelogram    construction     gluing    in    a     cylinder;
$b_1'=b_1''=g-2$.  The  points  $z+\gamma$,  $w+\gamma$  in   the
parallelogram  construction  are  chosen  to  be  symmetric  with
respect to  the  hyperelliptic  involution on $S'$. For
fixed data there are precisely $g-1$ surfaces $S$ that can
be built.
\item[ii)] $(\overline{(g_1-1)},\ \overline{(g_1-1)})
\to(\overline{g_2-1},\overline{g_2-1})\to$

The multiplicity  is  $2$;   $S'_i\in\cH^{hyp}(g_i-1,g_i-1)$.
where $g_1+g_2=g-1$.
We apply two creating a pair of holes constructions gluing the
components directly. The points $z_i+\gamma$, $w_i+\gamma$ in the
creating a pair of holes constructions  are  chosen  to  be  symmetric  with
respect    to    the    hyperelliptic    involution    on    each
$S'_i$, $i=1,2$. For the fixed  data  the  number  of
surfaces  $S$  that  can  be built is $g_1\cdot  g_2$  for
$g_1\neq g_2$ and $g_1\cdot g_2/2$ for $g_1= g_2$.
\end{itemize}
Surfaces can be related by a stratum interchange or
by a $\gamma\to-\gamma$ symmetry.
\end{lemma}

The proof is analogous to the proof of Lemma~\ref{lm:slit:hyp}.

\subsection{Constants for the Hyperelliptic Connected Components}

The  computation   of   the   constants   is   based  essentially
on~\eqref{eq:sad:conn:to:itself:const}.   However, since the flat
surfaces in the hyperelliptic components have extra symmetry
we  have  to  make  necessary adjustments.

\subsubsection{Hyperelliptic Component $\cH_1^{hyp}(2g-2)$}

The  admissible  assignments for this component are described  in
Lemma~\ref{lm:hyp:zero:to:itself:2g:minus2}. We always  have  the
$\gamma\to-\gamma$ symmetry  for  all three assignments. We never
have any stratum interchange symmetry.

The    only    modification   which   we   need   to   make    in
equation~\eqref{eq:sad:conn:to:itself:const}   is   as   follows.
Applying the parallelogram  constructions in this case we need to
choose the points $z_i+\gamma$, $w_i+\gamma$ to  be symmetric with
respect to hyperelliptic involution on each $S_i$ (see
Lemma~\ref{lm:hyp:zero:to:itself:2g:minus2}). Note  that for this
component  we   always   have   $b'_k=b''_k$  in  any  assignment
containing a creating a pair of holes construction. Thus we  have to replace
each  factor  $(b'_k+1)(b''_k+1)$  by  the  corresponding  factor
$(b'_k+1)$. We get the following

\begin{formula}
\label{f:toitself:hyp:2gm2}
For  almost  all  flat  surfaces  in  a  hyperelliptic
connected component $\cH_1^{hyp}(2g-2)$,  $g\ge  3$,  there are
only  three possible types  of  configurations  of saddle
connections joining zero to itself. The constants in  the
quadratic asymptotics for the  number of saddle connections
of each type are presented by the following formulae:

i) Assignment $(\overline{g-2},\ \overline{g-2})\to$
$$ 
c=\frac{g-1}{2}\cdot
\frac{\Vol(\cH_1^{hyp}(g-2,g-2))}{\Vol(\cH_1^{hyp}(2g-2))}
$$ 

ii) Assignment $(\overline{(g-2)+(g-2)})\Rightarrow$
$$ 
c=\frac{1}{2}\cdot\frac{2g-3}{2g-2} \cdot
\frac{\Vol(\cH_1^{hyp}(2g-4))}{\Vol(\cH_1^{hyp}(2g-2))}
$$ 

iii) Assignment
$(\overline{(g_1-1)+(g_1-1)})\to(\overline{g_2-1},\
\overline{g_2-1})\to$
$$ 
c=\frac{(2g_1-1)g_2}{4}\cdot
\frac{(2g_1-1)!\,(2g_2)!}{(2g-2)!} \cdot
\frac{\Vol(\cH_1^{hyp}(2g_1-2))\Vol(\cH_1^{hyp}(g_2-1,g_2-1))}
{\Vol(\cH_1^{hyp}(2g-2))}
$$ 
where $g_1+g_2=g-1$.
\end{formula}

\begin{example}{\bf Stratum} $\boldsymbol{\cH(2)}.$
This  stratum  is connected; it coincides with its  hyperelliptic
component. Flat surfaces from this stratum have a  single zero of
order  $2$.  It  is  easy  to  see  that  in  genus   $g=2$  only
multiplicity one is realizable.

\medskip
\noindent $\boldsymbol{ \to(\overline{0},\overline{0})\to}$\\
\noindent
In this  case is a  single saddle connection returning with angle
$3\pi$ so that there is no cylinder. After collapsing this saddle
connection  we  get a torus with  two  marked points to which  we
apply  the   parallelogram   construction,   gluing  the  circles
directly. Here we have
$$ 
c= \frac{2-1}{2}\cdot
\frac{\Vol(\cH_1(\torusemptyset))}{\Vol(\cH_1(2))} =
\frac{10}{3}\cdot\frac{1}{\zeta(2)} \approx 2.026
$$ 

\medskip
\noindent $\boldsymbol{ \Rightarrow (\overline{0}) \Rightarrow
}$\\ \noindent
The other possibility is that  the  saddle  connection returns at
angle  $\pi$.  Then   there  is  a  whole  cylinder  filled  with
homologous closed  geodesics.  In  particular,  there  is another
saddle connection returning  at the angle $\pi$ homologous to the
initial one. After collapsing the cylinder we get a torus  with a
single marked point. We apply the figure eight construction.
   %
$$ 
c= \frac{1}{4}\cdot
%
%
\frac{\Vol(\cH_1(\torusemptyset))}{\Vol(\cH_1(2))} =
\frac{5}{3}\cdot\frac{1}{\zeta(2)} \approx 1.013
$$ 
\end{example}

\begin{example} {\bf Component} $\boldsymbol{\cH^{hyp}(4)}.$
In genus $g=3$ all  three  possible  assignments  (see
Lemma~\ref{lm:hyp:zero:to:itself:2g:minus2})     are      already
admissible.

\medskip
\noindent $\boldsymbol{ \to (\overline{1},\overline{1}) \to }$\\
\noindent
The first represents the  multiplicity  $1$  case of a
saddle connection returning at angle $5\pi$.
   %
$$ 
c= \cfrac{3-1}{2}\cdot
%
%
\frac{\Vol(\cH(1,1))}{\Vol(\cH^{hyp}(4))} =
\frac{224}{27} \cdot \frac{1}{\zeta(2)} \approx 5.044
$$ 

\medskip
\noindent $\boldsymbol{ \Rightarrow (\overline{1+1}) \Rightarrow
}$\\ \noindent
If there is   a  cylinder, by symmetry, the  spacing  of the
angles
between the boundaries of the cylinders  is  $4\pi$.  Then  after
degeneration, we get a  single zero of order $2$ on a  surface of
genus   $2$.   We  apply  the  figure  eight  construction   with
$a'_1=a''_1=1$, and glue in a cylinder.
   %
We get
$$ 
c=\cfrac{1}{2}\cdot\cfrac{3}{4}\cdot
%
%
\frac{\Vol(\cH(2))}{\Vol(\cH^{hyp}(4))} =
\frac{7}{2} \cdot \frac{1}{\zeta(2)} \approx 2.128
$$ 

\medskip
\noindent $\boldsymbol{ \to (\overline{0+0}) \to
(\overline{0},\overline{0}) \to }$\\ \noindent
In the multiplicity $p=2$ case we  have  two  homologous  curves,
both returning to $z_0$  at  angles $3\pi$. After degeneration we
get two tori, one with two  marked points, the other with one. We
perform  the  figure  eight  construction  on  one torus and  the
creating a pair of holes construction  on  the  other, gluing the boundaries
directly to each other.
   %
   %
$$ 
c= \frac{(2\cdot 1-1)\cdot 1}{4}\cdot \frac{(2\cdot 1-1)!\,(2\cdot
1)!}{(6-1)!}\cdot
%
%
\frac{\Vol(\cH(\torusemptyset))^2}{\Vol(\cH^{hyp}(4))} =
\frac{70}{27}\cdot\frac{1}{\zeta(2)} \approx 1.576
$$ 
\end{example}

\subsubsection{Hyperelliptic Component $\cH^{hyp}(g-1,g-1)$}

The  admissible  assignments for this component are described  by
Lemma~\ref{lm:hyp:zero:to:itself:gm1:gm1}.  We  always  have  the
$\gamma\to-\gamma$ symmetry  for  both  assignments.  We have the
stratum interchange symmetry in the second assignment if and only
if $g_1=g_2$.

The   only   modification   which   we  need  to  make   in   the
equation~\eqref{eq:sad:conn:to:itself:const}   is   as   follows.
Applying the parallelogram  constructions in this case we need to
choose the points $z_i+\gamma$, $w_i+\gamma$ to  be symmetric with
respect to hyperelliptic involution on each $(S_i,\omega_i)$ (see
Lemma~\ref{lm:hyp:zero:to:itself:gm1:gm1}).  This  means that  we
have  to   replace   each   factor   $(b'_k+1)(b''_k+1)$  by  the
corresponding factor $(b'_k+1)$. Thus we get the following

\begin{formula}
\label{f:toitself:hyp:gm1:gm1}
For  almost  all  flat  surfaces  in  the  hyperelliptic
connected component $\cH^{hyp}(g-1,g-1)$, $g\ge  3$,  there
are only two possible types of configurations  of  saddle
connections joining any of two zeroes to itself.  The
constants in the quadratic asymptotics for  the number of
saddle connections of each  type are presented  by the
following formulae:

i) Assignment $(\overline{g-2}, \overline{g-2})\Rightarrow$
$$ 
c=\frac{(g-1)}{(2g-1)}\cdot
\frac{\Vol(\cH_1^{hyp}(g-2,g-2))}{\Vol(\cH_1^{hyp}(g-1,g-1))}
$$ 
\medskip

ii) Assignment $(\overline{(g_1-1)},\ \overline{(g_1-1)})
                \to(\overline{g_2-1},\overline{g_2-1})\to$
$$ 
c=\frac{g_1g_2}{2|\Gamma|}\cdot \frac{(2g_1)!\,(2g_2)!}{(2g-1)!}
\cdot
\frac{\Vol(\cH_1^{hyp}(g_1-1,g_1-1))\Vol(\cH_1^{hyp}(g_2-1,g_2-1))}
{\Vol(\cH_1^{hyp}(g-1,g-1))}.
$$ 
where  $g_1+g_2=g-1$, and
$$
|\Gamma|=\begin{cases}1&\text{when } g_1\neq g_2\\
2&\text{when } g_1= g_2\end{cases}
$$
\end{formula}

\begin{example} {\bf Stratum} $\boldsymbol{\cH(1,1)}.$
In genus $g=2$  the stratum $\cH(1,1)$ is the principal stratum.
It is connected, and it coincides  with  the hyperelliptic
component $\cH(1,1)=\cH^{hyp}(1,1)$.

Since $g=g_1+g_2+1$, and $g_i\ge 1$, we see that multiplicity two
does not occur in genus $g=2$. The value
$$
c=\frac{1}{3}\cdot\frac{\Vol(\cH_1^{hyp}(\torusemptyset))}{\Vol(\cH_1^{hyp}(1,1))}=
\cfrac{5}{2}
$$
of   the   constant   in   multiplicity  one  was   computed   in
section~\ref{s:distinct:principal}.
\end{example}

\begin{example} {\bf Component} $\boldsymbol{\cH^{hyp}(2,2)}.$

\medskip
\noindent $\boldsymbol{ \Rightarrow (\overline{1},\overline{1})
\Rightarrow }$\\ \noindent
   %
   %
In the multiplicity $1$ case there is a cylinder returning to the
other zero.
   %
   %
Here we get
$$ 
c=\cfrac{2}{5}\cdot
\frac{\Vol(\cH_1(1,1))}{\Vol(\cH_1^{hyp}(2,2))}=
\frac{14}{3}\cdot\frac{1}{\zeta(2)}\approx 2.837
$$ 

\medskip
\noindent $\boldsymbol{ \to (\overline{0},\overline{0}) \to (\overline{0},\overline{0}) \to }$\\
\noindent
In the multiplicity $2$  case  we have homologous curves, neither
of  which bounds  a  cylinder, one returning  to  each zero.  The
degenerating surfaces are tori each with  $2$  marked  points  on
each. We  perform  the  parallelogram  construction,  gluing  the
circles pairwise directly. Since $g_1=g_2$  we  get  not only the
$\gamma\to-\gamma$ but also the stratum interchange, $|\Gamma|=2$.
   %
$$ 
c=
\frac{1\cdot 1}{2\cdot 2}\cdot \frac{(2\cdot 1)!\, (2\cdot
1)!}{(2\cdot 3-1)!}\cdot
\frac{\Vol(\cH_1(\torusemptyset))^2}{\Vol(\cH_1^{hyp}(2,2))}=
\frac{35}{24}\cdot\frac{1}{\zeta(2)}\approx 0.8866
$$ 
\end{example}

\subsection{Connected Component $\cH^{nonhyp}(g-1,g-1)$; Even Genus $g$}

When the genus $g\ge 4$ is even, the stratum $\cH(g-1,g-1)$ has
two connected components: the hyperelliptic one,
$\cH^{hyp}(g-1,g-1)$, which we considered in the previous
section, and the nonhyperelliptic connected component
$\cH^{nonhyp}(g-1,g-1)$. The description of admissible
constructions for the connected component
$\cH^{nonhyp}(g-1,g-1)$ for even $g\ge 4$ is completely
analogous to the case of connected strata. However, in the
following two assignments
$$
(\overline{g-2}, \overline{g-2})\Rightarrow
$$
and
$$ (\overline{(g_1-1)},\ \overline{(g_1-1)})
\to(\overline{g_2-1},\ \overline{g_2-1})\to $$
we have interference with the hyperelliptic components. Thus
in these two cases we have to modify the general
formula~\ref{f:toitself:connected} by subtracting the
correction term corresponding to the constructions leading to
hyperelliptic flat surfaces, see
formula~\ref{f:toitself:hyp:gm1:gm1}.

\begin{formula}
\label{f:nonhyp:to:itself}
For  almost  all  flat  surfaces  in  a  nonhyperelliptic
connected component $\cH^{nonhyp}(g-1,g-1)$ for even $g \ge
4$, the constants in the quadratic asymptotics for  the
number of saddle connections of each of the following two
types are presented by the following formulae:

i) Assignment $(\overline{g-2}, \overline{g-2})\Rightarrow$
$$ 
c=\frac{(g-1)^2}{(2g-1)}\cdot \frac{\Vol(\cH_1(g-2,g-2))}
{\Vol(\cH_1^{nonhyp}(g-1,g-1))} -\frac{(g-1)}{(2g-1)}\cdot
\frac{\Vol(\cH_1^{hyp}(g-2,g-2))}{\Vol(\cH_1^{nonhyp}(g-1,g-1))}
$$ 
\medskip

ii) Assignment $(\overline{(g_1-1)},\ \overline{(g_1-1)})
                \to(\overline{g_2-1},\ \overline{g_2-1})\to$
\begin{multline*}
c=\frac{(g_1)^2 (g_2)^2}{2|\Gamma|}\cdot
\frac{(2g_1)!\,(2g_2)!}{(2g-1)!} \cdot
\frac{\Vol(\cH_1(g_1-1,g_1-1))\Vol(\cH_1(g_2-1,g_2-1))}
{\Vol(\cH_1^{nonhyp}(g-1,g-1))}\ - \\ - \frac{g_1
g_2}{2|\Gamma|}\cdot \frac{(2g_1)!\,(2g_2)!}{(2g-1)!} \cdot
\frac{\Vol(\cH_1^{hyp}(g_1-1,g_1-1))\Vol(\cH_1^{hyp}(g_2-1,g_2-1))}
{\Vol(\cH_1^{nonhyp}(g-1,g-1))}
\end{multline*}
where  $g_1+g_2=g-1$, and
$$ |\Gamma|=\begin{cases}1&\text{when } g_1\neq g_2\\
2&\text{when } g_1= g_2\end{cases} $$

iii) The constants for all other assignments are given by
equation~\eqref{eq:sad:conn:to:itself:const}, with
$\alpha=(g-1,g-1)$ where $\Vol(\cH(g-1,g-1))$ in the denominator
of
the rightmost fraction should be replaced by
$\Vol(\cH_1^{nonhyp}(g-1,g-1))$.
\end{formula}

As an example of this formula
we present the list of possible generic degenerations of a flat
surface
$S\in\cH^{nonhyp}(3,3)$ in Appendix~\ref{a:to:itself}.

\begin{remark}
\label{rm:nonhyp:to:itself}
Note, that one can represent the constant corresponding
to the two exceptional assignments
in the formula above as
$$
c=c_{regular}-c_{hyp}\cdot\cfrac{\Vol(\cH^{hyp}(g-1,g-1))}{\Vol(\cH^{nonhyp}(g-1,g-1))}.
$$
Here $c_{hyp}$ in the correctional term is the constant for the
corresponding assignment for the {\it hyperelliptic} component
$\cH^{hyp}(g-1,g-1)$, see Formula~\ref{f:toitself:hyp:gm1:gm1}.
\end{remark}

\subsection{Nonhyperelliptic Components with a Spin Structure}

Let now $\alpha$ be even  and  let $S$ belong to one
of the    nonhyperelliptic    components $\cH_1^{even}(\alpha)$
or $\cH^{odd}(\alpha)$. The admissible assignments for such
flat surfaces are described by the general
Lemmas~\ref{lm:spin:in:8:constr},
\ref{lm:spin:in:even:par:constr},
\ref{lm:spin:in:odd:par:constr}.

When $\alpha$ is one of the
two partitions $(2g-2), (g-1,g-1)$ there are special
assignments where there might be interference with
hyperelliptic components, see
Lemmas~\ref{lm:hyp:zero:to:itself:2g:minus2},
and~\ref{lm:hyp:zero:to:itself:gm1:gm1}.

In the formula below we use a function $\delta(\alpha',\phi')$ which
is equal
to $1$ when the stratum $\cH(\alpha')$ contains a hyperelliptic
component, and when, moreover, this hyperelliptic component has
parity
$\phi'$ of the spin structure. This function is equal to zero
otherwise (see~\eqref{eq:delta:hyp} for the explicit definition).
By convention we let the volume of several ``missing'' components
in small genera and the volume of nonexisting hyperelliptic
components
be equal to zero (see~\eqref{eq:dummy:volumes} for the
complete list).
\begin{formula}
\label{f:toitself:evenodd}
Let all the integers in $\alpha$ be even, and let $\alpha$ be
different
from any of the partitions $(2g-2), (g-1, g-1)$.

For  almost  all  flat  surfaces  in  the
connected component $\cH^{\phi}(\alpha)$ the constants in the
quadratic asymptotics for  the number of saddle connections
joining some zero to itself are presented by the following
formulae:

--- If the assignment contains only figure eight
constructions
\begin{multline*}
 c=\frac{1}{|\Gamma_-|\cdot|\Gamma|} \cdot
\prod_{m\in\alpha} \left(\frac{o(m)!}{\prod_{i=1}^p
o_i(m)!}\right)
\cdot \prod_{\substack{1\le i \le p\\a_i\in\alpha}} o_i(a_i)
\cdot \prod_{1\le i \le p} (a_i+1)\\
   %
   %
\cdot \frac{1}{2^{p-1}}
\cdot \frac{ \prod_{i=1}^p
(\tfrac{d_i}{2}-1)! }
     { (\frac{d}{2}-2)! }
\cdot \frac{1}{ \Vol(\cH^\phi_1(\alpha)) }\cdot\\
   %
   %
\cdot \sum_{\substack{\phi'_1,...,\phi'_p\in\{even,\, odd\}\\
  (\phi'_1+\dots+\phi'_p)+\qquad \\
  \quad + (a'_1+\dot+a'_p)+\qquad\\
  \qquad+p+1\equiv\phi\,
(mod\, 2)}} \prod_{i=1}^p
\Big(\Vol(\cH^{\phi_i}_1(\alpha'_i)) +
\delta(\alpha'_i,\phi'_i)\Vol(\cH^{hyp}_1(\alpha'_i))\Big)
\end{multline*}
\bigskip

--- If the assignment contains at least one parallelogram
construction with even $b'_i$
\begin{multline*}
c=\frac{1}{|\Gamma_-|}\cdot\frac{1}{|\Gamma|} \cdot
\prod_{m\in\alpha} \left(\frac{o(m)!}{\prod_{i=1}^p
o_i(m)!}\right) \cdot \\
\cdot \prod_{\substack{1\le i \le p\\z_i=w_i\\a_i\in\alpha}}
o_i(a_i) \cdot \prod_{\substack{1\le k \le p\\z_k\neq
w_k\\b'_k\neq b''_k\\b'_k\in\alpha}} o_k(b'_k) \cdot
\prod_{\substack{1\le k \le p\\z_k\neq w_k\\b'_k\neq
b''_k\\b''_k\in\alpha}}o_k(b_k'')\cdot \prod_{\substack{1\le
k \le p\\z_k\neq w_k\\b'_k=b''_k\in\alpha}} o_k(b'_k)
(o_k(b'_k)-1) \cdot \\
\cdot \prod_{\substack{1\le i \le p\\z_i=w_i}} (a_i+1)  \cdot
\prod_{\substack{1\le k \le p\\z_k\neq w_k}}
(b'_k+1)(b''_k+1)  \cdot \frac{1}{2^{p-1}}\cdot \frac{
\prod_{i=1}^p (\tfrac{d_i}{2}-1)! }
     { (\frac{d}{2}-2)! } \cdot\\
\cdot \frac{1}{ \Vol(\cH^\phi_1(\alpha)) }\cdot
\sum_{\substack{\phi'_1,...,\phi'_p\in\{even,\, odd\}\\
  (\phi'_1+\dots+\phi'_p)\equiv\phi\,
(mod\, 2)}} \prod_{i=1}^p
\Big(\Vol(\cH^{\phi_i}_1(\alpha'_i)) +
\delta(\alpha'_i,\phi'_i)\Vol(\cH^{hyp}_1(\alpha'_i))\Big)
\end{multline*}
\medskip

--- If the assignment contains at least one parallelogram
construction with odd $b'_i$
\begin{multline*}
c=\frac{1}{|\Gamma_-|}\cdot\frac{1}{|\Gamma|} \cdot
\prod_{m\in\alpha} \left(\frac{o(m)!}{\prod_{i=1}^p
o_i(m)!}\right) \cdot \\
\cdot \prod_{\substack{1\le i \le p\\z_i=w_i\\a_i\in\alpha}}
o_i(a_i) \cdot \prod_{\substack{1\le k \le p\\z_k\neq
w_k\\b'_k\neq b''_k\\b'_k\in\alpha}} o_k(b'_k) \cdot
\prod_{\substack{1\le k \le p\\z_k\neq w_k\\b'_k\neq
b''_k\\b''_k\in\alpha}}o_k(b_k'')\cdot \prod_{\substack{1\le
k \le p\\z_k\neq w_k\\b'_k=b''_k\in\alpha}} o_k(b'_k)
(o_k(b'_k)-1) \cdot \\
\cdot \prod_{\substack{1\le i \le p\\z_i=w_i}} (a_i+1)  \cdot
\prod_{\substack{1\le k \le p\\z_k\neq w_k}}
(b'_k+1)(b''_k+1) \cdot\\
\cdot\cfrac{1}{2}\, \cdot \frac{1}{2^{p-1}}\cdot \frac{
\prod_{i=1}^p (\tfrac{d_i}{2}-1)! }
     { (\frac{d}{2}-2)! } \cdot
\frac{ \prod_{i=1}^p \Vol( \cH_1(\alpha'_i)) }
     { \Vol(\cH_1^{\phi}(\alpha)) }
\qquad
\end{multline*}
\end{formula}
\begin{proof}
The formula above is obtained by elementary adjustment of
formula~\eqref{eq:sad:conn:to:itself:const}.

First suppose that the assignment under consideration does
not involve any parallelogram   constructions with odd  zeroes.
Then by Lemma~\ref{lm:spin:in:8:constr}          and
Lemma~\ref{lm:spin:in:even:par:constr} the  spin structure of
$S$ is  determined  by  the  spin  structures  of
$S'_i$ and by the parities of $a'_i$. Thus the
combinatorial constant is exactly the same as in
equation~\eqref{eq:sad:conn:to:itself:const}.
The numerator of the volume term  is      the      sum over those
products of
$\Vol(\cH^{even}(\alpha_i))$, $\Vol(\cH^{odd}(\alpha_j))$,
and $\Vol(\cH^{hyp}(\alpha_i))$ (when appropriate)
that satisfy the relevant relation
between spin structures of component surfaces and  the spin
structure
of the resulting surface, see
Lemma~\ref{lm:spin:in:8:constr} and
Lemma~\ref {lm:spin:in:even:par:constr} correspondingly.
In the denominator of the volume term we have
the   volume   of the component under consideration.

If there is a parallelogram  construction with odd $b_k'$
involved in the assignment, then by
Lemma~\ref{lm:spin:in:odd:par:constr} we have to multiply $c$
by $1/2$  to account for the fact that half of the
constructions produce flat surfaces with even spin structure
and another half of the constructions produce a flat surface
with  odd spin structure. Thus we need to divide by $2$ the
combinatorial factor
in equation~\eqref{eq:sad:conn:to:itself:const}.
In this case the numerator of the volume term is same as
in~\eqref{eq:sad:conn:to:itself:const}: it  is  a product of
$\Vol(\cH_1(\alpha_i))$, where we take the total volume of
all connected components in the stratum. In the denominator of the
volume term we again have the   volume   of the component under
consideration.
\end{proof}

We complete this section considering the exceptional strata.

\begin{formula}
Let $\alpha=(2g-2)$, or let $\alpha= (g-1, g-1)$ where $g$ is odd.
Let $g\ge 4$, or let $g=3$ and $\phi=odd$.

For  almost  all  flat  surfaces  in  the
connected component $\cH^{\phi}(\alpha)$ the constants in the
quadratic asymptotics for  the number of saddle connections
joining some zero to itself are presented by the following
formulae:

--- If the assignment is different from any of the
assignments listed in
Lemmas~\ref{lm:hyp:zero:to:itself:2g:minus2},
and~\ref{lm:hyp:zero:to:itself:gm1:gm1}, then the
constant is the same as in Formula~\ref{f:toitself:evenodd}.

--- For any of the assignments listed in one of the
Lemmas~\ref{lm:hyp:zero:to:itself:2g:minus2},
and~\ref{lm:hyp:zero:to:itself:gm1:gm1} the constant
is equal to
$$
c=c_{regular}-c_{hyp}\cdot\cfrac{\Vol(\cH^{hyp}(\alpha))}{\Vol(\cH^{\phi}(\alpha))}.
$$
Here $c_{regular}$ is given by Formula~\ref{f:toitself:evenodd},
and $c_{hyp}$ is the constant for the
corresponding assignment for the {\it hyperelliptic} component
$\cH^{hyp}(\alpha)$, see Formulae~\ref{f:toitself:hyp:2gm2}
and~\ref{f:toitself:hyp:gm1:gm1}.
\end{formula}

\begin{proof}
For nonexceptional assignments the proof is exactly the same
as for Formula~\ref{f:toitself:evenodd}. For the exceptional
assignments, which may lead to the flat surfaces from
hyperelliptic
components one makes appropriate adjustment, see the proof of
Formula~\ref{f:nonhyp:to:itself} and
Remark~\ref{rm:nonhyp:to:itself}
after this proof.
\end{proof}

\begin{example}{\bf Component} $\boldsymbol{\cH^{odd}(4)}.$
Recall  that  the stratum $\cH(4)$ has only one  nonhyperelliptic
connected  component  $\cH^{odd}(4)$.  Every   flat   surface  in
$\cH(4)$  having  even parity of the spin  structure  belongs  to
$\cH^{hyp}(4)$. We  have  $d=\dim_\reals  \cH(4)  = 12$. Consider
possible assignments for the flat  surfaces  from  the  component
$\cH^{odd}(4)$. We start with assignments of multiplicity one.

\medskip
\noindent $\boldsymbol{ \to (\overline{2},\overline{0}) \to }$\\
\noindent
We  have  a  single  saddle  connection  that does  not  bound  a
cylinder. It returns at    angle $3\pi$. After degeneration we
get a  surface  of  genus $2$ with a zero of order $2$ and with a
marked point. We apply the parallelogram  construction gluing the
circles directly.  Since  $b'_1=2$,  and  $b''_1=0$  are even, by
lemma~\ref{lm:spin:in:even:par:constr}  the parity  of  the  spin
structure of the resulting  surface  coincides with the parity of
the spin structure of a  surface  from $\cH(2)$, so it is  really
odd, see~\eqref{eq:spin:2g:minus2}. Since $b'_1\neq b''_2$ we  do
not  have  a  $\gamma\to-\gamma$  symmetry.
Thus,
$$ 
c=3\cdot \frac{(5-1)!}{(6-2)!} \cdot
\frac{\Vol(\cH(2))}{\Vol(\cH^{odd}(4))} =
\frac{81}{8}\cdot\frac{1}{\zeta(2)} \approx 6.155
$$ 

\medskip
\noindent $\boldsymbol{ \to (\overline{1},\overline{1}) \to }$\\
\noindent
Here we again have a single saddle connection that does not bound
a  cylinder,  but now  it  returns  at    angle  $5\pi$. After
degeneration we get a surface of  genus $2$ with a pair of simple
zeroes.  We  apply  the  parallelogram  construction  gluing  the
circles  directly.   Since   $b'_1=1$,   and   $b''_1=1$  we  get
$N=(b'_1+1)(b''_1+1)=4$ surfaces, but only half of  them have odd
parity       of        the       spin       structure,        see
lemma~\ref{lm:spin:in:odd:par:constr}.    Now     we    have    a
$\gamma\to-\gamma$    symmetry,    $|\Gamma_-|=2$.
   %
$$ 
c=1\cdot \frac{(5-1)!}{(6-2)!}\cdot
\frac{\Vol(\cH(1,1))}{\Vol(\cH^{odd}(4))} =
3\cdot \frac{1}{\zeta(2)} \approx 1.824
$$ 

\medskip
\noindent $\boldsymbol{ \Rightarrow  (\overline{0+2}) \Rightarrow
}$\\ \noindent
Suppose now that the saddle  connection  returns  at angle $\pi$.
Then  it  bounds   a   cylinder  filled  with  homologous  closed
geodesics. In  particular,  there  is  another  saddle connection
returning at the angle $\pi$ homologous to the initial one. Since
the surface is not hyperelliptic,  the  angle  between the saddle
connections  on  the opposite sides of the  cylinder  is  $6\pi$.
After degeneration we get a single zero of order $2$ on a surface
of genus $2$. We apply  the  figure  eight construction, breaking
the zero  of  order $2$ into a zero of order $2$ and one of order
$0$. This  corresponds to spacings of  $2\pi$ and $6\pi$,  or, in
our notation to  the choice $a'_1=0$,  $a''_1=2$. We glue  in  an
intervening  cylinder.   By  lemma~\ref{lm:spin:in:8:constr}  the
parity of the spin structure of the resulting flat surface equals
$1+\phi(S'_1)+(0+1)=1$      since       $\phi(S'_1)=1$,
see~\eqref{eq:spin:2g:minus2}. Since $a'_1\neq a''_1$  we  do not
have the $\gamma\to-\gamma$ symmetry.
   %
$$ 
c=3\cdot \frac{(4-1)!}{(6-2)!} \cdot
\frac{\Vol(\cH(2))}{\Vol(\cH^{odd}(4))}=
\frac{81}{32}\cdot\frac{1}{\zeta(2)} \approx 1.539
$$ 

\noindent $\boldsymbol{ \Rightarrow  (\overline{0+0}) \to (\overline{0+0}) \Rightarrow }$\\
\noindent
We  now  consider multiplicity  $2$.  We have  a  cylinder and  a
homologous curve that  does not bound a cylinder. The degenerated
surfaces $S_1',S_2'$ are  a  pair of  tori,  each with a  marked
point. We  apply the figure eight  construction giving a  pair of
circles on  each torus. We glue one pair  of circles directly and
put  an  intervening cylinder  between  the  other  pair.  Though
$\alpha'_1=\alpha'_2$ we  do  not  have  a  stratum  interchange,
$|\Gamma|=1$, but  we have the $\gamma\to-\gamma$ symmetry.
   %
$$ 
c=\frac{1}{2}\cdot\frac{1}{2^{2-1}}\cdot \frac{(2-1)!\cdot
(2-1)!}{(6-2)!}\cdot
\frac{\Vol(\cH(\torusemptyset))^2}{\Vol(\cH^{odd}(4))} =
\frac{15}{32} \cdot\frac{1}{\zeta(2)} \approx 0.2850
$$ 
\end{example}

\begin{example} {\bf Component} $\boldsymbol{\cH^{odd}(2,2)}.$
\label{ex:hodd22}

\noindent $\boldsymbol{ \to (2,\overline{0},\overline{0}) \to }$\\
\noindent
 In the multiplicity
$1$ case with no cylinder, the degenerating surface  is genus $2$
with a pair of marked points  and a zero of order $2$. We perform
the creating a pair of holes construction  on  the marked points, gluing the
circles   directly.    Since    $b'_1=b''_1=0$    is   even,   by
lemma~\ref{lm:spin:in:even:par:constr}  the parity  of  the  spin
structure   equals   $\phi(S')$,   which   is   odd    for
$S'\in\cH(2)$.
   %
$$ 
c=1\cdot\frac{(4+2-1)!}{(7-2)!} \cdot
\frac{\Vol(\cH_1(2))}{\Vol(\cH_1^{odd}(2,2))}=
6\cdot\frac{1}{\zeta(2)}\approx 3.648
$$ 

\noindent $\boldsymbol{ \Rightarrow (2,\overline{0+0}) \Rightarrow
}$\\ \noindent
In the case  of  a cylinder with multiplicity  $1$,  if the other
side returns to the same zero,  we get a genus $2$ surface with a
single marked point  and  a double  zero.  We perform the  figure
eight construction  at the marked  point, gluing in a cylinder to
the   two    circles.    We   have   $a'_1=a''_1=0$.   Thus,   by
lemma~\ref{lm:spin:in:8:constr}
$\phi(S)=1+\phi(S')+(0+1)=1+1+1=1\ \mod  2$. Now $|\Gamma_-|=2$.
   %
$$ 
c=1 \cdot \frac{(4+1-1)!}{(7-2)!} \cdot
\frac{\Vol(\cH_1(2))}{\Vol(\cH_1^{odd}(2,2))}=
\frac{6}{5}\cdot\frac{1}{\zeta(2)}\approx 0.7295
$$ 

\noindent $\boldsymbol{ \Rightarrow (\overline{1},\overline{1})
\Rightarrow }$\\ \noindent
If  the  other side of the  cylinder  returns to the other  zero,
after  degeneration  we  get  $2$ simple zeroes. We  perform  the
parallelogram   construction   gluing   in   a   cylinder.    Now
$b'_1=b''_1=1$  is  odd.  Thus  we  get   $N=(b'_1+1)(b''_1+1)=4$
surfaces,  but  by   lemma~\ref{lm:spin:in:odd:par:constr}   only
$N/2=2$ of them have odd parity of the spin structure.
There is also the $\gamma\to-\gamma$ symmetry.
   %
$$ 
c=2\cdot\frac{(5-1)!}{(7-2)!}\cdot
\frac{\Vol(\cH_1(1,1))}{\Vol(\cH_1^{odd}(2,2))}=
\frac{32}{15}\cdot\frac{1}{\zeta(2)}\approx 1.297
$$ 

\noindent $\boldsymbol{ \Rightarrow (\overline{0+0}) \Rightarrow
(\overline{0+0}) \Rightarrow
}$\\ \noindent
In the multiplicity $2$ case we have a pair of cylinders. The two
sides of each cylinder return to different zeroes.  The result is
a pair of tori,  each with a single marked point. We  perform the
figure eight construction  creating two pairs of circles to which
we glue intervening cylinders.  We  have $|\Gamma|=|\Gamma_-|=2$.
   %
$$ 
c=\cfrac{1}{2}\cdot\cfrac{1}{2^{2-1}}
\cdot\frac{(2-1)!\,(2-1)!}{(7-2)!}\cdot
\frac{\Vol(\cH_1(\torusemptyset))\Vol(\cH_1(\torusemptyset))}{\Vol(\cH_1^{odd}(2,2))}=
\frac{1}{6}\cdot\frac{1}{\zeta(2)}\approx 0.1013
$$ 
\end{example}

\bigskip
\centerline{\sc{Acknowledgments}}

The third author wants to thank MPI f\"ur
Mathematik at  Bonn and IHES at Bures-sur-Yvette for their
hospitality while  preparation of this  paper.

%
%

\newpage
\appendix
\addcontentsline{toc}{part}
{Appendix. Values of the Constants}
{\Large\bf Appendix. Values of the Constants}

\section{Saddle Connections Joining Distinct Zeroes}

\subsection{Connected Strata in Genus 4}
In the tables  below we consider  all connected strata  in  genus
$g=4$.  For  each   stratum  we  present  the  complete  list  of
all    possible topological configurations  of   geodesic   saddle
connections.  To  make calculations traceable  we  present  for
every configuration the orders  of  possible  symmetries  $\Gamma$
and  $\Gamma_-$,  the combinatorial  constant  $M$,  and  the
result: the Siegel---Veech constant  in quadratic  asymptotics  for
the number of  saddle  connections (closed geodesics) of that
type. The notations were introduced in
Section~\ref{ss:slit:construction} and
Section~\ref{ss:admissible:constructions}

\scriptsize
$$ 
\begin{array}{|c|c|c|c|c|c|}
\multicolumn{6}{c}{}\\
[-\halfbls]
\multicolumn{6}{c}{\text{Stratum }\cH(5, 1).}\\
\hline &&&&&\\
\text{Degeneration pattern} & |\Gamma_-| & |\Gamma| & M  & c & c \text{ approx.} \\
[-\halfbls]  &&&&&\\ \hline  &&&&& \\ [-\halfbls]
(1+5) \succ  & 1 & 1 & 7 & \cfrac{4311167}{373248} & 11.5504\\
[-\halfbls]  &&&&&\\ \hline  &&&&& \\ [-\halfbls]
(0+0) \succ (0+4) \succ  & 1 & 1 & 5 & \cfrac{38125}{93312} & 0.408576\\
[-\halfbls]  &&&&&\\ \hline  &&&&& \\ [-\halfbls]
(0+2) \succ (0+2) \succ  & 1 & 2 & \cfrac{9}{2} & \cfrac{21}{512} & 0.0410156\\
[-\halfbls]  &&&&&\\ \hline
%
%
\multicolumn{6}{c}{}\\
\multicolumn{6}{c}{}\\
[-\halfbls]\multicolumn{6}{c}{\text{Stratum }\cH(4, 1, 1).}\\
\hline &&&&&\\  
[-\halfbls]
(1+1, 4) \succ  & 2 & 1 & 3 & \cfrac{2403}{616} & 3.90097\\
[-\halfbls]  &&&&&\\ \hline  &&&&& \\ [-\halfbls]
(1+4, 1) \succ  & 1 & 1 & 12 & \cfrac{186624}{9625} & 19.3895\\
[-\halfbls]  &&&&&\\ \hline  &&&&& \\ [-\halfbls]
(0+0) \succ (0+0, 4) \succ  & 2 & 1 & 1 & \cfrac{61}{616} & 0.099026\\
[-\halfbls]  &&&&&\\ \hline  &&&&& \\ [-\halfbls]
(0+0) \succ (0+3, 1) \succ  & 1 & 1 & 8 & \cfrac{1024}{1925} & 0.531948\\
[-\halfbls]  &&&&&\\ \hline  &&&&& \\ [-\halfbls]
(0+2) \succ (0+1, 1) \succ  & 1 & 1 & 12 & \cfrac{108}{1375} & 0.0785455\\
[-\halfbls]  &&&&&\\ \hline
%
%
\multicolumn{6}{c}{}\\
\multicolumn{6}{c}{}\\
[-\halfbls]\multicolumn{6}{c}{\text{Stratum }\cH(3, 2, 1).}\\
\hline &&&&&\\
[-\halfbls]
(1+2, 3) \succ  & 1 & 1 & 4 & \cfrac{368}{63} & 5.84127\\
[-\halfbls]  &&&&&\\ \hline  &&&&& \\ [-\halfbls]
(1+3, 2) \succ  & 1 & 1 & 5 & \cfrac{55625}{7168} & 7.76018\\
[-\halfbls]  &&&&&\\ \hline  &&&&& \\ [-\halfbls]
(2+3, 1) \succ  & 1 & 1 & 6 & \cfrac{81}{7} & 11.5714\\
[-\halfbls]  &&&&&\\ \hline  &&&&& \\ [-\halfbls]
(0+0) \succ (0+2, 2) \succ  & 1 & 1 & 3 & \cfrac{765}{3584} & 0.213449\\
[-\halfbls]  &&&&&\\ \hline  &&&&& \\ [-\halfbls]
(0+0) \succ (0+1, 3) \succ  & 1 & 1 & 2 & \cfrac{10}{63} & 0.15873\\
[-\halfbls]  &&&&&\\ \hline  &&&&& \\ [-\halfbls]
(0+0) \succ (1+2, 1) \succ  & 1 & 1 & 4 & \cfrac{20}{63} & 0.31746\\
[-\halfbls]  &&&&&\\ \hline  &&&&& \\ [-\halfbls]
(0+0, 2) \succ (0+2) \succ  & 1 & 1 & 3 & \cfrac{27}{1024} & 0.0263672\\
[-\halfbls]  &&&&&\\ \hline  &&&&& \\ [-\halfbls]
(0+2) \succ (1+0, 1) \succ  & 1 & 1 & 6 & \cfrac{3}{64} & 0.046875\\
[-\halfbls]  &&&&&\\ \hline  &&&&& \\ [-\halfbls]
(1+1) \succ (0+1, 1) \succ  & 1 & 1 & 6 & \cfrac{3}{64} & 0.046875\\
[-\halfbls]  &&&&&\\ \hline  &&&&& \\ [-\halfbls]
(0+0) \succ (0+0) \succ (0+1, 1) \succ  & 1 & 1 & 2 & \cfrac{5}{288} & 0.0173611\\
[-\halfbls]  &&&&&\\ \hline
\end{array}
$$ 

\pagebreak 

\scriptsize
$$ 
\begin{array}{|c|c|c|c|c|c|}
\multicolumn{6}{c}{}\\
[-\halfbls]
\multicolumn{6}{c}{\text{ Configurations of geodesic saddle connections.}}\\
\multicolumn{6}{c}{\text{Stratum }\cH(3, 1, 1, 1).}\\
\multicolumn{6}{c}{}\\
\hline &&&&&\\
\text{Degeneration pattern} & |\Gamma_-| & |\Gamma| & M  & c & c \text{ approx.} \\
[-\halfbls]  &&&&&\\ \hline  &&&&& \\ [-\halfbls]
(1+1, 3, 1) \succ  & 2 & 1 & 9 & \cfrac{729}{62} & 11.7581\\
[-\halfbls]  &&&&&\\ \hline  &&&&& \\ [-\halfbls]
(1+3, 1, 1) \succ  & 1 & 1 & 15 & \cfrac{185625}{7936} & 23.3902\\
[-\halfbls]  &&&&&\\ \hline  &&&&& \\ [-\halfbls]
(0+0) \succ (0+0, 3, 1) \succ  & 2 & 1 & 3 & \cfrac{15}{62} & 0.241935\\
[-\halfbls]  &&&&&\\ \hline  &&&&& \\ [-\halfbls]
(0+0) \succ (0+2, 1, 1) \succ  & 1 & 1 & 9 & \cfrac{2025}{3968} & 0.510333\\
[-\halfbls]  &&&&&\\ \hline  &&&&& \\ [-\halfbls]
(0+2) \succ (0+0, 1, 1) \succ  & 1 & 1 & 9 & \cfrac{405}{7936} & 0.0510333\\
[-\halfbls]  &&&&&\\ \hline  &&&&& \\ [-\halfbls]
(0+1, 1) \succ (0+1, 1) \succ  & 1 & 2 & 12 & \cfrac{3}{62} & 0.0483871\\
[-\halfbls]  &&&&&\\ \hline
%
%
\multicolumn{6}{c}{}\\
\multicolumn{6}{c}{}\\
[-\halfbls]\multicolumn{6}{c}{\text{Stratum }\cH(2, 2, 1, 1).}\\
\multicolumn{6}{c}{}\\
\hline &&&&&\\
[-\halfbls]
(1+1, 2, 2) \succ  & 2 & 1 & 3 & \cfrac{4101}{1048} & 3.91317\\
[-\halfbls]  &&&&&\\ \hline  &&&&& \\ [-\halfbls]
(1+2, 2, 1) \succ  & 1 & 1 & 16 & \cfrac{3072}{131} & 23.4504\\
[-\halfbls]  &&&&&\\ \hline  &&&&& \\ [-\halfbls]
(2+2, 1, 1) \succ  & 2 & 1 & 5 & \cfrac{6875}{786} & 8.74682\\
[-\halfbls]  &&&&&\\ \hline  &&&&& \\ [-\halfbls]
(0+0) \succ (0+0, 2, 2) \succ  & 2 & 1 & 1 & \cfrac{85}{1048} & 0.0811069\\
[-\halfbls]  &&&&&\\ \hline  &&&&& \\ [-\halfbls]
(0+0) \succ (0+1, 2, 1) \succ  & 1 & 1 & 8 & \cfrac{200}{393} & 0.508906\\
[-\halfbls]  &&&&&\\ \hline  &&&&& \\ [-\halfbls]
(0+0) \succ (1+1, 1, 1) \succ  & 2 & 1 & 3 & \cfrac{25}{131} & 0.19084\\
[-\halfbls]  &&&&&\\ \hline  &&&&& \\ [-\halfbls]
(0+0, 2) \succ (0+0, 2) \succ  & 2 & 2 & 1 & \cfrac{3}{524} & 0.00572519\\
[-\halfbls]  &&&&&\\ \hline  &&&&& \\ [-\halfbls]
(0+0, 2) \succ (0+1, 1) \succ  & 1 & 1 & 8 & \cfrac{16}{393} & 0.0407125\\
[-\halfbls]  &&&&&\\ \hline  &&&&& \\ [-\halfbls]
(1+1) \succ (0+0, 1, 1) \succ  & 2 & 1 & 3 & \cfrac{5}{262} & 0.019084\\
[-\halfbls]  &&&&&\\ \hline  &&&&& \\ [-\halfbls]
(0+1, 1) \succ (1+0, 1) \succ  & 2 & 1 & 8 & \cfrac{128}{3537} & 0.0361889\\
[-\halfbls]  &&&&&\\ \hline  &&&&& \\ [-\halfbls]
(0+0) \succ (0+0) \succ (0+0, 1, 1) \succ  & 2 & 1 & 1 & \cfrac{25}{3537} & 0.00706814\\
[-\halfbls]  &&&&&\\ \hline
%
%
\multicolumn{6}{c}{}\\
\multicolumn{6}{c}{}\\
[-\halfbls]\multicolumn{6}{c}{\text{Stratum }\cH(2, 1, 1, 1, 1).}\\
\multicolumn{6}{c}{}\\
\hline &&&&&\\
[-\halfbls]
(1+1, 2, 1, 1) \succ  & 2 & 1 & 18 & \cfrac{1179}{50} & 23.58\\
[-\halfbls]  &&&&&\\ \hline  &&&&& \\ [-\halfbls]
(1+2, 1, 1, 1) \succ  & 1 & 1 & 16 & \cfrac{15872}{675} & 23.5141\\
[-\halfbls]  &&&&&\\ \hline  &&&&& \\ [-\halfbls]
(0+0) \succ (0+0, 2, 1, 1) \succ  & 2 & 1 & 6 & \cfrac{2}{5} & 0.4\\
[-\halfbls]  &&&&&\\ \hline  &&&&& \\ [-\halfbls]
(0+0) \succ (0+1, 1, 1, 1) \succ  & 1 & 1 & 8 & \cfrac{56}{135} & 0.414815\\
[-\halfbls]  &&&&&\\ \hline  &&&&& \\ [-\halfbls]
(0+0, 2) \succ (0+0, 1, 1) \succ  & 2 & 1 & 6 & \cfrac{1}{50} & 0.02\\
[-\halfbls]  &&&&&\\ \hline  &&&&& \\ [-\halfbls]
(0+0, 1, 1) \succ (0+1, 1) \succ  & 1 & 1 & 24 & \cfrac{16}{225} & 0.0711111\\
[-\halfbls]  &&&&&\\ \hline
\end{array}
$$


\newpage

$ \ \phantom{=} $

\begin{table}[!hbt]
\normalsize
%
\caption{
Principal stratum $\cH(1,\dots,1)$.  Approximate  values
of the constants $c$  for  saddle connections of multiplicity two
joining a pair of distinct zeroes.}
\label{tab:c:distinct:principal:mult2:approx}
$$
\cH(\overbrace{1,\dots,1}^{2g_1-2},0+0)\succ\cH(\overbrace{1,\dots,1}^{2g_2-2},0+0)\succ
\qquad\qquad g_1+g_2=g
$$

%

\scriptsize
$$ 
\begin{array}{|l|c|c|c|c|c|c|}
\hline &&&&&&\\
& g_2=1 & g_2=2 & g_2=3 & g_2=4 & g_2=5 & g_2=6 \\
[-\halfbls] &&&&&& \\ \hline &&&&&& \\ [-\halfbls]
 g_1=1 & 0.6250 & 0.8571 & 0.8355 & 0.8393 & 0.8426 & 0.8435\\
 [-\halfbls] &&&&&& \\ \hline &&&&&& \\ [-\halfbls]
 g_1=2 & 0.8571 & 0.07958 & 0.07193 & 0.04160 & 0.02708 & 0.01896\\
 [-\halfbls] &&&&&& \\ \hline &&&&&& \\ [-\halfbls]
 g_1=3 & 0.8355 & 0.07193 & 0.009358 & 0.007018 & 0.003195 & 0.001651\\
 [-\halfbls] &&&&&& \\ \hline &&&&&& \\ [-\halfbls]
 g_1=4 & 0.8393 & 0.04160 & 0.007018 & 0.9204\cdot 10^{-3} & 0.6185\cdot 10^{-3} & 0.2454\cdot 10^{-3}\\
 [-\halfbls] &&&&&& \\ \hline &&&&&& \\ [-\halfbls]
 g_1=5 & 0.8426 & 0.02708 & 0.003195 & 0.6185\cdot 10^{-3} & 0.7978\cdot 10^{-4} & 0.5011\cdot 10^{-4}\\
 [-\halfbls] &&&&&& \\ \hline &&&&&& \\ [-\halfbls]
 g_1=6 & 0.8435 & 0.01896 & 0.001651 & 0.2454\cdot 10^{-3} & 0.5011\cdot 10^{-4} & 0.6383\cdot 10^{-5}\\
 [-\halfbls] &&&&&& \\ \hline &&&&&& \\ [-\halfbls]
 g_1=7 & 0.8430 & 0.01398 & 0.9350\cdot 10^{-3} & 0.1100\cdot 10^{-3} & 0.1822\cdot 10^{-4} & 0.3840\cdot 10^{-5}\\
 [-\halfbls] &&&&&& \\ \hline &&&&&& \\ [-\halfbls]
 g_1=8 & 0.8418 & 0.01072 & 0.5675\cdot 10^{-3} & 0.5415\cdot 10^{-4} & 0.7420\cdot 10^{-5} & 0.1315\cdot 10^{-5}\\
 [-\halfbls] &&&&&& \\ \hline &&&&&& \\ [-\halfbls]
 g_1=9 & 0.8406 & 0.008475 & 0.3638\cdot 10^{-3} & 0.2872\cdot 10^{-4} & 0.3309\cdot 10^{-5} & 0.5001\cdot 10^{-6}\\
 [-\halfbls] &&&&&& \\ \hline &&&&&& \\ [-\halfbls]
 g_1=10 & 0.8393 & 0.006864 & 0.2438\cdot 10^{-3} & 0.1618\cdot 10^{-4} & 0.1590\cdot 10^{-5} & 0.2072\cdot 10^{-6}\\
 [-\halfbls] &&&&&& \\ \hline &&&&&& \\ [-\halfbls]
 g_1=11 & 0.8382 & 0.005670 & 0.1694\cdot 10^{-3} & 0.9585\cdot 10^{-5} & 0.8122\cdot 10^{-6} & 0.9226\cdot 10^{-7}\\
 [-\halfbls] &&&&&& \\ \hline &&&&&& \\ [-\halfbls]
 g_1=12 & 0.8372 & 0.004762 & 0.1213\cdot 10^{-3} & 0.5921\cdot 10^{-5} & 0.4372\cdot 10^{-6} & 0.4366\cdot 10^{-7}\\
 [-\halfbls] &&&&&& \\ \hline &&&&&& \\ [-\halfbls]
 g_1=13 & 0.8363 & 0.004056 & 0.8909\cdot 10^{-4} & 0.3790\cdot 10^{-5} & 0.2461\cdot 10^{-6} & 0.2177\cdot 10^{-7}\\
 [-\halfbls] &&&&&& \\ \hline &&&&&& \\ [-\halfbls]
 g_1=14 & 0.8354 & 0.003495 & 0.6691\cdot 10^{-4} & 0.2503\cdot 10^{-5} & 0.1440\cdot 10^{-6} & 0.1136\cdot 10^{-7}\\
 [-\halfbls] &&&&&& \\ \hline &&&&&& \\ [-\halfbls]
 g_1=15 & 0.8347 & 0.003043 & 0.5122\cdot 10^{-4} & 0.1697\cdot 10^{-5} & 0.8712\cdot 10^{-7} & 0.6173\cdot 10^{-8}\\
 [-\halfbls] &&&&&& \\ \hline &&&&&& \\ [-\halfbls]
 g_1=16 & 0.8340 & 0.002673 & 0.3986\cdot 10^{-4} & 0.1179\cdot 10^{-5} & 0.5431\cdot 10^{-7} & 0.3474\cdot 10^{-8}\\
 [-\halfbls] &&&&&& \\ \hline &&&&&& \\ [-\halfbls]
 g_1=17 & 0.8334 & 0.002367 & 0.3149\cdot 10^{-4} & 0.8358\cdot 10^{-6} & 0.3476\cdot 10^{-7} & 0.2017\cdot 10^{-8}\\
 [-\halfbls] &&&&&& \\ \hline &&&&&& \\ [-\halfbls]
 g_1=18 & 0.8329 & 0.002110 & 0.2520\cdot 10^{-4} & 0.6038\cdot 10^{-6} & 0.2278\cdot 10^{-7} & 0.1204\cdot 10^{-8}\\
 [-\halfbls] &&&&&& \\ \hline &&&&&& \\ [-\halfbls]
 g_1=19 & 0.8324 & 0.001893 & 0.2041\cdot 10^{-4} & 0.4436\cdot 10^{-6} & 0.1525\cdot 10^{-7} & 0.7378\cdot 10^{-9}\\
 [-\halfbls] &&&&&& \\ \hline &&&&&& \\ [-\halfbls]
 g_1=20 & 0.8319 & 0.001708 & 0.1670\cdot 10^{-4} & 0.3308\cdot 10^{-6} & 0.1041\cdot 10^{-7} & 0.4624\cdot 10^{-9}\\
 [-\halfbls] &&&&&& \\ \hline &&&&&& \\ [-\halfbls]
 g_1=21 & 0.8315 & 0.001549 & 0.1380\cdot 10^{-4} & 0.2501\cdot 10^{-6} & 0.7226\cdot 10^{-8} & 0.2960\cdot 10^{-9}\\
 [-\halfbls] &&&&&& \\ \hline &&&&&& \\ [-\halfbls]
 g_1=22 & 0.8312 & 0.001411 & 0.1150\cdot 10^{-4} & 0.1915\cdot 10^{-6} & 0.5098\cdot 10^{-8} & 0.1931\cdot 10^{-9}\\
 [-\halfbls] &&&&&& \\ \hline &&&&&& \\ [-\halfbls]
 g_1=23 & 0.8308 & 0.001290 & 0.9663\cdot 10^{-5} & 0.1482\cdot 10^{-6} & 0.3650\cdot 10^{-8} & 0.1282\cdot 10^{-9}\\
 [-\halfbls] &&&&&& \\ \hline &&&&&& \\ [-\halfbls]
 g_1=24 & 0.8305 & 0.001185 & 0.8177\cdot 10^{-5} & 0.1160\cdot 10^{-6} & 0.2648\cdot 10^{-8} & 0.8649\cdot 10^{-10}\\
 [-\halfbls] &&&&&& \\ \hline &&&&&& \\ [-\halfbls]
 g_1=25 & 0.8302 & 0.001091 & 0.6966\cdot 10^{-5} & 0.9163\cdot 10^{-7} & 0.1946\cdot 10^{-8} & 0.5923\cdot 10^{-10}\\
 [-\halfbls] &&&&&& \\ \hline &&&&&& \\ [-\halfbls]
 g_1=26 & 0.8299 & 0.001009 & 0.5971\cdot 10^{-5} & 0.7304\cdot 10^{-7} & 0.1446\cdot 10^{-8} & 0.4113\cdot 10^{-10}\\
 [-\halfbls] &&&&&& \\ \hline &&&&&& \\ [-\halfbls]
 g_1=27 & 0.8297 & 0.9353\cdot 10^{-3} & 0.5147\cdot 10^{-5} & 0.5870\cdot 10^{-7} & 0.1086\cdot 10^{-8} & 0.2893\cdot 10^{-10}\\
 [-\halfbls] &&&&&& \\ \hline &&&&&& \\ [-\halfbls]
 g_1=28 & 0.8294 & 0.8695\cdot 10^{-3} & 0.4461\cdot 10^{-5} & 0.4754\cdot 10^{-7} & 0.8236\cdot 10^{-9} & 0.2060\cdot 10^{-10}\\
 [-\halfbls] &&&&&& \\ \hline &&&&&& \\ [-\halfbls]
 g_1=29 & 0.8292 & 0.8104\cdot 10^{-3} & 0.3885\cdot 10^{-5} & 0.3878\cdot 10^{-7} & 0.6305\cdot 10^{-9} & 0.1483\cdot 10^{-10}\\
 [-\halfbls] &&&&&& \\ \hline &&&&&& \\ [-\halfbls]
 g_1=30 & 0.8290 & 0.7571\cdot 10^{-3} & 0.3400\cdot 10^{-5} & 0.3184\cdot 10^{-7} & 0.4869\cdot 10^{-9} & 0.1079\cdot 10^{-10}\\
[-\halfbls] &&&&&&\\ \hline
\end{array}
$$ 
\end{table}

\newpage

$ \ \phantom{=} $

\begin{table}[!hbt]
\caption{Connected  component  $\cH^{hyp}(g-1,g-1)$;  approximate
values  of  the constants for saddle connections of  multiplicity
two joining distinct zeroes; $g_1+g_2=g$}
\label{tab:c:dist:hyp:mult2}
\normalsize

$$
\cH_1^{hyp}\big((\overline{g_1-1})+(\overline{g_1-1})\big) \succ
\cH_1^{hyp}(\overline{(g_2-1)}+\overline{(g_2-1)}) \succ
$$

\scriptsize
$$ 
\!
\begin{array}{|l|c|c|c|c|c|c|c|c|c|c|}
\hline &&&&&&&&&&\\
& g_2=1 & g_2=2 & g_2=3 & g_2=4 & g_2=5 & g_2=6 & g_2=7 & g_2=8 & g_2=9 &\hspace*{-2truept} g_2\!=\!10\hspace*{-3truept}\\
[-\halfbls] &&&&&&&&&& \\ \hline &&&&&&&&&& \\ [-\halfbls]
 g_1=1 & 0.6250 & 1.9687 & 2.9297 & 4.1122 & 5.5112 & 7.1245 & 8.9513 & 10.991 & 13.244 & 15.709\\
 [-\halfbls] &&&&&&&&&& \\ \hline &&&&&&&&&& \\ [-\halfbls]
 g_1=2 & 1.9687 & 1.1074 & 2.7191 & 3.3679 & 4.1334 & 5.0053 & 5.9792 & 7.0526 & 8.2244 & 9.4935\\
 [-\halfbls] &&&&&&&&&& \\ \hline &&&&&&&&&& \\ [-\halfbls]
 g_1=3 & 2.9297 & 2.7191 & 1.4729 & 3.3412 & 3.8411 & 4.4225 & 5.0749 & 5.7932 & 6.5745 & 7.4168\\
 [-\halfbls] &&&&&&&&&& \\ \hline &&&&&&&&&& \\ [-\halfbls]
 g_1=4 & 4.1122 & 3.3679 & 3.3412 & 1.7750 & 3.8803 & 4.2917 & 4.7662 & 5.2948 & 5.8725 & 6.4961\\
 [-\halfbls] &&&&&&&&&& \\ \hline &&&&&&&&&& \\ [-\halfbls]
 g_1=5 & 5.5112 & 4.1334 & 3.8411 & 3.8803 & 2.0372 & 4.3611 & 4.7134 & 5.1173 & 5.5653 & 6.0526\\
 [-\halfbls] &&&&&&&&&& \\ \hline &&&&&&&&&& \\ [-\halfbls]
 g_1=6 & 7.1245 & 5.0053 & 4.4225 & 4.2917 & 4.3611 & 2.2714 & 4.7983 & 5.1083 & 5.4619 & 5.8527\\
 [-\halfbls] &&&&&&&&&& \\ \hline &&&&&&&&&& \\ [-\halfbls]
 g_1=7 & 8.9513 & 5.9792 & 5.0749 & 4.7662 & 4.7134 & 4.7983 & 2.4849 & 5.2016 & 5.4798 & 5.7955\\
 [-\halfbls] &&&&&&&&&& \\ \hline &&&&&&&&&& \\ [-\halfbls]
 g_1=8 & 10.991 & 7.0526 & 5.7932 & 5.2948 & 5.1173 & 5.1083 & 5.2016 & 2.6821 & 5.5776 & 5.8309\\
 [-\halfbls] &&&&&&&&&& \\ \hline &&&&&&&&&& \\ [-\halfbls]
 g_1=9 & 13.244 & 8.2244 & 6.5745 & 5.8725 & 5.5653 & 5.4619 & 5.4798 & 5.5776 & 2.8663 & 5.9309\\
 [-\halfbls] &&&&&&&&&& \\ \hline &&&&&&&&&& \\ [-\halfbls]
 g_1=10 & 15.709 & 9.4935 & 7.4168 & 6.4961 & 6.0526 & 5.8527 & 5.7955 & 5.8309 & 5.9309 & 3.0396\\
 [-\halfbls] &&&&&&&&&& \\ \hline &&&&&&&&&& \\ [-\halfbls]
 g_1=11 & 18.387 & 10.860 & 8.3192 & 7.1637 & 6.5763 & 6.2764 & 6.1434 & 6.1170 & 6.1642 & 6.2652\\
 [-\halfbls] &&&&&&&&&& \\ \hline &&&&&&&&&& \\ [-\halfbls]
 g_1=12 & 21.277 & 12.322 & 9.2807 & 7.8740 & 7.1344 & 6.7305 & 6.5198 & 6.4315 & 6.4266 & 6.4821\\
 [-\halfbls] &&&&&&&&&& \\ \hline &&&&&&&&&& \\ [-\halfbls]
 g_1=13 & 24.379 & 13.881 & 10.301 & 8.6259 & 7.7255 & 7.2128 & 6.9221 & 6.7710 & 6.7143 & 6.7250\\
 [-\halfbls] &&&&&&&&&& \\ \hline &&&&&&&&&& \\ [-\halfbls]
 g_1=14 & 27.694 & 15.536 & 11.379 & 9.4190 & 8.3487 & 7.7223 & 7.3487 & 7.1334 & 7.0243 & 6.9906\\
 [-\halfbls] &&&&&&&&&& \\ \hline &&&&&&&&&& \\ [-\halfbls]
 g_1=15 & 31.221 & 17.288 & 12.516 & 10.253 & 9.0032 & 8.2577 & 7.7983 & 7.5170 & 7.3547 & 7.2764\\
 [-\halfbls] &&&&&&&&&& \\ \hline &&&&&&&&&& \\ [-\halfbls]
 g_1=16 & 34.961 & 19.135 & 13.710 & 11.127 & 9.6885 & 8.8184 & 8.2697 & 7.9205 & 7.7040 & 7.5806\\
 [-\halfbls] &&&&&&&&&& \\ \hline &&&&&&&&&& \\ [-\halfbls]
 g_1=17 & 38.913 & 21.078 & 14.962 & 12.041 & 10.404 & 9.4038 & 8.7624 & 8.3431 & 8.0710 & 7.9019\\
 [-\halfbls] &&&&&&&&&& \\ \hline &&&&&&&&&& \\ [-\halfbls]
 g_1=18 & 43.077 & 23.116 & 16.271 & 12.994 & 11.150 & 10.014 & 9.2757 & 8.7839 & 8.4548 & 8.2390\\
 [-\halfbls] &&&&&&&&&& \\ \hline &&&&&&&&&& \\ [-\halfbls]
 g_1=19 & 47.453 & 25.251 & 17.638 & 13.987 & 11.926 & 10.647 & 9.8091 & 9.2425 & 8.8547 & 8.5913\\
 [-\halfbls] &&&&&&&&&& \\ \hline &&&&&&&&&& \\ [-\halfbls]
 g_1=20 & 52.042 & 27.481 & 19.062 & 15.020 & 12.731 & 11.304 & 10.362 & 9.7183 & 9.2702 & 8.9580\\
 [-\halfbls] &&&&&&&&&& \\ \hline &&&&&&&&&& \\ [-\halfbls]
 g_1=21 & 56.842 & 29.807 & 20.543 & 16.092 & 13.566 & 11.985 & 10.935 & 10.211 & 9.7007 & 9.3385\\
 [-\halfbls] &&&&&&&&&& \\ \hline &&&&&&&&&& \\ [-\halfbls]
 g_1=22 & 61.856 & 32.229 & 22.081 & 17.203 & 14.430 & 12.689 & 11.527 & 10.720 & 10.146 & 9.7325\\
 [-\halfbls] &&&&&&&&&& \\ \hline &&&&&&&&&& \\ [-\halfbls]
 g_1=23 & 67.081 & 34.746 & 23.677 & 18.354 & 15.323 & 13.416 & 12.138 & 11.246 & 10.606 & 10.140\\
 [-\halfbls] &&&&&&&&&& \\ \hline &&&&&&&&&& \\ [-\halfbls]
 g_1=24 & 72.518 & 37.359 & 25.329 & 19.544 & 16.245 & 14.166 & 12.768 & 11.788 & 11.079 & 10.559\\
 [-\halfbls] &&&&&&&&&& \\ \hline &&&&&&&&&& \\ [-\halfbls]
 g_1=25 & 78.168 & 40.067 & 27.039 & 20.772 & 17.197 & 14.939 & 13.417 & 12.345 & 11.567 & 10.992\\
 [-\halfbls] &&&&&&&&&& \\ \hline &&&&&&&&&& \\ [-\halfbls]
 g_1=26 & 84.030 & 42.871 & 28.806 & 22.040 & 18.177 & 15.735 & 14.085 & 12.919 & 12.069 & 11.436\\
 [-\halfbls] &&&&&&&&&& \\ \hline &&&&&&&&&& \\ [-\halfbls]
 g_1=27 & 90.104 & 45.770 & 30.630 & 23.347 & 19.186 & 16.553 & 14.771 & 13.508 & 12.584 & 11.893\\
 [-\halfbls] &&&&&&&&&& \\ \hline &&&&&&&&&& \\ [-\halfbls]
 g_1=28 & 96.391 & 48.766 & 32.510 & 24.692 & 20.224 & 17.394 & 15.475 & 14.113 & 13.113 & 12.361\\
 [-\halfbls] &&&&&&&&&& \\ \hline &&&&&&&&&& \\ [-\halfbls]
 g_1=29 & 102.89 & 51.856 & 34.448 & 26.077 & 21.291 & 18.257 & 16.198 & 14.733 & 13.655 & 12.842\\
 [-\halfbls] &&&&&&&&&& \\ \hline &&&&&&&&&& \\ [-\halfbls]
 g_1=30 & 109.60 & 55.042 & 36.443 & 27.500 & 22.387 & 19.143 & 16.939 & 15.369 & 14.210 & 13.334\\ [-\halfbls]
&&&&&&&&&&\\ \hline
\end{array}
$$ 
\end{table}

\newpage

\normalsize

\section{Saddle Connections Joining a Zero to Itself}
\label{a:to:itself}
\subsection{Connected Strata in Genus 4}

In the tables below we consider all connected strata in genus
$g=4$.
For each stratum we present the complete list of admissible
constructions (= all possible topological configurations
of homologous closed geodesics). To make calculations traceable
we present for every configuration the orders of possible
symmetries $\Gamma$ and $\Gamma_-$,
the combinatorial constant $M$, and the result: the constant
in quadratic asymptotics for the number of saddle connections
(closed geodesics) of this type.

\bigskip
\bigskip

\scriptsize
$$ 
\begin{array}{|c|c|c|c|c|c|}
\multicolumn{6}{c}{}\\
[-\halfbls]\multicolumn{6}{c}{\text{Stratum }\cH(5, 1).}\\
\multicolumn{6}{c}{\text{ Configurations of closed geodesics.}}\\
\multicolumn{6}{c}{}\\
\hline &&&&&\\
\text{Degeneration pattern} & |\Gamma_-| & |\Gamma| & M & c\cdot\zeta(2) & c\cdot\zeta(2) \text{ approx.} \\
[-\halfbls]  &&&&&\\ \hline  &&&&& \\ [-\halfbls]
\Rightarrow(\boldsymbol{\overline{0}},\boldsymbol{\overline{4}} )\Rightarrow & 1 & 1 & 5 & \cfrac{38125}{15552} & 2.45145\\
[-\halfbls]  &&&&&\\ \hline  &&&&& \\ [-\halfbls]
\to(\boldsymbol{\overline{0}},\boldsymbol{\overline{3}} ; 1)\to & 1 & 1 & 4 & \cfrac{2240}{243} & 9.21811\\
[-\halfbls]  &&&&&\\ \hline  &&&&& \\ [-\halfbls]
\Rightarrow(\boldsymbol{\overline{0+3}} ; 1)\Rightarrow & 1 & 1 & 4 & \cfrac{320}{243} & 1.31687\\
[-\halfbls]  &&&&&\\ \hline  &&&&& \\ [-\halfbls]
\Rightarrow(\boldsymbol{\overline{1+2}} ; 1)\Rightarrow & 1 & 1 & 4 & \cfrac{320}{243} & 1.31687\\
[-\halfbls]  &&&&&\\ \hline  &&&&& \\ [-\halfbls]
\to(\boldsymbol{\overline{1}},\boldsymbol{\overline{2}} ; 1)\to & 1 & 1 & 6 & \cfrac{175}{18} & 9.72222\\
[-\halfbls]  &&&&&\\ \hline  &&&&& \\ [-\halfbls]
\to(\boldsymbol{\overline{0+0}} )\Rightarrow(\boldsymbol{\overline{0}},\boldsymbol{\overline{2}} )\to & 1 & 1 & 3 & \cfrac{35}{288} & 0.121528\\
[-\halfbls]  &&&&&\\ \hline  &&&&& \\ [-\halfbls]
\to(\boldsymbol{\overline{0}},\boldsymbol{\overline{0}} )\Rightarrow(\boldsymbol{\overline{0+2}} )\to & 1 & 1 & 3 & \cfrac{35}{576} & 0.0607639\\
[-\halfbls]  &&&&&\\ \hline  &&&&& \\ [-\halfbls]
\to(\boldsymbol{\overline{0}},\boldsymbol{\overline{0}} )\Rightarrow(\boldsymbol{\overline{1+1}} )\to & 1 & 1 & 3 & \cfrac{35}{576} & 0.0607639\\
[-\halfbls]  &&&&&\\ \hline  &&&&& \\ [-\halfbls]
\to(\boldsymbol{\overline{0}},\boldsymbol{\overline{0}} )\Rightarrow(\boldsymbol{\overline{2+0}} )\to & 1 & 1 & 3 & \cfrac{35}{576} & 0.0607639\\
[-\halfbls]  &&&&&\\ \hline  &&&&& \\ [-\halfbls]
\to(\boldsymbol{\overline{0+0}} )\to(\boldsymbol{\overline{0}},\boldsymbol{\overline{1}} ; 1)\to & 1 & 1 & 2 & \cfrac{175}{486} & 0.360082\\
[-\halfbls]  &&&&&\\ \hline  &&&&& \\ [-\halfbls]
\to(\boldsymbol{\overline{0}},\boldsymbol{\overline{0}} )\to(\boldsymbol{\overline{0+1}} ; 1)\to & 1 & 1 & 2 & \cfrac{35}{243} & 0.144033\\
[-\halfbls]  &&&&&\\ \hline  &&&&& \\ [-\halfbls]
\to(\boldsymbol{\overline{0+0}} )\Rightarrow(\boldsymbol{\overline{0+1}} ; 1)\to & 1 & 1 & 2 & \cfrac{35}{486} & 0.0720165\\
[-\halfbls]  &&&&&\\ \hline  &&&&& \\ [-\halfbls]
\to(\boldsymbol{\overline{0+0}} )\Rightarrow(\boldsymbol{\overline{1+0}} ; 1)\to & 1 & 1 & 2 & \cfrac{35}{486} & 0.0720165\\
[-\halfbls]  &&&&&\\ \hline  &&&&& \\ [-\halfbls]
\to(\boldsymbol{\overline{0+0}} )\to(\boldsymbol{\overline{0+0}} )\Rightarrow(\boldsymbol{\overline{0}},\boldsymbol{\overline{0}} )\to & 1 & 1 & 1 & \cfrac{175}{7776} & 0.0225051\\
[-\halfbls]  &&&&&\\ \hline
\end{array}
$$ 

\newpage 

$$ \phantom{space} $$

$$ 
\begin{array}{|c|c|c|c|c|c|}
\multicolumn{6}{c}{}\\
[-\halfbls]\multicolumn{6}{c}{\text{Stratum }\cH(4, 1, 1).}\\
\multicolumn{6}{c}{\text{ Configurations of closed geodesics.}}\\
\multicolumn{6}{c}{}\\
\hline &&&&&\\
\text{Degeneration pattern} & |\Gamma_-| & |\Gamma| & M & c\cdot\zeta(2) & c\cdot\zeta(2) \text{ approx.} \\
[-\halfbls]  &&&&&\\ \hline  &&&&& \\ [-\halfbls]
\Rightarrow(\boldsymbol{\overline{0}},\boldsymbol{\overline{0}} ; 4)\Rightarrow & 2 & 1 & 1 & \cfrac{61}{88} & 0.693182\\
[-\halfbls]  &&&&&\\ \hline  &&&&& \\ [-\halfbls]
\Rightarrow(\boldsymbol{\overline{0}},\boldsymbol{\overline{3}} ; 1)\Rightarrow & 1 & 1 & 8 & \cfrac{1024}{275} & 3.72364\\
[-\halfbls]  &&&&&\\ \hline  &&&&& \\ [-\halfbls]
\to(\boldsymbol{\overline{0}},\boldsymbol{\overline{2}} ; 1, 1)\to & 1 & 1 & 3 & \cfrac{432}{55} & 7.85455\\
[-\halfbls]  &&&&&\\ \hline  &&&&& \\ [-\halfbls]
\Rightarrow(\boldsymbol{\overline{0+2}} ; 1, 1)\Rightarrow & 1 & 1 & 3 & \cfrac{54}{55} & 0.981818\\
[-\halfbls]  &&&&&\\ \hline  &&&&& \\ [-\halfbls]
\Rightarrow(\boldsymbol{\overline{1+1}} ; 1, 1)\Rightarrow & 2 & 1 & \cfrac{3}{2} & \cfrac{27}{55} & 0.490909\\
[-\halfbls]  &&&&&\\ \hline  &&&&& \\ [-\halfbls]
\to(\boldsymbol{\overline{1}},\boldsymbol{\overline{1}} ; 1, 1)\to & 2 & 1 & 2 & \cfrac{224}{55} & 4.07273\\
[-\halfbls]  &&&&&\\ \hline  &&&&& \\ [-\halfbls]
\to(\boldsymbol{\overline{0}},\boldsymbol{\overline{0}} )\Rightarrow(\boldsymbol{\overline{0}},\boldsymbol{\overline{2}} )\to & 1 & 1 & 6 & \cfrac{27}{275} & 0.0981818\\
[-\halfbls]  &&&&&\\ \hline  &&&&& \\ [-\halfbls]
\Rightarrow(\boldsymbol{\overline{0}},\boldsymbol{\overline{0}} )\Rightarrow(\boldsymbol{\overline{0+2}} )\Rightarrow & 1 & 1 & 6 & \cfrac{27}{1100} & 0.0245455\\
[-\halfbls]  &&&&&\\ \hline  &&&&& \\ [-\halfbls]
\Rightarrow(\boldsymbol{\overline{0}},\boldsymbol{\overline{0}} )\Rightarrow(\boldsymbol{\overline{1+1}} )\Rightarrow & 2 & 1 & 3 & \cfrac{27}{2200} & 0.0122727\\
[-\halfbls]  &&&&&\\ \hline  &&&&& \\ [-\halfbls]
\to(\boldsymbol{\overline{0+0}} )\to(\boldsymbol{\overline{0}},\boldsymbol{\overline{0}} ; 1, 1)\to & 2 & 1 & \cfrac{1}{2} & \cfrac{6}{55} & 0.109091\\
[-\halfbls]  &&&&&\\ \hline  &&&&& \\ [-\halfbls]
\to(\boldsymbol{\overline{0}},\boldsymbol{\overline{0}} )\to(\boldsymbol{\overline{0+0}} ; 1, 1)\to & 2 & 1 & \cfrac{1}{2} & \cfrac{2}{55} & 0.0363636\\
[-\halfbls]  &&&&&\\ \hline  &&&&& \\ [-\halfbls]
\to(\boldsymbol{\overline{0+0}} )\Rightarrow(\boldsymbol{\overline{0+0}} ; 1, 1)\to & 1 & 1 & 1 & \cfrac{2}{55} & 0.0363636\\
[-\halfbls]  &&&&&\\ \hline  &&&&& \\ [-\halfbls]
\to(\boldsymbol{\overline{0+0}} )\Rightarrow(\boldsymbol{\overline{0}},\boldsymbol{\overline{1}} ; 1)\to & 1 & 1 & 4 & \cfrac{8}{55} & 0.145455\\
[-\halfbls]  &&&&&\\ \hline  &&&&& \\ [-\halfbls]
\to(\boldsymbol{\overline{0}},\boldsymbol{\overline{0}} )\Rightarrow(\boldsymbol{\overline{0+1}} ; 1)\to & 1 & 1 & 4 & \cfrac{16}{275} & 0.0581818\\
[-\halfbls]  &&&&&\\ \hline  &&&&& \\ [-\halfbls]
\to(\boldsymbol{\overline{0}},\boldsymbol{\overline{0}} )\Rightarrow(\boldsymbol{\overline{1+0}} ; 1)\to & 1 & 1 & 4 & \cfrac{16}{275} & 0.0581818\\
[-\halfbls]  &&&&&\\ \hline  &&&&& \\ [-\halfbls]
\to(\boldsymbol{\overline{0+0}} )\to(\boldsymbol{\overline{0}},\boldsymbol{\overline{0}} )\Rightarrow(\boldsymbol{\overline{0}},\boldsymbol{\overline{0}} )\to & 2 & 1 & 1 & \cfrac{1}{110} & 0.00909091\\
[-\halfbls]  &&&&&\\ \hline  &&&&& \\ [-\halfbls]
\to(\boldsymbol{\overline{0+0}} )\Rightarrow(\boldsymbol{\overline{0}},\boldsymbol{\overline{0}} )\Rightarrow(\boldsymbol{\overline{0+0}} )\to & 2 & 1 & 1 & \cfrac{1}{220} & 0.00454545\\
[-\halfbls]  &&&&&\\ \hline
\end{array}
$$ 

\newpage 

$$ \phantom{space} $$

$$ 
\begin{array}{|c|c|c|c|c|c|}
\multicolumn{6}{c}{}\\
[-\halfbls]\multicolumn{6}{c}{\text{Stratum }\cH(3, 2, 1).}\\
\multicolumn{6}{c}{\text{ Configurations of closed geodesics.}}\\
\multicolumn{6}{c}{}\\
\hline &&&&&\\
\text{Degeneration pattern} & |\Gamma_-| & |\Gamma| & M & c\cdot\zeta(2) & c\cdot\zeta(2) \text{ approx.} \\
[-\halfbls]  &&&&&\\ \hline  &&&&& \\ [-\halfbls]
\Rightarrow(\boldsymbol{\overline{0}},\boldsymbol{\overline{2}} ; 2)\Rightarrow & 1 & 1 & 3 & \cfrac{765}{512} & 1.49414\\
[-\halfbls]  &&&&&\\ \hline  &&&&& \\ [-\halfbls]
\to(\boldsymbol{\overline{0}},\boldsymbol{\overline{0}} ; 3, 1)\to & 2 & 1 & \cfrac{1}{2} & \cfrac{20}{9} & 2.22222\\
[-\halfbls]  &&&&&\\ \hline  &&&&& \\ [-\halfbls]
\Rightarrow(\boldsymbol{\overline{0+0}} ; 3, 1)\Rightarrow & 2 & 1 & \cfrac{1}{2} & \cfrac{5}{18} & 0.277778\\
[-\halfbls]  &&&&&\\ \hline  &&&&& \\ [-\halfbls]
\Rightarrow(\boldsymbol{\overline{0}},\boldsymbol{\overline{1}} ; 3)\Rightarrow & 1 & 1 & 2 & \cfrac{10}{9} & 1.11111\\
[-\halfbls]  &&&&&\\ \hline  &&&&& \\ [-\halfbls]
\to(\boldsymbol{\overline{0}},\boldsymbol{\overline{1}} ; 2, 1)\to & 1 & 1 & 2 & \cfrac{25}{4} & 6.25\\
[-\halfbls]  &&&&&\\ \hline  &&&&& \\ [-\halfbls]
\Rightarrow(\boldsymbol{\overline{0+1}} ; 2, 1)\Rightarrow & 1 & 1 & 2 & \cfrac{25}{32} & 0.78125\\
[-\halfbls]  &&&&&\\ \hline  &&&&& \\ [-\halfbls]
\Rightarrow(\boldsymbol{\overline{1}},\boldsymbol{\overline{2}} ; 1)\Rightarrow & 1 & 1 & 6 & \cfrac{75}{32} & 2.34375\\
[-\halfbls]  &&&&&\\ \hline  &&&&& \\ [-\halfbls]
\to(\boldsymbol{\overline{0+0}} )\Rightarrow(\boldsymbol{\overline{0}},\boldsymbol{\overline{0}} ; 2)\to & 1 & 1 & 1 & \cfrac{25}{512} & 0.0488281\\
[-\halfbls]  &&&&&\\ \hline  &&&&& \\ [-\halfbls]
\to(\boldsymbol{\overline{0}},\boldsymbol{\overline{0}} )\Rightarrow(\boldsymbol{\overline{0+0}} ; 2)\to & 1 & 1 & 1 & \cfrac{5}{256} & 0.0195313\\
[-\halfbls]  &&&&&\\ \hline  &&&&& \\ [-\halfbls]
\to(\boldsymbol{\overline{0}},\boldsymbol{\overline{0}} )\Rightarrow(\boldsymbol{\overline{2}},\boldsymbol{\overline{0}} )\to & 1 & 1 & 3 & \cfrac{15}{256} & 0.0585938\\
[-\halfbls]  &&&&&\\ \hline  &&&&& \\ [-\halfbls]
\Rightarrow(\boldsymbol{\overline{0+0}} )\Rightarrow(\boldsymbol{\overline{0}},\boldsymbol{\overline{2}} )\Rightarrow & 1 & 1 & 3 & \cfrac{15}{512} & 0.0292969\\
[-\halfbls]  &&&&&\\ \hline  &&&&& \\ [-\halfbls]
\to(\boldsymbol{\overline{0}},\boldsymbol{\overline{0}} )\to(\boldsymbol{\overline{0}},\boldsymbol{\overline{1}} ; 1)\to & 1 & 1 & 2 & \cfrac{25}{144} & 0.173611\\
[-\halfbls]  &&&&&\\ \hline  &&&&& \\ [-\halfbls]
\to(\boldsymbol{\overline{0+0}} )\Rightarrow(\boldsymbol{\overline{1}},\boldsymbol{\overline{0}} ; 1)\to & 1 & 1 & 2 & \cfrac{25}{288} & 0.0868056\\
[-\halfbls]  &&&&&\\ \hline  &&&&& \\ [-\halfbls]
\to(\boldsymbol{\overline{0}},\boldsymbol{\overline{0}} )\Rightarrow(\boldsymbol{\overline{1}},\boldsymbol{\overline{1}} )\to & 1 & 1 & 4 & \cfrac{5}{72} & 0.0694444\\
[-\halfbls]  &&&&&\\ \hline  &&&&& \\ [-\halfbls]
\Rightarrow(\boldsymbol{\overline{0+0}} )\Rightarrow(\boldsymbol{\overline{0+1}} ; 1)\Rightarrow & 1 & 1 & 2 & \cfrac{5}{288} & 0.0173611\\
[-\halfbls]  &&&&&\\ \hline  &&&&& \\ [-\halfbls]
\to(\boldsymbol{\overline{0}},\boldsymbol{\overline{0}} )\to(\boldsymbol{\overline{0+0}} )\Rightarrow(\boldsymbol{\overline{0}},\boldsymbol{\overline{0}} )\to & 1 & 1 & 1 & \cfrac{25}{2304} & 0.0108507\\
[-\halfbls]  &&&&&\\ \hline  &&&&& \\ [-\halfbls]
\to(\boldsymbol{\overline{0+0}} )\Rightarrow(\boldsymbol{\overline{0+0}} )\Rightarrow(\boldsymbol{\overline{0}},\boldsymbol{\overline{0}} )\to & 1 & 1 & 1 & \cfrac{25}{4608} & 0.00542535\\
[-\halfbls]  &&&&&\\ \hline
\end{array}
$$ 

\newpage 

$$ \phantom{space} $$

$$ 
\begin{array}{|c|c|c|c|c|c|}
\multicolumn{6}{c}{}\\
[-\halfbls]\multicolumn{6}{c}{\text{Stratum }\cH(3, 1, 1, 1).}\\
\multicolumn{6}{c}{\text{ Configurations of closed geodesics.}}\\
\multicolumn{6}{c}{}\\
\hline &&&&&\\
\text{Degeneration pattern} & |\Gamma_-| & |\Gamma| & M & c\cdot\zeta(2) & c\cdot\zeta(2) \text{ approx.} \\
[-\halfbls]  &&&&&\\ \hline  &&&&& \\ [-\halfbls]
\Rightarrow(\boldsymbol{\overline{0}},\boldsymbol{\overline{0}} ; 3, 1)\Rightarrow & 2 & 1 & 3 & \cfrac{60}{31} & 1.93548\\
[-\halfbls]  &&&&&\\ \hline  &&&&& \\ [-\halfbls]
\Rightarrow(\boldsymbol{\overline{0}},\boldsymbol{\overline{2}} ; 1, 1)\Rightarrow & 1 & 1 & 9 & \cfrac{2025}{496} & 4.08266\\
[-\halfbls]  &&&&&\\ \hline  &&&&& \\ [-\halfbls]
\to(\boldsymbol{\overline{0}},\boldsymbol{\overline{1}} ; 1, 1, 1)\to & 1 & 1 & 2 & \cfrac{1575}{248} & 6.35081\\
[-\halfbls]  &&&&&\\ \hline  &&&&& \\ [-\halfbls]
\Rightarrow(\boldsymbol{\overline{0+1}} ; 1, 1, 1)\Rightarrow & 1 & 1 & 2 & \cfrac{175}{248} & 0.705645\\
[-\halfbls]  &&&&&\\ \hline  &&&&& \\ [-\halfbls]
\Rightarrow(\boldsymbol{\overline{0}},\boldsymbol{\overline{0}} )\Rightarrow(\boldsymbol{\overline{0}},\boldsymbol{\overline{2}} )\Rightarrow & 1 & 1 & 18 & \cfrac{405}{7936} & 0.0510333\\
[-\halfbls]  &&&&&\\ \hline  &&&&& \\ [-\halfbls]
\to(\boldsymbol{\overline{0+0}} )\Rightarrow(\boldsymbol{\overline{0}},\boldsymbol{\overline{0}} ; 1, 1)\to & 1 & 1 & 3 & \cfrac{225}{1984} & 0.113407\\
[-\halfbls]  &&&&&\\ \hline  &&&&& \\ [-\halfbls]
\to(\boldsymbol{\overline{0}},\boldsymbol{\overline{0}} )\Rightarrow(\boldsymbol{\overline{0+0}} ; 1, 1)\to & 1 & 1 & 3 & \cfrac{75}{1984} & 0.0378024\\
[-\halfbls]  &&&&&\\ \hline  &&&&& \\ [-\halfbls]
\to(\boldsymbol{\overline{0}},\boldsymbol{\overline{0}} )\Rightarrow(\boldsymbol{\overline{0}},\boldsymbol{\overline{1}} ; 1)\to & 1 & 1 & 12 & \cfrac{75}{496} & 0.15121\\
[-\halfbls]  &&&&&\\ \hline  &&&&& \\ [-\halfbls]
\Rightarrow(\boldsymbol{\overline{0}},\boldsymbol{\overline{0}} )\Rightarrow(\boldsymbol{\overline{0+1}} ; 1)\Rightarrow & 1 & 1 & 12 & \cfrac{15}{496} & 0.0302419\\
[-\halfbls]  &&&&&\\ \hline  &&&&& \\ [-\halfbls]
\to(\boldsymbol{\overline{0+0}} )\Rightarrow(\boldsymbol{\overline{0}},\boldsymbol{\overline{0}} )\Rightarrow(\boldsymbol{\overline{0}},\boldsymbol{\overline{0}} )\to & 1 & 1 & 6 & \cfrac{75}{7936} & 0.0094506\\
[-\halfbls]  &&&&&\\ \hline
\end{array}
$$ 

\newpage 

$$ \phantom{space} $$

$$ 
\begin{array}{|c|c|c|c|c|c|}
\multicolumn{6}{c}{}\\
[-\halfbls]\multicolumn{6}{c}{\text{Stratum }\cH(2, 2, 1, 1).}\\
\multicolumn{6}{c}{\text{ Configurations of closed geodesics.}}\\
\multicolumn{6}{c}{}\\
\hline &&&&&\\
\text{Degeneration pattern} & |\Gamma_-| & |\Gamma| & M & c\cdot\zeta(2) & c\cdot\zeta(2) \text{ approx.} \\
[-\halfbls]  &&&&&\\ \hline  &&&&& \\ [-\halfbls]
\Rightarrow(\boldsymbol{\overline{0}},\boldsymbol{\overline{0}} ; 2, 2)\Rightarrow & 2 & 1 & 1 & \cfrac{85}{131} & 0.648855\\
[-\halfbls]  &&&&&\\ \hline  &&&&& \\ [-\halfbls]
\to(\boldsymbol{\overline{0}},\boldsymbol{\overline{0}} ; 2, 1, 1)\to & 2 & 1 & 1 & \cfrac{600}{131} & 4.58015\\
[-\halfbls]  &&&&&\\ \hline  &&&&& \\ [-\halfbls]
\Rightarrow(\boldsymbol{\overline{0+0}} ; 2, 1, 1)\Rightarrow & 2 & 1 & 1 & \cfrac{200}{393} & 0.508906\\
[-\halfbls]  &&&&&\\ \hline  &&&&& \\ [-\halfbls]
\Rightarrow(\boldsymbol{\overline{0}},\boldsymbol{\overline{1}} ; 2, 1)\Rightarrow & 1 & 1 & 8 & \cfrac{1600}{393} & 4.07125\\
[-\halfbls]  &&&&&\\ \hline  &&&&& \\ [-\halfbls]
\Rightarrow(\boldsymbol{\overline{1}},\boldsymbol{\overline{1}} ; 1, 1)\Rightarrow & 2 & 1 & 4 & \cfrac{5600}{3537} & 1.58326\\
[-\halfbls]  &&&&&\\ \hline  &&&&& \\ [-\halfbls]
\to(\boldsymbol{\overline{0}},\boldsymbol{\overline{0}} )\Rightarrow(\boldsymbol{\overline{0}},\boldsymbol{\overline{0}} ; 2)\to & 1 & 1 & 4 & \cfrac{25}{393} & 0.0636132\\
[-\halfbls]  &&&&&\\ \hline  &&&&& \\ [-\halfbls]
\Rightarrow(\boldsymbol{\overline{0+0}} )\Rightarrow(\boldsymbol{\overline{0}},\boldsymbol{\overline{0}} ; 2)\Rightarrow & 2 & 1 & 2 & \cfrac{25}{1572} & 0.0159033\\
[-\halfbls]  &&&&&\\ \hline  &&&&& \\ [-\halfbls]
\Rightarrow(\boldsymbol{\overline{0}},\boldsymbol{\overline{0}} )\Rightarrow(\boldsymbol{\overline{0+0}} ; 2)\Rightarrow & 2 & 1 & 2 & \cfrac{5}{786} & 0.00636132\\
[-\halfbls]  &&&&&\\ \hline  &&&&& \\ [-\halfbls]
\to(\boldsymbol{\overline{0}},\boldsymbol{\overline{0}} )\to(\boldsymbol{\overline{0}},\boldsymbol{\overline{0}} ; 1, 1)\to & 2 & 1 & 1 & \cfrac{100}{1179} & 0.0848176\\
[-\halfbls]  &&&&&\\ \hline  &&&&& \\ [-\halfbls]
\to(\boldsymbol{\overline{0}},\boldsymbol{\overline{0}} )\Rightarrow(\boldsymbol{\overline{1}},\boldsymbol{\overline{0}} ; 1)\to & 1 & 1 & 8 & \cfrac{400}{3537} & 0.11309\\
[-\halfbls]  &&&&&\\ \hline  &&&&& \\ [-\halfbls]
\Rightarrow(\boldsymbol{\overline{0+0}} )\Rightarrow(\boldsymbol{\overline{0+0}} ; 1, 1)\Rightarrow & 2 & 1 & 1 & \cfrac{25}{3537} & 0.00706814\\
[-\halfbls]  &&&&&\\ \hline  &&&&& \\ [-\halfbls]
\Rightarrow(\boldsymbol{\overline{0+0}} )\Rightarrow(\boldsymbol{\overline{0}},\boldsymbol{\overline{1}} ; 1)\Rightarrow & 1 & 1 & 8 & \cfrac{200}{3537} & 0.0565451\\
[-\halfbls]  &&&&&\\ \hline  &&&&& \\ [-\halfbls]
\Rightarrow(\boldsymbol{\overline{0}},\boldsymbol{\overline{0}} )\Rightarrow(\boldsymbol{\overline{1}},\boldsymbol{\overline{1}} )\Rightarrow & 2 & 1 & 8 & \cfrac{80}{3537} & 0.022618\\
[-\halfbls]  &&&&&\\ \hline  &&&&& \\ [-\halfbls]
\to(\boldsymbol{\overline{0}},\boldsymbol{\overline{0}} )\to(\boldsymbol{\overline{0}},\boldsymbol{\overline{0}} )\Rightarrow(\boldsymbol{\overline{0}},\boldsymbol{\overline{0}} )\to & 2 & 1 & 2 & \cfrac{25}{3537} & 0.00706814\\
[-\halfbls]  &&&&&\\ \hline  &&&&& \\ [-\halfbls]
\to(\boldsymbol{\overline{0}},\boldsymbol{\overline{0}} )\Rightarrow(\boldsymbol{\overline{0+0}} )\Rightarrow(\boldsymbol{\overline{0}},\boldsymbol{\overline{0}} )\to & 2 & 1 & 2 & \cfrac{25}{7074} & 0.00353407\\
[-\halfbls]  &&&&&\\ \hline  &&&&& \\ [-\halfbls]
\Rightarrow(\boldsymbol{\overline{0+0}} )\Rightarrow(\boldsymbol{\overline{0+0}} )\Rightarrow(\boldsymbol{\overline{0}},\boldsymbol{\overline{0}} )\Rightarrow & 2 & 1 & 2 & \cfrac{25}{14148} & 0.00176703\\
[-\halfbls]  &&&&&\\ \hline
\end{array}
$$ 

\newpage 

$$ \phantom{space} $$

$$ 
\begin{array}{|c|c|c|c|c|c|}
\multicolumn{6}{c}{}\\
[-\halfbls]\multicolumn{6}{c}{\text{Stratum }\cH(2, 1, 1, 1, 1).}\\
\multicolumn{6}{c}{\text{ Configurations of closed geodesics.}}\\
\multicolumn{6}{c}{}\\
\hline &&&&&\\
\text{Degeneration pattern} & |\Gamma_-| & |\Gamma| & M & c\cdot\zeta(2) & c\cdot\zeta(2) \text{ approx.} \\
[-\halfbls]  &&&&&\\ \hline  &&&&& \\ [-\halfbls]
\Rightarrow(\boldsymbol{\overline{0}},\boldsymbol{\overline{0}} ; 2, 1, 1)\Rightarrow & 2 & 1 & 6 & \cfrac{18}{5} & 3.6\\
[-\halfbls]  &&&&&\\ \hline  &&&&& \\ [-\halfbls]
\to(\boldsymbol{\overline{0}},\boldsymbol{\overline{0}} ; 1, 1, 1, 1)\to & 2 & 1 & \cfrac{1}{2} & \cfrac{7}{3} & 2.33333\\
[-\halfbls]  &&&&&\\ \hline  &&&&& \\ [-\halfbls]
\Rightarrow(\boldsymbol{\overline{0+0}} ; 1, 1, 1, 1)\Rightarrow & 2 & 1 & \cfrac{1}{2} & \cfrac{7}{30} & 0.233333\\
[-\halfbls]  &&&&&\\ \hline  &&&&& \\ [-\halfbls]
\Rightarrow(\boldsymbol{\overline{0}},\boldsymbol{\overline{1}} ; 1, 1, 1)\Rightarrow & 1 & 1 & 8 & \cfrac{56}{15} & 3.73333\\
[-\halfbls]  &&&&&\\ \hline  &&&&& \\ [-\halfbls]
\Rightarrow(\boldsymbol{\overline{0}},\boldsymbol{\overline{0}} )\Rightarrow(\boldsymbol{\overline{0}},\boldsymbol{\overline{0}} ; 2)\Rightarrow & 2 & 1 & 12 & \cfrac{1}{40} & 0.025\\
[-\halfbls]  &&&&&\\ \hline  &&&&& \\ [-\halfbls]
\to(\boldsymbol{\overline{0}},\boldsymbol{\overline{0}} )\Rightarrow(\boldsymbol{\overline{0}},\boldsymbol{\overline{0}} ; 1, 1)\to & 1 & 1 & 12 & \cfrac{2}{15} & 0.133333\\
[-\halfbls]  &&&&&\\ \hline  &&&&& \\ [-\halfbls]
\Rightarrow(\boldsymbol{\overline{0+0}} )\Rightarrow(\boldsymbol{\overline{0}},\boldsymbol{\overline{0}} ; 1, 1)\Rightarrow & 2 & 1 & 6 & \cfrac{1}{30} & 0.0333333\\
[-\halfbls]  &&&&&\\ \hline  &&&&& \\ [-\halfbls]
\Rightarrow(\boldsymbol{\overline{0}},\boldsymbol{\overline{0}} )\Rightarrow(\boldsymbol{\overline{0+0}} ; 1, 1)\Rightarrow & 2 & 1 & 6 & \cfrac{1}{90} & 0.0111111\\
[-\halfbls]  &&&&&\\ \hline  &&&&& \\ [-\halfbls]
\Rightarrow(\boldsymbol{\overline{0}},\boldsymbol{\overline{0}} )\Rightarrow(\boldsymbol{\overline{0}},\boldsymbol{\overline{1}} ; 1)\Rightarrow & 1 & 1 & 48 & \cfrac{4}{45} & 0.0888889\\
[-\halfbls]  &&&&&\\ \hline  &&&&& \\ [-\halfbls]
\to(\boldsymbol{\overline{0}},\boldsymbol{\overline{0}} )\Rightarrow(\boldsymbol{\overline{0}},\boldsymbol{\overline{0}} )\Rightarrow(\boldsymbol{\overline{0}},\boldsymbol{\overline{0}} )\to & 2 & 1 & 12 & \cfrac{1}{180} & 0.00555556\\
[-\halfbls]  &&&&&\\ \hline  &&&&& \\ [-\halfbls]
\Rightarrow(\boldsymbol{\overline{0+0}} )\Rightarrow(\boldsymbol{\overline{0}},\boldsymbol{\overline{0}} )\Rightarrow(\boldsymbol{\overline{0}},\boldsymbol{\overline{0}} )\Rightarrow & 2 & 1 & 12 & \cfrac{1}{360} & 0.00277778\\
[-\halfbls]  &&&&&\\ \hline
\end{array}
$$ 

\vspace*{1truecm}

$$ 
\begin{array}{|c|c|c|c|c|c|}
\multicolumn{6}{c}{}\\
[-\halfbls]\multicolumn{6}{c}{\text{Stratum }\cH(1, 1, 1, 1, 1, 1).}\\
\multicolumn{6}{c}{\text{ Configurations of closed geodesics.}}\\
\multicolumn{6}{c}{}\\
\hline &&&&&\\
\text{Degeneration pattern} & |\Gamma_-| & |\Gamma| & M & c\cdot\zeta(2) & c\cdot\zeta(2) \text{ approx.} \\
[-\halfbls]  &&&&&\\ \hline  &&&&& \\ [-\halfbls]
\Rightarrow(\boldsymbol{\overline{0}},\boldsymbol{\overline{0}} ; 1, 1, 1, 1)\Rightarrow & 2 & 1 & 15 & \cfrac{3150}{377} & 8.35544\\
[-\halfbls]  &&&&&\\ \hline  &&&&& \\ [-\halfbls]
\Rightarrow(\boldsymbol{\overline{0}},\boldsymbol{\overline{0}} )\Rightarrow(\boldsymbol{\overline{0}},\boldsymbol{\overline{0}} ; 1, 1)\Rightarrow & 2 & 1 & 180 & \cfrac{90}{377} & 0.238727\\
[-\halfbls]  &&&&&\\ \hline  &&&&& \\ [-\halfbls]
\Rightarrow(\boldsymbol{\overline{0}},\boldsymbol{\overline{0}} )\Rightarrow(\boldsymbol{\overline{0}},\boldsymbol{\overline{0}} )\Rightarrow(\boldsymbol{\overline{0}},\boldsymbol{\overline{0}} )\Rightarrow & 2 & 3 & 120 & \cfrac{5}{754} & 0.0066313\\
[-\halfbls]  &&&&&\\ \hline
\end{array}
$$ 

\newpage


\normalsize

\addcontentsline{toc}{subsection}
{B.2.\hspace*{2.6truemm} Closed Geodesics for Nonconnected Strata in Genus 4}
{B.2. {\bf Closed Geodesics for Nonconnected Strata in Genus 4}}
%
%
%
\scriptsize
$$ 
\begin{array}{|c|c|c|c|c|}
\multicolumn{5}{c}{}\\
[-\halfbls]
\multicolumn{5}{c}{\text{Component }\cH^{hyp}(6).}\\
\multicolumn{5}{c}{}\\
\hline &&&&\\
\text{Degeneration pattern} & |\Gamma_-| & |\Gamma| & c\cdot\zeta(2) & c\cdot\zeta(2) \text{ approx.} \\
[-\halfbls]  &&&&\\ \hline  &&&& \\ [-\halfbls]
\to(\boldsymbol{\overline{2}},\boldsymbol{\overline{2}} )\to & 2 & 1 & \cfrac{384}{25} & 15.36\\
[-\halfbls]  &&&&\\ \hline  &&&& \\ [-\halfbls]
\Rightarrow(\boldsymbol{\overline{2+2}} )\Rightarrow & 2 & 1 & \cfrac{6}{1} & 6.\\
[-\halfbls]  &&&&\\ \hline  &&&& \\ [-\halfbls]
\to(\boldsymbol{\overline{0}},\boldsymbol{\overline{0}} )\to(\boldsymbol{\overline{1+1}} )\to & 2 & 1 & \cfrac{84}{25} & 3.36\\
[-\halfbls]  &&&&\\ \hline  &&&& \\ [-\halfbls]
\to(\boldsymbol{\overline{0+0}} )\to(\boldsymbol{\overline{1}},\boldsymbol{\overline{1}} )\to & 2 & 1 & \cfrac{896}{225} & 3.98222\\
[-\halfbls]  &&&&\\ \hline
%
%
\multicolumn{5}{c}{}\\
\multicolumn{5}{c}{}\\
[-\halfbls]\multicolumn{5}{c}{\text{Component }\cH^{even}(6).}\\
\multicolumn{5}{c}{}\\
\hline &&&&\\
[-\halfbls]
\to(\boldsymbol{\overline{0}},\boldsymbol{\overline{4}} )\to & 1 & 1 & \cfrac{30375}{4096} & 7.41577\\
[-\halfbls]  &&&&\\ \hline  &&&& \\ [-\halfbls]
\Rightarrow(\boldsymbol{\overline{0+4}} )\Rightarrow & 1 & 1 & \cfrac{10125}{8192} & 1.23596\\
[-\halfbls]  &&&&\\ \hline  &&&& \\ [-\halfbls]
\Rightarrow(\boldsymbol{\overline{1+3}} )\Rightarrow & 1 & 1 & \cfrac{875}{256} & 3.41797\\
[-\halfbls]  &&&&\\ \hline  &&&& \\ [-\halfbls]
\to(\boldsymbol{\overline{2}},\boldsymbol{\overline{2}} )\to & 2 & 1 & \cfrac{405}{128} & 3.16406\\
[-\halfbls]  &&&&\\ \hline  &&&& \\ [-\halfbls]
\to(\boldsymbol{\overline{1}},\boldsymbol{\overline{3}} )\to & 1 & 1 & \cfrac{15}{1} & 15.\\
[-\halfbls]  &&&&\\ \hline  &&&& \\ [-\halfbls]
\to(\boldsymbol{\overline{0+0}} )\to(\boldsymbol{\overline{0}},\boldsymbol{\overline{2}} )\to & 1 & 1 & \cfrac{2835}{2048} & 1.38428\\
[-\halfbls]  &&&&\\ \hline  &&&& \\ [-\halfbls]
\to(\boldsymbol{\overline{0}},\boldsymbol{\overline{0}} )\to(\boldsymbol{\overline{0+2}} )\to & 1 & 1 & \cfrac{2835}{4096} & 0.692139\\
[-\halfbls]  &&&&\\ \hline  &&&& \\ [-\halfbls]
\to(\boldsymbol{\overline{0+0}} )\Rightarrow(\boldsymbol{\overline{1+1}} )\to & 1 & 1 & \cfrac{2835}{8192} & 0.346069\\
[-\halfbls]  &&&&\\ \hline
%
%
\multicolumn{5}{c}{}\\
\multicolumn{5}{c}{}\\
[-\halfbls]\multicolumn{5}{c}{\text{Component }\cH^{odd}(6).}\\
\multicolumn{5}{c}{}\\
\hline &&&&\\
[-\halfbls]
\to(\boldsymbol{\overline{0}},\boldsymbol{\overline{4}} )\to & 1 & 1 & \cfrac{350}{27} & 12.963\\
[-\halfbls]  &&&&\\ \hline  &&&& \\ [-\halfbls]
\Rightarrow(\boldsymbol{\overline{0+4}} )\Rightarrow & 1 & 1 & \cfrac{175}{81} & 2.16049\\
[-\halfbls]  &&&&\\ \hline  &&&& \\ [-\halfbls]
\Rightarrow(\boldsymbol{\overline{1+3}} )\Rightarrow & 1 & 1 & \cfrac{25}{32} & 0.78125\\
[-\halfbls]  &&&&\\ \hline  &&&& \\ [-\halfbls]
\Rightarrow(\boldsymbol{\overline{2+2}} )\Rightarrow & 2 & 1 & \cfrac{175}{162} & 1.08025\\
[-\halfbls]  &&&&\\ \hline  &&&& \\ [-\halfbls]
\to(\boldsymbol{\overline{2}},\boldsymbol{\overline{2}} )\to & 2 & 1 & \cfrac{105}{16} & 6.5625\\
[-\halfbls]  &&&&\\ \hline  &&&& \\ [-\halfbls]
\to(\boldsymbol{\overline{1}},\boldsymbol{\overline{3}} )\to & 1 & 1 & \cfrac{256}{27} & 9.48148\\
[-\halfbls]  &&&&\\ \hline  &&&& \\ [-\halfbls]
\to(\boldsymbol{\overline{0+0}} )\Rightarrow(\boldsymbol{\overline{0+2}} )\to & 1 & 1 & \cfrac{7}{32} & 0.21875\\
[-\halfbls]  &&&&\\ \hline  &&&& \\ [-\halfbls]
\to(\boldsymbol{\overline{0+0}} )\Rightarrow(\boldsymbol{\overline{2+0}} )\to & 1 & 1 & \cfrac{7}{32} & 0.21875\\
[-\halfbls]  &&&&\\ \hline  &&&& \\ [-\halfbls]
\to(\boldsymbol{\overline{0+0}} )\to(\boldsymbol{\overline{1}},\boldsymbol{\overline{1}} )\to & 2 & 1 & \cfrac{7}{27} & 0.259259\\
[-\halfbls]  &&&&\\ \hline  &&&& \\ [-\halfbls]
\to(\boldsymbol{\overline{0+0}} )\to(\boldsymbol{\overline{0+0}} )\to(\boldsymbol{\overline{0}},\boldsymbol{\overline{0}} )\to & 2 & 1 & \cfrac{35}{432} & 0.0810185\\
[-\halfbls]  &&&&\\ \hline  &&&& \\ [-\halfbls]
\to(\boldsymbol{\overline{0+0}} )\to(\boldsymbol{\overline{0+0}} )\Rightarrow(\boldsymbol{\overline{0+0}} )\to & 2 & 1 & \cfrac{35}{864} & 0.0405093\\
[-\halfbls]  &&&&\\ \hline
\end{array}
$$ 

\newpage 

\scriptsize
$$ 
\begin{array}{|c|c|c|c|c|}
\multicolumn{5}{c}{}\\
[-\halfbls]
\multicolumn{5}{c}{\text{ Configurations of closed geodesics.}}\\
\multicolumn{5}{c}{}\\
\multicolumn{5}{c}{\text{Component }\cH^{hyp}(3, 3).}\\
\multicolumn{5}{c}{}\\
\hline &&&&\\
\text{Degeneration pattern} & |\Gamma_-| & |\Gamma| & c\cdot\zeta(2) & c\cdot\zeta(2) \text{ approx.} \\
[-\halfbls]  &&&&\\ \hline  &&&& \\ [-\halfbls]
\Rightarrow(\boldsymbol{\overline{2}},\boldsymbol{\overline{2}} )\Rightarrow & 2 & 1 & \cfrac{15}{2} & 7.5\\
[-\halfbls]  &&&&\\ \hline  &&&& \\ [-\halfbls]
\to(\boldsymbol{\overline{0}},\boldsymbol{\overline{0}} )\to(\boldsymbol{\overline{1}},\boldsymbol{\overline{1}} )\to & 2 & 1 & \cfrac{35}{9} & 3.88889\\
[-\halfbls]  &&&&\\ \hline
%
%
\multicolumn{5}{c}{}\\
\multicolumn{5}{c}{}\\
\multicolumn{5}{c}{\text{Component }\cH^{nonhyp}(3, 3).}\\
\multicolumn{5}{c}{}\\
\hline &&&&\\
\text{Degeneration pattern} & |\Gamma_-| & |\Gamma| & c\cdot\zeta(2) & c\cdot\zeta(2) \text{ approx.} \\
[-\halfbls]  &&&&\\ \hline  &&&& \\ [-\halfbls]
\Rightarrow(\boldsymbol{\overline{2}},\boldsymbol{\overline{2}} )\Rightarrow & 2 & 1 & \cfrac{3699}{1120} & 3.30268\\
[-\halfbls]  &&&&\\ \hline  &&&& \\ [-\halfbls]
\to(\boldsymbol{\overline{0}},\boldsymbol{\overline{1}} ; 3)\to & 1 & 1 & \cfrac{64}{5} & 12.8\\
[-\halfbls]  &&&&\\ \hline  &&&& \\ [-\halfbls]
\Rightarrow(\boldsymbol{\overline{0+1}} ; 3)\Rightarrow & 1 & 1 & \cfrac{64}{35} & 1.82857\\
[-\halfbls]  &&&&\\ \hline  &&&& \\ [-\halfbls]
\to(\boldsymbol{\overline{0+0}} )\Rightarrow(\boldsymbol{\overline{2}},\boldsymbol{\overline{0}} )\to & 1 & 1 & \cfrac{27}{80} & 0.3375\\
[-\halfbls]  &&&&\\ \hline  &&&& \\ [-\halfbls]
\to(\boldsymbol{\overline{0}},\boldsymbol{\overline{0}} )\to(\boldsymbol{\overline{1}},\boldsymbol{\overline{1}} )\to & 2 & 1 & \cfrac{1}{5} & 0.2\\
[-\halfbls]  &&&&\\ \hline  &&&& \\ [-\halfbls]
\to(\boldsymbol{\overline{0}},\boldsymbol{\overline{0}} )\to(\boldsymbol{\overline{0+0}} )\Rightarrow(\boldsymbol{\overline{0+0}} )\to & 2 & 1 & \cfrac{1}{32} & 0.03125\\
[-\halfbls]  &&&&\\ \hline
%
%
\multicolumn{5}{c}{}\\
\multicolumn{5}{c}{}\\
\multicolumn{5}{c}{\text{Component }\cH^{even}(4, 2).}\\
\multicolumn{5}{c}{}\\
\hline &&&&\\
\text{Degeneration pattern} & |\Gamma_-| & |\Gamma| & c\cdot\zeta(2) & c\cdot\zeta(2) \text{ approx.} \\
[-\halfbls]  &&&&\\ \hline  &&&& \\ [-\halfbls]
\to(\boldsymbol{\overline{0}},\boldsymbol{\overline{0}} ; 4)\to & 2 & 1 & \cfrac{4}{3} & 1.33333\\
[-\halfbls]  &&&&\\ \hline  &&&& \\ [-\halfbls]
\Rightarrow(\boldsymbol{\overline{0+0}} ; 4)\Rightarrow & 2 & 1 & \cfrac{4}{21} & 0.190476\\
[-\halfbls]  &&&&\\ \hline  &&&& \\ [-\halfbls]
\to(\boldsymbol{\overline{0}},\boldsymbol{\overline{2}} ; 2)\to & 1 & 1 & \cfrac{256}{45} & 5.68889\\
[-\halfbls]  &&&&\\ \hline  &&&& \\ [-\halfbls]
\Rightarrow(\boldsymbol{\overline{0+2}} ; 2)\Rightarrow & 1 & 1 & \cfrac{256}{315} & 0.812698\\
[-\halfbls]  &&&&\\ \hline  &&&& \\ [-\halfbls]
\Rightarrow(\boldsymbol{\overline{1+1}} ; 2)\Rightarrow & 2 & 1 & \cfrac{8}{9} & 0.888889\\
[-\halfbls]  &&&&\\ \hline  &&&& \\ [-\halfbls]
\Rightarrow(\boldsymbol{\overline{1}},\boldsymbol{\overline{3}} )\Rightarrow & 1 & 1 & \cfrac{32768}{8505} & 3.85279\\
[-\halfbls]  &&&&\\ \hline  &&&& \\ [-\halfbls]
\to(\boldsymbol{\overline{1}},\boldsymbol{\overline{1}} ; 2)\to & 2 & 1 & \cfrac{128}{27} & 4.74074\\
[-\halfbls]  &&&&\\ \hline  &&&& \\ [-\halfbls]
\to(\boldsymbol{\overline{0+0}} )\to(\boldsymbol{\overline{0}},\boldsymbol{\overline{0}} ; 2)\to & 2 & 1 & \cfrac{8}{27} & 0.296296\\
[-\halfbls]  &&&&\\ \hline  &&&& \\ [-\halfbls]
\to(\boldsymbol{\overline{0}},\boldsymbol{\overline{0}} )\to(\boldsymbol{\overline{0+0}} ; 2)\to & 2 & 1 & \cfrac{16}{135} & 0.118519\\
[-\halfbls]  &&&&\\ \hline  &&&& \\ [-\halfbls]
\to(\boldsymbol{\overline{0}},\boldsymbol{\overline{0}} )\to(\boldsymbol{\overline{0}},\boldsymbol{\overline{2}} )\to & 1 & 1 & \cfrac{32}{45} & 0.711111\\
[-\halfbls]  &&&&\\ \hline  &&&& \\ [-\halfbls]
\Rightarrow(\boldsymbol{\overline{0+0}} )\Rightarrow(\boldsymbol{\overline{1+1}} )\Rightarrow & 2 & 1 & \cfrac{2}{45} & 0.0444444\\
[-\halfbls]  &&&&\\ \hline  &&&& \\ [-\halfbls]
\to(\boldsymbol{\overline{0+0}} )\Rightarrow(\boldsymbol{\overline{1}},\boldsymbol{\overline{1}} )\to & 1 & 1 & \cfrac{256}{1215} & 0.2107\\
[-\halfbls]  &&&&\\ \hline
\end{array}
$$ 

\newpage 

\scriptsize
$$ 
\begin{array}{|c|c|c|c|c|}
\multicolumn{5}{c}{}\\
[-\halfbls]
\multicolumn{5}{c}{\text{Component }\cH^{odd}(4, 2).}\\
\multicolumn{5}{c}{}\\
[-\halfbls]\hline &&&&\\
\text{Degeneration pattern} & |\Gamma_-| & |\Gamma| & c\cdot\zeta(2) & c\cdot\zeta(2) \text{ approx.} \\
[-\halfbls]  &&&&\\ \hline  &&&& \\ [-\halfbls]
\to(\boldsymbol{\overline{0}},\boldsymbol{\overline{0}} ; 4)\to & 2 & 1 & \cfrac{49}{18} & 2.72222\\
[-\halfbls]  &&&&\\ \hline  &&&& \\ [-\halfbls]
\Rightarrow(\boldsymbol{\overline{0+0}} ; 4)\Rightarrow & 2 & 1 & \cfrac{7}{18} & 0.388889\\
[-\halfbls]  &&&&\\ \hline  &&&& \\ [-\halfbls]
\to(\boldsymbol{\overline{0}},\boldsymbol{\overline{2}} ; 2)\to & 1 & 1 & \cfrac{147}{16} & 9.1875\\
[-\halfbls]  &&&&\\ \hline  &&&& \\ [-\halfbls]
\Rightarrow(\boldsymbol{\overline{0+2}} ; 2)\Rightarrow & 1 & 1 & \cfrac{21}{16} & 1.3125\\
[-\halfbls]  &&&&\\ \hline  &&&& \\ [-\halfbls]
\Rightarrow(\boldsymbol{\overline{1+1}} ; 2)\Rightarrow & 2 & 1 & \cfrac{3}{10} & 0.3\\
[-\halfbls]  &&&&\\ \hline  &&&& \\ [-\halfbls]
\Rightarrow(\boldsymbol{\overline{1}},\boldsymbol{\overline{3}} )\Rightarrow & 1 & 1 & \cfrac{128}{45} & 2.84444\\
[-\halfbls]  &&&&\\ \hline  &&&& \\ [-\halfbls]
\to(\boldsymbol{\overline{1}},\boldsymbol{\overline{1}} ; 2)\to & 2 & 1 & \cfrac{7}{2} & 3.5\\
[-\halfbls]  &&&&\\ \hline  &&&& \\ [-\halfbls]
\to(\boldsymbol{\overline{0+0}} )\Rightarrow(\boldsymbol{\overline{0+0}} ; 2)\to & 1 & 1 & \cfrac{7}{80} & 0.0875\\
[-\halfbls]  &&&&\\ \hline  &&&& \\ [-\halfbls]
\Rightarrow(\boldsymbol{\overline{0+0}} )\Rightarrow(\boldsymbol{\overline{0+2}} )\Rightarrow & 1 & 1 & \cfrac{21}{320} & 0.065625\\
[-\halfbls]  &&&&\\ \hline  &&&& \\ [-\halfbls]
\to(\boldsymbol{\overline{0+0}} )\Rightarrow(\boldsymbol{\overline{1}},\boldsymbol{\overline{1}} )\to & 1 & 1 & \cfrac{7}{45} & 0.155556\\
[-\halfbls]  &&&&\\ \hline  &&&& \\ [-\halfbls]
\to(\boldsymbol{\overline{0+0}} )\to(\boldsymbol{\overline{0}},\boldsymbol{\overline{0}} )\to(\boldsymbol{\overline{0}},\boldsymbol{\overline{0}} )\to & 2 & 1 & \cfrac{7}{144} & 0.0486111\\
[-\halfbls]  &&&&\\ \hline  &&&& \\ [-\halfbls]
\to(\boldsymbol{\overline{0+0}} )\Rightarrow(\boldsymbol{\overline{0+0}} )\Rightarrow(\boldsymbol{\overline{0+0}} )\to & 2 & 1 & \cfrac{7}{576} & 0.0121528\\
[-\halfbls]  &&&&\\ \hline
%
%
\multicolumn{5}{c}{}\\
\multicolumn{5}{c}{}\\
[-\halfbls]\multicolumn{5}{c}{\text{Component }\cH^{even}(2, 2, 2).}\\
\multicolumn{5}{c}{}\\
[-\halfbls]\hline &&&&\\
[-\halfbls]
\to(\boldsymbol{\overline{0}},\boldsymbol{\overline{0}} ; 2, 2)\to & 2 & 1 & \cfrac{180}{37} & 4.86486\\
[-\halfbls]  &&&&\\ \hline  &&&& \\ [-\halfbls]
\Rightarrow(\boldsymbol{\overline{0+0}} ; 2, 2)\Rightarrow & 2 & 1 & \cfrac{45}{74} & 0.608108\\
[-\halfbls]  &&&&\\ \hline  &&&& \\ [-\halfbls]
\Rightarrow(\boldsymbol{\overline{1}},\boldsymbol{\overline{1}} ; 2)\Rightarrow & 2 & 1 & \cfrac{225}{37} & 6.08108\\
[-\halfbls]  &&&&\\ \hline  &&&& \\ [-\halfbls]
\to(\boldsymbol{\overline{0}},\boldsymbol{\overline{0}} )\to(\boldsymbol{\overline{0}},\boldsymbol{\overline{0}} ; 2)\to & 2 & 1 & \cfrac{225}{296} & 0.760135\\
[-\halfbls]  &&&&\\ \hline  &&&& \\ [-\halfbls]
\Rightarrow(\boldsymbol{\overline{0+0}} )\Rightarrow(\boldsymbol{\overline{1}},\boldsymbol{\overline{1}} )\Rightarrow & 2 & 1 & \cfrac{5}{37} & 0.135135\\
[-\halfbls]  &&&&\\ \hline
%
%
\multicolumn{5}{c}{}\\
\multicolumn{5}{c}{}\\
[-\halfbls]\multicolumn{5}{c}{\text{Component }\cH^{odd}(2, 2, 2).}\\
\multicolumn{5}{c}{}\\
[-\halfbls]\hline &&&&\\
[-\halfbls]
\to(\boldsymbol{\overline{0}},\boldsymbol{\overline{0}} ; 2, 2)\to & 2 & 1 & \cfrac{252}{31} & 8.12903\\
[-\halfbls]  &&&&\\ \hline  &&&& \\ [-\halfbls]
\Rightarrow(\boldsymbol{\overline{0+0}} ; 2, 2)\Rightarrow & 2 & 1 & \cfrac{63}{62} & 1.01613\\
[-\halfbls]  &&&&\\ \hline  &&&& \\ [-\halfbls]
\Rightarrow(\boldsymbol{\overline{1}},\boldsymbol{\overline{1}} ; 2)\Rightarrow & 2 & 1 & \cfrac{144}{31} & 4.64516\\
[-\halfbls]  &&&&\\ \hline  &&&& \\ [-\halfbls]
\Rightarrow(\boldsymbol{\overline{0+0}} )\Rightarrow(\boldsymbol{\overline{0+0}} ; 2)\Rightarrow & 2 & 1 & \cfrac{9}{155} & 0.0580645\\
[-\halfbls]  &&&&\\ \hline  &&&& \\ [-\halfbls]
\Rightarrow(\boldsymbol{\overline{0+0}} )\Rightarrow(\boldsymbol{\overline{1}},\boldsymbol{\overline{1}} )\Rightarrow & 2 & 1 & \cfrac{16}{155} & 0.103226\\
[-\halfbls]  &&&&\\ \hline  &&&& \\ [-\halfbls]
\to(\boldsymbol{\overline{0}},\boldsymbol{\overline{0}} )\to(\boldsymbol{\overline{0}},\boldsymbol{\overline{0}} )\to(\boldsymbol{\overline{0}},\boldsymbol{\overline{0}} )\to & 2 & 3 & \cfrac{4}{93} & 0.0430108\\
[-\halfbls]  &&&&\\ \hline  &&&& \\ [-\halfbls]
\Rightarrow(\boldsymbol{\overline{0+0}} )\Rightarrow(\boldsymbol{\overline{0+0}} )\Rightarrow(\boldsymbol{\overline{0+0}} )\Rightarrow & 2 & 3 & \cfrac{1}{186} & 0.00537634\\
[-\halfbls]  &&&&\\ \hline
\end{array}
$$ 


\newpage

\end{document}